\documentclass{book}
\input{amssymb.sty}
\usepackage{makeidx}
\makeindex

\def\Ind#1#2{#1\setbox0=\hbox{$#1x$}\kern\wd0\hbox to 0pt{\hss$#1\mid$\hss}
\lower.9\ht0\hbox to 0pt{\hss$#1\smile$\hss}\kern\wd0}
\def\dnf{\mathop{\mathpalette\Ind{}}}
\def\Notind#1#2{#1\setbox0=\hbox{$#1x$}\kern\wd0\hbox to 0pt{\mathchardef
\nn=12854\hss$#1\nn$\kern1.4\wd0\hss}\hbox to
0pt{\hss$#1\mid$\hss}\lower.9\ht0 \hbox to
0pt{\hss$#1\smile$\hss}\kern\wd0}
\def\df{\mathop{\mathpalette\Notind{}}}

\def\ind#1#2#3{{  { #1 {{\dnf}^s_{#3}} #2 } }}

\def\stp{\mathop{\rm stp}\nolimits}
\def\CB{\mathop{\rm CB}\nolimits}
\def\Gal{\mathop{\rm Gal}\nolimits}
\def\VS{\mathop{\rm VS}\nolimits}
\def\Mult{\mathop{\rm Mult}\nolimits}
\def\Cb{\mathop{\rm Cb}\nolimits}
\def\wt{\mathop{\rm wt}\nolimits}
\def\prwt{\mathop{\rm prwt}\nolimits}
\def\Int{\mathop{\rm Int}\nolimits}

\def\alg{\mathop{\rm alg}\nolimits}
\def\rep{\mathop{\rm rep}\nolimits}

\def\GL{\mathop{\rm GL}\nolimits}
\def\Res{\mathop{\rm Res}\nolimits}
\def\dom{\mathop{\rm dom}\nolimits}
\def\rad{\mathop{\rm rad}\nolimits}
\def\Aut{\mathop{\rm Aut}\nolimits}
\def\Inv{\mathop{\rm Inv}\nolimits}
\def\RED{\mathop{\rm RED}\nolimits}

\def\eq{\mathop{\rm eq}\nolimits}
\def\div{\mathop{\rm div}\nolimits}
\def\acl{\mathop{\rm acl}\nolimits}
\def\dcl{\mathop{\rm dcl}\nolimits}
\def\tp{\mathop{\rm tp}\nolimits}
\def\trdeg{\mathop{\rm trdeg}\nolimits}
\def\rk{\mathop{\rm rk}\nolimits}
\def\max{\mathop{\rm max}\nolimits}

\def\res{\mathop{\rm res}\nolimits}
\def\ACVF{\mathop{\rm ACVF}\nolimits}

\def\Th{\mathop{\rm Th}\nolimits}

\def\qed{\hfill$\Box$}
\def\char{\mathop{\rm char}\nolimits}
\def\phi{\varphi}
\def\epsilon{\varepsilon}
\def\RM{\mathop{\rm RM}\nolimits}
\def\DM{\mathop{\rm DM}\nolimits}
\def\min{\mathop{\rm min}\nolimits}
\def\B{{\cal B}}
\def\L{{\cal L}}

\def\U{{\cal U}}

\def\bF{\bar{F}}
\def\br{\bar{\rho}}

\def\GG{{\cal G}}
\def\Aa{\mathbb A}
\def\Qq{\mathbb Q}
\def\Nn{\mathbb N}
 \def\Cc{{\mathbb C}}
  \def\Rr{{\mathbb R}}
 \def\Zz{{\mathbb Z}}
 \def\inv{{^{-1}}}
 
\def\F{{\cal F}}
\def\QED{\qed \medskip}

\def\meet{\cap}
\def\union{\cup}

\def\to{\rightarrow}

\def\acl{{\rm acl}}

\def\tensor{\otimes}

\def\dom{{\rm dom}}

\def\lceil{\ulcorner}
\def\rceil{\urcorner}

\def\Hct{\mathop{\rm Hct}\nolimits}
\def\Sub{\mathop{\rm Sub}\nolimits}
\def\wt{\mathop{\rm wt}\nolimits}
\def\prwt{\mathop{\rm prwt}\nolimits}
\def\Int{\mathop{\rm Int}\nolimits}

\def\rep{\mathop{\rm rep}\nolimits}

\def\GL{\mathop{\rm GL}\nolimits}
\def\Res{\mathop{\rm Res}\nolimits}
\def\dom{\mathop{\rm dom}\nolimits}
\def\rad{\mathop{\rm rad}\nolimits}
\def\Aut{\mathop{\rm Aut}\nolimits}

\def\eq{\mathop{\rm eq}\nolimits}
\def\div{\mathop{\rm div}\nolimits}
\def\acl{\mathop{\rm acl}\nolimits}
\def\dcl{\mathop{\rm dcl}\nolimits}
\def\tp{\mathop{\rm tp}\nolimits}
\def\trdeg{\mathop{\rm trdeg}\nolimits}
\def\rk{\mathop{\rm rk}\nolimits}
\def\max{\mathop{\rm max}\nolimits}

\def\red{\mathop{\rm red}\nolimits}
\def\res{\mathop{\rm res}\nolimits}

\def\Th{\mathop{\rm Th}\nolimits}

\def\qed{\hfill$\Box$}
\def\char{\mathop{\rm char}\nolimits}
\def\phi{\varphi}
\def\epsilon{\varepsilon}
\def\RM{\mathop{\rm RM}\nolimits}
\def\proves{\vdash}
\def\ct{\mathop{\rm ct}\nolimits}
\def\min{\mathop{\rm min}\nolimits}
\def\SL{\mathop{\rm SL}\nolimits}
\def\th{\mathop{\rm th}\nolimits}
\def\red{\mathop{\rm red}\nolimits}
\def\St{\mathop{\rm St}\nolimits}
\def\st{\mathop{\rm st}\nolimits}

\def\Ast{A^{\rm st}}
\def\Bst{B^{\rm st}}

\def\Dst{D^{\rm st}}
\def\Mst{M^{\rm st}}
\def\ast{a^{\rm st}}

\def\dst{d^{\rm st}}
\def\xst{x^{\rm st}}

\def\domind{{\dnf}^{d}}

 \def\g{\gamma}
 \def\k{{k}}

\def\meet{\cap}
\def\union{\cup}

\def\to{\rightarrow}

\def\acl{{\rm acl}}

\def\dom{{\rm dom}}

\def\Id{\mathop{\rm id}\nolimits}
\def\Max{{\rm Max\,}}

\def\lceil{\ulcorner}
\def\rceil{\urcorner}

\def\G{\Gamma}
\def\Keq{K^{\eq}}

\def\M{{\cal M}}

\def\Ueq{\U^{\eq}}
\def\B{{\cal B}}
\def\L{{\cal L}}

\def\U{{\cal U}}

\def\S{{\cal S}}
\def\T{{\cal T}}

\def\TB{\mathop{\rm TB}\nolimits}

\newtheorem{defi}{Definition}[chapter]
\newtheorem{theorem}[defi]{Theorem}
\newtheorem{definition}[defi]{Definition}
\newtheorem{lemma}[defi]{Lemma}
\newtheorem{proposition}[defi]{Proposition}

\newtheorem{corollary}[defi]{Corollary}

\newtheorem{remark}[defi]{Remark}
\newtheorem{example}[defi]{Example}

\begin{document}

\title{Stable domination \\ and \\ independence in algebraically closed valued fields}
\author{Deirdre Haskell
\and Ehud Hrushovski 
\and Dugald Macpherson}

\maketitle


\tableofcontents

%
%

\section{Preface}

Valuations are among the fundamental structures of number theory and of algebraic geometry.
This was recognized early by model theorists, with gratifying results:  Robinson's description
\cite{rob}
of algebraically closed valued fields as the model completion of the theory of valued fields,
the Ax-Kochen, Ershov study of Henselian fields of large residue characteristic with the application to Artin's conjecture \cite{ak1,ak2,ak3,er}, work of Denef and others on integration,
work of Macintyre, Delon,  Prestel, Roquette, Kuhlmann and   others 
on $p$-adic fields and positive characteristic.   The model theory of valued fields is 
thus one
of the most established and deepest areas of the subject.  

However, precisely because of the complexity of valued fields, much of the work centers on quantifier 
elimination and basic properties of formulas.  
Few tools are available for a more structural model-theoretic analysis.  This contrasts
with the situation for the classical model complete theories, of algebraically closed and real closed fields, 
where stability theory and o-minimality make possible a study of the category of definable sets.
Consider for instance the statement that  fields interpretable over ${\mathbb C}$ are finite or algebraically closed. 
 Quantifier elimination by itself is of little use in proving this statement.  
One uses instead the notion of   $\omega$-stability; it is preserved under interpretation,
implies a chain condition on definable subgroups,
and   by a theorem of Macintyre, $\omega$-stable fields are algebraically closed.    
With more analysis, using notions such as generic types, one can show 
that indeed every interpretable field is finite or definably isomorphic to
$ {\mathbb C}$ itself.  
This method can be extended to differential and difference  fields.  Using a combination
of such methods and of ideas of manifolds and Lie groups in a definable setting, 
Pillay was able to prove similar results for  
fields definable over  ${\mathbb R}$ or  $ {\mathbb Q_p}$.  But just a step beyond,
a description of interpretable fields seems out of reach of  the classical
methods.    When
$p$-adic or valuative geometry enters in an essential way, an intrinsic analog
of the notion of generic type becomes necessary.

For another example, take the notion of connectedness.  In stability, no topology is given in advance, but one manages to define stationarity of types or connectedness of definable groups, by looking at the type space.  In o-minimality, natural topologies on definable sets exist, and
connectedness, defined in terms of {\em definable paths}, is a central notion.   In valued fields,
the valuation topology is analogous to the o-minimal order topology, and one has the linearly
ordered value group that may serve as the domain of a path; but every continuous definable map from the value group to the field is 
constant, and a model-theoretic definition of connectedness is missing.   
The lack of such structural
model-theoretic understanding of valued fields is a central  obstacle to a wider interaction of model theory with geometry in general.

It is this gap that the present monograph is intended to address.   We suggest an approach
with two components.  We identify
a certain subset of the type space, the set of {\em stably dominated types}, that
behaves in many ways like the types in a stable theory.  Since these types are not literally
stable, it is necessary to first describe abstractly an extension of stability theory that includes
them.    Secondly, we note the existence of o-minimal families of
stably dominated types, and show, at least over sufficiently rich bases, that any type can 
be viewed as a limit of such a family.   This requires imaginaries in a concrete form,
serving as canonical bases of stably dominated types; and a theory of definable maps from $\G$
into such sets of imaginaries.   Thus, whereas type spaces work best in stable theories, and definable 
sets and maps in o-minimal theories, we suggest here an approach mixing the two.  
  As both the method and the intended applications depend heavily
on imaginary elements, we develop some techniques for dealing with them, including
prime models that often allow a canonical passage from imaginary to real bases.

We work throughout with the model completion, algebraically closed valued fields. 
 This is analogous to Weil's program of understanding geometry first at the level of algebraic closure.  One can hope
that the geometry of other valued fields could also, with additional work, be elucidated by this approach.  As an example of this viewpoint, consider the known 
elimination of quantifiers
for Henselian fields of residue characteristic zero, relative to the value group and 
residue field.  This was originally derived as an independent theorem.   But it is also 
an immediate consequence of a  
fact about algebraically closed valued fields of characteristic zero, namely that over any subfield $F$,
 any $F$-definable set is definably isomorphic to pullbacks of definable subsets
of the residue field and value group.  See \cite{hk} for more details on this short argument.
It is also noted there that definable sets in $F^n$ can be fibered over the residue
field and value group in fibers that are ACVF-definable; and a similar statement 
for definable bijections between them can be made.    Another example is
the elimination of imaginaries  for $\Qq_p$, proved in \cite{hmartin} using the ACVF methods
of invariant types.
 We illustrate the phenomenon  in Chapter~\ref{henselian}
using  the theory of $\Cc((t))$.

 We restrict ourselves in this monograph to laying the foundations of this approach.  
Definable maps from $\G$ into imaginaries were described in \cite{hhm};
We   concentrate here on the  theory of stable domination.  
Only future work based on these foundations can show to what extent they are successful.
We note here that Pillay has defined an analogue of stable domination, {\em compact domination}, that 
appears to be useful for thinking about  o-minimal groups (cf. 
\cite{hpp}).  Some progress has been made with the analysis of definable
 groups using the present methods; a   notion of {\em metastability}  has been abstracted 
 (Definition~\ref{metastable}), and results obtained for  Abelian groups in a general 
 metastable setting,  and linear groups interpretable in ACVF.  
See \cite{hrush2}.  
  And in very recent work of one of the authors with Loeser, connections
 with the Berkovich theory of rigid analytic spaces are beginning to emerge.

\bigskip

{\em Acknowledgements.} This work has benefitted from discussions over ten years with
 many mathematicians, and it is hard to single out individuals. The research of Haskell was supported by NSF and NSERC  grants,
 and that of Hrushovski 
 was partially conducted whilst he was a Clay Mathematics Institute Prize Fellow, and was also partially supported by the Israel 
Science Foundation. 
 The work of Macpherson was supported by EPSRC and LMS grants. All three authors were supported by  NSF funding for parts of a 
semester at the 
Mathematical Sciences Research Institute, Berkeley.

%
%
\chapter{Introduction}

As developed in \cite{shelah}, stability theory is based on the notion of an  {\em invariant} type,
  more specifically a {\em definable type}, and the closely related theory  of {\em independence of substructures}.   
We will review the   definitions in Chapter~\ref{stabilityintro} below;
suffice it to recall here that an (absolutely) invariant type gives a recipe yielding, for any
substructure $A$ of any model of $T$, a type $p|A$, in a way that respects
elementary maps between substructures; in general one relativizes to a set $C$
of parameters, and considers only $A$ containing $C$.
  Stability arose in response to questions in pure model theory, but 
has also provided effective tools for the analysis of algebraic and geometric structures.  
 The theories  of algebraically and differentially closed fields are stable, and the stability-theoretic analysis of types in these theories provides considerable information about
algebraic and differential-algebraic varieties.  The model companion of the theory of fields with an automorphism is not quite stable, but satisfies the related hypothesis of simplicity;
in an adapted form, the theory of independence remains valid and has served well  in applications to difference fields and definable sets over them.    On the other hand,
such tools have  played a rather limited role, so far,  in o-minimality and its applications to real geometry.  

Where do valued fields lie?  Classically, local fields are viewed as closely analogous to the real numbers.    We take a  
``geometric'' point of view however,  in the sense of Weil, and adopt the model completion   as the 
setting for our study.  This is Robinson's theory ACVF of algebraically closed valued fields.  We will view valued fields  
as substructures of models of ACVF.  Moreover, we admit other substructures involving imaginary elements, notably codes for lattices; these have been classified in \cite{hhm}.    
This will be essential not only for increasing the strength of the statements, but even for formulating our basic definitions.   

A glance at ACVF reveals  immediately a stable part, the residue field $k$; and an o-minimal part,
the value group $\G$.  Both are stably embedded, and have the induced structure of an   algebraically closed field, and an 
 ordered divisible abelian group, respectively.  But they amount between them to  a
  small part  of the theory.  For instance, over the uncountable field $\Qq_p$, the residue
field has only finitely many definable points, and both $k$ and $\G$ are countable in the
model $\Qq_p^a$.  As observed by Thomas Scanlon \cite{scan}, ACVF is not stable over $\G$, 
in the sense of  \cite{shelah2}.    
  
We seek to show nevertheless that 
stability-theoretic ideas can play a significant role in the description of valued fields.  To this end we undertake two logically independent  but mutually motivating endeavors.  In Part I we introduce  an extension of stability theory.  We consider 
  theories that have a stable  part,   
  define the notion of a  {\em stably dominated} type, and study its properties.   
The idea is that a type can be controlled by a very small part, lying in the stable part;
by analogy, (but it is more than an analogy), a power series is controlled, with respect to the question of invertibility for 
 instance, by its constant coefficient.  
Given a large model $\U$ and
a set of parameters $C$ from $\U$, we  define $\St_C$ to be a many-sorted structure whose sorts 
are the $C$-definable stably embedded stable subsets of the universe. 
The basic relations of $\St_C$ are those given by
$C$-definable relations of $\U$. Then 
$\St_C(A)$ (the stable part of $A$) is the definable closure of $A$ in $\St_C$. 
We  write $A\domind_C B$ if $\St_C(A) \dnf \St_C(B)$ in the
stable structure $\St_C$ and $\tp(B/C\St_C(A))\proves \tp(B/CA)$, and  say that $\tp(A/C)$ is
{\em stably dominated} if,  for all $B$, whenever $\St_C(A)\dnf \St_C(B)$, we have $A\domind_C B$.
In this case $\tp(A / \acl(C))$ 
lifts uniquely to an $\Aut(\U/ \acl(C))$-invariant type $p$.    Base-change results (under an extra assumption
 of existence of invariant extensions of types) show that 
if $p$ is also $\Aut(\U/ \acl(C'))$-invariant then $p| C'$ is stably dominated; hence, under this assumption,
stable domination is in fact
a property of this invariant type, and not of the particular base set.  We formulate
a general notion of domination-equivalence of invariant types (\ref{domination-equivalence}).
 In these terms,  an invariant type is stably dominated iff it is domination-equivalent
 to a type of elements in a stable part $\S_C$.

Essentially the whole forking calculus becomes available for stably dominated types.
Properties such as definability, symmetry, transitivity, characterization in terms
of dividing, lift easily from
$\St_C$ to $\domind$.  Others, notably  the descent  part of base change,  require more work and in fact 
an additional assumption:   that 
 for any algebraically closed substructure $C \subseteq M \models T$, 
any type $p$ over $C$ extends to an $\Aut(M/C)$-invariant type $p'$ over $M$.  

We isolate a further property of definable types in stable theories.    Two functions are said to have the same {\em germ}
relative to an invariant type $p$  if they agree generically on $p$.  In the o-minimal context,
an example of this is the germ at $\infty$ of a function on $\Rr$.  
  Moving from the function to the germ one is able to abstract
 away from the artifacts of a particular definition.   In 
 stability, this is an essential substitute for a topology.
 For instance, if $f$ is a function into a sort $D$, one shows that the germ is internal to $D$;
 this need not be the case for a code for the function itself.  In many stable applications,
 the strength of this procedure 
depends on the ability to reconstruct a representative of the germ from the germ alone.  
 We say that a germ is {\em strong} if this is the case.  
 
It is easy to see the importance of strong germs for the coding of imaginaries.  One
wants to code a function; as a first approximation, code the germ of the function;
if the code is strong, one has succeeded in coding at least a (generic) piece of the function in
question.  If it is not, one seems to have nothing at all.  

 We show 
that germs of stably dominated invariant types are always strong.   The proof depends on a combinatorial lemma saying that finite set functions on pairs, with a certain triviality property on triangles, arise from a function on singletons; in this sense it evokes
 a kind of primitive 2-cohomology, rather as the fundamental combinatorial lemma behind simplicity has a feel of 2-homology.  Curiously, both can be proved using the fundamental lemma of stability.

In \cite{hrush2} is is shown that   stable domination works well with definable groups.   
A group $G$ is called generically metastable if it has a 
 translation invariant stably dominated definable type.  In this case there exists a unique
  translation invariant definable type; and   the stable domination can be witnessed by a   definable homomorphism
$h: G \to H$ onto a connected stable definable group.  Conversely, given such a homomorphism $h$, $G$  is generically metastable iff the fiber of $h$ above
a generic element of $H$ is a complete type.  Equivalently, for any definable subset $R$
of $G$, the set $Y$ of elements $y \in H$ such that $h \inv(y) $ is neither contained in, nor disjoint from $R$ is a small set; no finite union of translates of $Y$ covers $H$.
We show this in Theorem~\ref{groupdom}, again using strong germs.

The general theory is  at present developed locally, at the level of a single type.   
It is necessary to say when we expect it to be meaningful globally.  The condition cannot be that every type be stably dominated; this would imply stability.  Instead we would like
to say that uniformly definable families of stably dominated types capture, in some sense, all types. 
Consider  theories with a distinguished predicate $\G$, that we assume to be linearly 
ordered so as to sharply distinguish it from the stable part.    We define 
 a theory to be {\em metastable over $\G$} (Definition~\ref{metastable})
  if every type over an algebraically closed set 
 extends to an invariant type, and, over sufficiently rich base sets, every type
 falls into a $\G$-parameterized family of stably dominated types.     We show that this 
 notion is preserved under passage to imaginary sorts.  
 
The proviso of
 ``sufficiently rich base set'' is familiar from stability, where the primary
 domination results are valid only over sufficiently saturated models; a great deal of
 more technical work is then needed to obtain some of them over arbitrary base.
 The saturation requirement (over ``small'' base sets) is effective since types over
 a model are always based on a small set.   In the metastable context,    more
 global conditions incompatible with stability are preferred.  This will be discussed
 for ACVF below.   
 
 For some purposes, extensions of the base are harmless and the theory can be used directly.
 This is so for results asserting the existence of a canonical definable set or relation
 of some kind, since a posteriori the object in
 question is defined without extra parameters.  This occurred   in the
 classification of maps from $\G$ in \cite{hhm}.  Another instance is
 in \cite{hrush2}, where under  certain finiteness of rank assumptions, it is shown that 
  a metastable Abelian group is an extension of a group
 interpretable over $\G$ by a   definable direct limit of generically metastable groups.

In Part II we study ACVF.  
  This is a $C$-minimal theory, in the 
sense of  \cite{ms}, \cite{hm}:   there exists a uniformly
definable  family of equivalence relations, linearly ordered by refinement; their
classes are referred to as (ultrametric) {\em balls}; and 
any definable set (in 1-space) is a Boolean combination of balls.   
 In strongly minimal and o-minimal contexts, one often argues by induction on dimension,
fibering an $n$-dimensional set over an $n-1$-dimensional set with $1$-dimensional
fibers, thus reducing many questions to the one-dimensional case over parameters.
This can also be done in the $C$-minimal context.  Let us call this procedure ``d\'evissage".

A difficulty arises:    many such arguments require canonical 
 parameters, not available in the field sort alone.     And  
certainly all our notions, from algebraic closure to stable embeddedness, must
be understood with imaginaries.   The imaginary sorts of ACVF were given
concrete form in \cite{hhm}:  the spaces $S_n$ of $n$-dimensional lattices, 
and certain spaces $T_n$, fibered over $S_n$ with fibers isomorphic to finite dimensional
vector spaces over the residue field.  But though concrete, these are not in any
sense one-dimensional;  attempting to reduce complexity by 
induction on the number of coordinates only leads to subsets of $S_n$, which is hardly simpler than $(S_n)^m$.  

Luckily, $S_n$ itself admits a sequence of fibrations $S_n =X_N \to X_{N-1} \to \ldots \to X_0$,
with $X_0$ a point and such that the fibers of $X_{i+1} \to X_i$ are $o$- or $C$-minimal. This uses the transitive action of
the solvable group of upper triangular matrices on $S_n$; 
see  the paragraph following Proposition~\ref{getun}.    There is a similar
statement for $T_n$ (where strongly minimal fibers also occur.)  It follows that any 
definable set of imaginaries admits a sequence of fibrations with successive
fibers that are strongly, $o$- or $C$-minimal (``unary sets''), or finite.
 D\'evissage arguments are thus possible. 

One result obtained this way is the existence
of invariant extensions.  A type over a base set $C$ can only have
an invariant extension if it is {\em stationary}, i.e. implies a complete type over
$\acl(C)$.  We show that in ACVF, every stationary type over $C$ has an
$\Aut(\U/C)$-invariant extension.  For $C$-minimal 
sets (including strongly minimal and $o$-minimal ones), there is a standard
choice of invariant extension:  the extension avoiding balls of radius smaller
than necessary.  

But this does not suffice to set up an induction, since for 
finite sets there is no invariant extension at all.  Thus a minimal step of induction
consists of   {\em finite covers of $C$-minimal sets}, i.e. with 
sets $Y$ admitting a finite-to-one map $\pi: Y \to X$, with $X$ unary.
This is quite typical of ACVF, and   
 resembles algebraic geometry, where 
 d\'evissage can reduce as far as  {\em curves} 
but not to a single variable.    In  the o-minimal case, 
by contrast,   one can  
do induction on   ambient dimension, or the number of coordinates of a tuple;
this explains much of the more ``elementary'' feel of basic o-minimality vs. strong minimality.  

   The additional ingredient needed
to  obtain invariant extensions of types is   
  the
 {\em stationarity  lemma} from \cite{hhm}, implying that if $\pi$
  admits a section
 over a larger base, then it admits a section over $\acl(C)$.  See Proposition~\ref{endofhhm}. 
For the theory  ACF over a perfect field, stationarity corresponds to the notion of a regular extension, and the stationarity lemma to the existence of a geometric notion of
irreducibility of varieties.  
It is instructive to recall the proof for ACVF.   Given a finite cover $\pi: Y \to X$ as above,
a section $s$ of $Y$ will have a strong germ with respect to the canonical
invariant extension of any  type of $X$.  
Generic types of closed balls are stably dominated; for these,   by 
the results of Part I, all functions have strong germs.  Other types are viewed
as limits of definable maps from $\G$ into the space of generics of closed balls.
For instance if $\tilde{b}$ is an open ball, consider the family of closed sub-balls 
$b$ of $\tilde{b}$; these can be indexed by their radius $\gamma \in \Gamma$ the moment one 
fixes a point in $\tilde{b}$; by the stably dominated case, one has a section
of $\pi$ over each $b$.  The classification of definable maps from $\G$
(actually from finite covers of $\G$)  is then used
to glue them into a single section, over the original base.   This could be done abstractly
for $C$-minimal theories whose associated (local) linear orderings satisfy $\dcl(\G)=\acl(\G)$.
The proof of elimination of imaginaries itself has a similar structure; see a sketch
at the end of Chapter~\ref{Jindependence}.    

Another application of the unary decomposition is the existence of canonical resolutions,
or prime models.  
In the field sorts, 
ACVF has prime models trivially; the prime model over a nontrivially valued field $F$
is just the algebraic closure $F^{alg}$.  In the geometric sorts the situation becomes
more interesting.  The algebraic closure does not suffice, but we show that finitely generated
structures (or structures finitely generated over models) do admit canonical prime models.  
A key point is that the  prime model over a finitely generated structure $A$ 
add to $A$ no elements of the residue field or value group.  This is important 
in the theory of motivic integration; see the discussion of resolution in \cite{hk}.   
A further application of canonical resolution  is a quantifier-elimination 
for $\Cc((t))$ in the $\GG$-sorts, relative to the value group $\G$.  In essence   resolution is 
used to produce functions on imaginary sorts; in fact for any $\GG$-sort represented
as $X/E$ and any function $h$ on $X$ into the value group or residue field, there exists a function $H$ on $X/E$ such that
 $H(u)=h(x)$ for some $x\in X/E$.  
 
The construction of prime models combines the decomposition 
into unary sets with the idea of opacity.  An equivalence relation $E$ on $X$ is called {\em opaque} if any definable subset of $X$
 is a union of classes of $E$, up to a set contained in finitely many classes.  This is another
 manifestation of a recurring theme.    Given $f: X \to Y$
 and an ideal $I$ on $Y$, we say $X$ is dominated by $Y$ via $(f,I)$ if for any subset $R$
 of $X$, for $I$-almost every $y \in Y$, the fiber $f \inv (y)$ is contained in $R$ or is disjoint from it.  For stable domination, $Y$ is stable and $I$ is the forking ideal; for stationarity, $f$ has finite fibers, and  $I$ is the dual ideal to an invariant type;   for opacity, $I$ is the ideal of finite sets.
The equivalence relations associated with the  analyses of $S_n$ and $T_n$ above are opaque.
For an opaque equivalence
relation, all elements in a non-algebraic class have the same type (depending only on the
class); this gives a way to choose elements in such a non-algebraic class canonically up to isomorphism.    Algebraic classes are dealt with in another way.  

We now discuss the appropriate notion of a ``sufficiently rich'' structure.  In the
stable part,  saturation is the right requirement; this will not actually be felt
in the present work, since the stable part is $\aleph_1$-categorical and does
not really need saturation.  For the o-minimal part, a certain completeness condition
turns out to be useful; see Chapter 13.2.   It allows the description of the semi-group of invariant types up to domination-equivalence, and   a characterization of forking in ACVF over
very rich bases.  For the most part however neither of these play any role; the significant
condition is richness over the stable and the o-minimal parts.  Here we adopt
Kaplansky's maximally complete fields.  An   algebraically closed  valued 
is maximally complete if it has  no proper immmediate extensions.  It follows from
(\cite{kap1}, \cite{kap2}) that   any model of ACVF
embeds  in a maximally complete field, uniquely up to isomorphism.   Since we use all the geometric sorts, 
a `rich base' for us is a model of ACVF whose field part is maximally complete.  

Over such a base $C$, we prove first, using standard results on finite dimensional vector spaces over maximally complete fields, that any field extension $F$ is dominated by
its parts in the residue field $k(F)$ and the value group $\G_F$.     This kind of domination does not admit descent.   A stronger statement is that $F$ is dominated by the stable part 
{\em over} $C$ together with $\G_F$, so that the type of any element of $F^n$ over $C \union \G_F$ is stably dominated.  After an algebraic interpretation of this statement, 
it is deduced from the previous one by a perturbation argument.   Both  
these results are then  extended 
to imaginary elements.  

We interpret the last result as follows:  an arbitrary type 
lies in a family of stably dominated types, definably indexed by $\G$.  Note that $\k$
and $\G$ play asymmetric roles here.  Indeed, at first approximation, 
we develop what can be thought of as the model theory of $\k^\G$, rather than 
$\k \times \G$.  However $\k^\G$ is presented by a $\G$-indexed
system of {\em opaque} equivalence relations, each hiding the structure on the finer ones until a specific class is chosen.  
 This kind of phenomenon, with hidden forms of $\k^n$ given by 
finitely many nested equivalence relations, is familiar from stability theory;
the presence of a definable directed system of levels is new here.

Even for fields, the  stable domination in the stronger statement cannot be understood  without imaginaries. 
Consider a  field extension $F$ of $C$; for simplicity suppose the value group 
$\G_F$ of $F$ is generated over $\G_C$ by one element $\g$. 
  There is then a canonical vector space $V_\g$
over the residue field.  If $\g$ is viewed as a code for a closed ball $E_\g= \{x: v(x) \geq \g \}$, the elements
of $V_\g$ can be taken to be codes   for the maximal open sub-balls of  $E_\g$.  The vector space
$V_\g$ lies in the stable part of the theory, over $C(\g)$.
  We show
that $F$ is dominated over $C(\g)$ by elements of  $k(F) \cup V_\g(F)$.    Note that $k(F)$
may well be empty. 
 
Over arbitrary bases, invariant types orthogonal to the value group are shown to be dominated
by their stable part; this follows from existence of invariant extensions, and descent.

At this point, we have the metastability of ACVF. 
   We now seek to relate this still somewhat abstract picture more directly with the geometry of valued fields.   We characterize the stably dominated types as those invariant types that
   are orthogonal to the value group (Chapter 8.)  
    In Chapter 14, we describe geometrically the connection between a stably dominated type $P$ and the associated invariant type $p$, when $P$ is contained in an algebraic variety $V$.
In the case of ACF, the invariant  extension is obtained by avoiding all proper subvarieties.   In ACVF, the demand is not only to avoid but to stay as far away as possible
from any given subvariety.  See Theorem~\ref{maxmodapp}.  In ACF the same prescription 
yields the unique invariant type of any definable set; it is not necessary to pass through
types.  In ACVF 
 the picture for general definable sets is more complicated.  But for  a definable subgroup $G$ of $\GL_n(K)$, or for a definable affine homogeneous space, we show that 
 a
 translation invariant stably dominated type is unique if it exists, and that in this case
 it is again the type of elements   of maximal distance from any proper subvariety of
 the Zariski closure of $G$.    
 
 In chapter 15 the ideas are similar, but the focus is on canonical bases.  Any definable type,
 in general, has a smallest substructure over which it is defined.  
 In ACF, this is essentially
  the field of definition of the associated prime ideal.  
 We obtain a similar geometric description 
   for stably dominated types; the ideal of regular functions vanishing on the type is replaced by the $R$-module of functions taking small values on it.

 While presented here for {\em stably dominated types}, where the theory flows
 smoothly from the main ideas,
   within the text we  try to work with weaker hypotheses on the types when possible.     Over sufficiently rich base structures, all our results
 can be read off from the main domination results discussed above.  But over 
 smaller bases this is not always the case,  leading us to think that
  perhaps a general principle  remains  to be discovered.  An example is the Theorem of Chapter 10, that an indiscernible
 sequence whose canonical base (in an appropriate sense) is orthogonal to $\G$,
 is in fact an indiscernible set, and indeed a Morley sequence for a stably dominated type.
Others are phrased in the language of {\em independence of substructures}.  

   Classical stability theory yields
a notion of independence of two substructures $A,B$
over their intersection $C$, defined in many equivalent ways.  
One is directly connected to invariant types:    
If $A$ is generated by elements $a$,  then $A,B$ are independent over $C$
iff $\tp(a/B)$ has a $C$-invariant extension to any model.   Intuitively
 `$A$ is independent from
$B$ over $C$' should say that `$B$ provides 
as little as possible extra information about $A$, beyond what 
$C$ provides'. 
In other words, the locus of $A$ over $B$ is as large as possible
inside the locus of $A$ over $C$.  
A number of the above ideas lead to    notions of independence for substructures
of models of ACVF, i.e. for valued fields.

The simplest  notion, {\em sequential independence} (Chapter 8), depends on the choice of an ordered tuple $a$ of  
generators of $A$ over $C$.    Let $p$ be the invariant extension of $\tp(a/C)$
constructed above by d\'evissage.  
  We  say that $A$ is sequentially independent from $B$ over $C$,
  $A  \dnf^g_C B$, if $\tp(a/B) \subset p$.      In general, the notion 
depends on the order of the tuple, and is not symmetric.     

 A point in an 
irreducible variety is generic if it does not lie in any smaller dimensional 
variety over the same parameters, and in an algebraically closed field, 
$A$ is independent from $B$ over $C$  if every tuple from
 $A$ which is generic in 
a variety defined over $C$ remains generic in the same variety with the additional 
parameters from $C \cup B$. Since varieties are defined by polynomial equations, this is
equivalent to saying that for every $a\in A$, the ideal of polynomials 
which vanish on the $\tp(a/C)$ is the same as those which vanish on $\tp(a/C \cup B)$.
We extend both of these points of view to an algebraically closed valued field.

In this setting, the definable sets depend on the valuation as well, so we consider 
the set of polynomials which satisfy a valuation inequality on $\tp(a/C)$. This is no
longer an ideal, but naturally gives a collection of modules over the valuation ring,
 which we call $J(\tp(a/C))$. We define $A$ to be 
$J$-independent from $B$ over $C$ if $J(\tp(a/C))=J(\tp(a/C \cup B))$ for all tuples $a$ from $A$
(Chapter 15).
This definition does not depend on the order of the tuple $a$.

Our final notion of independence is defined here only for variables of the field
sort. We define $\tp^+(A/C)$ to be the positive quantifier-free type
of $A$ over $C$, and say that $A$ and $B$ are {\em modulus independent} over
$C$ if $\tp^+(AB/C)$ is determined by the full types of $A$ and $B$ separately (Chapter 14).
In a pure algebraically closed field, the positive type corresponds 
precisely to the ideal of polynomials which vanish on the type, so modulus 
independence is in that setting another way of stating non-forking.
 In an algebraically closed valued field,
we use modulus independence  as a step from sequential independence to $J$-independence.
In this setting,  the quantifier-free positive formulas refer to the maximum norm 
that a polynomial can take on a type.  
 
All these notions agree    on stably dominated types, and   have  the good properties
of independence for stable theories (see Theorem~\ref{summary}).   Under more general conditions, they diverge,
 and various properties can fail; for instance, 
for  non-stably dominated types, $J$-independence need not be 
symmetric in $A$ and $B$, nor need 
  $\tp(A/C)$ have a $J$-independent extension over $C \cup B$.   We give numerous
  examples to showing this.  
    We do show 
  however that  if $A$ and $B$ are fields, and $C\cap K\leq A$ with $\Gamma(C)=\Gamma(A)$, then
sequential independence over $C$ implies both modulus independence and $J$-independence.

In the final chapter we briefly illustrate the idea mentioned in the preface, 
that the methods of this monograph should be useful for valued fields beyond ACVF.
Theorem~\ref{vdf} asserts that Scanlon's model completion of valued differential fields is metastable.
This gives for the first time a language to pose structural questions about this
rich theory, and we hope it will be fruitful.  We also show the metastability
of the theory of $\Cc((t))$, and related theories.  Here we prove nothing
 anew, reducing all questions to properties of ACVF. 
  The property of metastability itself can only hold for a limited number
of theories valued fields, but the method of reduction to the algebraic closure
 is much more general.

%
%

\part{Stable domination}

\chapter{Some background on stability theory}\label{stabilityintro}

We give here a brief preview of stability theory, as it underpins stable domination. 
We also introduce some of the model-theoretic notation used later.
Familiarity with the basic notions of
logic (languages, formulas, structures, theories, types, compactness) is assumed,
but we explain the model theoretic notions beginning with saturation, algebraic closure, imaginaries.  
 We have in mind a reader who is 
familiar with o-minimality or some model theory of valued fields,
 but has not worked with notions from stability.   Sources include Shelah's
 {\em Classification Theory} as well as  books by Baldwin \cite{bal}, Buechler \cite{buechler},  Pillay \cite{pil} and Poizat
\cite{poi}. There is also 
a broader introduction by Hart \cite{hart}  intended partly for non-model theorists, and an introduction to stability 
theory intended for a wider audience in \cite{hrush3} . Most of the stability theoretic results below should 
be attributed to Shelah.  Our treatment will mostly follow Pillay \cite{pil}.

Stability theory is a large body of abstract model theory developed in the 1970s and 1980s by Shelah and others, but having 
its roots 
in Morley's 1965 Categoricity Theorem: if a complete theory in a countable language is categorical in some uncountable power, 
then it is categorical in all uncountable powers.   Shelah formulated a radical generalization
of Morley's theorem, weakening the categoricity assumption from one isomorphism type
to any number less than the set-theoretic maximum.  The conclusion is that
all models of the theory, in any power, are classifiable by a small tree of numerical dimensions.
This can be viewed as a description and analysis of all  complete  theories 
in which the large models are classifiable by a small family of numerical invariants.   This work  brought out an impressive
list of model-theoretic ideas and properties beyond stability itself: superstability, regular types,
domination and orthogonality, higher properties such as shallowness, NDOP and NOTOP.
These yield sharp tools for taking structures apart, explaining much in terms of 
simpler embedded structures, the regular types.  These achievements are analogous to Ax-Kochen, Ershov theorems in valued fields, where many properties of Henselian valued fields
(with appropriate additional assumptions) are reduced to the residue field and value group.
However, in stability no special parts of the theory are given in advance.  All notions are defined
in terms of an  abstractly defined notion of
 independence.
 
Among fundamental mathematical theories, a few are stable: algebraically closed fields (ACF), and more generally, 
separably closed fields (SCF);
differentially closed fields (the model companion of the theory of fields with a differential operator, with theory denoted DCF);
 and modules over any ring.  By recent work of Sela, the theory of
free groups is also stable.

Since the mid 1980s, model-theorists have noticed stability-theoretic phenomena  in a number of structures which formally are 
unstable but have significant mathematical interest. 
These include certain structures with {\em simple theory}, where the notion of independence is slightly less constrained: for example,
smoothly approximable structures, pseudofinite fields, algebraically closed fields with a generic automorphism (ACFA).

Many structures of algebraic geometry can be interpreted without quantifiers in fields  or in differential or difference fields.  The stability or simplicity of ACF and SCF, DCF, and ACFA 
thus    opens
the way for the use of stability theoretic methods in algebraic geometry.
   However many 
other algebraic geometric and number theoretic  constructions involve valuations, and these were not accessible 
up to now.

 In a different
 direction, 
a totally ordered structure $(M,<,\ldots)$ is {\em o-minimal} \index{o-minimal} if every parameter-definable subset of $M$ is a finite union 
of singletons and open intervals.  Because of the total ordering, any o-minimal theory is unstable.  
 Ideas from stability theory have  influenced the development of o-minimality, but
 for the most part o-minimal technique uses topologies and ideas closer to classical
 geometric ones, especially in the local (or definably compact) parts of the theory.
   In Part~II, we make occasional reference
 to o-minimality, since the value group of an algebraically closed valued field is o-minimal, but essentially no knowledge
 of o-minimality is needed.

The theory of stable domination, which we develop in Part~I, provides a new 
setting where stability theory has application, for unstable structures with a stable constituent. The motivating example, as 
developed in Part~II of this monograph, is the theory of algebraically closed valued fields.

The stability-theoretic features of ACVF as developed in Part~II are seen most strongly through stable domination. 
 However, we also prove several other results which have analogues for stable theories.
For example, in Theorem~\ref{densestone} we show that in ACVF, over an algebraically closed set definable types are dense
in the Stone space (in a stable theory, all types are definable). In Corollary~\ref{lascar}, we show that
 in  ACVF strong type and Lascar strong type agree. And in Proposition~\ref{indiscseq}, we show that any indiscernible sequence over the value group is indiscernible as a {\em set} over the value group (in a stable theory, any indiscernible sequence is an indiscernible set). 
Furthermore, we develop several different theories of independence, each of which has some properties of independence defined through non-forking in a stable theory.
Below we summarise some basic stability-theoretic facts and terminology used later, emphasising stability
ideas which lift to ACVF. Proofs are omitted, and the presentation is mostly taken from Pillay \cite{pil}.

\medskip

\section{Saturation, the universal domain, imaginaries.}

\medskip
For this chapter, we assume that $\L$ is a first order language, and that $T$ is a complete theory over $\L$. 
Our convention  will be to use a single variable $x$ (rather than $\bar{x}$) for a finite sequence of variables, and a single parameter $a$ for a finite sequence of parameters. We often omit union signs, writing for example
$ABcd$ for $A\cup B \cup\{c,d\}$, where $A,B$ are sets of parameters. Often we 
write $a\equiv_C b$ to mean that $\tp(a/C)=\tp(b/C)$. We shall use symbols $A,B,C$ for sets of parameters, and 
$M,N$ for models of an ambient theory (usually denoted $T$) which is understood from the context. When we say that, for
 some model $M$, the set $X\subseteq M^n$ 
is {\em definable} \index{definable}, we mean that it is {\em definable with parameters}, i.e. that it is 
the solution set in $M^n$ of some formula
$\phi(x_1,\ldots,x_n,a_1,\ldots,a_m)$ where $a_1,\ldots,a_m\in M$.

Recall that if $C\subseteq M\models T$ then the {\em algebraic closure of $C$},
\index{algebraic closure} written $\acl(C)$ consists of those $c\in M$ 
such that, for some 
formula $\phi(x)$ over $C$, $\phi(M)$ is finite and contains $c$. The {\em definable closure} 
$\dcl(C)$ \index{definable closure} is the union of the
 1-element $C$-definable sets.
Thus, $C\subseteq \dcl(C)\subseteq \acl(C)$, $\dcl(C)$ is definably closed (that is, 
$\dcl(\dcl(C))=\dcl(C)$), and $\acl(C)$ is algebraically closed.

In general, the language $\L$ is assumed to be multi-sorted. By `type' we shall always mean `complete type 
over some small base set', and if completeness is not assumed we say `partial type'.
For any set $C$ of 
parameters in a model of $T$, and any set ${ x}$ of variables, we write $S_{x}(C)$
for the set of types over $C$ in the  variables ${ x}$, and $S(C)$ when the variables are 
understood from the context.  This is more correct
for many-sorted theories, since the variables in $  { x}$ carry with them a specification of the 
relevant sort.   
The space $S_x(C)$ is the Stone space of the Boolean algebra of formulas in  free variables $x$ up to
 equivalence over $T$, and is a compact totally disconnected topological space.
When we have a particular sort in mind, we write $S_n(C)$ for $S_{x_1,\ldots,x_n}(C)$,
with $x_i$ of that sort.  (Note that the set of $n$-types of all sorts combined is not compact.)
 The basic open sets have the form
$[\phi]=\{p\in S_{x}(C): \phi\in p\}$, where $\phi$ is a formula in   free variables among $x$ over $C$. In 
particular, an {\em isolated type} \index{type!isolated} is a   type which contains a formula $\phi(x)$
(an {\em isolating formula}) such that for all $a$ (in an ambient model of $T$), if $\phi(a)$ holds, then $a$ realises $p$. 

The set of variables $x$ is usually assumed to be finite, but for general questions this is not needed. 
Types in infinitely many variables are called $*$-types in \cite{shelah}, and they are often convenient to use.   
In particular, in this monograph we often talk of $\tp(A/C)$, where $A,C$ are infinite sets.  
Implicitly, we have in mind some enumeration $(a_i:i<\lambda)$ of $A$, and are considering a type
in the variables $(x_i:i<\lambda)$, consisting of all formulas $\phi(x_{i_1},\ldots,x_{i_n})$ (for any 
$n<\omega$) over $A$  such that
$\phi(a_{i_1},\ldots,a_{i_n})$ holds.   Alternatively, choose a variable $x_a$ of the appropriate sort for each $a \in A$,
and use the tautological correspondence of variable with element.  This makes sense
as long as one discusses properties of types that are invariant under permutations of the
variables.

Following Morley, we present various model-theoretic ranks in terms of Cantor-Bendixson ranks 
\index{rank!Cantor-Bendixson} of type spaces. Recall that if $X$ is a compact totally disconnected topological 
space then the Cantor-Bendixson rank
$\CB_X(p)$ of $p\in X$ is defined by transfinite induction as follows: $\CB_X(p) \geq 0$ for all $p\in X$,
and $\CB_X(p)=\alpha$ if $p$ is isolated in the closed subspace $\{q\in X: \CB_X(q)\geq \alpha\}$. If $\CB_X(p)<\infty$ for all $p\in X$, then by topological compactness $\{\CB_X(p):p\in X\}$ has a maximal element $\alpha$, and $\{p\in X:\CB_X(p)=\alpha\}$ is finite; its size is denoted by CB-Mult$(X)$.

If $\lambda$ is an infinite cardinal, the  model $M \models T$ is {\em $\lambda$-saturated},
 \index{saturated} if, for all $C\subset M$ with $|C|<\lambda$, $M$ realises all types in
$S(C)$; it is {\em saturated} if it is $|M|$-saturated.
We shall work in a universal domain $\U$, which is assumed to be a large model of $T$.
Any parameter sets $A,B,C,\ldots$ which we mention are subsets of $\U$, and any models 
$M,N,\ldots\models T$ are assumed 
to be elementary substructures of $\U$.  We write $\phi(a)$ to mean:  $\U \models \phi(a)$.
We often consider types $p\in S(\U)$, and occasionally realisations $a$ of $p$, which in general are not in $\U$.
It is common  to assume that $\U$ is saturated, since this guarantees that any elementary map between 
subsets of $\U$ of size less than $|\U|$ extends to
 an element of $\Aut(\U)$.  For unstable theories (such as ACVF), existence of saturated models depends on
 set-theoretic assumptions, but a variety of set-theoretic tools (such as absoluteness)
 can be used to dispense with these a posteriori.  For most arguments
 it suffices to assume that  $\U$ is {\em sufficiently saturated}.
 Some cardinal $\kappa$ is fixed, and sets of size $<\kappa$ are called {\em small};
 on the other hand $\U$ is assumed to be  $\kappa$-saturated   and
 {\em $\kappa$-homogeneous}, i.e. two elements realizing
the same type over a small set  $A$ are conjugate under the automorphism groups.  These two
conditions can be achieved without any special set-theoretic assumptions or tools;
$\kappa$ is chosen safely above the cardinalities of any objects of interest.    

If $T$ is a multi-sorted theory, $S$ is a sort of $T$, and $M\models T$, write $S(M)$ for the set 
of elements of $M$ of sort $S$. We  say that $T$ has 
{\em elimination of imaginaries} \index{imaginaries!elimination of}
 if, for any $M\models T$, any collection $S_1,\ldots,S_k$ of sorts in $T$, and any $\emptyset$-definable 
equivalence relation $E$ on $S_1(M)\times\ldots \times S_k(M)$, there is a $\emptyset$-definable function 
$f$ from $S_1(M)\times \ldots \times S_k(M)$ into 
a product of sorts of $M$, such that
for any $a,b\in S_1(M) \times\ldots \times S_k(M)$, we have $Eab$ if and only if $f(a)=f(b)$. 
Given a complete theory $T$, it is possible to extend it to a complete theory $T^{\eq}$ over a 
language $L^{\eq}$ by adjoining, 
for each collection $S_1,\ldots,S_k$ of sorts and $\emptyset$-definable equivalence relation $E$ on 
$S_1\times \ldots \times S_k$, a sort $(S_1\times \ldots \times S_k)/E$, together
 with a function symbol for the natural map $a\mapsto a/E$. Any $M\models T$ can be canonically
 extended to a model of $T^{\eq}$, denoted $M^{\eq}$.
 The theory $T^{\eq}$ automatically has elimination of imaginaries.
For the purposes of stability theory, elimination of imaginaries is very helpful. Therefore, in the development
 of stability theory in this chapter we 
shall assume that $T$ has elimination of imaginaries (though not necessarily that $T$ is
 formally  a theory of form $(T')^{\eq}$). We shall refer to the sorts of $T^{\eq}$ as imaginary sorts, and to 
elements of them as {\em imaginaries}.\index{imaginaries}

Suppose that $D$ is a definable set in $M\models T$, defined say by the formula $\phi(x,a)$. There is a
$\emptyset$-definable equivalence relation $E_\phi(y_1,y_2)$, where $E_\phi(y_1,y_2)$ holds if and only if
$\forall x(\phi(x,y_1)\leftrightarrow \phi(x,y_2))$. Now 
$a/E_\phi$ is identifiable with an element of an imaginary sort; it is determined uniquely
(up to interdefinability over $\emptyset$) by $D$, and will often be referred to as a {\em code}
\index{imaginaries!code for}
 for $D$, and denoted $\lceil D\rceil$. We prefer to think of $\lceil D \rceil$ as a fixed object 
(e.g. as  a member of $\U^{\eq}$) rather than
as an equivalence class of $M$; for  viewed as an equivalence class it is formally a different set
(as is $D$ itself) in elementary extensions of $M$.

Under our assumption that $T$ eliminates imaginaries, $\lceil D\rceil$ can be regarded as a tuple in $M$; without elimination
 of imaginaries we would just know it to be in $M^{\eq}$. An automorphism of $M$ will fix $D$ setwise 
if and only if it fixes $\lceil D \rceil$. 

In a purely theoretical context, the device of moving to $M^{\eq}$ gives a soft way of 
dealing with imaginaries.  However when working with specific theories 
a more concrete description of a collection of sorts admitting elimination of imaginaries
is useful.  It turns out that o-minimal fields, separably and differentially closed fields, 
and existentially closed difference fields all admit elimination of imaginaries in the field
sort.  For  algebraically closed valued fields, the field sort, value group and residue field
do not exhaust the necessary imaginaries.    
The main purpose of \cite{hhm} was to identify a geometrically meaningful family  of sorts
which, when adjoined to the field sort, suffice for elimination of imaginaries.  See Chapter~\ref{acvfbackground}
 for a description.

\medskip
  
\section{Invariant types}  \label{invariant}

Let $\Aut(\U/C)$ denote the subgroup of $\Aut(\U)$ fixing $C$ pointwise.
Then $\Aut(\U/C)$ acts naturally on $S(\U)$: if $g\in \Aut(\U/C)$ then
$g(p)=\{\phi(x,g(a)): \phi(x,a)\in p\}$.

\begin{definition} \label{invariantdef}\rm \index{type!invariant} \index{invariant type}

The type $q\in S(\U)$ is {\em $\Aut(\U/C)$-invariant} if it
 is fixed by this action.
 \end{definition}
 
 If $q$ is $\Aut(\U/C)$-invariant, then for any formula $\phi(x,y)$ and $a_1,a_2\in \U$, if 
 $a_1\equiv_C a_2$ then $\phi(x,a_1)\in q$ if and only if $\phi(x,a_2)\in q$, and assuming saturation this statement 
is equivalent to invariance. 
This says that $p\in S(\U)$ is $\Aut(\U/C)$-invariant precisely if it does not {\em split} over $C$, in the sense 
of Shelah \cite{shelah}.
This  gives
 an alternative definition of {\em invariant type}, without reference to automorphisms and saturation assumptions.
 
 Suppose that  $p\in S(\U)$ is $\Aut(\U/C)$-invariant, and $e_1,e_2\in \U$ with $e_1\equiv_C e_2$. Let $a\models p|Ce_1e_2$. Then $e_1 \equiv_{Ca}e_2$. This observation will be used frequently without explicit mention.

Call a type $p$ over $\U$ {\em   invariant} if 
it is   $\Aut(\U/A)$ -invariant for
some  small $A$ (in the sense of the previous section.)  
 Any such $A$ is called a {\em base} \index{type!base for} for $p$.

  Let $\Inv_x$ denote the set of invariant types in the variable $x$.  
If $p \in \Inv_x$, $q \in \Inv_y$ are  invariant types, define $p \tensor q$  as follows: let $A$ be small such that $p$ is
$\Aut(\U/A)$ -invariant; let $d \models q |A$, and let $c \models p | Ad$.  Then 
$\tp(cd/A)$ does not depend on the choice of $c,d$; call it 
$(p \tensor q) | A$.  Clearly there exists a unique type $p \tensor q$ over
$\U$ whose restriction to any such $A$ is $(p \tensor q) | A$.   Moreover
$p \tensor q$ is $\Aut(\U/B)$-invariant whenever $p,q$ are.  
Thus $p \tensor q \in \Inv_{xy}$.  
This gives an associative
product, not in general commutative.  

We take this opportunity to observe that a notion of {\em domination} can be defined for invariant types in full generality.  In stable theories, it will agree with the usual notion.
 Write $P \equiv Q$ if the two partial types $P,Q$ imply each other.

\begin{definition}  \label{domination-equivalence}  \index{domination-equivalence}  
If $p,q$ are  invariant  types  over $\U$, call $p,q$ domination-equivalent if for
some     base $A$ for $p,q$  
and some $r \in S_{xy}(A)$ containing  $p(x) | A \union q(y) |A$, we have
$p(x) \union r(xy) \equiv q(y) \union r(xy)$ \end{definition}

   It is easy to 
see that domination-equivalence is an equivalence relation on invariant types, and
is a congruence with respect to $\tensor$.  The set of domination-equivalence classes
of invariant types thus becomes an associative, commutative semigroup.  Denote it 
 by $\overline{\Inv}(\U)$. 

We will use   domination-equivalence almost exclusively when one of these types lies in the stable part of
the theory; in this case the other will be said to be {\em stably dominated}.  

\section{Conditions equivalent to stability.}

\medskip

We give some definitions which yield notions equivalent to stability.

\begin{definition} \label{ranks} \rm
(i) The theory $T$ is {\em $\lambda$-stable} (where $\lambda$ is an infinite cardinal), if for all $C\subset \U$
with $|C|=\lambda$, $|S(C)|=\lambda$. It is {\em stable} if it is $\lambda$-stable for some infinite $\lambda$.
\index{stable!theory}

(ii) The formula $\phi(x,y)$ (possibly with extra parameters) is {\em unstable} if there are $a_i,b_i\in \U$ 
for $i<\omega$ such that for all $i,j<\omega$, $\phi(a_i,b_j)$ holds if and only if $i<j$. We say $\phi(x,y)$ 
is {\em stable} if it is not unstable.\index{stable!formula}

(iii) Let $p(x) \in S(B)$ and $C\subseteq B$. Then $p$ is {\em definable over $C$} \index{type!definable} 
if for every $\L$-formula $\phi(x,y)$ there is an
$\L(C)$-formula $\psi(y)$ (often denoted $(d_px)(\phi(x,y))$) such that, for all $b \in B$, $\phi(x,b)\in p$ if and only if
$\psi(b)$ holds. We say that $p$ is definable {\em almost over $C$} if it is definable over $\acl(C)$ 
{\em imaginaries included},  and that $p$ is {\em definable} if it is definable over $B$.   An isolated type $p$ over 
$B$ (isolated via a formula $\psi(x)$) is an example of
a definable type over $B$ in this sense, since $\phi(x,b)\in p$  if and only if $(\forall x)(\psi(x) \to \phi(x,b))$.  However, this definition does not extend to structures $B'$ containing $B$.  Definable
types are useful inasmuch as they are `extendible', i.e. the same definition scheme
defines a complete type $p|B'$ over any bigger base set $B'$. \index{type!extendible} In this monograph, we will
only consider definable types which are extendible, and indeed, think of the function $B' \mapsto (p|B')$
as being the definable type.   Definable types over models are always uniquely
extendible; much less trivially, in stable theories with elimination of imaginaries, definable
types over algebraically closed sets are uniquely extendible.  The uniqueness is the content
of the {\em finite equivalence relation theorem}, mentioned below.

(iv) Suppose that $\Delta$ is a finite set of formulas $\delta(x,y)$,
and let $S_\Delta(\U)$ denote the  space of complete $\Delta$-types over $\U$; here a $\Delta$-type
\index{type!deltatype@$\Delta$-type} is a maximal
(subject to consistency with $T$)
 set of $\Delta$-formulas, that is, Boolean combinations of formulas $\delta(x,a)$ for $\delta\in \Delta$ and $a \in \U$. 
 Then $S_\Delta(\U)$ is a compact totally disconnected topological space, where a basic open set is the 
collection of $\Delta$-types containing a fixed $\Delta$-formula. Let $\Phi(x)$ be a set of formulas over 
some $A\subset \U$ (so $|\Phi|<|\U|$). Then the {\em $\Delta$-rank} of $\Phi(x)$, denoted
$R_\Delta(\Phi(x))$, \index{rank!deltarank@$\Delta$-rank}
 is the Cantor-Bendixson rank of the subspace
$Y=\{q\in S_\Delta(\U):q \mbox{~is consistent with~} \Phi(x)\}$ of $S_\Delta(\U)$. If this is finite, then
 $\Mult_\Delta(\Phi)$
is the Cantor-Bendixson multiplicity of $\Phi$.
\end{definition}

\begin{theorem} \label{equivstable}
The following are equivalent.

(i) $T$ is stable.

(ii) Every formula $\phi(x,y)$ is stable.

(iii) For every model $M$, every 1-type over $M$ is definable.

(iv) For every model $M$, every type over $M$ is definable.

(v) For every $\Delta$ and $\Phi$ as in Definition~\ref{ranks}, $R_\Delta(\Phi)$ is finite.

(vi)  For every $\Delta$, the space $S_\Delta(M)$ has cardinality at most $|M|+|L|$.

(vii) For some $\kappa$, for any $M$, the space of types $S(M)$ over $M$ has
cardinality at most $|M|^{\kappa}$.
\end{theorem}

(vii)  is often the easiest way to verify stability, while (v) is most useful in proofs.

\begin{remark} \label{pos} \rm Suppose that $\phi(x,y)$ is a stable formula, $C\subset M$ with $M$ $|C|^+$-saturated, and $p\in S(M)$. Then 
 there
are $a_1,\ldots,a_n\in M$ with $a_{i}\models p|C\cup\{a_j:j<i\}$ for each $i=1,\ldots, n$, such that $(d_px)\phi(x,y)$ is equivalent to a {\em positive} Boolean combination of formulas $\phi(a_i,y)$.
See for example Lemma 1.2.2 of \cite{pil} for a proof.
\end{remark}

If $C\subseteq M\models T$ and $p\in S(M)$ is $C$-definable, then the {\em defining schema} for $p$
(the map $\phi(x,y)\mapsto (d_px)\phi(x,y)$) yields canonically an extension
 $p|B$ of $p$ over any $B\supseteq C$, and $p|B$
 is also $C$-definable. In particular, there is a canonical extension $p|\U$ to $\U$.
In the other direction, to restrict the parameter set, if $q$ is {\em any} type over $\U$, then $q|C$ denotes the restriction of $q$ to  $C$, so $q|C\in S(C)$.

Clearly, any $C$-definable type over $\U$ is $\Aut(\U/C)$-invariant, but other (undefinable) 
 examples of invariant types play a major role in Part~II. 
 
\medskip

\section{Independence and Forking.}

\medskip

Stable theories admit a unique theory of independence, that can be approached
in a number of ways.  One such approach, that we will not mention, is
 the Lascar-Poizat notion of {\em fundamental order}.  Another, essential in the
 generalization to simplicity, de-emphasizes uniqueness properties in favour of 
 their consequences in terms of amalgamation of independent triangles.  We 
 will focus on two aspects:  (1) the choice of a  canonical element of a type space;
 especially, if   $C\subseteq B \subset \U$ and $p\in S(C)$, a canonically chosen type  $q\in S(B)$ 
which extends $p$.   (2) 
For the best such $q$, the information content in $q$ should be as small as possible, and the locus of $q$ (i.e. the set
 of realisations in $\U$)
will in some sense be as large as possible.   On the other hand, the other extensions will 
be `small'; a reasonable technical notion of smallness is given by `dividing', below.  

\begin{definition} \rm
Let $C\subseteq B \subset \U$ and $\phi(x,b) \in S(B)$ (so $b\in B$). Then 
$\phi(x,b)$ {\em divides} over $C$  \index{divides}  if there is $k\in\omega$
and a sequence $(b_i:i\in \omega)$ of realisations of $\tp(b/C)$ such that
$\{\phi(x,b_i):i\in \omega\}$ is $k$-inconsistent, that is, any subset of it of size $k$ is inconsistent.
A partial type $\pi(x)$ over $B$ {\em forks} \index{forks} over $C$ if there are $n\in \omega$ and formulas 
$\phi_0(x),\ldots,\phi_n(x)$ such that $\pi(x)$ implies
$\bigvee_{i\leq n} \phi_i(x)$, and each $\phi_i(x)$ divides over $C$.
If $p\in S(B)$ does not fork over $C$ and $a$ realises $p$, we write $a\dnf_C B$, \index{independence!non-forking, $\dnf$}
and say that $a$ is independent from $B$ over $C$.
More generally, if $A,B,C$ are subsets of $\U$, we write $A\dnf_C B$ if for any finite tuple $a$ from $A$,
$a\dnf_C B\cup C$.
\end{definition}

\begin{theorem}
Let $T$ be stable, and $C\subseteq B$. Then

(i) if $p\in S(C)$, then $p$ has a non-forking extension over $B$,

(ii) if $q \in S(B)$, then $q$ does not fork over $C$ if and only if it is definable almost over $C$,

(iii) if $q\in S(B)$, then $q$ does not fork over $C$ if and only if $R_\Delta(q|C)=R_\Delta(q)$ for all finite sets 
$\Delta(x)$ of $\L$-formulas.
\end{theorem}

A {\em finite equivalence relation} \index{equivalence relation!finite} is a definable equivalence relation with 
finitely many classes. By the following theorem, in a stable theory the
non-forking extensions of a type are governed by finite equivalence relations.

\begin{theorem}[Finite Equivalence Relation Theorem]\label{FERT}
Assume $T$ is stable. Let $C\subseteq B\subset \U$, and $p_1,p_2\in S(B)$ be distinct non-forking extensions of $p\in S(C)$. Then there
is a $C$-definable finite equivalence relation $E$
such that $p_1(x) \cup p_2(y) \vdash \neg Exy$.
\end{theorem}

Under our assumption of elimination of imaginaries, the statement becomes simpler:  Let $C=\acl(C) \subseteq B \subset \U$,
and let $p_1,p_2\in S(B)$ be non-forking extensions of $p\in S(C)$.  Then $p_1=p_2$.

A type $p\in S(C)$ is {\em stationary} \index{type!stationary} if it has a unique non-forking extension over $\U$. By the
 Finite Equivalence Relation Theorem, in a stable theory,
if $C=\acl(C)\subset \U$, then every type over $C$ is stationary: indeed, by elimination of imaginaries, the
 equivalence classes of the equivalence relation given in Theorem~\ref{FERT} are coded  
by an element of $\U$, and this element will lie in $\acl(C)$. In particular, any type over a model is stationary.
If $p\in S(C)$ is stationary, and $C\subset B$, we write $p|B$ for the unique non-forking extension of $p$ to $B$.

If $C\subset \U$, we say that $a,b$ have the same {\em strong type} \index{type!strong} over $C$, written, 
$\stp(a/C)=\stp(b/C)$, if $a,b$ lie in the same class of any $C$-definable finite equivalence relation.  Equivalently,
with imaginaries:  $\tp(b/\acl(C))=\tp(a/\acl(C))$. 

For the sake of one result (Corollary~\ref{lascar}), we mention the following variant, which has been important in 
the study of simple theories: $a$ and $b$ have the  same {\em Lascar strong type} over $C$ \index{type!Lascar strong} 
if there is a sequence $M_1,\ldots,M_n$ of models containing $C$, and for each $i=1,\ldots, n$
some $f_i\in \Aut(\U/M_i)$, such that $b=f_n \circ \ldots \circ f_1(a)$. In stable theories, the notion of strong type and Lascar strong type coincide, but there are other (so far, all artificially constructed) theories in which they are known to differ. A major question is whether they agree in all simple theories.

\begin{remark} \label{forking/inv}\rm
Suppose $C\subseteq B$.

(i) If $\tp(a/B)$ has an $\Aut(\U/\acl(C))$-invariant extension over $\U$, then
$a\dnf_C B$. 

(ii) If $T$ is stable and $a\dnf_C B$, then $\tp(a/B)$ has a unique $\Aut(\U/\acl(C)$-invariant extension over $\U$.

Here, (i) follows almost immediately from the definition of forking, and (ii) is a consequence of \ref{FERT}. 
Under our ambient assumption of elimination of imaginaries, we can work with $\acl(C)$, not $\acl^{\eq}(C)$.
\end{remark}

Here is a standard list of basic properties of the non-forking relation $A\dnf_C B$ on sets in a  stable theory. 

(i) (Existence) If $p\in S(C)$ and $M$ is a model containing $C$, then $p$ has a non-forking extension $q\in S(M)$; if
$C=\acl(C)$ then $q$ is unique.

(ii) (Finite character) If $p\in S(B)$ and $B\supseteq C$, then $p$ does not fork over $C$ if and only if for every formula $\phi(x)\in p$, $\phi$ does not fork over $C$. In particular, $p$ does not fork over $C$ if and only if for all finite $B_0\subseteq B$, $p|CB_0$ does not fork over $C$.

(iii) (Transitivity) If $p\in S(B)$ and $C\subseteq D\subseteq B$ then $p$ does not fork over $C$ if and only if $p$ does not fork over $D$ and $p|D$ does not fork over $C$.

(iv) (Symmetry) If $C\subseteq B$, then $\tp(a/B)$ does not fork over $C$ if and only if for all finite tuples $b$ from $B$,
$\tp(b/Ca)$ does not fork over $C$.

(v) (Invariance) If $p\in S(B)$ does not fork over $C\subset B$, and $f\in \Aut(\U)$, then $f(p)\in S(f(B))$ does not fork over $f(C)$.

(vi) If $p\in S(C)$ and $B\supseteq C$ then $p$ has at most $2^{|T|}$ non-forking extensions over $B$. 

(vii) If $p\in S(B)$ then there is $C\subseteq B$ with $|C|\leq |T|$ such that $p$ does not fork over $C$.

(viii) If $p\in S(B)$ has finitely many realisations and does not fork over $C\subset B$, then $p|C$ has finitely many realisations.

\begin{remark} \label{localsymmetry} \rm
Symmetry of non-forking can also be viewed locally. Let $\delta(x,y)$ be a stable formula, and let $p(x)\in S_\delta(\U)$,
$q(y)\in S_\delta(\U)$. These are both definable -- this comes from the proof of Theorem~\ref{equivstable}.
Let $\psi(y)$ be the $\delta$-definition $d_px\delta(x,y)$ of $p$, and $\chi(x)$ be the $\delta$-definition
$d_qy\delta(x,y)$ of $q$. Then $\chi(x)\in p(x)$ if and only if $\psi(y) \in q(y)$ -- see \cite[Lemma I.2.8]{pil}.
\end{remark} 

In Part~I we make occasional use of the notion of {\em weight}. In a stable theory, $T$, the {\em preweight} \index{preweight}
of a type $p(x)=\tp(a/C)$
is the supremum of the set of cardinals $\kappa$ for which there is an $C$-independent set $\{b_i:i<\kappa\}$  such that
$a\df_C b_i$ for all $i$. The {\em weight} $\wt(p)$ \index{weight} is the supremum of
$\{\prwt(q): q\mbox{~is a non-forking extension of~} p\}$. 

\begin{lemma}  \label{morley-div} In a stable theory,
any type has weight bounded by the cardinality of the language; in a superstable theory,
any type in finitely many variables has finite weight. \end{lemma}

{\em Proof.} See Pillay \cite{pil} ch. 1,
 Section 4.4 or Section 5.6.3 of Buechler \cite{buechler}.  \qed

\medskip

 Suppose that the theory $T$ is stable. Let $p\in S_x(C)$ be stationary, and $\phi(x,y)$ be a formula. By definability of 
$p$, there is a $\phi$-definition of $p$, namely, a formula $\psi(y)$ of form $(d_p x)(\phi(x,y))$. 
The set defined by $\psi(y)$ has a code $c$, unique up to definable closure, and $c\in \U$ by elimination of imaginaries. 
The {\em canonical base} $\Cb(p)$ of $p$ \index{canonical base} is the definable closure of the set of all codes 
$c$ for sets defined by formulas
$(d_p x)(\phi(x,y))$ ($x$ fixed, $y$ and $\phi$ varying). If $p=\tp(a/C)$, we write $\Cb(a/C)$ for $\Cb(p)$.  This is the 
smallest set over which $p$ is defined.

\begin{lemma}\label{canbasebasic}
Assume $T$ is stable, let $C\subseteq B$, and let $p\in S(B)$ be stationary. Then

(i) $\Cb(p) \subseteq \dcl(B)$;

(ii) $\Cb(p)\subseteq \acl(C)$ if and only if $p$ does not fork over $C$;

(iii) $\Cb(p)\subseteq \dcl(C)$ if and only if $p$ does not fork over $C$ and $p|C$ is stationary. 
\end{lemma}

\medskip

\section{Totally transcendental theories and Morley rank.}

\medskip

The theory $T$ is {\em totally transcendental} \index{totally transcendental} if the Cantor-Bendixson rank of each type 
space $S_x(\U)$  is  ordinal-valued. If $\L$ is countable, this is equivalent to $T$ being $\omega$-stable. 
Let $\Phi(x)$ be a set of formulas, and $X=\{p\in S_x(\U): \Phi\subseteq p\}$ (a closed subspace of $S_x(\U)$). Then
the Morley rank $\RM(\Phi)$ of $\Phi$ is just $\CB(X)$. If $\RM(\Phi)<\infty$ (as holds if $T$ is totally transcendental), then
the Morley degree $\DM(\Phi)$ of $\Phi$ is CB-Mult$(X)$. We often write $\RM(a/C)$ for $\RM(\tp(a/C))$, and similarly for
 Morley degree. In a totally transcendental theory, forking is determined by Morley rank: if $C\subseteq B$, then $a\dnf_C B$ if and only if $\RM(a/C)=\RM(a/B)$.
	
	We give a slightly more explicit definition of Morley rank and degree for a definable set $D$. It is defined inductively, as follows. 

(i) $\RM(D)\geq 0$ if and only if $D\neq \emptyset$.

(ii) For a limit ordinal $\alpha$, $\RM(D)\geq \alpha$ if and only if, for all $\lambda<\alpha$, $\RM(D)\geq \lambda$.

(iii) For any ordinal $\alpha$, $\RM(D) \geq \alpha+1$ if and only if there are 
infinitely many pairwise disjoint definable subsets
$D_i$ ($i<\omega)$ of $D$ each of Morley rank at least $\alpha$.

If $\RM(D)\geq \alpha$ but $\RM(D)\not\geq \alpha+1$ we write $\RM(D)=\alpha$, and if $\RM(D)\geq \alpha$ for all ordinals $\alpha$ we put $\RM(D)=\infty$. Finally, if $\RM(D)=\alpha$, then $\DM(D)$ is the largest $n$ such that there are $n$ disjoint definable subsets of $D$ each of Morley rank $\alpha$.

If the definable set $X$ in $\U$ has Morley rank and degree both equal to 1, then 
$X$ is said to be {\em strongly minimal} \index{strongly minimal}. 
 A 1-sorted theory $T$ is strongly minimal if the set (in one variable) defined by $x=x$ is
strongly minimal.  Examples of  strongly minimal theories are the theories of pure sets, vector spaces,
and algebraically closed fields. In an algebraically closed valued field, the residue
 field is an algebraically closed field with no extra structure, so is strongly minimal. Any  strongly minimal theory is 
$\aleph_1$-categorical, and conversely, according to Baldwin and Lachlan \cite{baldlach}, any model of 
an $\aleph_1$-categorical theory is prime over a certain strongly minimal set  and the parameters used to define it.

\section{Prime models.}

\medskip

Prime models play a central role in Shelah's classification theory. We summarise the basic facts here. In Chapter~\ref{primemodels},
a slight generalisation of the notion is developed, partly in the ACVF context.

\begin{definition}
Let
$T$ be a complete theory, and $C\subseteq M\models T$.

(i) $M$ is  {\em prime over $C$} \index{prime} if, for every $N\models T$ with $C\subseteq N$, there is an elementary embedding
$f:M\rightarrow N$ which is the identity on $C$. 

(ii) $M$ is {\em atomic over $C$} \index{atomic} if for any finite tuple $a$ of elements of $M$, $\tp(a/C)$ is
isolated.

(iii) $M$ is {\em minimal over $C$} \index{minimal} if there is no proper elementary substructure of $M$ which contains $C$. 
\end{definition}

An easy consequence of the Omitting Types Theorem is that if $\L$ and $S(\emptyset)$ are
 countable then $T$ has an atomic model over $\emptyset$. More generally, we have the following 
(see e.g. \cite[4.2.10]{mark}).
 
 \begin{theorem}
 If $T$ is a complete theory with infinite models over a countable language, then the following are equivalent.
 
 (i) $T$ has a prime model.
 
 (ii) $T$ has an atomic model.
 
 (iii) For all $n$, the set of isolated $n$-types is dense in $S_n(\emptyset)$.
 \end{theorem}

In many algebraic theories, prime models exist and have an algebraic characterisation. For  the theories of
algebraically closed fields and real closed fields,  the prime model over $C$ is exactly
the  (field-theoretic) algebraic closure of $C$, respectively the real closure.   For differentially closed fields it is the differential closure; uniqueness and 
non-minimality 
of the differential closure was one of the early applications of stability to this area.
For  $p$-adically closed fields a similar characterization holds when $C$ is a field, but
not when imaginary elements are allowed.  

\begin{theorem}
Let $T$ be $\omega$-stable or o-minimal, and $C\subset \U$. 
Then there is a (unique up to isomorphism over $C$) prime model $M$ of $T$ over $C$, and $M$ is also atomic over $C$.
\end{theorem}

In the above theorem, the $\omega$-stable case is due to Shelah, and the o-minimal case due to Pillay and Steinhorn.  

\medskip

\section{Indiscernibles, Morley sequences.}

\medskip

If $(I,<)$ is a totally ordered set, then the sequence $(a_i:i\in I)$ of distinct elements of $\U$ is {\em (order)-indiscernible
over $C$} \index{indiscernible} if, for any $n\in \omega$, and $i_1,\ldots,i_n, j_1,\ldots,j_n \in I$ with
$i_1<\ldots<i_n$ and $j_1<\ldots<j_n$, $\tp(a_{i_1}\ldots a_{i_n}/A)=\tp(a_{j_1}\ldots a_{j_n}/C)$. 
We say that $\{a_i:i\in I\}$
is an {\em indiscernible set} over $C$ if for any distinct $i_1,\ldots, i_n \in I$ and distinct $j_1,\ldots,j_n\in I$,
$\tp(a_{i_1}\ldots a_{i_n}/C)=\tp(a_{j_1}\ldots a_{j_n}/C)$ (so the order doesn't matter). 

A standard application of Ramsey's theorem, compactness, and the saturation of $\U$ ensures that for any ordered set $(I,<)$
which is `small' relative to $\U$ and any small parameter set $C$, there is in $\U$ an indiscernible sequence over 
$C$ of order type $I$. If the theory $T$ is stable, any infinite indiscernible sequence is an indiscernible set.

Assume now that $T$ is stable, and let $p(x)\in S(C)$ be stationary. A {\em Morley sequence} \index{Morley sequence} 
for $p$ of length $\lambda$
is a sequence $(a_i:i<\lambda)$ such that for each $i<\lambda$, $a_i$ realises $p|C\cup\{a_j:j<i\}$. Such a sequence  is
independent over $C$ and is an indiscernible set over $C$. Furthermore, if $C=\acl(C)$  any infinite indiscernible set over
 $C$ which is independent over $C$
is a Morley sequence for some strong type over $C$.

Below, and elsewhere, if $J\subseteq I$, we  write $a_J$ for the subsequence $(a_j:j\in J)$.

\begin{lemma} \label{Morleybuechler}
Let $\{a_i:i\in I\}$ be an infinite indiscernible set over $C$, and suppose that for any finite disjoint $J,J'\subset I$, 
$\acl(Ca_J) \cap \acl(Ca_{J'})=\acl(C)$. Then $\{a_i:i\in I\}$ is a Morley sequence 
over $C$. 
\end{lemma}

{\em Proof.} This follows easily from Lemma 5.1.17 of Buechler \cite{buechler}, which states that the indiscernible
 sequence $I$ is a  Morley sequence over $C\cup J$ for any infinite $J\subset I$. To apply this, let $i_1,\ldots,i_n\in I$ be distinct,
and let $J_1,J_2$ be infinite disjoint subsets of $I\setminus \{i_1,\ldots,i_n\}$. Then $a_{i_1}\dnf_{CJ_j}a_{i_2}\ldots a_{i_n}$
for $j=1,2$,
so $\Cb(a_1/Ca_{i_2}\ldots a_{i_n}) \subseteq \acl(CJ_1 \cap CJ_2)=\acl(C)$.
\qed

\medskip

\section{Stably embedded sets.}

\medskip

A $C$-definable set $D$ in $\U$ is {\em stably embedded} \index{stably embedded}  if, for any definable set $E$ and 
$r>0$, $E \cap D^r$ is definable over $C\cup D$. If instead we worked in a small model $M$, and $C,D$ were from 
$M^{\eq}$, we would say that $D$ is stably embedded if for any definable $E$ in $M^{\eq}$ and any $r$,
$E\cap D^r$ is definable over $C\cup D$ {\em uniformly} in the parameters defining $E$; that is, for any formula $\phi(x,y)$ there is a formula $\psi(x,z)$ such that for all $a$ there is a sequence $d$ from $D$ such that
$$\{x\in D^r:\models\phi(x,a)\}=\{x\in D^r:\models \psi(x,d)\}.$$ 
For more on stably embedded sets, see the Appendix of \cite{ch}. We mention in particular that if 
$D$ is a $C$-definable stably embedded set and $\U$ is saturated, then any permutation
of $D$ which is elementary over $C$ extends to an element of $\Aut(\U/C)$.
Basic examples of stably embedded sets include the field of constants in a differentially closed field, and
 the fixed field in a model of ACFA.
By definability of types, in a stable theory all definable sets are stably embedded.

If  $D$ is $C$-definable, then we say that $D$ is {\em stable} \index{stable!set} if the structure with 
domain $D$, when equipped with all the $C$-definable relations, is stable. We emphasise the distinction between saying 
that a definable set $D$ is stable, and that a formula $\phi(x,y)$ is a stable formula
(as in Definition~\ref{ranks}).

In this monograph, working over a set $C$ of parameters,
 $C$-definable sets which are both stable and stably embedded play a crucial role. Such sets are in some sources just called 
{\em stable}. \index{stable, stably embedded}

\begin{lemma}\label{critstabemb}
Let $D$ be a $C$-definable set. Then the following are equivalent.

(i) $D$ is stable and stably embedded.

(ii)
Any formula $\phi(x_1\ldots x_n,y)$ which implies $D(x_1) \wedge \ldots \wedge D(x_n)$ is stable 
(viewed as a formula $\phi(x,y)$, where $x=(x_1,\ldots,x_n)$). 

(iii) If $\lambda>|T|+|C|$ with $\lambda=\lambda^{\aleph_0}$ and  $B\supseteq C$ with $|B|=\lambda$,  then there are at most
$\lambda$ 1-types over $B$ realised in $D$.

\end{lemma}

In particular, by (ii), if $D$  is $C$-definable, and $C\subseteq B$,
then $D$ is stable and stably embedded over $C$ if and only if it is stable and stably embedded over $B$, so there is no need 
to mention the parameter set.  
Hence any sort in $\St_C$ is also in $\St_{C'}$
for some finite $C' \subseteq C$; thus:

\begin{corollary}  \label{kappa-st} Let $S$ be a sort of $\St_C$.   The  language of $S$ is obtained from
a language of size $\leq |T|$ by the addition of constants.  \end{corollary}

If $D,E$ are definable sets in $\U$, then $E$ is said to be {\em $D$-internal} \index{internal} if there is a finite set $A$ 
of parameters such that $E\subset \dcl(D\cup A)$. It is immediate that any definable set which is internal to a stable and 
stably embedded set is stable and stably embedded. In ACVF, by Proposition~\ref{263}, any stable and stably embedded 
set is internal to the residue field.

%
%

\chapter{Definition and basic properties of $\St_C$} \label{stableindependence}

Given any theory $T$, we consider the stable, stably embedded definable
sets, including imaginary sorts.  The induced structure $\St_{\emptyset}$ on these sets can be viewed as 
the maximal stable theory   interpretable
without parameters in $T$.   We define a class of types of $T$, called {\em stably dominated},
that are not necessarily stable, but have a stable part; and such that any interaction
between a realization of the type and any other set $B$ is preceded by an interaction of
the stable part with $B$.  See below for a precise definition.  

We show that the theory of independence on the stable part lifts to the stably dominated
types.  The same is true of a slightly less well-known theory of independence based on
dcl, rather than acl, and we describe this refinement too.   We develop the basic properties
of stably dominated types.  

At this level of generality, elimination of quantifiers and of imaginaries can be, and will be assumed.

Let $\U$ be a sufficiently saturated model of a complete theory $T$ 
with elimination of imaginaries. Throughout the monograph, $A$, $B$, $C$, $D$
will be arbitrary subsets of the universe. Symbols $a,b,c$ 
denote possibly infinite tuples, sometimes regarded as sets, sometimes as 
tuples, and indeed sometimes $A,B$ are implicitly regarded as tuples. We occasionally write $a\subset X$ to mean that the set enumerated by $a$ is a subset of $X$.
We write $A\equiv_C B$ to mean that, under some enumeration of $A$, $B$ which is fixed but usually not specified,
$\tp(A/C)=\tp(B/C)$. By the saturation assumption, this is equivalent to $A$ and $B$, in the given enumerations, 
being in the same orbit of $\Aut(\U/C)$
(the group of automophisms of $\U$ which fix $C$ pointwise); sometimes we just say that $A$ and $B$ 
are {\em conjugate}\index{conjugate} over $C$.

\begin{definition} \rm Let $C$ be any small set  of parameters. Write $\St_C$ \index{St} for the multi-sorted 
structure $\langle D_i, R_{j}\rangle_{i\in I, j\in J}$ 
whose sorts $D_i$ are the $C$-definable, stable, stably embedded subsets
of $\U$. For each finite set of sorts $D_i$, all the $C$-definable relations
on their union are included as $\emptyset$-definable relations $R_j$.  For any 
$A\subset \U$, write $\St_C(A):= \St_C \cap \dcl(CA)$. 
We often write $A^{\rm st}$ for $\St_C(A)$ when the base $C$ is 
unambiguous.
\end{definition}

We begin with some basic observations about $\St_C$.

\begin{lemma}
The stucture $\St_C$ is stable.
\end{lemma}

{\em Proof.} It suffices to observe that if $D_1$ and $D_2$ are $C$-definable stable and stably embedded sets, then 
so is $D_1 \cup D_2$. By Lemma~\ref{critstabemb}, it is sufficient to show that if $C'\supseteq C$ and $|C'|=\lambda$, where $\lambda=\lambda^{\aleph_0}>|T|+|C|$, then there are at most $\lambda$ 1-types over $C'$ realised in $D_1 \cup D_2$. This is immediate, as by the same lemma there are at most $\lambda$ 1-types realised in each of $D_1$ and $D_2$.\qed

\medskip

We shall use the symbol $\dnf$ for 
non-forking (for subsets of $\St_C$) in the usual sense of stability theory;
the context should indicate that this takes place in $\St_C$. In particular,
 $A\dnf_C B$
means that $A,B$ are independent (over $C$) in $\St_C$, unless
 we explicitly state that the independence is in another structure.
 Notice that, if
$A, B \subset \St_C$, although we do not expect $A = \St_C(A)$, still 
$A\dnf_C B$ if and only if $\St_C(A)\dnf_C \St_C(B)$.

We begin with some remarks on how the structure $\St_C$ varies when the parameter set $C$ increases. As noted after Lemma~\ref{critstabemb}, the property of a definable set $D$ being stable and stably embedded does not depend on the choice of defining parameters. Thus, if $D$ is a sort of $\St_C$ and $B\supseteq C$, then $D$ is also a sort of $\St_B$. 
In general, if $C\subseteq B$ then $\St_B$ has more sorts than $\St_C$.  Notice that $\St_{\acl(C)}$ has essentially the same domain as $\St_C$: if $D_1$ is a sort  of $\St_{\acl(C)}$ whose conjugates over $C$ are $D_1,\ldots,D_r$, then
$D_1 \cup \ldots \cup D_r$ is a sort of $\St_C$. Each element of $\St_{\acl(C)}$ has a code over $C$ in $\St_C$.

\begin{lemma} \label{3.21}
Let $D_a$ be a stable, stably embedded $a$-definable set. Then $D_a$ can be defined by a formula
$\phi(x,a)$ with $\phi(x,y)$ stable.
\end{lemma}

{\em Proof.} Let $p=\tp(a/\emptyset)$. Suppose that $D_a$ is defined by the formula $\psi(x,a)$. 

{\em Claim.} There do not exist elements $a_i,b_i$ (for $i\in \omega$) such that $p(a_i)$ holds for each $i$
and $b_i\in D_{a_j}$ if and only if $i>j$.

{\em Proof of Claim.} Suppose such $a_i,b_i$ exist. Then $b_i \in D_{a_0}$ for all $i>0$. Thus, 
$\psi(b_i,a_0)\wedge \psi(b_i,a_j)$ holds if and only if $i>j$ (for $i,j>0$). It follows that the formula
$\psi(x,a_0) \wedge\psi(x,y)$ is unstable. This contradicts the fact that $D_{a_0}$ is stable and stably embedded.

Given the claim, it follows by compactness that there is a formula $\rho(y)\in p$ such that there do not exist
$a_i,b_i$ (for $i\in \omega$) such that $\rho(a_j)$ holds for all $j$ and $\psi(b_i,a_j)$ if and only if $i>j$.
Thus, the formula $\phi(x,y)=\psi(x,y)\wedge \rho(y)$ is stable, and $\phi(x,a)$ defines the same set as $\psi(x,a)$. \qed

\begin{lemma}\label{3.22}
Let $D_b$ be  a stable stably embedded $Cb$-definable set. Let $p$ be a $C$-definable type, and assume that
$a\models p|Cb$ and $a\in D_b$. Then there is a $C$-definable stable stably embedded set $D'$ with $a\in D'$.
\end{lemma}

{\em Proof.} We may suppose that $C=\acl(C)$. For if such a set $D'$ may be found over $\acl(C)$, then the union $D''$ of the $C$-conjugates of $D'$ works over $C$.

By Lemma~\ref{3.21}, there is a stable formula $\delta(x,y)$ such that $D_b$ is defined by $\delta(x,b)$. Put $q:=\tp(b/C)$.
 Let $\chi(x)$ be the formula $d_q y \delta(x,y)$, a formula over $C$. By forking symmetry for stable formulas 
(Remark~\ref{localsymmetry}), $\chi(a)$ holds. 
Also, by Remark~\ref{pos}, $\chi(x)$ is a (finite) positive Boolean combination of formulas $\delta(x,b_i)$ 
where $b_i\equiv_C b$ for each $i$. Each $\delta(x,b_i)$ defines a stable and stably embedded set, and
 hence so does $\chi(x)$. \qed

\bigskip

We shall say that $\St_C$ satisfies the condition $(*_C)$ if:
whenever $b\in \St_C$, and $D_b$ is a stable and stably embedded $Cb$-definable set, there is a formula $\psi(y)\in \tp(b/C)$ such that whenever $\psi(b')$ holds, $D_{b'}$ is stable and stably embedded.

\begin{lemma}\label{*C0}
Suppose that $(*_C)$ holds, and $C\subseteq C'\subset \St_C$. Then every element of $\St_{C'}$ is interdefinable over $C$ with an element of $\St_C$.
\end{lemma}

{\em Proof.} Let $a\in \St_{C'}$. Then there is $b\in \St_C$ and a stable and stably embedded $Cb$-definable set 
$D_b$ such that $a\in D_b$. Let $\psi(y)$ be the formula provided by $(*_C)$. As $b\in \St_C$, we may suppose that $\psi(y)$ defines a stable stably embedded set. Put 
$E:=\{x: \exists y(\psi(y)\wedge x\in D_y)\}$, a $C$-definable set. Then $a\in E$, and by the type-counting criterion in 
Lemma~\ref{critstabemb}, $E$ is stable and stably embedded. \qed

\begin{remark}\label{*C1}\rm The condition $(*_C)$ holds for each $C$ if $\St_B$ is $\St_{\emptyset}$-internal for each $B$. 
For suppose $E$ is a sort of $\St_C$, that $b\in E$, and that $D_b$ is a stable stably embedded $Cb$-definable set. There is a stable stably embedded  $\emptyset$-definable set $D^*$ and some $e$-definable surjection
$f_e: D^*\rightarrow D_b$. Now let $\psi(y)$ be the formula
$$y\in E \wedge\exists u (f_u\mbox{~is a surjection~} D^*\rightarrow D_y).$$
If $\psi(y)$ then $D_y$ is $D^*$-internal, so is stable and stably embedded.

As observed later in Remark~\ref{*C} below, $(*_C)$ therefore holds in ACVF.
\end{remark}

In  ACVF, a family of $\k$-vector spaces  $\VS_{k,C}$ is defined over any base set $C$;
see Chapter~\ref{acvfbackground} below.  (This modifies the notation of
\cite[Section 2.6]{hhm}, where we wrote $\Int_{k,C}$ for $\VS_{k,C}$).  
  $\VS_{k,C}$ is a subset of  $\St_C$; while formally it is a proper subset, 
any $C$-definable set in   $\St_C$ is in $C$-definable bijection with
an element of   $\VS_{k,C}$, so they can be viewed as the same.  

We shall frequently in the text write statements like $\tp(A/C)\vdash \tp(A/CB)$. This means that for any $A'$, if $A'\equiv_C A$ then $A'\equiv_{CB} A$. Observe that $\tp(A/C)\vdash \tp(A/CB)$ is equivalent to
$\tp(B/C)\vdash \tp(B/CA)$. Indeed, suppose $\tp(A/C)\vdash \tp(A/CB)$, and let $B'\equiv_C B$. There is $g\in \Aut(\U/C)$ with $g(B')=B$. Then $g(A)\equiv_C A$, so $g(A)\equiv_{CB}A$. Thus there is $h\in \Aut(\U/CB)$ with $hg(A)=A$. Then $hg(B')=B$, so $B'\equiv_{CA} B$. This symmetry is used often without comment.

\begin{remark}\label{symmetry} \rm 
The fact that $\St_C$ is  stably embedded implies immediately that for any 
sets $A,B$, $\tp(B/C\Bst)\vdash \tp(B/C\Bst\Ast)$ (equivalently, that
$\tp(\Ast / C\Bst) \vdash \tp(\Ast / CB)$).  To see this, let $B'\equiv_{CB^{\st}} B$. We must show
$B'\equiv_{CB^{\st}A^{\st}}B$. Let 
$\varphi(x,a) \in \tp(B/C\Ast\Bst)$ 
with $a\in \Ast$ and $\phi(x,y)$ a formula over $CB^{\st}$, and consider $\{y:\varphi(B,y)\mbox{~holds}\}$. This
is a subset of $\St_C$, defined with parameters from $C \cup B$. Since $\St_C$ is
stably embedded, it must be definable by some formula $\psi(x,b')$, say, with parameters $b'$ from $\St_C$. The
set of such parameters $b'$ (an element of $\Ueq=\U$) is definable from $B$, hence is in 
$\dcl(B)\cap\St_C=\Bst=B'^{\st}$. It follows
that $\{y:\phi(B,y)\mbox{~holds}\}$ is $C\Bst$-definable, so equals
$\{y:\phi(B',y)\mbox{~holds}\}$. Thus, $\phi(B',a)$ holds, so $\phi(x,a)\in \tp(B'/C\Bst\Ast)$, as required.

By the last paragraph we have in particular the following.  Suppose
$g$ is an automorphism of $\U$ fixing $C\Bst$. Then there is an automorphism
$h$ fixing $C\Ast\Bst$ such that $h(B)=g(B)$.
\end{remark}

\begin{lemma} \label{symmetric}
(i) 
$\tp(B/C\Ast)\vdash \tp(B/CA)$ if and only if $\tp(A/C\Bst)\vdash \tp(A/CB)$.

(ii)  $\tp(B/C\Ast)\vdash \tp(B/CA)$
if and only if $\tp(\Bst/C\Ast)\cup \tp(B/C) \proves \tp(B/CA)$.
\end{lemma}

{\em Proof.} (i) Assume  
$\tp(B/C\Ast)\vdash \tp(B/CA)$. Suppose $g$ is an automorphism fixing $C\Bst$.
By Remark~\ref{symmetry}, there is an automorphism $h$ fixing $C\Ast\Bst$ such that 
$h|_B=g|_B$. As $h^{-1}g$ fixes $B$, $g(A) \equiv_{CB} h(A)$. By our
assumption and as $h$ fixes $C\Ast$, $B\equiv_{CA} h^{-1}(B)$, or equivalently, $A\equiv_{CB} h(A)$. Hence
$A \equiv_{CB} g(A)$, as required. The other direction is by symmetry.

(ii) The right-to-left direction is immediate as $\tp(B/C\Ast)\vdash \tp(B/C)$. For the left-to-right direction, suppose
$g$ is an automorphism fixing $C$ and such that $\Ast\Bst \equiv_{C} \Ast g(\Bst)$. 
Since also  $\Ast\Bst \equiv_{C} g(\Ast) g(\Bst)$, there is an automorphism
$h$ fixing $g(\Bst)$ such that $h(g(\Ast))= \Ast$. By Remark~\ref{symmetry}, 
there is automorphism $k$ fixing $\Ast g(\Bst)$ such that $kh(g(B))=g(B)$.
Now $khg$ fixes $\Ast$ and maps $B$ to $g(B)$. Thus $B\equiv_{C\Ast} g(B)$, so by the assumption, it
follows that $B\equiv_{CA} g(B)$, as required. \qed

\begin{definition}\label{stableind} \index{stable domination} \index{independence!domination$\domind$} \rm
Define $A\domind_C B$ to hold if $A^{\st}\dnf_C B^{\st}$ and 
$\tp(B/C\Ast)\proves \tp(B/CA)$.
We say that $\tp(A/C)$ is {\em stably dominated} if, whenever 
$B\subset \U$  and $\Ast\dnf_C \Bst$, we have  $A \domind_C B$.
\end{definition}

\begin{remark}\label{downtofinite}\rm
It follows from Proposition~\ref{bits}(iii) below that if $a$ is an infinite tuple, then $\tp(a/C)$ is stably dominated if and only if $\tp(a'/C)$ is stably dominated 
for every finite subtuple $a'$ of $a$.
\end{remark}

Symmetry of $\domind$ follows immediately from Lemma~\ref{symmetric}(i) and symmetry of stable non-forking 
applied in $\St_C$.

\begin{lemma}\label{symmetric2}
 $A\domind_C B$ if and only if $B\domind_C A$.
\end{lemma}

The definition of stable domination is analogous to that of compact domination introduced in \cite{hpp}, in the following sense:
$\tp(a/C)$ is stably dominated precisely if, for any definable set $D$, if $\dst:=\St_C(\lceil D\rceil)$, then
either
$\{x\models p: \xst \dnf_C \dst\}\subseteq D$, or 
$\{x\models p: \xst \dnf_C \dst\}\cap D = \emptyset$.
 For compact domination, there is 
a similar definition, except that $\ast$ is a finite tuple, $\St_C$ is a compact topological space,  `small' 
means `of measure zero', and the uniformity condition on fibres is slightly weaker.

In many situations only a finite part of $\ast$ is required for domination: that is, there is a $C$-definable map $f$
to $\St_C$ such that for any $B$, if $f(a)\dnf_C \Bst$ then $\tp(B/Cf(a))\vdash \tp(B/Ca)$. It will
emerge that this is the situation 
in ACVF. In particular, the generic type of the valuation ring is stably dominated by that of the residue field: essentially, the 
type of a generic element of the valuation ring is dominated by its residue.

\medskip
The following lemma will yield that stably dominated types have 
definable (so also invariant) extensions.

\begin{lemma} \label{beth}
Let $M$ be a sufficiently saturated structure, and $C$ a small subset of 
$M$. Let $p(x,y)$ be a type over $M$ with $x,y$ possibly infinite tuples, 
and let $q(x)$ be the restriction of $p$ to the $x$-variables. Suppose that 
$p(x,y)$ is the unique extension over $M$ of $p(x,y)|_C \cup q(x)$, and that 
$q$ is $C$-definable. Then $p$ is $C$-definable.
\end{lemma}

{\em Proof.}
Suppose $\phi(x,y,b)\in p$. Then there are formulas $\chi(x,b,c)\in q$ and $\psi(x,y,c)\in p|C$ such that
$\chi(x,b,c) \wedge \psi(x,y,c)\vdash \phi(x,y,b)$. Then the formula $\rho(x,b,c)$, namely 
$\chi(x,b,c)\wedge \forall y(\psi(x,y,c)\rightarrow \phi(x,y,b))$, is in $q$. Thus, as $q$ is $C$-definable, 
for the  formula $d_q x\rho(x, y)$ over $C$ we have $d_q x\rho(x,b)$.  Thus,
$\phi(x,y,b)\in p$ if and only if for one of at most $|L(C)|$-many formulas $d_q x \rho$ over $C$, we have
$d_q x\rho(x,b)$. Likewise, $\neg \phi(x,y,b)\in p$ is equivalent to a disjunction of at most $|L(C)|$-many formulas 
over $C$ about $b$. Definability of $p$ now follows by compactness. \qed

\medskip

\begin{proposition}\label{definable}
Let $M$ be a model, and $C\subset M$.

(i) Suppose that $\tp(A/C)$ is stably dominated. Then $\tp(A/C)$ is extendable to an
$\acl(C)$-definable (so $\Aut(M/\acl(C))$-invariant) type over $M$, such that if $A'$ 
realises this type
then $A'\domind_C M$.

(ii) Suppose  $\tp(A/C)$ is stably dominated and $M$ is saturated, with $|A|, |C|<|M|$. Then
$\tp(A/\acl(C))$ has a unique  $\Aut(M/\acl(C))$-invariant extension. \index{invariant extension}
\end{proposition}

{\em Proof.} 
(i) We may suppose that $\Ast \dnf_C \Mst$. Then, by the hypothesis, 
$\tp(\Ast A/M)$ is the unique extension of $\tp(\Ast A/C)\cup \tp(\Ast/M)$ 
to $M$.
Since $\tp(\Ast/M^{\st})$ is definable over $\acl(C)$ 
in the stable structure $\St_C$, it
follows from Lemma~\ref{beth} that $\tp(A/M)$ is definable over $\acl(C)$. In particular,
$\tp(A/M)$ is $\Aut(M/\acl(C))$-invariant.

(ii) 
Now suppose 
$A''\equiv_{\acl(C)} A'$ and $\tp(A'/M)$, $\tp(A''/M)$ are both 
 $\Aut(M/\acl(C))$-invariant extensions of $\tp(A/C)$. Then 
$\tp((A')^{\rm st}/\Mst)$ and $\tp((A'')^{\rm st}/\Mst)$
are both invariant extensions of $\tp(\Ast/\acl(C))$ in $\St_C$: indeed, any automorphism
in $\Aut(\St_C/\acl(C))$ extends to an automorphism of $M$ over $\acl(C)$ (saturation of $M$ and stable embeddedness), so fixes $\tp(A'/M)$ and $\tp(A''/M)$, and hence fixes
$\tp((A')^{\rm st}/\Mst)$ and $\tp((A'')^{\rm st}/\Mst)$.
Hence $\tp((A')^{\rm st}/\Mst)$ and 
$\tp((A'')^{\rm st}/\Mst)$ are equal, 
as  invariant extensions of a type over an algebraically closed base  are 
unique in a saturated model of a stable theory (see,
for example Remark~\ref{forking/inv}). In particular, $(A')^{\rm st}\dnf_C \Mst$. By stable domination,
$\tp(A'/C\Mst)\vdash \tp(A'/CM)$. Thus, applying 
Lemma~\ref{symmetric}(ii) with $A',M$ replacing $B,A$,  $\tp(A'/M)=\tp(A''/M)$.
\qed

\begin{remark} \label{C''}\rm
(i) In the proof of (i), $\tp(A\Ast/\Mst)$ is definable over $C'':=\acl(C)\cap \dcl(CA)$. It follows
 that $\tp(A/C)$ extends to a $C''$-definable type over $M$. Likewise, in (ii), $\tp(A/C'')$ has a 
unique $\Aut(\U/C'')$-invariant extension.

(ii) By the proof of Proposition~\ref{definable}, if $C\subseteq B$ and $\tp(C/B)$ has an $\Aut(\U/\acl(C))$-invariant 
extension over $\U$,
then $\St_C(A) \domind_C \St_C(B)$.
\end{remark}

 We will say that  $\Aut(\U/C)$-invariant type $p$ is {\em stably
dominated} if $p|C$ is stably dominated.   \index{Stably dominated! for invariant types}

In ACVF there are `many' stably dominated types, enough in some sense to control
the behaviour of all types, if $\G$-indexed families are taken into account.  We 
  expect similar facts to hold in other theories of intended application, such as the  rigid analytic expansions of 
Lipshitz \cite{l}. Further examples are given in Chapter 16. One cannot however expect {\em all}
types to be stably dominated, even over a model, unless the entire theory is   stable.

\begin{corollary}\label{stable}
Suppose that for each singleton $a$ and every model   $M$,
 $\tp(a/M)$ is stably dominated. Then the theory $T$ is stable.
\end{corollary}

{\em Proof.} This follows immediately from Proposition~\ref{definable} (i), since one of the standard characterisations
 of stability (Theorem~\ref{equivstable}(iv)) is that every 1-type  over a model is definable. \qed

\bigskip

Next,  we consider the relationship between indiscernible sequences and their trace in $\St_C$. Recall that
if $(a_i:i\in I)$ is an indiscernible sequence, and $J \subset I$, we write
$a_J:=\{a_j:j \in J\}$.

\begin{proposition} \label{stab0.7} \index{indiscernible}

(i)  Let $(a_i:i \in I)$ be an infinite indiscernible sequence over $C$,
with $\acl(Ca_J) \meet \acl(Ca_{J'}) = \acl(C)$ for any finite disjoint $J,J' \subset I$.
Then $(\St_{C}(a_i):i \in I)$ forms a Morley sequence in $\St_C$.

(ii) Let $q$ be an $\Aut({\cal U}/C)$-invariant type, and suppose that the sequence $(a_i:i\in I)$ satisfies that for all $i\in I$, 
$a_i \models q| C\cup\{a_j:j<i\}$. Then $(a_i:i \in I)$ is $C$-indiscernible,
and  $\acl(Ca_J) \meet \acl(Ca_{J'}) = \acl(C)$ for
any finite disjoint $J,J' \subset I$.
\end{proposition}

{\em Proof.} (i) Clearly, $(\St_C(a_i):i\in I)$ is an indiscernible sequence 
(over $C$) in $\St_C$, so by stability is an indiscernible set. That it is a 
Morley sequence now follows, for example, from Lemma~\ref{Morleybuechler}.

(ii) The indiscernibility is straightforward. For the algebraic closure condition, suppose 
for a contradiction that there is $b\in (\acl(Ca_J)\cap \acl(Ca_{J'}))\setminus \acl(C)$, 
with $J\cap J'=\emptyset$. We may suppose that $J \cup J'$ is minimal 
subject to this (for $b$).
Let $n:=\max\{J\cup J'\}$, with $n\in J'$, say. Let $M$ be
 a sufficiently saturated model containing $C\cup\{a_i:i<n\}$. We may suppose that $a_n\models q|M$. Then 
$b\in\acl(Ca_J)\subset M$, so $b\in M$. Also, the orbit of $b$ under
 $\Aut (M/C\cup a_{J'\setminus \{n\}})$
is infinite, so 
there is $f\in \Aut(M/C \cup a_{J'\setminus\{n\}})$ such that $f(b) \not\in \acl(Ca_{J'})$. Such
$f$
 is not
 elementary over $a_n$. This 
contradicts the invariance of $q$. \qed

 \begin{proposition} \label{stab-divide}  Let $p$ be an $\Aut(\U/C)$-invariant type,
 with $p|C$ stably dominated.  Then for any formula $\phi$ over $\U$,
 $\phi \in p$ if and only if for some $\theta \in p|C$,  $\theta \wedge \neg \phi$ 
divides over $C$.
 \end{proposition}
 
{\em Proof.}  Write $\phi=\phi(x,b)$.   If $\theta \wedge \neg \phi$ divides over $C$,
let $(b_j: j \in J)$ be an infinite sequence with $b=b_1$ and such that 
$\{ \theta(x) \wedge \neg \phi(x,b_j) \}$ is $k$-inconsistent. 
 Let $a \models p | C(b_j: j \in J)$.  Then $\theta(a)$ holds, 
so necessarily $\neg \phi(a,b_j)$ fails for some $j$, i.e. $\phi(a,b_j)$ holds.
Thus $\phi(x,b_j) \in p$.   By  $\Aut(\U/C)$-invariance, $\phi \in p$.  

Conversely, assume $\phi \in p$.  
Let $I$ be an index set, $|I| > |C|+|L|$, and let
$(b_i:i \in I)$ be an indiscernible sequence over $C$, such that  $b=b_1$ and
$(\St_{C}(b_i):i \in I)$ forms a Morley sequence over $C$ in $\St_C$. If $a \models p|C$, then
by \ref{morley-div} $\St_C(a)$ is independent from some $\St_C(b_i)$ for some $i$; hence 
$a \models p | Cb_i$.  Thus   $p|C \union \{ \neg \phi(x,b_i): i \in I \}$
is inconsistent, and the proposition follows by compactness.   \qed

\bigskip

 We will 
be concerned with issues of independence just in $\St_C$ itself, and of stationarity 
when moving from $C$ to its algebraic closure. 
The rest of this chapter is primarily concerned with these stationarity issues.
We therefore introduce the
following notion of `stationary' independence. We often use notation $\St_b$ rather than $\St_C$, 
particularly when the base is varying.

\begin{definition}\label{stationary}\index{stationary independence|(} \index{independence!stationary, $\dnf^s$} \rm
 Suppose at least one of $X,Y$ is a subset of $\St_b$. We 
shall write $\ind{X}{Y}{b}$ if $\St_b(X)$ and $\St_b(Y)$ are independent over $b$ in 
$\St_b$, and in addition,

(s1) for any $a \in \acl(b) \meet \dcl(Xb)$, $\tp(a/b) \proves \tp(a/Yb)$.    
\end{definition}

If $\dcl(b)=\acl(b)$ then (s1) has no content. So
in particular, $\ind{X}{Y}{\acl(b)}$ means only 
that $ \St_b(X),  \St_b(Y)$ are independent over $ \acl(b)$ in the 
stable structure $\St_{\acl(b)}$, and hence also in $\St_b$. 
(Recall that $\St_b$ has essentially the same domain as $\St_{\acl(b)}$.)
The `$s$' in $\dnf^s$ stands for `stationary'.

\begin{proposition}\label{station}
The following conditions are equivalent to (s1) (without any extra assumptions on $b,X,Y$, such as independence).

(s2) For any $a\in \dcl(Xb)$, $\tp(a/b) \vdash \tp(a/\dcl(Yb)\cap \acl(b))$.

(s3) For any $a\in \acl(b) \cap \dcl(Xb)$, $\tp(a/b)\vdash \tp(a/\dcl(Yb)\cap \acl(b))$.

\noindent Assume now that $X\subset \St_b$ and $X$ and $\St_b(Y)$ are independent 
 in $\St_b$. Then (s1)--(s3) are equivalent also to the following conditions. 

(s4) For any $a\in X$, $\tp(a/b)$ is stationary as far as $\St_b(Y)$; that is, 
if $a'\equiv_b a$ and $a'$ and $\St_b(Y)$ are independent in $\St_b$, then
$a'\equiv_{\St_b(Y)} a$.

(s5) For any $a\in X$ and  all finite sets of stable formulas $\Delta$, $\tp(a/b)$ and $\tp(a/\St_b(Y))$ 
have the same $\Delta$-rank and 
$\Delta$-multiplicity.

(s6) For any $a\in X$, any stable formula $\phi$, and any $\phi$-type $q$ consistent with
$\tp(a/\St_b(Y))$ and definable over $\acl(b)$, the disjunction of all 
$\St_b(Y)$-conjugates of the $q$-definition of $\phi$ is defined over $b$.
\end{proposition}

{\em Proof.}
$(s1)\Rightarrow (s2)$ Suppose $\sigma\in \Aut(\U/b)$.  Let $a\in \dcl(Xb)$, and assume  that
$\{e_1,\ldots,e_r\}$ is some complete set of conjugates over $b$, with $e_1\in \dcl(Yb)$. 
Suppose $\phi(a,e_1)$ holds ($\phi$ over $b$). We must show $\phi(\sigma(a),e_1)$ holds.  So 
suppose
without loss that 
$$\{e_1,\ldots,e_s\}:=\{x\in \{e_1,\ldots,e_r\}: \phi(a,x)\}.$$
Let $e:=\lceil \{e_1,\ldots,e_s\}\rceil$. Then $e\in \acl(b) \cap \dcl(Xb)$. Hence
$\sigma(e) \equiv_{Yb} e$. In particular $e_1$ is in the set coded by $\sigma(e)$, so
$\phi(\sigma(a),e_1)$. 

$(s2)\Rightarrow (s3)$.  This is immediate.

$(s3)\Rightarrow (s1)$. Let the conjugates of $a$ over $b$ be $a=a_1,\ldots,a_r$,
 and those over $Yb$
be (without loss) $a_1,\ldots,a_s$. Suppose $a'\equiv_b a$. Suppose $\phi(a,y)$ holds,
 where $y$ is a tuple from $Y$ and
$\phi$ is over $b$. Since $\lceil\{a_1,\ldots,a_s\}\rceil \in \dcl(Yb) \cap \acl(b)$, 
by $(s3)$, $a$ and $a'$ have the same type over $b\lceil\{a_1,\ldots,a_s\}\rceil $, so
$a'\in \{a_1,\ldots,a_s\}$. 
Since
$\phi(a,y)$, also $\phi(a_i,y)$ for $i=1,\ldots,s$. Hence $\phi(a',y)$ holds, as required.

Under the extra conditions, 
the equivalence of (s1) with (s4)--(s6) is an exercise in stability theory.\qed

\medskip

Note that (s2) is symmetric between $X$ and $Y$. For assume (s2) and let
$a'\in \dcl(Yb)\cap \acl(b)$. Then by (s2), $\tp(a'/b)\vdash \tp(a'/b\dcl(X))$, which is (s1) with $X$ and $Y$ reversed. Since (s1) implies (s2), (s2) also holds with $X$ and $Y$ reversed.
It follows that the 
condition $\ind{X}{Y}{b}$ is also symmetric.

 We now have a sequence of easy lemmas giving basic properties of the stationarity condition. 
The first one is a slight  extension of  Proposition~\ref{definable}.

\begin{lemma}\label{stab0.1}  \index{invariant type!and stationary independence}
Suppose $a$ is a tuple in $\St_b$.  

(i) $a \dnf_{\acl(b)} \St_b(Y)$ (in $\St_b$) if and only if  $\tp(a/Yb)$ extends to
an $\Aut({\cal U}/\acl(b))$-invariant type.   

(ii) $\ind{a}{Y}{b}$ if and only if $\tp(a/Yb)$ is the unique extension over $Yb$ of
$\tp(a/b)$ which extends to an $\Aut({\cal U}/\acl(b))$-invariant type. 

(iii) If $p$ is an $\Aut(\U/b)$-invariant type and $a\models p|bY$, then $\ind{a}{Y}{b}$.
\end{lemma}

{\em Proof.} (i) If $\tp(a/Yb)$ extends to an $\Aut({\cal U}/\acl(b))$-invariant 
type $p$, then $p | \St_b$ is $\Aut(\St_b / \acl(b))$-invariant, and extends 
$\tp(a/\St_b(Y))$; hence $a,\St_b(Y)$ are independent in $\St_b$ (as noted in the 
proof of Proposition~\ref{definable}).  Conversely, assume 
$a \dnf_{\acl(b)} \St_b(Y)$. 
Then, as $a\subset \St_b$,  $\tp(a/\St_b(Y))$ extends to an $\Aut(\St_b / \acl(b))$-invariant 
type $p'$.  By stable embeddedness, since $p'$ is a type in a sort of $\St_b$, $p'$ generates a complete 
type $p$ over ${\cal U}$; so clearly $p$ is $\Aut({\cal U}/\acl(b))$-invariant. 
Similarly, $\tp(a/\St_b(Y))$ implies $\tp(a/bY)$, so $p$ extends $\tp(a/bY)$.
 
For (ii), it suffices as above to work within $\St_b$, 
replacing $Y$ by $\St_b(Y)$. Both conditions imply $a\dnf_b \St_b(Y)$.  The statement then is
that $\ind{a}{\St_b(Y)}{b}$ if and only if $\tp(a/b)$ is   stationary as far 
as $\St_b(Y)$. This is just the equivalence 
$(s1) \Leftrightarrow (s4)$ of Proposition~\ref{station}.

(iii) Independence in $\St_b$ is as above. 
Apply Proposition~\ref{station}(s3) for the stationarity condition.
\qed

\bigskip

The next lemma will be extended in Proposition~\ref{simplifications}.

\begin{lemma}\label{stab0.2}
Let $c \subset  \St_b$.  Then  $\St_b(ac) =   \dcl(\St_b(a)c) \meet \St_b$.
\end{lemma}

{\em Proof.}  Let   $E = \dcl(\St_b(a),c) \meet \St_b$.  Clearly $E \subseteq  \St_b(ac)$.  
We must show $\St_b(ac) \subseteq E$. By Remark~\ref{symmetry},
  $\tp(a/\St_b(a))$ implies $\tp(a/\St_b)$, so 
$\tp(a/E) $ implies $\tp(a/\St_b)$; as $c \in E$, it follows that $\tp(ac/E)$ implies 
$\tp(ac/\St_b)$ and
in particular $\tp(ac/E)$ implies $\tp(ac / \St_b(ac))$.  
Hence, $\tp(\St_b(ac)/E)\vdash \tp(\St_b(ac)/bac)$. As $\St_b(ac)\subseteq \dcl(bac)$,
it follows that 
$$\St_b(ac)\subseteq \dcl(E) \cap \St_b =E.$$ \qed

The next technical proposition yields the natural properties of $\dnf^s$ given in 
Corollaries~\ref{stab0.4} and \ref{stab0.5}.

\begin{proposition}\label{stab0.3}
  Let $b$ be a tuple from a set $Z$, and $c\subset \St_b$. Then the following are equivalent.

(i)  $\ind{X}{c}{Z}$.

(ii)  $\ind{\St_b(XZ)} {c} { \St_b(Z)}$.

(iii)  $\ind{\St_b(XZ)} {c} { \St_b(Z)}$ holds in the many-sorted stable structure $\St_b$.

\noindent In case $Z \subset \St_b$, (i)--(iii) are also equivalent to each of : 

(iv)  $\ind{\St_b(X)}{c}{Z}$

(v) $\ind{\St_b(X)}{c}{Z}$ in the structure $\St_b$.
\end{proposition}

{\em Proof.} Notice that $\St_{\St_b(Z)}(\St_b(XZ)) \supseteq \St_b(XZ)$.

(i) $\Rightarrow$ (iii). Assume (i).   By Lemma~\ref{stab0.1}(i) and the symmetry of 
Definition~\ref{stationary},  $\tp(c/XZ)$ extends to an 
$\Aut({\cal U}/\acl(Z))$-invariant type $p$.  By stable embeddedness, for any 
$W$ the restriction map $\Aut({\cal U}/ W) \to \Aut(\St_b / \St_b(W))$ 
is surjective. Hence, as $\St_b(\acl(Z)) \subseteq \acl(\St_b(Z))$, the type   
$p' := p | \St_b$ is $\Aut(\St_b/\acl(\St_b(Z)))$-invariant. Also, 
$p'$ extends $\tp(c/\St_b(XZ))$. Thus, in the structure $\St_b$, we have 
$\St_b(XZ)\dnf_{\acl(\St_b(Z))} c$. To obtain (iii), we show that the 
symmetric form of condition (s1) of Definition~\ref{stationary} holds. 
If $d \in  \acl(\St_b(Z)) \meet \dcl(c\St_b(Z))$,
then $d \in \acl(Z) \meet \dcl(cZ)$, so by (i), $\tp(d/Z)$ implies 
$\tp(d/ZX)$.  Now by stable embeddedness, $\tp(d/\St_b(Z))$ implies
$\tp(d/Z)$.  Hence, $\tp(d/\St_b(Z))$ implies $\tp(d/ZX)$, which implies
 $\tp(d/\St_{\St_b(Z)}(\St_b(XZ)) )$, as required.

(ii) $\Leftrightarrow$ (iii) The structure $\St_{\St_b(Z)}$ may be larger than 
$\St_b$ (it may have more sorts), but in the sorts of $\St_b$, it has precisely 
the structure of $\St_b$ enriched with constants for $\St_b(Z)$.  So for the 
subsets $\St_b(XZ)$ and $c$ of $\St_b$, independence and stationarity
over $\St_b(Z)$ have the same sense in the two structures.

(ii) $\Rightarrow$ (i). Assume (ii).  Then by Lemma~\ref{stab0.1}(i),
$\tp(c/\St_b(XZ))$ extends to an $\Aut({\cal U}/\acl(\St_b(Z)))$-invariant 
type $q$; a fortiori, $q$ is $\Aut({\cal U}/\acl(Z))$-invariant.
As $\tp(c/\St_b(XZ))$ implies $\tp(c/XZ)$ (stable embeddedness), $q$ 
extends $\tp(c/XZ)$.  Hence, by Lemma~\ref{stab0.1}(i), 
we have $\ind{X}{c}{\acl(Z)}$. To prove (i), we now use (twice) the second
part of Lemma~\ref{stab0.1}. Indeed, if $q'$ is another 
$\Aut({\cal U}/\acl(Z))$-invariant extension of $\tp(c/Z)$, 
then (as in the proof of (i)$\Rightarrow$ (ii)) $q' | \St_b$
is an $\Aut(\St_b/\acl(\St_b(Z))$-invariant extension of $\tp(c/\St_b(Z))$,
so $q' | \St_b(XZ) = q| \St_b(XZ)$.
Hence, since $c\in \St_b$,  $q' |XZ=q|XZ$ by stable embeddedness.

Finally, if $Z\subset \St_b$, then $\St_b(Z)=\dcl(Z)$, and the equivalences
 (ii) $\Leftrightarrow$ (iv) and (iii) $\Leftrightarrow$ (v) are immediate, using Lemma~\ref{stab0.2}. \qed

\begin{corollary}[Transitivity] \label{stab0.4}
Assume that either $ac \subset \St_b$ or $d \subset \St_b$.  Then $\ind{ac}{d}{b}$
if and only if $\ind{c}{d}{b}$ and  $\ind{a}{d}{bc}$.  
\end{corollary}

{\em Proof.} First, observe that if $ac \in \St_b$, then in all the relations, 
by definition, $d$ can be replaced by $\St_b(d)$. Thus, we may assume $d \subset \St_b$.
By Proposition~\ref{stab0.3} (i) $\Leftrightarrow$ (iii), we have to show, in $\St_b$, that 
$\ind{\St_b(abc)}{d}{b}$ holds if and only if $\ind{\St_b(bc)}{d}{b}$
and $\ind{\St_b(abc)}{d}{\St_b(bc)}$.  This is just transitivity in the stable 
structure $\St_b$ (for the stationarity condition, it is easiest 
to use (s4)).  \qed

\begin{corollary} \label{stab0.5}
Assume $ab \subset \St_C$ with $\ind{a}{b}{C}$ and  $\ind{ab}{e}{C}$.  Then $\ind{a}{b}{Ce}$.
\end{corollary}

{\em Proof.} We use Corollary~\ref{stab0.4} repeatedly. From $\ind{ab}{e}{C}$,
this gives
$\ind{a}{e}{Cb}$, and hence $\ind{a}{be}{Cb}$. By \ref{stab0.4} again, since 
$\ind{a}{b}{C}$, we have $\ind{a}{be}{C}$, and hence $\ind{a}{b}{Ce}$. \qed

\bigskip

We introduce next a notion of {\em $\dcl$-canonical base} appropriate for the relation
$\ind{X}{Y}{Z}$. See Section 2.3 above, and \cite{pil}, for basics of canonical bases in stable theories.

\begin{definition} \rm \index{canonical base!dcl@$\dcl$}
Let $Y \subset \St_b$.  Write $\Cb(X/Y;b)$ for the smallest 
$Z = \dcl(Zb) \subset \dcl(Yb)$ such that  $\ind{X}{Y}{Z }$ holds.
\end{definition}

\noindent
Note that by Proposition~\ref{stab0.3}, as $Y\subset \St_b$ and  $Z=\dcl(Zb)\subset \St_b$ then 
$\ind{X}{Y}{Z }$ holds if and only if 
$\ind{\St_b(X)}{Y}{Z}$ in the stable structure $\St_b$. Such $\dcl$-canonical
bases exist in stable structures; this is clearest perhaps 
from characterisation $(s4)$ in Proposition~\ref{station}, together with Section 4 (especially Example 4.3) of \cite{pillnorm}. Thus,
$\Cb(X/Y;b)$ always exists (and is unique).

\begin{lemma}\label{stab0.6}
 Assume $d \subset \St_b$, $\ind{b'}{d}{b}$, and $d' \in \St_{b'}$ where $d'$
 is a finite tuple.  Assume $d' \in \dcl(bb'd)$.
Then there exists a finite tuple $f \in F := \Cb(b'd'/bd;b)$ such that $F = \dcl(bf)$.
\end{lemma}

{\em Proof.} We have $d'\in \dcl(bb'F)$: indeed, $\ind{b'd'}{d}{F}$,
so $\ind{d'}{d}{Fb'}$ (Corollary~\ref{stab0.4}), and $d'\in \dcl(Fdb')$ (as $d'\in \dcl(bb'd)$), so $d'\in \dcl(Fb')$.
 Choose $f$ to be any finite tuple 
in $F$ such that
$d'\in \dcl(bb'f)$. Then as $\ind{b'}{d}{b}$ and $f\in \dcl(bd)$, we have
$\ind{b'}{d}{bf}$ (\ref{stab0.4} again), so
$\ind{b'd'}{d}{bf}$. By minimality of $F$, it follows that $F =\dcl(bf)$. \qed

\begin{lemma}\label{canbase}
Let $b,d \in \St_{\emptyset}$, and let $a$ be arbitrary. Assume $\ind{a}{b}{d}$
 and $\ind{a}{d}{b}$.
Let $f:=\Cb(a/b;\emptyset)$. Then $f\subseteq \dcl(d)$.
\end{lemma}

{\em Proof.} Let $f':=\Cb(a/bd;\emptyset)$. By the minimality in the definition of
 canonical base, 
$f'\subseteq \dcl(d)$
(as $\ind{a}{bd}{d}$). Similarly, $f'\subseteq \dcl(b)$. But we have $\ind{a}{bd}{f'}$, so
$\ind{a}{b}{f'}$. Hence by minimality of $\dcl(f)$ it follows that $f\subseteq \dcl(f')$. 
Thus, $f\subseteq \dcl(d)$.  \qed
\index{stationary independence|)}


\begin{definition}\label{*type} \rm
By a {\em $*$-type} \index{type!aatype@$*$-type}we mean a type in variables $(x_i)_{i \in I}$, where $I$ is a 
possibly infinite index set.  A {\em $*$-function} \index{function!aafunction@$*$-function} is a tuple $(f_i)_{i \in I}$ of 
definable functions.  Definable functions 
always depend on finitely many variables, but may have other dummy variables. 
A $*$-function is {\em $C$-definable} if all its component functions are $C$-definable. We emphasise that in this chapter types have generally been infinitary. The purpose of this definition is rather to express how $*$-functions are decomposed into definable functions.

We say that a $*$-type  $p$ over $C$ is {\em stably dominated  via a $C$-definable $*$-function $f$} 
\index{stable domination!via a $*$-function} if, whenever $a \models p | C$,  we have $f(a) \subset \St_C$, and for any 
$e$, if $\ind{f(a)}{e}{\acl(C)}$ (equivalently $f(a)\dnf_{\acl(C)} \St_C(e)$) then 
$\tp(e/Cf(a)) \vdash \tp(e/Ca)$. 
\end{definition} 

It is clear that a $*$-type over $C$ is stably dominated  if and only if it 
is stably dominated  by some $C$-definable $*$-function.

\begin{proposition}\label{stab1}   Let $c$ enumerate $\acl(C)$.  Let $f$ be a $C$-definable 
$*$-function from $\tp(a/C)$ to $\St_C$.  Then  $\tp(a/C)$ is stably dominated via $f$ 
 if and only if both of the following hold:

 (i) $\tp(a/Cf(a)) $ implies $\tp(a/ cf(a))$,

 (ii) $\tp(a/\acl(C))$ is stably dominated via $f$.
\end{proposition}

{\em Proof.} Assume that $\tp(a/C)$ is stably dominated via $f$. So for any $e$,
if $\St_C(f(a))\dnf_{C} \St_C(e)$ then  $\tp(e/Cf(a))\vdash \tp(e/Ca)$. Putting 
$e=c$, we get $\tp(c/Cf(a))\vdash \tp(c/Ca)$, which is equivalent to (i). 

To see (ii), suppose $\St_c(f(a))\dnf_{c} \St_c(e)$. Then $\St_C(f(a))\dnf_c \St_C(e)$, so
$\St_C(f(a))\dnf_{c} \St_C(ec)$. Thus, $\tp(ec/Cf(a))\vdash \tp(ec/Ca)$. 
In particular, $\tp(e/cf(a))\vdash \tp(e/ca)$, which gives (ii).

Conversely, assume that (i) and (ii) hold. Suppose $\St_C(f(a))\dnf_{c} \St_C(e)$. 
Then by (ii), $\tp(e/cf(a))$ implies $\tp(e/ca)$.   Thus $\tp(a/cf(a))$ implies 
$\tp(a/cf(a)e)$. But by (i), $\tp(a/Cf(a))$ implies $\tp(a/cf(a))$, so it implies 
$\tp(a/cf(a)e)$ and in particular $\tp(a/Cf(a)e)$, as required. \qed

\begin{remark}\label{partacl} \rm
The above proof also shows that if $f$ is a $C$-definable $*$-function from $\tp(a/C)$ to $\St_C$
 and $\tp(a/C)$ is stably dominated by $f$, then for any $C''$ with $C\subseteq C''\subseteq \acl(C)$, 
$\tp(a/C'')$ is stably dominated by $f$.
\end{remark}

\begin{corollary}\label{stab1.1} \index{stable domination!over algebraic closure of base}
(i) Let $f$ be a $C$-definable $*$-function, and suppose that
 $\tp(a/\acl(C))$ is stably dominated via $f$. Then $\tp(a/C)$ is 
stably dominated via $f'$, where $f'(a)$ enumerates 
$f(a) \cup (\acl(C) \cap \dcl(Ca))$.  

(ii) Suppose that $\tp(a/\acl(C))$ is stably dominated 
via $F$, where $F$ is an $\acl(C)$-definable $*$-function. Then $\tp(a/\acl(C))$ is 
stably dominated  via a $C$-definable $*$-function.

(iii) $\tp(a/C)$ is stably dominated if and only if 
$\tp(a/\acl(C))$ is stably dominated.
\end{corollary}

{\em Proof.} Recall that the domains of $\St_C$ and $\St_{\acl(C)}$ are the same, 
though there are more relations on $\St_{\acl(C)}$. 

(i) In Proposition~\ref{stab1}, condition (i) is equivalent to the condition
$\acl(C)\cap \dcl(Ca)=\acl(C) \cap \dcl(Cf(a))$. This certainly holds if $f(a)$ contains
$\acl(C) \cap \dcl(Ca)$. Thus, (i) follows from Proposition~\ref{stab1}.

(ii)  We may write $F(x)=F^*(d,x)$, where $F^*$ is $C$-definable, and $d\in \acl(C)$. 
Now define $f(x)$ to be a code for the $*$-function $d\mapsto F^*(d,x)$ (more precisely, a sequence of codes for the component functions): so 
$f(x)$ is a code for a $*$-function
on the finite $C$-definable set of conjugates of $d$ over $C$. Then $\tp(a/\acl(C))$ is stably dominated  by 
$f$, which is $C$-definable.

(iii) By (i) and (ii), stable domination over $\acl(C)$ implies stable domination 
over $C$. The other direction is immediate from \ref{stab1}. \qed

\begin{proposition} \label{bits}
(i) Suppose $\tp(a/C)$ is stably dominated, but $\St_C(a)=C$. Then $a\subset \dcl(C)$.

(ii) Suppose $B\subset \St_C$ and $A$ is dominated by $B$ over $C$ (that is, whenever 
$B\dnf_C \Dst$ we have $\tp(D/CB)\vdash \tp(D/CA)$). Then $\tp(A/C)$ is 
stably dominated.

(iii) Suppose $\tp(A/C)$ is stably dominated, 
and $a'\subset \dcl(CA)$. Then $\tp(a'/C)$ is stably dominated.
\end{proposition}

{\em Proof.} (i) Suppose that $a\not\subset \dcl(C)$ and choose a distinct conjugate 
$a'$ of $a$ over $C$. Then $\St_C(a')=C$, so $\St_C(a)\dnf_C \St_C(a')$. Also, 
$a\equiv_{\St_C(a)} a'$,
so by stable domination $a \equiv_{Ca} a'$, which is a contradiction.

(ii) We may suppose $\Ast\dnf_C \Dst$, and must show that 
$\tp(D/C\Ast)\vdash \tp(D/CA)$. So suppose $D\equiv_{C\Ast} D'$. It is possible to
 find $D''\equiv_{\acl(C)} D'$ with $D'' \equiv_{CA} D$. (Indeed, suppose $\{e_1,\ldots,e_r\}$ 
is an orbit over $C$, and $\tp(D/CA)$ implies that
$\phi(D,e_i)$ holds if and only if $i=1$; then $e_1\in CA^{\st}$, so $\phi(D',e_1)$ holds; thus the 
two conditions on $D''$ are consistent.) Observe that we may replace $B$ by any $B'\equiv_{CA} B$, without affecting the assumption. 
By stable embeddedness, if $B'\equiv_{C\Ast} B$ then $B'\equiv_{CA} B$. Thus, we may choose $B$ so that
$B\dnf_{C\Ast} (D''D')^{\st}$ in $\St_C$. Since $\Ast\dnf_C D''^{\st}$ and $\Ast\dnf_C D'^{\st}$, 
by transitivity we have $B\dnf_C D''^{\st}$ and $B\dnf_C D'^{\st}$. 
Thus $D''\equiv_{CB} D'$, so $D''\equiv_{CA} D'$ by the domination assumption. It 
follows that $D\equiv_{CA} D'$.

(iii) We may apply (ii), since $\tp(a'/C)$ is dominated by $\Ast$. \qed

\begin{remark} \rm It would be possible to localise the notion of stable domination by working with
just some of the sorts of $\St_C$. Let ${\cal F_C}$ be a collection of $C$-definable stable stably embedded
sets, and let $\St_C^{{\cal F}}$ be the many-sorted structure whose sorts are the members of ${\cal F}$, with 
the induced $C$-definable structure. There will be a corresponding notion of ${\cal F}$-stable domination, defined 
as above but with $\St_C^{{\cal F}}$ replacing $\St_C$. To develop a satisfactory theory, one would need certain 
closure properties of ${\cal F}$. We do not pursue this here.
\end{remark}

\chapter{Invariant types and change of base}\label{invarianttypes}

The property of a type being stably dominated is quite
sensitive to the base set of parameters. In this chapter, we will
consider conditions under which stable domination is preserved under
increasing or decreasing the base. Notice that it is not in general true 
that if $\tp(a/C)$ is stably dominated and $C\subset B$ then $\tp(a/B)$ is 
stably dominated. For example, in ACVF, in the notation of Part~II (or \cite{hhm}),
 let $C=\emptyset$, and let $a$ be a 
generic element of the valuation ring. Let $s$ be a code for the open ball $B_{<1}(a)$, and identify $s$ with $\res(a)$. Then 
$\tp(a/C)$ is stably dominated by the residue map, but $\tp(a/Cs)$ is not stably dominated; in particular, $\tp(a/Cs)$ is not stably dominated by the residue map, since $s=\res(a)$ is already in the base $Cs$, so independence in $\St_{Cs}$ from $s$ has no content. (This will 
become clearer in Part~II; essentially, by Lemma~\ref{eqorthog} and Corollary~\ref{equivorthstabdom},
the generic type of a closed ball is stably
 dominated, but the generic type of an open ball is not.)

More generally, suppose every extension of $\tp(a/C)$ is stably dominated. Let $b=\St_C(a)$. Then $\tp(a/b)$ is stably dominated.
But $\St_b(a)=b$, so $\tp(a/b)$ is stably dominated by $b$. Thus, by Proposition~\ref{bits}(i), $a\in \dcl(b)$. 
Thus, under this assumption, $a$ must be a tuple from $\St_C$.

 Thus, some further 
hypothesis is needed, and the one we use is that $\tp(a/B)$ has an 
$\Aut(\U/C)$-invariant extension to $\U$. 
This gives quite an easy result (Proposition~\ref{dropstabdom} below)
for the case of increasing the parameter set. To show that stable
domination is preserved under decreasing the parameter set is 
rather more difficult. We give an easy special case in Lemma~\ref{stab1.2},
and our most general result in Theorem~\ref{stab7}. The hypothesis in \ref{stab7} that $\tp(B/C)$
is invariant is stronger than it might be; it would be good to investigate the 
weakest assumptions under which some version of  Theorem~\ref{stab7} could
be proved.

First, we give the easy direction of base change.

\begin{proposition} \label{dropstabdom}
Suppose that $C\subseteq B \subset \U$ and that $p$ is an $\Aut(\U/\acl(C))$-invariant
 type over $\U$. 
Suppose also that $p|C$ is stably dominated. Then $p|B$ is stably dominated.
\end{proposition}

{\em Proof.} Let $A \models p|B$. Suppose that $\St_B(A) \dnf_B \St_B(D)$ 
(in $\St_B$) and that $D \equiv_{\St_B(A)} D'$. We must show that $D\equiv_{BA} D'$. 

As $A\models p|B$, $\St_C(A)\dnf_C \St_C(B)$ (Lemma~\ref{stab0.1}). Also, $\St_C(A)\dnf_B \St_C(BD)$
 in $\St_B$,
so $\St_C(A)\dnf_{\St_C(B)} \St_C(BD)$ in the structure $\St_C$.
 Hence, by transitivity,
 $\St_C(A)\dnf_C \St_C(BD)$ in $\St_C$.
The assumption $D\equiv_{\St_B(A)} D'$ gives $BD\equiv_{\St_C(A)} BD'$ (for $B \subset \St_B(A)$);
 so by stable 
domination
of $p|C$, $BD\equiv_{CA} BD'$. Hence $D\equiv_{BA} D'$, as required. \qed

\bigskip

We next tackle the harder (downwards) direction of base change, first dealing with an easy special case.

\begin{lemma} \label{stab1.2}
Let $g$ be a $C$-definable $*$-function, and let $p,q $ be $\Aut({\cal U}/C)$-invariant types. 
Assume that when  $b \models q |C$, $p|Cb$ is stably dominated  via $g$.  
Then $p|C$ is stably dominated  via $g$.
\end{lemma}

{\em Proof.} 
First observe that by Lemma~\ref{3.22}, if $a\models p|C$ then the entries $a_i$ of $g(a)$ lie in $\St_C$; for any such entry satisfies an $\acl(C)$-definable type $q_i$, and if $b\models q|C$ and $a_i\models p|Cb$, then $a_i$ lies in an $\acl(C)b$-definable stable stably embedded set, and so lies in one definable over $\acl(C)$ and hence in one defined over $C$.

Since $p$ is $\Aut(\U/C)$-invariant, $\dcl(Ca) \cap \acl(C)=\dcl(C)$. Thus,
by Corollary~\ref{stab1.1}(i), it suffices to prove that if $c$ is an enumeration of $\acl(C)$, then $p|c$ is stably dominated by $g$. 
 Let $a\models p$ and suppose that for some $e$ we have $\ind{g(a)}{e}{c}$. 
We must show that $\tp(e/cg(a))$ implies $\tp(e/ca)$. So suppose that $e'\equiv_{cg(a)} e$.
Let $b\models q|caee'$. Then 
$\ind{b}{g(a)}{ce}$ (via Lemma~\ref{stab0.1} (iii) to obtain first 
$\ind{\St_{ce}(b)}{g(a)}{ce}$).
As 
$\ind{e}{g(a)}{c}$, by transitivity, $\ind{be}{g(a)}{c}$, 
and so $\ind{e}{g(a)}{cb}$. Likewise,
$\ind{e'}{g(a)}{cb}$. By the choice of $b$, $e'\equiv_{cb} e$. It follows that
$e'\equiv_{cbg(a)} e$ (use Proposition~\ref{station} (s4) and 
an automorphism, as $g(a)\in \St_C$). 
 By domination over $Cb$ and Remark~\ref{partacl}, $\tp(e/cbg(a))$ implies $\tp(e/cba)$.  
Thus, $e'\equiv_{cba} e$, so $e' \equiv_{ca} e$. \qed

\bigskip

In Theorem~\ref{stab7} we substantially extend this result, removing the 
assumption that  stable domination over $p|Cb$ is via a fixed $C$-definable function.

\begin{lemma}\label{stab2}
 Let $f_0,f_1,f_2$ be $C$-definable $*$-functions, all with range in $\St_C$, and 
let $p$ be a type over $C$ stably dominated via $f_1$ and via $f_2$. 
 Suppose that for $a \models p$,
$\ind {f_1(a)} {f_2(a)} {Cf_0(a)}$.  Also suppose $(*)$ that
$\acl(C) \meet \dcl(Cf_0(a)) = \acl(C) \meet \dcl(Ca) $. 
Then $p$ is stably dominated via $f_0$.
\end{lemma}

{\em Proof.}  We prove the result over $c$, an enumeration of $\acl(C)$. 
For suppose the result is proved in this case. The assumption $(*)$
ensures that $\tp(a/Cf_0(a))\vdash \tp(a/c)$. Hence, it will follow from 
Proposition~\ref{stab1} that $p$ is stably dominated via $f_0$ over $C$. 

Clearly, $\ind {f_1(a)} {f_2(a)} {cf_0(a)}$, and by 
Proposition~\ref{stab1}, any extension of $p$ over $c$ is 
stably dominated  via $f_1$ and via $f_2$. Thus, we may now replace $C$ by $c$.

Let $a \models p$, and let $a_i = f_i(a)$ for $i=0,1,2$.  So $\ind{a_1}{a_2}{ca_0}$.

\medskip

{\em Claim.}  For any $d$, if $\ind{d}{a_0}{c}$ then $\ind{d}{a_1}{c}$.  

{\em Proof of Claim.}  Since $a_1,a_0$ lie in $\St_c$ which is stably embedded,
it suffices to prove the claim for $d \subset \St_c$. We now work with 
stable independence $\dnf$ in the structure $\St_c$. We may
assume $d\dnf_{ca_0a_1} a_2$.  As $a_1\dnf_{ca_0} a_2$, by transitivity 
and symmetry we have $da_1\dnf_{ca_0} a_2$, so $d\dnf_{ca_0} a_2$.  But $d\dnf_c a_0$, so 
$d\dnf_c a_2$, and hence, as $c$ is algebraically closed, $\ind{d}{a_2}{c}$.  By the domination assumption,
$\tp(d/ca_2)$ implies $\tp(d/ca)$.  In particular $\tp(d/ca_2)$ implies 
$\tp(d/ca_2a_1)$. Choose $d'\equiv_{ca_2} d$ with $d'\dnf_{ca_2} a_1$.
As $d\dnf_c a_2$, also $d'\dnf_c a_2$. By transitivity, 
$d'\dnf_c a_1a_2$, so $\ind{d'}{a_1a_2}{c}$. Hence also  $\ind{d}{a_1a_2}{c}$,
and in particular $\ind{d}{a_1}{c}$.

\medskip

Now to prove the lemma, assume  $\ind{e}{a_0}{c}$.   Then $\ind{e}{a_1}{c}$ by the claim, so,
as $\tp(a/c)$ is stably dominated by $a_1$, $\tp(e/ca_1) \proves \tp(e/ ca)$.  
Also, since $\tp(e/ca_0)$ implies $\ind{e}{a_1}{c}$,
it implies $\tp(e/ca_1)$, e.g. by Lemma~\ref{stab0.1}(ii).  Thus $\tp(e/ca_0)$ implies $\tp(e/ca)$. \qed

\begin{remark} \label{stab2.1}\rm
If $p$ extends to an $\Aut(\U/C)$-invariant type and $a\models p|\U$, then every elementary permutation of $\acl(C)$ which is the identity on $C$ is elementary over $a$.  In particular, $\acl(C) \cap \dcl(Ca)=\dcl(C)$. 
It follows that our assumption $(*)$ automatically holds in this case.
\end{remark}

 The next lemma gives  a useful symmetry condition for invariant extensions. Without the stable domination assumption, any o-minimal structure $\U$ would provide a counterexample: let $p$ and $q$ both be the 
 $\Aut(\U)$-invariant type `$x<\U$'.

\begin{lemma}\label{stabrem3}
Let $p,q$ be $\Aut({\cal U}/C)$-invariant types, with $p|C$ stably dominated.  
Let $a \models p |C$, $b \models q| Ca$.  Then $a \models p | Cb$.
\end{lemma}

{\em Proof.} Since $p|C$ is stably dominated, we may suppose that $a\in \St_C$. 
By Lemma~\ref{stab0.1}(ii) it suffices to show $\ind{a}{b}{C}$.
By forking symmetry and Proposition~\ref{station} (s4), it suffices to show
$\ind{\St_C(b)}{a}{C}$. This follows from the assumption by 
Lemma~\ref{stab0.1}(ii) again. \qed

\begin{lemma}\label{stab4}  Let $((b_i,d_i): i \in I)$ be an indiscernible 
sequence, with 
$d_i \subset \St_{b_i}$ for each $i \in I$.    Fix $i \in I$
and let $J$ be an infinite subset of $\{i' \in I: i' <i \}$. Put
 $b_J := \{b_j: j \in J\}, d_J := \{d_j : j \in J\}$. Then  
$\ind{d_i}{b_{<i}d_{<i}} {b_{J}b_i d_J}$. 
\end{lemma}

{\em Proof.} We apply criterion (s5) in Proposition~\ref{station}. we may suppose that $d_i$ is a finite tuple. 
 Let  $\Delta$ be a finite set of formulas $\phi(x, y, z)$, 
$y=(y_1,\ldots,y_n)$, such that $\phi(x,y,b_j)$
implies $x  \in \St_{b_j}$.  Then  $\phi(x,y,b_j)$ is a stable formula.  
Let $\Delta_i = \{ \phi(x,y,b_i): \phi(x,y,z) \in \Delta \}$.
We have to show that (for any such $\Delta$)
$$\rk_{\Delta_i } \tp(d_i / d_{<i} b_{\leq i}) 
                  = \rk_{\Delta_i }\tp(d_i / d_{J} b_{J}b_i) 
$$
and similarly for multiplicity. 
Now $\leq$ is clear. For the $\geq$-inequality,  let 
$\psi \in  \tp(d_i / d_{<i} b_{\leq i})$ have minimal $\Delta_i$ -rank and 
multiplicity. Then $\psi$ uses  finitely many constants from $d_{<i}$ and from 
$b_{<i}$ (as well as $b_i$). Since $J$ is infinite, and by indiscernibility, 
these can be replaced by constants from $d_J,b_J$ (keeping $b_i$) so that the 
resulting formula $\psi'$ has the same $\Delta_i$-rank.  This gives 
the $\geq$ inequality and proves the lemma. \qed

\bigskip

In the next lemma, we use the clumsy symbols $i,J,J'$ for indices, instead of 
just 0,1,2. The reason is that in the proof 
of Lemma~\ref{stab6}, this notation will fit with that of Lemma~\ref{stab4} 
(which will yield condition (A3) below). In applications, $J$, $J'$ are sets of 
indices (much as in Lemma~\ref{stab4}).
The final two clauses of (A1) can be regarded as the first two with $Ji$ in place of
 $J$ (respectively $J'i$ in place of $J'$). And the final clause of (A2) is the first 
with $J\cup J'$ in place of $J$.

\begin{lemma}\label{stab4.1}  

Assume we are given $b_J,b_{J'},b_i,d_{J},d_{J'},d_i$, with  
$d_\nu \subset \St_{b_\nu}$ for each $\nu\in \{i,J,J'\}$, satisfying the following independence conditions.

$$    \ind{d_J} {b_{J'}} {b_{J}} ,\ \  \ind{d_{J'}} {b_J} {b_{J'}},
\ \
 \ind{d_Jd_i}{b_{J'}}{b_Jb_i}, \ \ \ind{d_{J'}d_i}{b_J}{b_{J'}b_i}\ \ \eqno(A1) 
$$

$$    \ind{b_i}{\St_{b_J}(b_{J}b_{J'}d_{J}d_{J'})} {b_J},\ \ 
\ind{b_i}{\St_{b_{J'}}(b_{J}b_{J'}d_{J}d_{J'})} {b_{J'}},\ \ 
\ind{b_i}{\St_{b_Jb_{J'}}(b_Jb_{J'}d_Jd_{J'})}{b_Jb_{J'}}\ \ \eqno(A2)
$$

$$   \ind{  d_i }  {d_{ J' } }  { b_{ J}b_{ J'}b_id_J}, \ \   
             \ind{  d_i }  {d_{ J} }  { b_{ J}b_{ J'}b_id_{J'}} \ \ \eqno(A3)
$$ 
Then the following hold.

   (i) $$ \ind{ b_id_i }  {d_{ J' } }  { b_{ J} b_{J'}d_J}. $$

(ii)  Let $f $ enumerate  $\Cb(b_id_i / b_Jd_J;b_J)$, and let $f'$ 
enumerate
$\Cb(b_id_i / b_{J'}d_{J'};b_{J'})$.  Then $\dcl(f b_Jb_{J'}) = 
\dcl(f' b_Jb_{J'})$.    
\end{lemma}

{\em Proof.}  (i)  By (A3), $\ind{ b_id_i }  {d_{ J' }} { b_Jb_{J'}b_id_J} $. 
 By forking symmetry,
$ \ind {d_{ J' } } { b_id_i }  { b_Jb_{J'}b_id_J} $.  But by
 assumption (A2) (last clause), 
$\ind {d_{ J' } } { b_i  }  { b_{ J}b_{ J'} d_J} $.  By 
transitivity, $ \ind {d_{ J' } } { b_id_i }  { b_{ J}b_{ J'} d_J} $.

(ii)  We have $\ind{b_id_i} {d_J} {b_J f}$, by definition of $f$.  Thus by 
transitivity,

$$ 
    (1) \ \ \ \ \ \ \ \ \ \ \  {\ind{d_i}{d_J}{b_Jb_if} }  .          
$$
By (A1), 
$$
    (1.1)\  \  \ind{d_id_J}{b_{J'}}{b_Jb_i}.
$$
We have $f \subset \dcl(b_Jd_J)$, so 
$\ind{fd_id_J}{b_{J'}}{b_Jb_i}$, and hence  
$\ind{d_id_J}{b_{J'}} {b_{J}b_if}$.   By Corollary~\ref{stab0.5},               
$$ 
   (2)  \ \ \ \ \ \ \ \ \ \ \  {\ind{d_i}{d_J}{b_{J'}b_Jb_if}} . 
$$
By (A3), 
$$ 
    (3')   \ \ \ \ \ \ \ \ \ \ \  \ind{d_i}{d_{J'}} {b_{J} b_{J'}b_i d_J}.
$$
As $f \subset \dcl(b_Jd_J)$, by (2), $(3')$,  and transitivity, 
$$ 
   (3)   \ \ \ \ \ \ \ \ \ \ \  {\ind{d_i}{d_{J'}} {b_{J}b_{J'} b_i f}}. 
$$
Now  by (A2) (last clause) 
$\ind{b_i}{d_{J}d_{J'}}{b_Jb_{J'}}$, so $\ind{b_i}{d_{J'}f}{b_{J}b_{J'}}$, and 
thus $\ind{b_i}{d_{J'}}{b_{J'}b_Jf}$. 
From this,  (3), and transitivity, we obtain 

$$  
    (4')  \ \ \ \ \ \ \ \ \ \ \    \ind{b_id_i}  {d_{J'}} {b_{J}b_{J'}f} .
$$
Let $e = \St_{b_{J'}}(b_J f)$.  Also,  
  $\ind{b_J f}{d_{J'}}{b_{J'}e}$ (for example, use Proposition~\ref{stab0.3} (iii) $\Rightarrow$ (i),
 with $X=b_Jf, b=b_{J'}, Z=\St_{b_J'}(b_Jf), c=d_{J'}$).
  Together with  ($4'$), transitivity gives 
$$ 
     (4)   \ \ \ \ \ \ \ \ \ \ \  \ind{b_id_i}  {d_{J'}} {b_{J'}e} .       
$$
Now by (A3), we have
$$  
    (5')   \ \ \ \ \ \ \ \ \ \ \  \ind{d_i}{d_{J}}{b_{J}b_{J'}b_i d_{J'}}  .   
$$
As $e   \subset \dcl(b_{J'}b_Jd_J)$,

$$   
     (5)   \ \ \ \ \ \ \ \ \ \ \  \ind{d_i}{e}{b_{J}b_{J'}b_i d_{J'}}  .   
$$
Dually to  (1.1), we have $\ind{d_id_{J'}}{b_{J}}{b_{J'}b_i}$, so 
$$  
    (6)   \ \ \ \ \ \ \ \ \ \ \   \ind{d_i}{b_J}{b_{J'} b_i d_{J'}}.
$$
By (5), (6) and transitivity,  we obtain
$$ 
    (7) \ \ \ \ \ \ \ \ \ \ \ \    \ind{d_i}{e}{b_{J'}b_i d_{J'}} .  
$$
By (A2), we have  $\ind{b_i}{\St_{b_{J'}}(b_{J'}b_{J}d_{J}d_{J'})} {b_{J'}}$. Thus,
$\ind{b_i}{d_{J'}e}{b_{J'}}$. Note that $e\subset \St_{b_{J'}}$, so Corollary~\ref{stab0.4} is applicable. 
Hence by transitivity $\ind{b_i}{e}{b_{J'}d_{J'}}$, so by transitivity and (7),
$$ 
         (8)  \ \ \ \ \ \ \ \ \ \ \ \    \ind{b_id_i}{e}{b_{J'} d_{J'}}   .
$$
Now $f'$ enumerates the canonical base of $b_id_i/{b_{J'} d_{J'}}$. It
follows from (4), (8), and Lemma~\ref{canbase} (applied over $b_{J'}$)
 that $f' \subseteq \dcl(b_{J'}e) \subseteq \dcl(b_{J}b_{J'} f)$.  Dually,
$f \subseteq \dcl(b_{J}b_{J'} f'$). Thus, $\dcl(fb_Jb_{J'})=\dcl(f'b_Jb_{J'})$.  \qed

\bigskip

For the purposes of the lemma below, let  $p,q$ be 
$\Aut({\cal U}/C)$-invariant   $*$-types, and $F$ a $C$-definable $*$-function.
We will be interested in the property  
\begin{eqnarray*}
    \Pi^*(p,q,F): & &  \mbox{ if } b \models q|C , b' \models q|Cb, a \models p| Cbb' \\
         & & \mbox{ then } \dcl(Cbb'F(a,b))= \dcl(Cbb'F(a,b')). 
\end{eqnarray*}

If $q$ is an $\Aut({\cal U}/C)$-invariant $*$-type, let 
$q^{\omega} = \tp((b_0,b_1,\ldots,b_\omega)/C)$ where $b_\alpha \models q | \dcl(Cb_\beta:\beta<\alpha)$.
Then $q^{\omega}$ is also an $\Aut({\cal U}/C)$-invariant $*$-type.

\begin{lemma}\label{stab6}
 Let $p,q$ be $\Aut({\cal U}/C)$-invariant $*$-types,
 $F=(F_\lambda)_{\lambda\in \Lambda}$ a $*$-definable function, with 
each $F_\lambda$ $C$-definable. 
Assume that  if $b \models q | C$, then $p|Cb$ is stably dominated  via $F(x,b)$.  
Then there exists 
 a $C$-definable $*$-function $f$  
such that:

(i)   if $b \models q^{\omega} | C$, then $p|Cb$ is stably dominated over 
$Cb$ via 
$f(x,b)$;

(ii)   $\Pi^*(p,q^{\omega},f)$.

\noindent
If already $\Pi^*(p,q ,F)$, then we may choose $f$ to have form $f=(f_\lambda)_{\lambda \in\Lambda}$, so that 
for each $\lambda \in \Lambda$, $f_\lambda$ is a $C$-definable function and 
$\Pi^*(p,q^{\omega},f_\lambda)$ holds.
\end{lemma}

{\em Proof.}  
 For convenience we now work over $C$, so  drop all
reference to $C$. Let $a\models p|\emptyset$. Let $I:=\omega+\omega+1$ and let $(b_i: i \in I)$ 
be an indiscernible sequence over $a$, such that 
$b_i \models q | ab_{<i}$ for each $i$ (so obtained as in Proposition~\ref{stab0.7}). 
Let $d_i(\lambda):=F_\lambda(a,b_i)$ for each $i,\lambda$, and put
$d_i=(d_i(\lambda):\lambda\in \Lambda)$. For each $\lambda\in \Lambda$,
$(b_id_i(\lambda):i\in I)$ is an indiscernible sequence over $a$.

Let $J=\omega=(0,1,2,\ldots)$, 
$J'=(\omega+n:n\in \omega)=(\omega,\omega+1,\omega+2,\ldots)$, and $i=\omega 2:=\omega+\omega$.
As in Lemma~\ref{stab4}, let $b_J=\{b_j:j\in J\}$, and similarly 
$b_{J'}$, $d_J$, $d_{J'}$. Also, let ${b^*}_1 = (b_0,  b_1,\ldots,b_\omega)$, 
${b^*}_2 = (b_{\omega+1},b_{\omega+2},\ldots, b_{ \omega 2})$,
${d^*}_1 =  (d_0,d_1,\ldots,d_\omega)$,  
${d^*}_2 = (d_{\omega+1},b_{\omega+2},\ldots, d_{\omega 2})$.
The hypotheses of Lemma~\ref{stab4} hold. 
We shall verify that the hypotheses of Lemma~\ref{stab4.1} also hold.
 
A1 is   easily seen via Lemma~\ref{stabrem3}, which ensures
$\tp(a/b_I) = p | b_I$.   For instance $\ind{d_J}{b_{J'}}{b_J}$ follows from $\ind{a}{b_{J'}}{b_J}$.
Condition A2 follows from the fact that $b_i \models q | b_{<i}a$,
so $b_i \models q |\St_{b_J}(b_{<i}d_{<i})$, using Lemma~\ref{stab0.1} 
and invariance of $q$.  
The first clause of (A3) follows from Lemma~\ref{stab4}, which yields
$\ind{d_{\omega 2}}{b_Jb_{J'}d_Jd_{J'}}{b_Jb_{\omega 2}d_J}$, and hence
$\ind{d_{\omega 2}}{d_{J'}}{b_Jb_{J'}b_{\omega 2}d_J}$. The second 
clause follows similarly.

Now $ {b^*}_1 \models q^\omega$. This is an $\Aut({\cal U})$-invariant 
type, and ${b^*}_2 \models  q^{\omega} |  {b^*}_1$.

Now $\Cb(b_{\omega 2 }d_{\omega 2} / b_{J}d_{J} ; b_J)\subseteq \dcl(ab_J)$.
Thus, there is a $*$-function $f$ such that
$f(a,{b^*}_1)$ enumerates $\Cb(b_{\omega 2}d_{\omega 2}/ b_{J}d_{J} ; b_{J})$.
We may suppose $f$ is a sequence of $*$-functions
$(\hat{f}_\lambda)$, where $\hat{f}_\lambda(a,{b^*}_1)$ enumerates
$\Cb(b_{\omega 2}d_{\omega 2}(\lambda) / b_{J}d_{J} ; b_{J})$.
Here, the $d_{\omega 2}(\lambda)$ form a directed system of finite tuples with limit $d_{\omega 2}$, and each function $\hat{f}_\lambda$ depends on $a$ and finitely many of 
the variables corresponding to $b_1,b_2,\ldots$; notice that $b_\omega$ is 
not used. Likewise, by indiscernibility, $f(a,{b^*}_2)$  enumerates 
$\Cb(b_{\omega 2}d_{\omega 2}/b_{J'}d_{J'};b_{J'})$.
It follows by Lemma~\ref{stab4.1}(ii) that
$$
   \dcl(f(a,{b^*}_1), b_J,b_{J'})=\dcl(f(a,{b^*}_2),b_J,b_{J'}).
$$ 
This yields (ii) of the lemma.

Next, we obtain (i). By definition of canonical basis, 
$\ind{b_{\omega 2}d_{\omega 2}}{b_Jd_J}{f(a,{b^*}_1)}$, so
$\ind{d_{\omega 2}}{d_J}{b_{\omega 2}b_Jf(a,{b^*}_1)}$.
It  follows from Lemma~\ref{stab2} and Remark~\ref{stab2.1}, applied 
over $C':=\dcl(b_Jb_{\omega 2})$, that $p|C'$ is stably dominated 
 by $f(a,{b^*}_1)$. Hence, since $f(a,{b^*}_1)$ depends only on $ab_J$ and
$b_\omega,b_{\omega 2}$ have the same type over $b_J$, $p|{b^*}_1$ is stably dominated by $f(a,{b^*}_1)$, as required.

For the final assertion, we shall apply Lemma~\ref{stab0.6}. We have 
$\ind{b_{\omega 2}}{d_J}{b_J}$, as in the proof of (A1). 
Also, as $\Pi^*(p,q,F)$,  we have 
$\dcl(b_1b_{\omega 2}d_1)=\dcl(b_1b_{\omega 2}d_{\omega 2})$. Hence 
$d_{\omega 2}(\lambda)\in \dcl(b_1b_{\omega 2}d_1)\subset \dcl(b_Jb_{\omega 2}d_J)$ for each $\lambda$. 
It follows from \ref{stab0.6} that there is a finite tuple 
$f[\lambda]\in \hat{f}_\lambda(a,{b^*}_1)= 
\Cb(b_{\omega 2}d_{\omega 2}(\lambda)/b_Jd_J;b_J)$ 
such that $\hat{f}_\lambda(a,{b^*}_1)=\dcl(b_Jf[\lambda])$. 

There is a $\emptyset$-definable function $f_\lambda$
so that $f[\lambda]=f_\lambda(a,b_J)$. Thus, the function $f$ defined 
earlier can be viewed as having components
$f_\lambda$ for $\lambda\in \Lambda$.  We now apply Lemma~\ref{stab4.1}(ii)
with $b_{\omega 2}$ as $b_i$ and $f_\lambda(a,b_J)\in d_{\omega 2}$ as $d_i$. 
This yields that $\dcl(b_Jb_{J'}f_\lambda(a,b_J))=
\dcl(b_Jb_{J'}f_\lambda(a,b_{J'}))$.
Hence, $\dcl({b^*}_1{b^*}_2f_\lambda(a,b_J))= \dcl({b^*}_1{b^*}_2f_\lambda(a,b_{J'}))$.
\qed

\begin{theorem}\label{stab7} \index{invariant type!and stable domination} 
Let $p,q$ be $\Aut({\cal U}/C)$-invariant   $*$-types.  
Assume that  whenever $b \models q | C$, the type   $p|Cb$ is
stably dominated.  Then $p|C$ is stably dominated.
\end{theorem}

{\em Proof.}  The strategy is to show that $p|Cb$ is stably dominated via a $C$-definable function $g$, and then to apply Lemma~\ref{stab1.2}. 
As a first approximation of the proof, let $a\models p|Cb$, and put $d=f_b(a)$, where $f_b$ is a $Cb$-definable function witnessing that $p|Cb$ is stably dominated. Define an equivalence relation
$E$ on $\tp(bd/C)$, putting $(b,d)E(b',d')$ if and only if, for $a\models p|Cbb'$, $f_b(a)=d \Leftrightarrow f_{b'}(a)=d'$.
Then put $g(a)=\lceil (b,d)/E \rceil$. 
The details are rather more intricate, and seem to require the previous  lemma. The technical problems are:
reducing to finite tuples which lie in $\St_C$; identifying the right equivalence relation, and proof of transitivity; and proof that  a code for the equivalence class lies in $\St_C$.

 Note at the outset that the symmetry Lemma \ref{stabrem3} is true of $p,q$:  if 
 $b_0 \models q$, then $p|Cb_0$ is stably dominated, so \ref{stabrem3} can be applied
 over $b_0$.  Let $a \models p|Cb_0, b \models q | Cb_0a$.  Then $a \models 
 p | Cb_0b$.  In particular, $a \models p | Cb$ and $b \models q | Ca$; that is, symmetry holds
 already over $C$. Likewise if $b_\lambda$ is a subtuple of $b$ and $q_\lambda$ is the restriction of $q$ to the corresponding variables, then $a\models p|Cb_\lambda$ and $b_\lambda\models q_\lambda|Ca$, so again symmetry holds.
 
The goal is to obtain a $C$-definable $*$-function $g$ such that
for $b\models q|C$, the type $p|Cb$ is stably dominated via $g$. We then apply
Lemma~\ref{stab1.2}.

We work over $C$, so omit all reference to $C$. There is a $\emptyset$-definable $*$-function $F$
 such that, if $b\models q|\emptyset$, then the type
$p|b$ is stably dominated via $F(x,b)$.
We apply Lemma~\ref{stab6} twice, first to obtain a function $f'$
satisfying \ref{stab6}(i) and (ii), and then, to use $f'$ in place of $F$
to obtain a new function $f$ satisfying also the last  assertion of \ref{stab6}. That is,
we obtain an $\Aut({\cal U}/\emptyset)$-invariant $*$-type $q^{\omega}$ (which we shall
continue to write as $q$), and $f=(f_\lambda)_{\lambda\in \Lambda}$, $f_\lambda$ 
a definable function, such that $p|b$ is stably dominated via $f(x,b)$ for any 
$b\models q|\emptyset$, and  $\Pi^*(p,q,f_\lambda)$ for each $\lambda \in \Lambda$. It 
follows that for each $\lambda$ there is a type $q_\lambda$ (the restriction of 
the $*$-type $q$ to a finite subset of the variables) such that $f_\lambda$ 
depends only on the $p$-variables and $q_\lambda$-variables, and such that the 
following holds: if $b_{\lambda}\models q_\lambda|\emptyset$,
$b_{\lambda}'\models q_\lambda|b_{\lambda}$, and 
$a\models p|b_{\lambda}b'_{\lambda}$, then 
$$
   \dcl(b_{\lambda}b'_{\lambda}f_\lambda(a,b_{\lambda}))=
                     \dcl(b_{\lambda}b'_{\lambda} f_\lambda(a,b'_{\lambda})).
$$
This means that there is a $\emptyset$-definable function $h_\lambda$ with inverse 
${h_\lambda}^{-1}$ such  that if 
$b_\lambda \models q_\lambda|\emptyset$, $b'_\lambda \models q_\lambda|b_\lambda$, 
$a \models p| b_\lambda b'_\lambda$,  then
$$ 
  (*) \ \ \ \   h_\lambda(b_\lambda,b'_\lambda,f_\lambda(a,b_\lambda))=  
       f_\lambda(a,b'_\lambda) \mbox{~and~}
   {h_\lambda}^{-1}(b_\lambda,b'_\lambda,f_\lambda(a,b'_\lambda))=
                                                       f_\lambda(a,b_\lambda).
$$

Let $b \models q|\emptyset$, let $b_\lambda$ be its restriction to the 
$q_\lambda$-variables, and suppose $a \models p | b$. Put 
$d_\lambda:=f_\lambda(a,b_\lambda)$, and
$d:=(d_\lambda)_{\lambda\in \Lambda}$.  
Let $Q$ be the set of realisations of $q|\emptyset$, $Q_\lambda$ the set of realisations of 
$q_\lambda|\emptyset$, $R$ be the set of realisations of $\tp(bd/\emptyset)$, and 
$R_\lambda$ be the set of realisations of $\tp(b_\lambda d_\lambda/\emptyset)$.

If $b'_\lambda \models q_\lambda | b_\lambda$
and $a \models p| b_\lambda b'_\lambda$, then 
$\ind{f_\lambda(a,b'_\lambda)}{b_\lambda}{b'_\lambda}$, since
$\tp(f_\lambda(a,b'_\lambda) / b'_\lambda b_\lambda)$ extends to 
an $\Aut(\U/\emptyset)$-invariant
type.   Thus since $R_\lambda$ is a complete type,  we can add to $(*)$ the following:

$(*1)$  if $(b_\lambda,d_\lambda) \in R_\lambda$, $(b_1)_\lambda \models q_\lambda | (b_\lambda,d_\lambda)$ then
$\ind{h_\lambda(b_\lambda,(b_1)_\lambda,d_\lambda)}{b_\lambda}{(b_1)_\lambda}$

Since $p$ is $\Aut(\U/b)$-invariant, $\tp(a/b)$ implies $\tp(a/\acl(b))$, so
$\tp(d/b)$ implies $\tp(d/\acl(b))$.  Thus, for any $X$,  $\ind{d}{X}{b}$ if and only if 
$\ind{d}{X}{\acl(b)}$.  The same holds with $b,d$ replaced by 
$b_\lambda,d_\lambda$.

The domain of $h_\lambda(b_\lambda,b'_\lambda,x)$ contains
$\{d_\lambda:b_\lambda d_\lambda\in R_\lambda \mbox{~and~} 
\ind{d_\lambda}{b'_\lambda}{b_\lambda}\}$. Also, as $\tp(a/b_\lambda b'_\lambda)$ has an
$\Aut(\U/b'_\lambda)$-invariant extension, 
$\ind{h_\lambda(b_\lambda,b'_\lambda,d)}{b_\lambda}{b'_\lambda}$, so we have:

 $(**)$
the image of $h_\lambda(b_\lambda,b'_\lambda,x)$
(and the domain of $h_\lambda^{-1}(b_\lambda,b'_\lambda,y)$) contains
$$\{y: b'_\lambda y\in R_\lambda \mbox{~and~} 
\ind{y}{b_\lambda}{b'_\lambda}\}.$$

For each $\lambda \in \Lambda$, define an equivalence relation 
$E_\lambda$ on $R_\lambda$:  
$(b_\lambda,d_\lambda) E_\lambda (b'_\lambda,d'_\lambda)$ if and only if, whenever 
$b''_\lambda \models q_\lambda | b_\lambda d_\lambda b'_\lambda d'_\lambda$, 
we have
$h_\lambda(b_\lambda,b''_\lambda,d_\lambda) = 
h_\lambda(b'_\lambda,b''_\lambda,d'_\lambda)$.
Until Claim 8 below, we are concerned only with 
$q_\lambda,b_\lambda,d_\lambda$, so for ease of notation we drop the 
subscript $\lambda$  until then
(except for $E_\lambda, R_\lambda, f_\lambda, h_\lambda$).

\medskip

{\em Claim 1.}  Assume $bd,b'd' \in R_\lambda,  b_1 \models q_\lambda | bd$ and  
$  b_1 \models q_\lambda | b'd'$.  Then 
$(b,d) E_\lambda (b',d')$ if and only if 
$h_\lambda(b,b_1,d) = h_\lambda(b',b_1,d')$.

{\em Proof of Claim.}

 Let $d_1 = h_\lambda(b,b_1,d) $, $d'_1= h_\lambda(b',b_1,d')$.  
 Let $b_2 \models q_\lambda | bb'b_1dd'$.  Let 
$d_2= h_\lambda(b_1,b_2,d_1)$, $d'_2= h_\lambda(b_1,b_2,d'_1)$.   
So $d_1=  h_\lambda^{-1}(b_1,b_2,d_2),
d'_1=  h_\lambda^{-1}(b_1,b_2,d'_2)$.  Thus, $d_1=d'_1$ if and only 
if $d_2=d'_2$.

Thus, to prove the claim, we must show   
$(b,d) E_\lambda (b',d')$ if and only if   $d_2=d_2'$. To show this, 
by the definition of $E_\lambda$, is suffices to prove that
$d_2 = h_\lambda(b,b_2,d)$ and $ d'_2= h_\lambda(b',b_2,d')$.

We shall prove this last assertion for $d_2$.  Let 
$a^* \models p | bb_1b_2$, and let $d^*=f_\lambda(a,b)$.  Then 
$d^* \in \St_{b}$, and $\ind{d^*}{b_1b_2}{b}$.  Thus 
$\tp(d^*bb_1b_2/\emptyset) = \tp(dbb_1b_2/\emptyset)$.  So there exists $a' \models p|bb_1b_2$
with $f_\lambda(a',b)=d$.  By $(*)$, 
$$f_\lambda(a',b_1) = h_\lambda(b,b_1,f_\lambda(a',b)) =h_\lambda(b,b_1,d)= d_1,$$
and then   
$$h_\lambda(b,b_2,d)=h_\lambda(b,b_2,f_\lambda(a',b))=f_\lambda(a',b_2) = 
h_\lambda(b_1,b_2,f_\lambda(a',b_1))=h_\lambda(b_1,b_2,d_1) = d_2.$$
This and a similar argument for $d_2'$ yield the claim.

\medskip

{\em Claim 2.}  $E_\lambda$ is a definable equivalence relation 
(i.e. the intersection with $R_\lambda^2$ of one).

{\em Proof of Claim.}  By Claim 1, if $bd, b'd'\in R_\lambda$, then $(b,d) E_\lambda (b',d')$ if and only if
there is $b_1$ such that 
$$  
    b_1\models q_\lambda|bd,~~ b_1 \models q_\lambda|b'd' \mbox{~and~} 
           h_\lambda(b,b_1,d)= h_\lambda(b',b_1,d').
$$
Also, $\neg(b,d)  E_\lambda (b',d')$ if and only if there is $b_1$ such that
$$
   b_1 \models q_\lambda|bd, b_1 \models q_\lambda|b'd' \mbox{~and~} 
         h_\lambda(b,b_1,d) \neq h_\lambda(b',b_1,d').
$$
Thus both $E_\lambda$ and $R_\lambda^2 \setminus E_\lambda$ 
are $\infty$-definable over $\emptyset$.  The claim follows by compactness. (In more detail, 
by compactness, there is a definable relation $E_\lambda'$ such that $E_\lambda'\cap R_\lambda^2=E_\lambda$; and in general, if an equivalence relation on a complete type is induced by some definable relation, then it is induced by some definable {\em equivalence} relation on some definable set containing the complete type.)

\medskip

{\em Claim 3.} Let $(b,d) \in R_\lambda$, $b' \in Q_\lambda$, and 
suppose $\ind{d}{b'}{b}$.  
Then there exists $d'$ with $(b,d) E_\lambda (b',d')$.

{\em Proof of Claim.}  Let $b_1 \models q_\lambda | bb'd$, 
$d_1 = h_\lambda(b,b_1,d)$.   By (*1) we have
$\ind{d_1}{b}{b_1}$.  Also $\ind{b_1}{d}{bb'}$ by genericity,
so   $\ind{b'b_1}{d}{b}$ by transitivity.  Thus $\ind{b'}{d}{bb_1}$, so 
$\ind{b'}{d_1}{bb_1}$.  But then $\ind{b'}{d_1}{b_1}$
(by transitivity, since $\ind{d_1}{b}{b_1}$.)  So by $(**)$,
$d_1 $ is in the range of $h_\lambda(b',b_1,x)$. Hence there is
 $d' := h_\lambda^{-1}(b',b_1,h_\lambda(b,b_1,d))$. 
Now $b_1\models q_\lambda|bd$ and $b_1\models q_\lambda|b'd'$:
to see the latter, note that since $\ind{b'}{d_1}{b_1}$, 
if $a' \models p | b'b_1$ then $\tp(f_\lambda(a',b_1) / b'b_1) = \tp(d_1 / b'b_1)$.
Thus there exists $a \models p | b'b_1$ with $f_\lambda(a,b_1) = d_1$.  Now
  $d'=h_\lambda^{-1}(b',b_1, f_\lambda(a,b_1))=f_\lambda(a,b')$,
and $b'\models q_\lambda|b_1a$, so $b_1\models q_\lambda|b'd'$.
  Hence, by Claim 1, 
  $(b,d) E_\lambda (b',d')$.  

\medskip

{\em Claim 4.}  Let $(b,d),(b,d') \in R_\lambda$, with $(b,d) E_\lambda (b,d')$.  
Then $d=d'$.

{\em Proof of Claim.}  Let $b' \models q_\lambda | bdd'$.  Then by Claim 1,
$h_\lambda(b,b',d)=h_\lambda(b,b',d')=d^*$, say.  So 
$d=h_\lambda^{-1}(b,b',d^*)=d'$.

\medskip

Now let $(b_i:   i < |T|^+)$ be an indiscernible sequence, with 
$b_i \models q_\lambda |b_{<i}$ for each $i$. 

\medskip

{\em Claim 5.}  Let $(b^*,d^*) \in R_\lambda$. Then there is $i<|T|^+$ and
some $d_i$ such that $(b_i,d_i)\in R_\lambda$ and $(b^*,d^*)E_\lambda (b_i,d_i)$.

{\em Proof of Claim.}
  Let $b  \models q_\lambda | b_{<|T|^+}b^*d^*$.  
Then $\ind{b}{d^*}{b^*}$.  By Claim 3, $(b^*,d^*)E_\lambda(b,d)$ for 
some $d$. Now $\{\St_{b}(b_i): i < |T|^+ \}$ forms a Morley sequence in 
$\St_{b}$, by Proposition~\ref{stab0.7} (i) and (ii).
So $\ind{b_i}{d}{b}$ for some $i$, by properties of preweight (see
 Chapter 2).  By Claim 3 again, 
there exists $d_i$ with $(b_i,d_i) E_\lambda (b,d)$.

\medskip

{\em Claim 6.}  $R_\lambda/E_\lambda \subset \St_{\emptyset}$.  

{\em Proof of Claim.} 
By Claim 5, any element of $R_\lambda/E_\lambda$ has 
the form  $(b_i,d)/E_\lambda \in \St_{b_i}$, for some $i$. Hence,
$R_\lambda/E_\lambda \subset 
\St_{\{b_i: i < |T|^+\}}$. Hence there is a $\emptyset$-definable subset of $\St_{\{b_i: i < |T|^+\}}$  which contains $R_\lambda/E_\lambda$. As it must be stable and stably embedded,  the claim follows.

\medskip

{\em Claim 7.}  Let $b \models q_\lambda | \emptyset$.  Let $a \models p |b$.  
Let $d = f_\lambda(a,b)$.  
Then $(b,d)/E_\lambda \in \dcl(a)$. 
 
{\em Proof of Claim.}  Let $(b',d')$ be $a$-conjugate to $(b,d)$.  
Then $d' = f_\lambda(a,b')$.  Let
 $b'' \models q_\lambda | bb'a$.   Then 
$h_\lambda(b,b'',d) = f(a,b'') = h_\lambda(b',b'',d')$, so 
$(b,d) E_\lambda (b',d')$.

\medskip

Now let $a, b_\lambda,d_\lambda$ be as in $a,b,d$ of Claim 7 (so we revert to subscript notation), and put
$g_\lambda(a):=\lceil (b_\lambda,d_\lambda)/E_\lambda\rceil$. The following claim completes 
the proof 
of the proposition.

\medskip

{\em Claim 8.}  $p|\emptyset$ is stably dominated  via  $g=(g_\lambda: 
\lambda \in \Lambda)$.

{\em Proof of Claim.}  As $a\models p|b$, symmetry yields  $b \models q |a$.  By Claim 4, 
$f_\lambda(a,b_\lambda)$ is the unique $d_\lambda$ such that  
$g_\lambda(a) = \lceil (b_\lambda,d_\lambda)/E\rceil$.  Thus
$f_\lambda(a,b_\lambda) \in \dcl(b_\lambda,g_\lambda(a))$.   As $p|b$ is stably dominated 
via $f$, it is also stably dominated via $g$.  By Lemma~\ref{stab1.2}, $p|\emptyset$ is stably dominated.
 \qed

\begin{corollary} \label{st-dom-eq-0} Let $T$ be a theory,  $\U$ be a universal domain for $T$.  
 Let $\G$ be a 0-definable stably embedded set with  a 0-definable
linear ordering, such that every type $\tp(\g/E)$ (with $\g \in \G, E \subset \G)$ extends to an $E$ - invariant type.  
 Let $f$ be a 0-definable function.  

(i)  Assume 
$\tp(a/B)$ extends to an $\Aut(\U/B)$-invariant type.  Then so does $\tp(f(a)/B)$.  

(ii)  Assume   $\tp(a/ B,\G(Ba))$ is stably dominated.
  Then so is $\tp(f(a)/ B,\G(B,f(a)))$.

\end{corollary}

{\em Proof}    
(i)  Let $\tp(a'/\U)$ be an $\Aut(\U/B)$-invariant
extension of $\tp(a/B)$.  Then clearly $\tp(f(a')/\U)$ is an $\Aut(\U/B)$-invariant
extension of $\tp(f(a)/B)$.  

(ii) 
Let $\g$ enumerate $\G(Ba)$, and $\g'$ enumerate $\G(B, f(a))$.
    Since $\tp(a/B,\g)$ is stably dominated, by Proposition~\ref{bits} so is $\tp(f(a)/ B, \g)$.
But $\tp(f(a) / B,\g') \vdash \tp(f(a) / B, \G)$ by stable embeddedness,
and in particular $\tp(f(a) / B,\g') \vdash \tp(f(a)/ B, \g)$.  Let $q$
be an $\Aut(\U/ B,\g')$-invariant extension of $\tp(\g / B,\g')$; this exists by stble embeddedness of $\G$.
Then $\tp(\g / B(f(a))) = q | B(f(a))$, since both equal the unique extension of
$\tp(\g / B(\g'))$ to $B(f(a))$.  By Theorem~\ref{stab7}, $\tp(f(a)/B,\g')$
is stably dominated.   

\qed

Let $T$ be a theory, with universal domain $\U$.  Let   ${\mathcal S}$ a collection
 of sorts, and among them  
let $\G$ be a stably embedded sort with a 0-definable
linear ordering.   

\begin{definition} \label{metastable} \index{metastable}  We say that  ${\mathcal S}$ is
 {\em metastable}  over $\G$ if
for any finite product $D$ of sorts of $S$, and any small $C \leq \U$,   we have:

(i)  if $C = \acl(C)$, then for any $a \in D(\U)$, $\tp(a/C)$ extends to an
$\Aut(\U/C)$-invariant type;

(ii)  for some small $B$ with $C \subseteq B \leq \U$, for any $a \in D(\U)$, $\tp(a / B, \G(Ba))$ is stably dominated.

If  ${\mathcal S}$ consists of all sorts, we say $T$ is metastable over $\G$.
\end{definition} 

In Theorem~\ref{fulldom} we will show that algebraically closed valued fields are metastable.

\begin{corollary} \label{metastable-eq} 
Assume every sort of $T$ lies in the definable closure of ${\mathcal S}$, and
  ${\mathcal S}$ is {\em metastable} over $\G$.  Then $T$ is  metastable over $\G$.




\end{corollary}

{\em Proof.}   By assumption, any  finite sequence of elements   can be written as $f(a)$ for some 0-definable function
$f$ and some $a \in D$, where $D$ is a product of sorts in $ {\mathcal S}$.   The corollary is thus 
immediate from Corollary~\ref{st-dom-eq-0}.  \qed

\chapter{A combinatorial lemma}\label{combinatorial}
We give in this chapter three versions of a combinatorial lemma which may be of
 independent interest. 
The proof of the first version uses Neumann's Lemma (stated below), and the proof of the 
second 
uses basic combinatorial facts about stability. The second version has content in any
 model and gives an explicit bound 
of $n/2$.  Thirdly in Lemma~\ref{PQRs} we lift the result from finite sets to stable sets.
Our application of
Lemma~\ref{PQR} will be 
in the next chapter, for the existence of `strong' codes for germs of functions;
$R(a,b)$ will be the set of points on which two functions with the same germ disagree.

\begin{lemma}\label{PQR}
Let $M$ be an $\omega$-saturated structure with partial 1-types $S,Q$ over $\emptyset$, and with a 
$\emptyset$-definable relation $R_0$ inducing a relation $R\subseteq Q^2 \times S$.
Assume

(i) for any  $a,b\in Q$, $R(a,b):=\{z\in S:R(a,b,z)\}$ is a finite 
subset of $S$, and $R(a,b)=R(b,a)$;

(ii) for all  $a,b,c\in Q$, $R(a,c)\subseteq R(a,b) \cup R(b,c)$

(iii) for all $a\in Q$, $R(a,a)=\emptyset$.

\noindent
Then for all  $a,b\in Q$ we have $R(a,b) \subseteq \acl(a) \cup \acl(b)$.
\end{lemma}

{\em Proof.}  
By the saturation assumption, $R(a,b)$ is a finite 
subset of $S$ of bounded size.  Hence (i)-(iii) remain true in any elementary extension,
while the conclusion clearly descends.  So 
we may assume $M$ is $\omega$-homogeneous.  

We may also assume $\acl(\emptyset) \cap S=\emptyset$, 
and that $Q\not\subseteq \acl(\emptyset)$. For the lemma is trivial if $Q\subseteq \acl(\emptyset)$.
 If  $\acl(\emptyset) \cap S \neq \emptyset$, we may replace $S$ by $S':=S\setminus \acl(\emptyset)$ 
and $R$ by $R':=R\cap(Q^2\times S')$.
Then (i), (ii) still hold, and the conclusion for $R'$ implies that for $R$.

\medskip

{\em Claim 1.} For any finite  $E,C\subset M$ and any finite tuple $d \in M$, there
 is an $\Aut(M/C)$-translate $d'$ of $d$ with $\acl(d'C) \cap \acl(EC)=\acl(C)$.

{\em Proof of Claim.} This follows by compactness from Neumann's Lemma \cite[Lemma 2.3]{neumann} \index{Neumann's Lemma}, 
which states that if $G$ is  a permutation group on an infinite set $X$
 with no finite orbits, then any finite set $F\subset X$ has some (and hence
 infinitely many) disjoint translates.\index{Neumann's Lemma}

\medskip

\def\QQ{{\mathbf Q}_2}
Now let $\QQ:=\{(a,b)\in Q^2: \acl(a) \cap \acl(b)=\acl(\emptyset)\}.$

\medskip

{\em Claim 2.} It suffices to prove the lemma for $(a,b)\in \QQ$.

{\em Proof of Claim.} Suppose the lemma holds for elements of $\QQ$, and 
let $a,b\in Q$ be distinct. By Claim 1, there is $c\in Q$ 
with $\acl(c) \cap \acl(a,b)=\acl(\emptyset)$. Hence $(a,c),(b,c)\in \QQ$, 
so $R(a,c)\subseteq \acl(a) \cup \acl(c)$ and $R(b,c) \subseteq \acl(b) \cup \acl(c)$. 
Thus, $$R(a,b) \subseteq R(a,c) \cup R(b,c) \subseteq \acl(a) \cup \acl(b) \cup \acl(c).$$ Also,
$R(a,b) \subseteq \acl(a,b)$, so $R(a,b) \cap \acl(c) =\acl(\emptyset)$.
 Thus, $R(a,b) \subseteq \acl(a) \cup \acl(b)$, as required.

\medskip

Since the size of $R(a,b)$ is bounded, we can proceed by reverse induction. So suppose that
$R(a,b) \subseteq \acl(a) \cup \acl(b)$ whenever $(a,b)\in \QQ$ and $|R(a,b)|>n$. We must
 prove $R(a,b)\subseteq \acl(a) \cup \acl(b)$ when $(a,b)\in \QQ$ and $|R(a,b)|=n$. Suppose for
 a contradiction that this is false for some $(a,b) \in \QQ$.

Put $m:=|\acl(b) \cap R(a,b)|$ and $m':=|\acl(a) \cap R(a,b)|$. 
We will
show in Claim 3 that 
$m\geq n/2$. A symmetrical argument will give $m'\geq n/2$. Since $(a,b) \in \QQ$ and $S\cap \acl(\emptyset)=\emptyset$, we have $\acl(a) \cap \acl(b)\cap S=\emptyset$. 
 By counting it will follow  from Claim 3 that 
$R(a,b) \subseteq \acl(a) \cup \acl(b)$,  contradicting our assumption on $a,b$.

\medskip

{\em Claim 3.} $m\geq n/2$.

{\em Proof of Claim.} By Claim 1 applied over $b$, there is $a'\equiv_b a$ 
with $\acl(a'b) \cap \acl(ab)=\acl(b)$. Then $$\acl(a) \cap \acl(a') \subseteq \acl(a)
 \cap\acl(ab) \cap \acl(a'b) \subseteq
\acl(a) \cap \acl(b)=\acl(\emptyset),$$ so $(a,a') \in \QQ$.

If $|R(a,a')|>n$, then by induction we have $R(a,a')\subseteq \acl(a) \cup \acl(a')$. 
Then 
$$R(a,b) \subseteq R(a,a') \cup R(a',b)\subseteq\acl(a)\cup \acl(a') \cup R(a',b).$$
But $R(a,b) \subseteq \acl(ab)$. So 
$$R(a,b) \subseteq (\acl(a) \cup \acl(a'b)) \cap \acl(ab)\subseteq \acl(a) \cup \acl(b),$$
again contradicting the assumption on $a,b$.

Thus, we may suppose  $|R(a,a')|\leq n$. However, by (ii), $R(a,a')$ contains the symmetric difference of 
$R(a,b), R(a',b)$. We have $|R(a,b)|=|R(a',b)|=n$, since $(a,b),(a',b)$ are conjugate. Also,
$$R(a,b) \cap R(a',b)\subseteq R(a,b) \cap (\acl(ab) \cap \acl(a'b)) \subseteq
 R(a,b) \cap \acl(b),$$
so $|R(a,b) \cap R(a',b)|\leq m$. Thus, the symmetric difference of
$R(a,b), R(a',b)$ has size at least $2(n-m)$. Hence $n\geq 2(n-m)$, so $m \geq n/2$, as required. \qed

\bigskip

In the second version below, if $T\subseteq Q\times S$ and $a\in Q$ then
$T(a):=\{y\in S: (a,y)\in T\}$. Likewise, if $a,b\in Q$ then 
$R(a,b):=\{y\in S: R(a,b,y)\}$.

\begin{lemma}\label{comb2}
Let $n \in {\mathbb N}$. Then there is an integer $n_0$ depending on $n$ with the following 
property.
Let $Q,S$ be disjoint sets, and $R \subseteq Q^2\times S$. Suppose that for all
 $a,b,c\in Q$, 
we 
have $R(a,b)=R(b,a) \subseteq R(a,c) \cup R(b,c)$, and  $|R(a,b)|\leq n$ and $R(a,a)=\emptyset$.
Then there is $S_0 \subseteq S$ with $|S_0|\leq n_0$ and
$T\subseteq Q \times S$, with   $|T(a)| \leq n/2$ for all $a \in Q$
such that, letting $S' = S \setminus S_0$, $R' = R \meet (Q^2 \times {S'})$,
$$T(a) \triangle T(b) \subseteq R'(a,b) \subseteq T(a) \union T(b)$$
 for all $a,b \in Q$.
Moreover, $S_0$ and $T$ are defined in the structure $(Q,S;R)$ by parameter-free 
first-order formulas,
 depending only on $n$. 
\end{lemma}
 
{\em Proof.}  We shall prove the result in a fixed (possibly infinite) structure $M=(Q,S;R)$. 
Compactness then implies that the
 formulas and bounds are uniform. 

\medskip

{\em Claim 1.}  Let $N=  n+2$.   Assume there are $a_i,b_i,c_i \in M$ for $i=0,\ldots,N$
 such that for all $j=0,\ldots, N$, $R(a_i,b_j,c_i)$ holds
whenever $i<j$.   Then $R(a_i,b_j,c_i)$ holds for some $i,j$ with $i\geq j$.

{\em Proof of Claim.}  Suppose otherwise. 
We first argue that $c_0,\ldots,c_N$ are 
distinct. To see this,
suppose that $c_i=c_j$ where $i<j$. Then $c_j\in R(a_i,b_j)$, and
 $c_j\not\in R(a_i,b_i)$. As $R(a_i,b_j) \subseteq R(a_i,b_i)
\cup R(b_i,b_j)$, it follows that $c_j\in R(b_i,b_j)$. Now 
$R(b_i,b_j) \subseteq R(a_j,b_i) \cup R(a_j,b_j)$, so
$c_j\in R(a_j,b_i)$ or $c_j\in R(a_j,b_j)$, and each of these contradicts our assumption. 

Whenever $j'<i<j$ we have
 $c_i \in R(b_j',b_j)$: indeed,
$c_i  \in R(a_i,b_j) \subseteq R(a_i,b_{j'}) \union R(b_j,b_{j'})$ but 
$c_i \notin R(a_i,b_{j'})$.  Thus $c_i \in R(b_0,b_N)$
for all $i$ with $0<i<N$.  As the $c_i$ are distinct, this contradicts the 
assumption that $|R(b_0,b_N)| \leq n$.

\medskip
 
Let $x=(x_1,x_2) $ be a variable ranging over $Q \times S$, let $y$ range over $Q$, and
let $\phi(x,y) =  R(x_1,y,x_2) $.  Then Claim 1 states that $\phi$  is a stable
formula.  Let $q_1(y),\ldots,q_m(y)$ be the  $\phi$-types of maximal $\{\phi\}$-rank 
in the variable $y$. 
For $j=1,\ldots,m$, put
$T_j =  \{(x_1,x_2) \in Q \times S: (d_{q_j}y)\phi(x_1x_2,y) \}$.  Let 
$T' = T_1\cup \ldots \cup T_m$. As $\phi$ is stable, the $q_i$ are definable $\phi$-types (see 
Remark~\ref{localsymmetry}), so the $T_i$ and $T'$ are definable. Observe that as $|T_j(a)|\leq n$ for all $a$, 
$|T'(x)|$ is bounded as $x$ varies through $Q$.

\medskip

{\em Claim 2.} If $a,b\in Q$ then  $R(a,b) \subseteq T'(a) \union T'(b)$.

{\em Proof.} It suffices to show that for each $j$, if $c \in R(a,b)$ then
$c \in T_j(a) \union T_j(b)$.  To see this, choose $d \models q_j | \{a,b,c\}$.
Since $R(a,b)\subseteq R(a,d)\cup R(b,d)$, $c\in R(a,d) \cup R(b,d)$. We suppose $c\in R(a,d)$.
Then $\phi(ac,d)$ holds, so by the choice of $d$, $(a,c)\in T_j$, so $c\in T_j(a)$. Similarly, 
if $c\in R(b,d)$ then $c\in T_j(b)$.

\medskip

{\em Claim 3.}  If $a,b\in Q$ then the symmetric difference $T'(a) \triangle T'(b) $ is a subset of $R(a,b)$. 

{\em Proof.} Let $c \in T'(a) \triangle T'(b)$, 
say with $c \in T_j(a) \setminus T'(b)$.  Let $d \models q_j | \{a,b,c\}$.
Then $\phi(ac,d)$ (as $c\in T_j(a)$) and $\neg\phi(bc,d)$ (as $c\not\in T_j(b)$). So
$c \in R(a,d)$ and 
$c \notin R(b,d)$.   But $R(a,d) \subseteq R(b,d) \union R(a,b)$.  So $c \in R(a,b)$.

\medskip

We now choose $S_0$ to be a finite $\emptyset$-definable subset of $S$, chosen so as to minimise
$\ell:=\Max\{|T'(x)\setminus S_0|:x\in Q\}$. 
Define $T:=T' \cap (Q\times (S\setminus S_0))$, so that $T(a) = T'(a) \setminus S_0$ for
 each $a$. Then $T$ is $\emptyset$-definable. 
By Claims 2 and 3,
$$T(a) \triangle T(b) \subseteq R'(a,b) \subseteq T(a) \union T(b)$$
 for all $a,b \in Q$.
To complete the proof,  we will show that 
$\ell \leq n/2$.
 
Let 
$$H_0 = \{w \subseteq S:  |w| \leq \ell, \mbox {~ for any~ } a'\in Q  \mbox{~with~ } |T(a')|=\ell,
 w \cap T(a') \neq \emptyset \}.$$
Let $H$ be the set of minimal elements of $H_0$ (under inclusion).  
Then $H$ is finite.  For otherwise, there exists   $H' \subseteq H$, $|H'| > \ell$,  
forming a $\Delta$-system,
that is, there is some $w_0 \subseteq S$ such that any two distinct elements of $H'$ have 
intersection $w_0$.
Any $T(a')$ with $|T(a')|=\ell$ must meet any element
of $H'$ nontrivially; but   it can  meet at most $ m$  of the 
disjoint sets $w \setminus w_0$;
so it must meet $w_0$ nontrivially.    But then $w_0 \in H_0$, contradicting the 
minimality in the definition of $H$. In particular, the $\emptyset$-definable set $ \bigcup H$
is finite (possibly empty).

There is $a\in Q$ such that $|T(a)| =\ell $ and $T(a)  \cap \bigcup H = \emptyset$; 
otherwise, $|T(a) \setminus \bigcup H | < \ell$
for all $a\in Q$, contradicting the minimality in the choice of $\ell$.  
Hence $T(a)\not\in H_0$, as otherwise $T(a)$ containes some member of $H$. So by definition of $H_0$
 there exists $a'\in Q$ such that $|T(a') | = \ell$ and $T(a) \cap T(a') = \emptyset$. 
 But then, by Claim 3, $R(a,a') \supseteq T(a) \union T(a')$, 
so $R(a,a') \geq 2\ell$.  It follows that $\ell\leq n/2$.  \qed

\begin{remark} \rm

1. In Lemma~\ref{comb2}, the formulas defining $T$ and $S_0$ can also be taken to be quantifier-free (in the structure $(Q,S;R)$), but with parameters. For $S_0$ this is immediate, as $S_0$ is finite, and for
$T_0$ it holds since each $T_j$ is quantifier-free definable (by a `majority rule'
definition for the average of the appropriate indiscernible sequence, see \cite{pil}.)

2. It is easy to deduce Lemma~\ref{PQR} from Lemma~\ref{comb2}, by compactness. 

\end{remark}

 \begin{lemma}\label{PQRs}
Let $M$ be an $\omega$-saturated structure,
$\Delta$ a finite set of stable formulas.    Let $Q$ be a partial type 
over $\emptyset$, $S'$ a $\emptyset$-definable set, and $X_0$ an 
$\emptyset$-definable relation   inducing a relation $X \subseteq Q^2 \times S'$.  

For $a,b \in Q$, let  $X(a,b):=\{z\in S':X(a,b,z)\}$.  Let $r \in {{\mathbb N}}$.  Assume:

(i) for any  $a,b\in Q$, $\rk_\Delta (X(a,b)) \leq r$.  
 
(ii)  for all  $a,b,c\in Q$, $X(a,c)\subseteq X(a,b) \cup X(b,c)$.     

(iii) for all $a,b \in Q$, $X(a,b)=X(b,a)$; and $X(a,a) = \emptyset$.   

\noindent
Then there exists a $\emptyset$-definable set $D_0$ inducing $D \subseteq Q \times S'$
such that 

(i')  for any $a \in Q$, $\rk_\Delta (D(a)) \leq  r$.  

(ii')  for  all  $a,b\in Q$ we have $X(a,b) \subseteq D(a) \cup D(b)$.
\end{lemma}

{\em Proof.}     We may assume $M$ is $\omega$-homogeneous.  
For $r=0$, the statement amounts to Lemma~\ref{PQR}.  
We proceed by induction on $r$.  
We may suppose there are $b,b'\in Q$ with $\rk_\Delta(X(b,b'))=r$.

Let $S$ be the set of all $\Delta$-types over $M$ of rank $r$, in the sort $S'$; we may
 regard $S$ as a set of canonical parameters of these types.
For any $b,b'\in Q$, among the $\Delta$-types consistent with $X(b,b')$, 
let $R(b,b')$ be the set of  $\Delta$-types over $M$ consistent with $X(b,b')$, 
and of  $\Delta$-rank $r$; so $R \subset Q^2\times S$, the setting of Lemma~\ref{comb2}.
  Then for some $n$, for all $b,b'$, we have
$|R(b,b')|\leq n$.  For any $b,b',b''\in Q$ it follows from (ii) that 
$R(b,b') \subseteq R(b,b'') \cup R(b',b'')$, and from (iii) that $R(a,b)=R(b,a)$
and $R(a,a) = \emptyset$.  

Let $S_0, T$ be as in  Lemma~\ref{comb2}.  By the last statement of this lemma,
$S_0$ is $\Aut(M)$-invariant, and $T(b)$ is $\Aut(M/b)$-invariant.  Each element
$q$ of the finite set $S_0 \union T(b)$ thus contains a formula $D^q_b$ defined over $b$, and
of $\Delta$-rank  $\rk_{\Delta}(q)=r$.  Taking the disjunction of these formulas,
we obtain, for each $b$ with $S_0\cup T(b)\neq \emptyset$,  a formula $D_b$ of $\Delta$-rank  $r$, such that $D_b$ lies in each
$q \in S_0 \union T(b)$.  
 Thus $D_b \vee D_{b'}$ lies in each $q \in S_0 \union T(b) \union T(b') \supseteq R(b,b')$. 
 Since $R(b,b')$ contains all $\Delta$-types consistent
with $X(b,b')$, it follows that $X'(b,b') := X(b,b') \setminus \{x:D_b(x) \vee D_{b'}(x)\}$ has
$\Delta$-rank $<r$.  By compactness, we can replace $\{D_b\}$ by a uniformly
definable family $\{D'_b: b \in Q\}$ with the same property. For each $b$, let
$D'(b):=\{x:D'_b(x) \mbox{~holds}\}$. 

Now by induction, there exists a uniformly definable  $D''(b)$ of $\Delta$-rank $<r$, 
such that  for  all  $a,b\in Q$ we have  $X'(a,b) \subseteq D''(a) \union D''(b)$.  

Let $D(a) = D'(a) \union D''(a)$; then clearly (i'), (ii') hold.   \qed

\chapter{Strong codes for germs}\label{strongcodes}

In this chapter we find variants in our context of another well-known result for stable theories: 
that given the germ of a function on a definable type,
one can define from that germ a particular function having the same germ. This phenomenon is important 
for example in group configuration arguments in stability theory: see for example Ch. 5, Definition 1.3 and the 
remarks following it in \cite{pil}. Versions of this were also important in
\cite{hhm}, especially in Section 3.3, where it is shown that in ACVF, the germ of a 
function on the generic type of a closed ball has a `strong code'
in the sense of the next definition. 
 The main result is Theorem~\ref{strongcode1}, which gives strong codes for functions on 
stably dominated types; the proof uses Lemma~\ref{PQR}. We also obtain strong codes in Proposition~\ref{strongcode2}, 
where stable domination is replaced by an assumption on the range of the function and the condition `(BS)'. The 
 chapter concludes with some applications of strong codes, such as the useful transitivity result Proposition~\ref{transitive}.
 
 In this section tuples are finite except where indicated, and types are in finitely many variables. 

\begin{definition} \rm Let $p$ be a $C$-definable type over $\U$. Let $\varphi(x,y,b)$ be a formula defining a function
 $f_b(x)$ whose domain contains all realisations of $p$.  The {\em germ of $f_b $ on $p$ }, \index{germ} or 
{\em $p$-germ of $f_b$}, is the equivalence class of $b$ under the equivalence relation $\sim$,
where $b\sim b'$ if the formula $f_b(x)=f_{b'}(x)$ is in $p$. Equivalently, $b\sim b'$ if and only if for 
any $a\models p|Cbb'$, $f_b(a)=f_{b'}(a)$. As $p$ is $C$-definable, $\sim$ is also $C$-definable, and the germ of 
$f_b$ on $p$ is a definable object. A code $e$ for the germ $b/\sim$ of $f_b$ on $p$ is 
{\em strong} over $C$ \index{code!strong} \index{strong code} if there is a $Ce$-definable function $g$ 
such that the formula $f_b(x)=g(x)$ 
is in $p$. Equivalently, the code $e$ is strong if for any $a\models p|Cb$, $f_b(a)\in \dcl(Cea)$.
\end{definition}

If $p$ is not definable, but is an $\Aut(\U/C)$-invariant type over $\U$, we still sometimes say that definable functions 
$f$ and $g$ have the same germ on $p$. This means that for $a\models p|C\lceil f\rceil\lceil g\rceil$, the formula 
$f(a)=g(a)$ lies in $p$. This gives an equivalence relation on any definable family of functions, but in general not a 
$C$-definable one.

\begin{lemma} \label{stepforstrongcode} \index{strong code!on a definable type}
Suppose that $p$ is a $C$-definable type over $\U$, and $f=f_b$ is a definable function on the set 
of realisations of $p|Cb$. Let $e$ be a code for the germ of $f$ on $p$. Suppose that for any
$b'\equiv_{Ce}b$, and any $a$ with $a\models p|Cb$ and $a\models p|Cb'$, we have
$f_b(a)=f_{b'}(a)$. Then $e$ is a strong code for the $p$-germ of $f$ (over $C$). 
\end{lemma}

{\em Proof.} There is a 
well-defined function $F$
on $p$, having the same $p$-germ as $f$, and  defined  by putting $F(a)=f_{b'}(a)$ for any $b'\equiv_{Ce} b$ 
with $a\models p|Cb'$. The function $F$ is $\Aut(\U/Ce)$-invariant, so there is a 
formula $\phi(x,y)$
over $Ce$ such that $\phi(a,F(a))$ for all $a\models p$, and so that $\phi$ determines 
 a function on realisations of $p$. Adjusting $\phi$, we may suppose it is the graph of a function whose domain contains all realisations of $p|Ce$. 
The function defined by $\phi$ is now $Ce$-definable, and has the same $p$-germ as $f$, 
as required.\qed

\bigskip

From the lemma we see easily that in a stable theory, strong codes exist. For in the notation of the lemma,
 suppose $a\models p|Cb$ and $a\models p|Cb'$. Pick $b''\equiv_{Ce} b$ with $b'' \dnf_{Ce} abb'$. 
Then $f_{b''}$ has the same germ on $p$ as $f_b$ and $f_{b'}$. By properties of non-forking,
 $a\models p|Cbb''$ and $a\models p|Cb'b''$, so $f_b(a)=f_{b''}(a)=f_{b'}(a)$.

\begin{theorem}\label{strongcode1} \index{strong code!on a stably dominated type}
Let $p_0$ be a stably dominated type over $C$ with an extension $p$ over 
$C':=\acl(C)$, and let $f$ be a definable function  whose domain contains the set of realisations of $p|\U$. 
Let $a\models p|\U$ and $C'':=\acl(C) \cap \dcl(Ca)$. 

(i) The code for the germ of $f$ on $p$ is strong  over $C''$.

(ii) Suppose that $f(a)\in \St_{Ca}$. Then the code 
for the $p$-germ of $f$ is in $\St_C$. 
\end{theorem}

{\em Proof.}  (i) 
We shall work over $C''$, since by
Remark~\ref{partacl}, $p|C''$ is stably dominated, and by Remark~\ref{C''}, $p$ is $C''$-definable. For convenience
suppose $C=C''$, and write $p$ for $p|C''$.

\medskip

{\em Claim 1.} We may assume that $p$ is a type consisting of elements in $\St_C$. 

{\em Proof of Claim.} For this,
 we may write the type $p$
as $p(x,y)$, so that if $p(a',a)$ then $a'$ is an enumeration of $\St_C(a)$ (so in general 
is infinite). Write 
$p_x$ for the restriction of $p$ to the $x$-variables.
Note that the notion of {\em  germ}
 of a function on a definable type with infinitely many variables also makes sense.

Write $f=f_b=f(b,x,y)$, where $b$ is the parameter defining $f$ (over $C$). Put 
$b\sim b'$ if $f_b$, $f_{b'}$
have the same $p$-germ.  Then $b\sim b'$ if and only if
$$\mbox{~for all~} c\models p_x|Cbb', \mbox{~for all~} d \mbox{~realising~}
 p(c,y)|C, f(b,c,d)=f(b',c,d).$$
For by stable domination, $cd\models p|Cbb'$ if and only if $c\models p_x|Cbb'$
and $p(c,d)$ holds.
There is a definable function $F_{bc}=F_{bc}(y)$ which agrees with $f(b,c,y)$ on realisations of
$p(c,y)$. Since $F_{bc}$ is $bc_0$-definable for some finite subtuple
 $c_0$ of $c$, we may suppose now that
$c$ (and so $x$) has finite length. Let $F_b$ be the function defined on realisations of $p_x|Cb$,
such that $F_b(x)=\lceil F_{bx}\rceil$. Now $b\sim b'$ if and only if 
the functions $F_b(x)$ and $F_{b'}(x)$
have the same $p_x$-germ. Furthermore, if the code of $b/\sim$, regarded as the germ of $F_b$ on $p_x$ is strong over $C$, then it is also strong over $C$ when regarded as the germ of $f_b$ on $p$. This yields the claim.

\medskip

Write $f=f_b$, and let $\Sigma$ be the sort of the type $p$. We may suppose that $f_b$ is total, by giving it some formal value  
when undefined. 
Let $e$ be a code for the 
$p$-germ of $f$ and $E:= \dcl(Ce)$.
Let $Q$ be the set of conjugates of $b$ over $E$, and $P$ be the set of realisations of
$p|E$. 
For $b,b'\models Q$, let
$$X_0(b,b'):=\{x\in \Sigma:    
f_b(x)\neq f_{b'}(x)\}.$$
Clearly:  $X_0(b,b') \subseteq X_0(b,b'') \cup X_0(b',b'')$ for any $b,b',b''\in Q$.

As $\St_C$ is  stably embedded, there is a finite set $\Delta$ of formulas $\phi(x,y)$
such that for any $b,b'\in Q$ the set $X_0(b,b')$ is equal to $\phi(x,d)$ 
for some tuple $d=d(b,b')$ from $\St_C$ and some $\phi(x,y) \in \Delta$.

Let $\rho$ be the $\Delta$-rank of $p$; there exists in $p$ a formula $D(x)$ over $C$
equivalent to a $\Delta$-formula of $\Delta$-rank $\rho$ and multiplicity $1$.
For any $\Delta$-formula $\psi(x)$, we have $\psi(x) \in p | \U$ iff 
$\rk_{\Delta}(\psi \wedge D) = \rho$. 

Put $$X(b,b')=X_0(b,b')\cap D$$
Observe that if $a\models p|Cbb'$ then $a\not\in X(b,b')$, since $f_b$ and 
$f_{b'}$ have the same $p$-germ.   Hence $\rk_\Delta(X(b,b'))< \rho$.

By Lemma ~\ref{PQRs}, uniformly in $b \in Q$  there is a
 definable set $X(b)$ of $\Delta$-rank $<\rho$ such that for all $b,b'\in Q$,
 $X(b,b') \subseteq X(b) \union X(b')$.      
 
It follows that $f_b,f_{b'}$ agree away from $X(b) \union X(b')$; in particular
they agree on any $c$ such that $c \models p|Cb$ and $c \models p | Cb'$.
By Lemma~\ref{stepforstrongcode},
there is an $E$-definable function on $P$ with the same germ on $P$ as $f_b$.

\medskip

(ii) Again, we may assume $C=\acl(C)$. Let $e$ be a code for the $p$-germ of $f$.
Under the assumption that $f(a)\in \St_{Ca}$ for 
all $a\in P$, we must show that 
 $e\in \St_C$. By replacing $f$ by the function defined from $e$ with the same germ,
 we may suppose that
$f=f_e$, and is $Ce$-definable. 

We first need the following claim.

\medskip\noindent
{\em Claim 2.} There is a natural number $n$ and elements 
$a_1,\ldots,a_n$ of $P$ such that $e\in\dcl(Ca_1\ldots a_nf(a_1)\ldots f(a_n))$.

\par\noindent
{\em Proof of Claim.}
Let $\kappa = \wt(\tp(\ast/C))$, which is bounded in terms of $|T|$ (see Section 2.3).
Choose an indiscernible sequence $(a_i:i<\kappa^+)$ of realisations
of $p$, with $a_\lambda\models p|C(a_\mu:\mu<\lambda)$. By Proposition~\ref{stab0.7}, 
$(a_i^{\st}:i<\kappa^+)$ is a Morley sequence over $C$.
Put
$d_i=f(a_i)$. Then $e\in\dcl(C,a_id_i:i<\kappa^+)$. For suppose 
$e'$ is conjugate to $e$ over $\{C,a_id_i:i<\kappa^+\}$. Then 
for each $i<\kappa^+$ we have $f(a_i)=f_{e'}(a_i)$ (where $f_{e'}$ 
is the function defined from $e'$ in the same way that $f$ is defined 
from $e$). By the weight assumption, there is $i<\kappa^+$ with 
$a_i^{\rm st} \dnf_C (ee')^{\rm st}$. Hence, for any $a\in P$ with 
$a \domind_C ee'$, we have $f(a)=f_{e'}(a)$. As $e$ is the code for 
the germ of $f$, $e=e'$.

In particular, there are $i_1,\ldots,i_n$ and 
a $C$-definable function $h$ with $e=h(a_{i_1},\ldots,a_{i_n},d_{i_1},\ldots,d_{i_n})$.

\medskip
Since $f(a)\in \St_{Ca}$ for every $a\in P$, there is a $C$-definable relation 
$R(x,y)$ so that for each $j$, $R(a_{i_j},y)$ defines a stable, stably embedded set
which contains $d_{i_j}$. As $p$ is a definable type, 
 the following defines a first-order formula
$S(z)$ over $C$ (where $d_pu \phi(u)$ means that $\phi(u)$ holds for $u\models p|\U$):
$$
  d_pu_1\ldots d_pu_n \exists y_1\ldots \exists y_n
                  \left(\bigwedge_{j=1}^n R(u_j,y_j) \wedge
                           z=h(u_1,\ldots,u_n,y_1,\ldots,y_n)\right).
$$
Certainly $S(e)$ holds, and it remains to check that $S(y)$ is 
stable and stably embedded. Now $S$ can also be  defined, over $Ca_{i_1}\ldots a_{i_n}$,
by the formula 
$$
    \exists y_1\ldots \exists y_n \left(\bigwedge_{j=1}^n R(a_{i_j},y_j)
                   \wedge z=h(a_{i_1},\ldots,a_{i_n},y_1,\ldots,y_n)\right).
$$
Thus $S$ is internal over the stable, stably embedded set 
$\bigcup_{j=1}^n R(a_{i_j},y_j)$ and hence is itself stable and 
stably embedded. 
\qed

\begin{definition} \rm 
Let $A=\dcl(Ca)$, and let $\kappa$ be  a cardinal.   We say that $\Ast$ is
 {\em $\kappa$- generated over $C$} \index{generated!@$\kappa$-generated}  if 
there is a set  $Y\subseteq \Ast$ such that  $|Y| \le \kappa$ and
$\Ast\subseteq \dcl(CY)$.
\end{definition}

We now work towards a variant of Theorem~\ref{strongcode1} (ii) for arbitrary 
definable types, under a certain additional assumption (BS) restricting 
the growth of $St_C(a)$ as $C$ grows.  It is possible that (BS) follows from metastability over
an o-minimal $\G$;
we show at all events later that it holds for ACVF.

\begin{lemma} \label{boundedm}  Let $\Ast$ be any definably closed subset of $\St_C$.
   Let $\kappa \geq |T|$.  Then $\Ast$ is $\kappa$-generated over $C$ if and only if 
there is no strictly ascending chain of length $\kappa^+$
of sets $B=\dcl(B) \cap \St_C$ between $C$ and $\Ast$.

\end{lemma}

{\em Proof.} 
 For the forward direction, assume $\Ast$ is  $\kappa$-generated generated over $C$ by
 the set $Y$, let $\lambda=\kappa^+$ 
and suppose
there is a  chain $(A_\alpha:\alpha<\lambda)$ of subsets of $\Ast$ ordered by inclusion, 
each $A_\alpha$ definably closed in $\St_C$, with $A_0=C$. By stability, 
$\tp(Y/\bigcup_i A_i)$ is definable, hence
 using Corollary~\ref{kappa-st} it
is definable over $A_\mu$ for 
some $\mu <\lambda$. Consider $c\in A_{\mu +1}$. Then $c=f(y)$ for some $y\in Y$ 
and $C$-definable function $f$.
As $\tp(Y/CA_{\mu +1})$ is $A_\mu$-definable, $\{x\in A_{\mu +1}: f(y) = x\}$ is
definable over $A_\mu$ and is a singleton. Since $A_\mu$ is definably closed,
$c\in A_\mu$ and hence  the chain stabilises.

For the other direction, suppose $\Ast$ is not boundedly generated and take
$\lambda = (2^{|T|})^+$. Construct a sequence $(c_\alpha:\alpha<\lambda)$ of elements of $\Ast$ 
with
$c_\alpha \not\in \dcl(Cc_\beta:\beta<\alpha)$, for each  $\alpha$. Then for each 
$\mu<\lambda$, put $A_\mu:=\dcl(Cc_\alpha:\alpha<\mu)$. This is a strictly 
ascending chain of length $\lambda$ of subsets of $\St_C$.

\qed

\begin{definition}\label{BS} \rm
We say that the theory $T$ has the bounded stabilising property (BS) \index{BS} if, for every $C$ and 
every $A=\dcl(Ca)$, where $a$ is a finite tuple, there is no strictly ascending chain of length $(2^{|T|})^+$
of sets $B=\dcl(B) \cap \St_C$ between $C$ and $\Ast$. 
 
\end{definition}

We will see in Proposition~\ref{ACVFBS} that property (BS) holds in ACVF.

\begin{proposition}\label{strongcode2} Assume $T$ has (BS), and let 
$p$ be a $C$-definable type over $\U$.  
Let $f$ be a definable function on $P$, the set of realisations of $p|C$, 
and suppose that $f(a)\in \St_{Ca}$ for all $a\in P$. Then 

(i) the code $e$ 
for the germ of $f$ on $p$ (over $C$) is strong, and 

(ii) $e\in \St_C$. 
\end{proposition}

{\em Proof.} (i) Suppose that $f$ is $b$-definable. 
For any ordered set $I$ and any set $B$ containing $b$, we can find an indiscernible 
sequence $\{a_i:i\in I\}$ over $CB$, with $a_i\models p|C\cup\{Ba_j:j<i\}$ for 
any $i\in I$. Now let $A:=C\cup\{a_i:i\in I\}$. For any $i\in I$ let
$f(a_i)=d_i$, and for any $I'\subseteq I$ let
$D_{I'}:= \dcl(A\cup\{d_i:i\in I'\})$. Since by assumption, 
$f(a_i)\in \St_{Ca_i}$, $D_{I'}\subseteq \St_A(B)$.
By choosing $I$ large, we may find a strictly increasing chain
$(I_\alpha:\alpha<(2^{|T|})^+)$ of initial subintervals of $I$. Then by (BS) (applied within $\St_A)$),
the sequence $(D_{I_\alpha}:\alpha<(2^{|T|})^+)$ is eventually 
constant, so there is $i\in I$ so that $d_i\in \dcl(Ad_j:j<i)$. Choosing $I$ 
to have no greatest or least element, it follows by indiscernibility that this holds for all $i\in I$.
By a similar argument we may also suppose that $d_i\in \dcl(Ad_j:j>i)$. Thus, we have the following, for some $n\in \omega$
and any $i_1 < \cdots <i_n$ from $I$:
\begin{eqnarray}
 d_{i_n}\in & \dcl(Ca_{i_1}\ldots a_{i_n}d_{i_1}\ldots d_{i_{n-1}}) \\
 d_{i_1}\in & \dcl(CAd_{i_2}\ldots d_{i_n}). 
\end{eqnarray}
To see (6.1) above, certainly there are $m,n$ such that for any $j_1,\ldots,j_m$ with $i_1<\ldots<i_n<j_1<\ldots<j_m$,
$d_{i_n}\in \dcl(Ca_{i_1}\ldots a_{i_n}d_{i_1}\ldots d_{i_{n-1}}a_{j_1}\ldots a_{j_m})$. Suppose
$\tp(d_{i_n}'/Ca_{i_1}\ldots a_{i_n}d_{i_1}\ldots d_{i_{n-1}})=\tp(d_{i_n}/Ca_{i_1}\ldots a_{i_n}d_{i_1}\ldots d_{i_{n-1}})$.
Choose $a_{j_1},\ldots,a_{j_m}$ so that for each $k$,
$a_{j_k}\models p|Ca_{i_1}\ldots a_{i_n}d_{i_1}\ldots d_{i_{n}}d_{i_n}'a_{j_1}\ldots a_{j_{k-1}}$.
Then 
$$\tp(d_{i_n}/Ca_{i_1}\ldots a_{i_n}d_{i_1}\ldots d_{i_{n-1}}a_{j_1}\ldots a_{j_{m}})=
\tp(d_{i_n}'/Ca_{i_1}\ldots a_{i_n}d_{i_1}\ldots d_{i_{n-1}}a_{j_1}\ldots a_{j_{m}}).$$
This forces $d_{i_n}=d'_{i_n}$.

We can now show  that
if $f'$ is a definable function conjugate to $f$ over $C$ and 
has the same germ on $p$, then for any $a$, if $a\models p|C\lceil f\rceil$ and $a\models p|C\lceil f'\rceil$ 
then $f(a)=f'(a)$. For argue as above with  $B=\dcl(C\lceil f\rceil \lceil f'\rceil a)$. Let $a^*$ be an enumeration of $A$. Then 
$\lceil f \rceil aa^* \equiv_C \lceil f'\rceil aa^*$, and in particular
$$
        aa^*f(a)f(a_{i_2})\ldots f(a_{i_n}) \equiv_C 
                aa^*f'(a)f'(a_{i_2})\ldots f'(a_{i_n}).
$$
As $f,f'$ have the same $p$-germ, $f(a_{i_j})=f'(a_{i_j})=d_j$, say, for each 
$j=2,\ldots,n$.  By (6.2) above (applied over $Cb$ in place of $B$),
we have  $f(a), f'(a)\in \dcl(Caa^*d_2\ldots d_n)$. 
Thus $f(a)=f'(a)$.
It follows by Lemma~\ref{stepforstrongcode} that $e$ is a strong code for the germ of 
$f$ over $C$.

(ii) It suffices to prove (under the present hypotheses) 
Claim 2 from the proof of Theorem~\ref{strongcode1}, since 
then that proof can be mimicked.
Notice that (6.1) implies that any automorphism fixing
$Ca_{i_1}\ldots a_{i_n},d_{i_1},\ldots d_{i_{n-1}}$ also fixes $f(a_{i_n})$.
It follows that
$e\in\dcl(Ca_{i_1}\ldots a_{i_{n-1}} d_{i_1} \ldots d_{i_{n-1}})$. This gives Claim 2.
\qed

\begin{remark} \label{undefinable} \rm
The proof of (i) yields also the following. Assume (BS), and let $p$ be an 
$\Aut(\U/C)$-invariant type, not 
necessarily definable. Suppose that $f$ is a $CB$-definable function such that if 
$a\models p|CB$ then $f(a) \in \St_{Ca}$. Suppose also that $D$ is a
parameter set containing $C$, 
and that for any automorphism $\sigma$ fixing $D$ pointwise, $f$ and $\sigma(f)$ have 
the same $p$-germ. Then there is a $D$-definable function with the same germ 
as $f$ on $p$.
\end{remark}

\begin{lemma}\label{acltodclnew}
Let $C\subseteq C''\subseteq \acl(C)$, and $C\subseteq B$. Then $\St_{C''}(B)=\dcl(\St_C(B)C'')$.
\end{lemma}

{\em Proof.} 
Clearly $\dcl(\St_{C}(B)C'')\subseteq \St_{C''}(B)$.

For the other direction, let $e\in \St_{C''}(B)$. We have $\St_C=\St_{C''}=\St_{\acl(C)}$,  so
$\St_C''(B)=\St_{C''}\cap \dcl(C''B)=\St_C\cap \dcl(C''B)$.
Thus there are finite $c''\in C''$ and $b\in B$ with $e\in \dcl(c''b)$. Let $e'$ be the (finite) set of conjugates of $ec''$ over $CB$. 
Then $e'\in \St_C(B)$. Also, $\tp(ec''/e')\vdash \tp(ec''/B)$, so $\tp(e/c''e')\vdash \tp(e/c''B)$. It follows that
$e\in \dcl(c''e')\subseteq \dcl(c'', \St_C(B)\subseteq \dcl(\St_C(B) C'')$. \qed

\medskip

We now give some applications of the existence of strong codes.

\begin{proposition}\label{simplifications} Assume $C\subseteq B$, $\tp(A/C)$ is stably 
dominated, and $A\domind_C B$. Let $\Ast=\St_C(A)$  Then the following hold.

(i) $\tp(A/B)$ is stably dominated.

(ii) $\St_C(AB)=\St_C(\Ast\Bst)$.

(iii) If $\tp(B/C)$ is stably dominated then $\St_B(A)=\dcl(B\Ast)\cap \St_B$.

(iv) If $\tp(B/C)$ is extendable to a definable type, and 
$T$ satisfies (BS), then $\St_B(A)=\dcl(B\Ast)\cap \St_B$.
\end{proposition}

{\em Proof.} (i) This is just Proposition~\ref{dropstabdom}.

(ii) Put $C'':=\dcl(Ca) \cap \acl(C)$. Let $d\in \St_C(AB)=\dcl(AB)\cap\St_C$. Then $d=f(a)$ for some $a\in A$ and 
 $B$-definable function $f$. By Proposition~\ref{bits}(iii), $\tp(a/C)$ is stably dominated.
Let $C'':= \acl(C)\cap \dcl(Ca)$.
By Theorem~\ref{strongcode1}, the germ of $f$ on the definable extension of
$\tp(a/\acl(C))$ is strongly coded with the code $e$ over $C''$ lying in $\St_C$. Now $e\in \St_{C''}(B)$, so by
Lemma~\ref{acltodclnew}, $e\in \dcl(\St_C(B) C'')$. Since the code is strong, there is a $C''e$-definable function $f'$ with the same germ as $f$,
and hence $d=f'(a)$ is definable over $\St_C(B)C''A$, so is $A\Bst$-definable. Now
write $d=g(b)$, where $b\in \Bst$ and $g$ is an $A$-definable function. By
stable embeddedness, $g$ is $\Ast$-definable, hence $d\in \dcl(\Ast\Bst)$.

(iii) Let $d\in \St_B(A)=\dcl(BA)\cap\St_B$. Then $d=g(b)$ for some $b\in B$ and
$A$-definable function $g$. 
As $d\in \St_B$, we can expand $b$ to a larger tuple from $B$ to arrange 
$g(b)\in \St_{Cb}$. By Theorem~\ref{strongcode1} (with the roles of 
$A$ and $B$ reversed) the code $e$ over $\dcl(Cb) \cap \acl(C)$ for the germ over $C^*:=\acl(C)\cap \dcl(Cb)$ of $g$ on the definable extension of
$\tp(b/\acl(C))$ is strong.  Hence $e\in \St_C$, by Theorem~\ref{strongcode1}(ii). 
Thus, arguing as in (ii),  there is an $\St_C(A) C^*$-definable function $g'$ with $g'(b)=d$, that 
is, $d$ is 
definable from $\Ast b$.

(iv) The proof is the same as (iii), using Proposition~\ref{strongcode2} instead
of Theorem~\ref{strongcode1}.
\qed

\bigskip

\begin{proposition}\label{transitive} \index{stable domination!transitivity}
Suppose $\tp(A/C)$ and $\tp(B/CA)$ are stably dominated. Then
$\tp(AB/C)$ is stably dominated.
\end{proposition}

{\em Proof.} By Corollary~\ref{stab1.1}(iii) and Proposition~\ref{dropstabdom}, we may assume that
$C=\acl(C)$. Throughout this argument, for any set $X$ we will write 
$X^{\rm st}:= \dcl(CX)\cap \St_C$ and $X^*:= \dcl(CAX)\cap\St_{CA}$. Notice
that $B^*\cap\St_C=(AB)^{\rm st}=(B^*)^{\rm st}$. Suppose 
$(AB)^{\rm st} \dnf_C \Dst$. We need to show that 
$\tp(D/C(AB)^{\rm st}) \proves \tp(D/CAB)$.

As in Proposition~\ref{simplifications} we can show that 
$D^*=\dcl(CA\Dst)\cap \St_{CA}$. For consider 
$d^*\in D^*$. There are $a\in A$ and a $CD$-definable function $h$ with 
$h(a)=d^*\in \St_{Ca}$. By Theorem~\ref{strongcode1}, the germ (over $C'':=\acl(C) \cap \dcl(Ca)$) of $h$ on 
the definable extension of $\tp(a/\acl(C))$ is 
strongly coded in $\St_C$. Hence, as in Proposition~\ref{simplifications}(ii), there is a function $H$ definable over 
$\St_C(D) C''$ with $H(a)=d^*$. Thus $d^*\in \dcl(CA\Dst)$.

By Remark~\ref{symmetry} over $\St_{CA}$, we know that 
$\tp(D/CAD^*)\proves \tp(D/CAD^*B^*)$, and hence by the above paragraph, that
$\tp(D/CAD^{\rm st})\proves \tp(D/CAD^{\rm st}B^*)$. Rewriting this
as $\tp(AB^*/C\Dst)\proves \tp(AB^*/CD)$ and recalling the hypothesis that
$(AB)^{\rm st} \dnf_C \Dst$, we have
\begin{eqnarray}
  AB^*\domind_C D.
\end{eqnarray}
This implies that $\tp(D/C(AB^*)^{\rm st}) \proves \tp(D/CAB^*)$; that is,
\begin{eqnarray}
  \tp(D/C(AB)^{\rm st}) \proves \tp(D/CAB^*).
\end{eqnarray}

It follows from (3) that $B^*\dnf_{CA} D^*$ (in $\St_{CA}$). For there is 
certainly some 
$B'$ with $\tp(B'/CA)=\tp(B/CA)$ for which this is true. Intersecting with
$\St_C$ and applying transitivity (as $\Ast\dnf_C \Dst$), we get $(AB')^{\rm st} \dnf_C \Dst$. Since we also know 
$(AB)^{\rm st} \dnf_C \Dst$, it follows that 
$\tp(\Dst (AB')^{\rm st}/C) = \tp(\Dst (AB)^{\rm st}/C)$.
By (3), $\tp(AB^*/C\Dst)\vdash \tp(AB^*/CD)$. Hence, as $AB\equiv_{C\Dst} AB'$, we have
$(AB)^* \equiv_{C\Dst} (AB')^*$. Clearly $AB^*=(AB)^*$ and $A(B')^*=(AB')^*$,  so $AB^* \equiv_{CD} A(B')^*$. Hence,
as
$D^* \dnf_{CA} (B')^*$, we have
$D^*\dnf_{CA}B^*$. 

By stable domination of $\tp(B/AC)$ it follows that
$$
  \tp(D/CAB^*)\proves \tp(D/CAB).
$$
Together with (4), this gives the required result.
\qed

\smallskip

\begin{corollary} \label{uptoacl2}
Suppose $\tp(A/C)$ is 
stably dominated. Then $\tp(\acl(CA)/C)$ is stably dominated.
\end{corollary}

{\em Proof.} 
This follows from Proposition~\ref{transitive}, as $\tp(\acl(CA)/CA)$ is clearly stably dominated.
\qed

\bigskip

We close this part  with an illustration of   stable domination and strong germs in the context of
groups; cf. \cite{hrush2}.

 \def\fg{{\frak g}}

Let $G$ be a definable group, and $p$  
 a definable type over $\U$ of elements of $G$ (i.e.  containing the formula $x\in G$.)
 
For $a \in G(\U)$, we define the {\em translate} $^a p$ to be $\tp(ag/\U)$, where $g \models p | \U$.
This gives an action of $G(\U)$ on the definable types.  

We are interested in translation invariant types.  
In this case, since $^g p =p$ for $g \in G$, it follows that any element $g$ of $G$ is a
product of two elements of $p$.    

\begin{theorem} \label{groupdom} Let $G$ be a definable group, $p$ a stably dominated definable type of $G$.
Assume $p$ is translation invariant.  
  Then there
exist definable stable groups   $\fg_i$, and  definable homomorphisms $g_i: G \to \fg_i$, such that
$p$ is stably dominated via  $g=(g_i : i \in I )$.   
\end{theorem}

{\em Proof}  Let $\theta(a)$ enumerate $ \St_C(a)$.  Then $p|C$ is stably dominated via $\theta$.
Consider the map 
$f_a$ defined by:
$$f_a(b) = \theta(ab).$$
  The $p$-germ is strong, and is in $\St_C(a)$,
so it factors through $\theta(a)$: $f_a = f'_{\theta{a}}$.    By stable embeddedness, it factors through
$\theta(b)$ too: let $c =   f_a(b) = f'_{\theta{a}}(b) $.   Since $c \in \St_C$,
$\tp(\theta(a),c / C,\theta(b)) \vdash \tp(\theta(a),c / C, b)$.
Thus $c \in \dcl(\theta(a),\theta(b))$; i.e. $\theta(ab) = c =
 F(\theta(a),\theta(b))$ for $a \models p$, $b \models p|Ca$.   Associativity   of the group operation on $G$
 immediately gives associativity of $F$  on independent triples of realizations of 
 $q=\tp(c/C)$, within the stable structure $\St_C$.  
Hence by the   group chunk theorem of \cite{h-uni} (alternatively see 
\cite{hrush2}, or Poizat's book \cite{poizat2}), 
there exists an inverse limit system of stable
groups $\fg_i$ with inverse limit $\fg$, such that $\theta(a)$ is a generic element of $\fg$.  
Now    $\theta$
is generically a homomorphism.  If  $a,b,c,d $ realize $p$ and $ab = cd$ then $\theta(a)\theta(b) = \theta(c)\theta(d)$.
(Let $e \models p|C(a,b,c,d)$; then $abe=cde$, and
$ \theta(abe) = \theta(a) \theta(be) = \theta(a) \theta(b) \theta(e)$; and similarly for $cde$.)  Using the fact that
 any element of $G$ is a product of two generics,
$\theta$ extends uniquely to a homomorphism $G \to \fg$.   
\qed

\part{Independence in ACVF}

\chapter{Some background on algebraically closed valued fields}\label{acvfbackground}


In this chapter we give a little background on the model theory of valued fields, emphasising algebraically closed
valued fields. This monograph depends heavily on results and methods from \cite{hhm}. We shall summarise both the main 
results  
from \cite{hhm} which we use, and some of the methods developed there which we shall exploit.

\medskip

\section{Background on valued fields.}

A valued field consists of a field $K$ together with a homomorphism $|-|$ from its multiplicative group to an ordered 
abelian group $\Gamma$, which satisfies the ultrametric inequality. We shall follow here the notation of \cite{hhm},
and view the value group $(\Gamma, <,.,1)$ multiplicatively, with identity 1. Abusing notation, we shall usually
 suppose  that it contains an additional formal element $0$. So $0<\gamma$ for all $\gamma \in \Gamma$,
and the axioms for a valuation are as follows (with $x,y\in K$):

(i) $|xy|=|x|.|y|$;

(ii) $|x+y|\leq \Max\{|x|,|y|\}$;

(iii) $|x|=0$ if and only if $x=0$.
\noindent
The valuation is {\em non-trivial} \index{valuation!non-trivial} if its range properly contains $\{0,1\}$. 

For most  arguments in this monograph we find this multiplicative notation more intuitive than the usual additive one. However, we do 
occasionally adopt additive 
notation (with value map denoted $v$) during more valuation-theoretic arguments. For example, Lemma~\ref{place} 
and the proof of Proposition~\ref{domovergamma}
 are written additively,
as are \ref{lem1} and \ref{lotsinv}. Also, viewed additively, the value group $\Gamma$ of a field is a vector space over ${\mathbb Q}$, 
and we sometimes write $\rk_{{\mathbb Q}}(\Gamma)$ for its vector space dimension, even when viewing it multiplicatively. 

If the value map $|-|:K\rightarrow \Gamma$ is surjective, we say that $K$ has {\em value group}\index{value group} 
$\Gamma$  (though formally the group has domain $\Gamma \setminus \{0\}$). Often the value map is implicit, and we 
just refer to the valued field $(K,\Gamma)$.
The {\em valuation ring} \index{valuation ring} of $K$ is $R:=\{x\in K: |x|\leq 1\}$. This is a local ring with unique maximal ideal
$\M=\{x:|x|<1\}$. The {\em residue field} \index{residue field} is  $k:=R/\M$, and there is a natural map $\res:R\rightarrow k$. Later, when talking about ACVF, we shall slightly adjust this notation,
 viewing $K,\Gamma, R, \M, k$ as definable objects in a large saturated model of ACVF. If $\gamma\in \Gamma$ we often use the notation
 $\gamma R=\{x\in K: |x|\leq \gamma\}$ and $\gamma \M=\{x\in K: |x|<\gamma\}$.

In many texts (for example Ribenboim \cite{rib}),  part of the definition of {\em valued field} requires that the value group 
(written additively) is archimedean, that is, embeds in $({\mathbb R},+)$. Such sources refer to our more general notion as a
{\em Krull valuation}. Since we work with saturated models, for which the value group will be non-archimedean, the more general setting is forced on us.

The most familiar valued fields are probably the $p$-adic fields ${\mathbb Q}_p$. Other examples (written additively) 
are the fields of rational functions $F(T)$ ($F$ any field): here, if $p(T),q(T)\in F[T]$ are coprime, then 
$v( p(T)/q(T))=\deg(p(T))-\deg(q(T))$.
Also, if $F$ is any field, and $\Gamma$ is any ordered abelian group (written additively), we may form the field of
generalised power series
$F((T))^\Gamma$  consisting of elements $\Sigma_{\gamma \in \Gamma}a_\gamma T^\gamma$ whose {\em support}~
$\{\gamma \in \Gamma: a_\gamma \neq 0\}$ is well-ordered. Addition and multiplication are defined as for power series, 
and we put $v(\Sigma_{\gamma \in \Gamma}a_\gamma T^\gamma)=\min \{\gamma\in \Gamma: a_\gamma \neq 0\}$. 
The field $F((T))^\Gamma$ has value group $\Gamma$ and residue field $F$. In the particular case when $\Gamma$ 
is isomorphic to $({\mathbb Z},+)$ we write $F((T))^\Gamma$ as $F((T))$, the field of Laurent series.

Given a valued field $(K,\Gamma)$ with residue field $k$, there are three possibilities for the pair $(\char(K), \char(k))$:
 $(0,0), (p,p)$, and $(0,p)$ (the `mixed characteristic' case). For ${\mathbb Q}_p$ we have $(0,p)$, and 
for generalised power series fields $F((T))^\Gamma$, the field has the same characteristic as its residue field $F$.

There are various natural ways to view a valued field model-theoretically, and these mostly give the same universe of 
interpretable sets. The simplest is to view the object as having two sorts $K$ and $\Gamma$, with the value map between 
them. The same structure can also be parsed just in the field sort, with a unary predicate for the valuation ring, or with a 
binary predicate interpreted as $|x|\leq |y|$.  Indeed, if the valuation ring $R$ is specified, and $U$ is its group of units, 
then the value group $\Gamma$ is isomorphic to $K^*/U$, 
so the valuation is determined up to an automorphism of the value group.
Another option is to view it as a pair  $(K,k\cup \{\infty\})$ with a  place  $\pi: K\rightarrow k\cup \{\infty\}$
where $\pi(x)=\res(x)$ if $x\in R$, and $\pi(x)=\infty$ for $x\not\in R$. 
Recall here that a {\em place} \index{place} is a map  $\pi:K\rightarrow F\cup\{\infty\}$ (where $K,F$ are fields) such that
 $\phi(x+y)=\phi(x)+\phi(y)$
and $\phi(xy)=\phi(x)\phi(y)$ whenever the expressions on the right hand side are defined, and such that $\phi(1)=1$.
 Any surjective place $\pi: K\rightarrow F \cup\{\infty\}$ determines the valuation ring (as $\{x\in K:\pi(x)\in F\}$), so
determines  a valuation on $K$ with residue field $F$, uniquely up to an automorphism of the value group.

Chevalley's Place Extension Theorem can be stated as: given a valued field $(K,\Gamma)$ and a field extension $L>K$, 
it is always possible to extend the valuation to $L$. The value group and residue field of $K$ will embed canonically in 
those of $L$. We say that an extension $K<L$ is {\em immediate} \index{immediate extension} if $K$ and $L$ have the 
same value group and the same residue field. The valued field $(K,\Gamma)$ is
{\em maximally complete}\index{maximally complete} if it has no proper immediate extensions. 

Generalised power series fields $F((T))^\Gamma$ are maximally complete, as is ${\mathbb Q_p}$. 
Every valued field $(K,\Gamma)$ has an immediate maximally complete extension. The key ingredient in the proof 
of this is the notion of {\em pseudo-convergent (p.c.) sequence}: \index{pseudo-convergent sequence} the sequence 
$(a_\gamma: \gamma < \alpha)$  ($\alpha$ an ordinal)
is {\em pseudo-convergent} if $|a_\nu-a_\mu|< |a_\lambda-a_\mu|$      whenever $\lambda<\mu<\nu$.  The element 
$a\in K$ is a {\em pseudo-limit} of $(a_\gamma:\gamma < \alpha)$ if there is $\lambda< \alpha$ such that for all 
$\mu, \nu <\alpha$ with $\lambda<\mu<\nu$, $|a_\mu-a_\nu|=|a-a_\mu|$.
The valued field $(K,\Gamma)$ is maximally complete if and only if every p.c. sequence from it has a pseudo-limit in $K$. 
It is straightforward, given a p.c. sequence in $(K,\Gamma)$ to adjoin a pseudo-limit in an immediate extension, and 
iteration of this procedure yields a maximally complete immediate extension.
A detailed study of maximally complete fields, with criteria
 for uniqueness of immediate maximally complete extensions, was undertaken by Kaplansky \cite{kap1,kap2}.

The valued field $(K,\Gamma)$ is {\em Henselian} \index{Henselian} if for any monic polynomial $f(X)\in R[X]$ and any simple
 root $\alpha \in K$ of $\res(f)$ (the reduction of $f$ modulo $\M$), there is $a\in R$
such that $f(a)=0$ and $\res(a)=\alpha$. This property is expressible by a first order axiom scheme. Every 
valued field $K$ has a {\em henselisation}, \index{henselisation}$K^h$. This is an immediate valued 
field extension of $K$ which is henselian,
is an algebraic field extension of $K$, 
and has the property that the valuation has a unique extension from $K^h$ to the algebraic closure (denoted $K^{\alg}$ in this text) 
of $K$. The henselisation 
is unique up to valued-fields isomorphism over $K$ 
(and the isomorphism is unique).
Any maximally complete valued field is henselian, but the converse is in general false.

\section{Some model theory of valued fields.}

The best known results in the model theory of valued fields are the Ax-Kochen/Ershov principles 
(\cite{ak1}, \cite{ak2}, \cite{ak3}, \cite{er}). This is a body of results which reduce
problems in the elementary theory of valued fields to that of the value group and residue field. One version states that
if two henselian valued fields of residue characteristic zero have elementarily equivalent value groups and residue fields,
 then the valued fields themselves are elementarily equivalent. The results and methods also give information, for example
 model completeness,
for $\Th({\mathbb Q_p})$.
The methods do not handle the general case of residue characteristic $p$, but yield, for example, that for any sentence $\phi$ 
in a language for valued fields, for sufficiently large $p$, $\phi$ holds in ${\mathbb Q_p}$  if and only if it holds in 
 ${\mathbb F_p}((T))$.

The model completeness for $\Th({\mathbb Q_p})$ was extended to a quantifier elimination by Macintyre in \cite{macint}, 
in the language of rings extended by
predicates $P_n$ (for each $n\geq 2$) interpreted by the set of $n^{\th}$ powers. The obvious analogy is the Tarski 
quantifier-elimination for real closed fields,
 when the order relation is adjoined to the language of rings. 
 
 AKE and quantifier elimination results have been greatly extended by Prestel and Roquette, Basarab, Pas, and F-V. Kuhlmann. 
The AKE reduction to value group and residue field
 is picked up, at the level of Grothendieck rings, in the very recent work of Hrushovski and Kazhdan \cite{hk}.

\section{Basics of ACVF}

Suppose now that $K$ is an algebraically closed field, equipped with a surjective and non-trivial 
valuation $|-|:K\rightarrow \Gamma$. It is immediate that the value group $\Gamma$ is divisible and that
 the residue field $k$ is algebraically closed. It follows that the valuation topology is not locally compact.

The easiest examples of algebraically closed valued fields to describe, in characteristics (0,0) and $(p,p)$,  are the 
generalised power series fields $F((T))^\Gamma$ where $\Gamma$ is a divisible ordered abelian group and $F$ is an 
algebraically closed field (of characteristic 0 and $p$ respectively). If  $\Gamma$ (written additively) is $({\mathbb Q},+)$, 
and $F$ has characteristic 0, then
 $F((T^\Gamma))$ is more familiar as the field of {\em Puiseux series}, \index{Puiseux series}that is $\bigcup_{n=1}^\infty F((T^{1/n}))$.
 Another familiar algebraically closed valued field is ${\mathbb C}_p$, the completion of the algebraic closure
 of ${\mathbb Q}_p$. It has characteristic 0, residue field ${\mathbb F}^{\alg}_p$, and value group (viewed additively) 
isomorphic to $({\mathbb Q},+)$.

The model theory of algebraically closed valued fields was initiated by Abraham Robinson well before the AKE principles, in
 \cite{rob}. He showed that any complete theory of non-trivially valued algebraically closed fields is determined by the pair $(\char(K),\char(k))$.

Robinson's results were stated in terms of  model-completeness, but with a little extra work yield the following. 
Part (iii) below is proved in \cite{hhm} (Theorem 2.1.1). The language $\L_{\div}$ is the language
 $(+,-,.,0,1,\div)$, where $\div$ is the binary predicate of $K$
interpreted by $\div(x,y)$  whenever $|y|\leq |x|$.

\begin{theorem}\label{robin}
Let $K$ be an algebraically closed valued field.

(i) The theory of $K$ has quantifier elimination in the language
$\L_{\div}$.

(ii) The theory of $K$ has quantifier elimination in a 2-sorted language
with a sort $K$ for the field (equipped with the language of rings), a
sort $\Gamma$ for the value group written multiplicatively (with the language $(<,.,0)$ with
usual conventions for $0$), and a value map $|- |:K\rightarrow \Gamma$
with $|0|=0$.

(iii) The theory of $K$ has quantifier elimination in a 3-sorted language
$\L_{\Gamma k}$
with the sorts and language of (ii) together with a sort $k$ for the
residue field, with the language of rings, and a map $\Res: K^2\rightarrow
k$ given by putting $\Res(x,y)$ equal to the residue of $xy^{-1}$ (and
taking value $0\in k$ if $|x|>|y|$). \end{theorem}

A key ingredient in proofs of such results is that if $K<L$ is a finite Galois extension, or if $L$ is the algebraic closure of $K$, and $|-|$  is a valuation on $K$, then the valuations on $L$ which extend $|-|$ are all conjugate under $\Gal(L/K)$.

There is  a partial converse of (i), due to Macintyre, McKenna and van den Dries \cite{mmd}: any non-trivially valued field whose theory has quantifier elimination in the language $\L_{\div}$ is algebraically closed.

If $(K,\Gamma)$ is a valued field, then an {\em open ball of radius $\gamma$} \index{ball} in $K$ is a set of the form
$B_{\gamma}(a):=\{x\in K:|x-a|<\gamma\}$, and a closed ball of radius $\gamma$ has form
$B_{\leq \gamma}(a):=\{x\in K:|x-a|\leq \gamma\}$ (where $a\in K$ and $\gamma\in \Gamma$). Here, we allow $\gamma=0$, so view field elements as  closed balls. Of course, in the valuation topology on $K$, both open balls and closed balls of non-zero radius are clopen. It follows easily from quantifier elimination ((i) above) that if $K$ is algebraically closed, then any parameter-definable subset of $K$ (i.e. one-variable set) is a Boolean combination of balls. Jan Holly's more precise statement, in terms of Swiss cheeses, will be given when we discuss 1-torsors. Quantifier elimination very easily yields the following. It also ensures that algebraically closed valued fields are {\em $C$-minimal}: this is a variant of o-minimality introduced in \cite{ms} and \cite{hm}, and extended in \cite{hk}.

\begin{proposition}\cite[Proposition 2.1.3]{hhm}\label{fullemb} 
(i) The value group $\Gamma$ of $K$ is o-minimal in the sense that every
$K$-definable subset of $\Gamma$ is a finite union of intervals.

(ii) The residue field $k$ is strongly minimal in the sense that any
$K$-definable subset of $k$ is finite or cofinite (uniformly in the
parameters).

(iii) $\Gamma$ is stably embedded in $K$.\index{stably embedded!value group}

(iv) If $A\subset K$ then the model-theoretic algebraic closure
$\acl(A)\cap K$ of $A$ in the field sort $K$ is equal to the
field-theoretic algebraic closure.

(v) If $S\subset k$ and $\alpha\in k$ and $\alpha\in \acl(S)$ (in the sense of
$K^{\eq}$), then $\alpha$ is in the field-theoretic algebraic closure of $S$ in
the sense of $k$.

(vi)  $k$ is stably embedded in $K$.\index{stably embedded!residue field}
\end{proposition}

\section{Imaginaries, and the ACVF sorts}

It is easily seen that the theory of algebraically closed valued fields does not have elimination of imaginaries just in the field sort, or even when the sorts $\Gamma$ and $k$ (which are $\emptyset$-definable quotients of subsets of $K$) are added. When the work in \cite{hhm} was begun, we had expected to prove elimination of imaginaries when sorts for open and closed balls are added, but this too proved false \cite[Proposition 3.5.1]{hhm}. 
The main  result of \cite{hhm} was the identification of certain sorts (from $T^{\eq}$) for which
algebraically closed valued fields do have elimination of imaginaries.

In addition to the sorts $K,k,\Gamma$ we shall describe below certain sorts $S_n$ and $T_n$ (for $n\geq 1$).
We write ACVF for the theory of algebraically closed fields with a non-trivial valuation
in a multisorted language $\L_{\GG}$ with the sorts $K,\Gamma,k, S_n,T_n$ (for $n\geq 1$), and we  call these
the
{\em geometric sorts} \index{geometric sorts}.
We denote by $\U$ a large sufficiently saturated model of ACVF in these sorts, so write $s\in \U$
  to mean that $s$ is a member of one of these sorts (in the large model). Just occasionally we will consider a type $p$ 
in say the field sort {\em over} $\U$, and might say that $a\in K$ realises $p$; 
that is, we sometimes regard $K, \Gamma$ etc. as sorts, and sometimes as the corresponding subsets of the large 
model $\U$.  Here $K$ is the sort for the  field with the usual ring language, $\Gamma\setminus\{0\}$ is the value 
group in the language of multiplicative ordered groups, and the   valuation is given as a norm
$|-|:K\rightarrow \Gamma \cup \{0\}$. As above, the ultrametric inequality has the form
$|x+y|\leq \Max\{|x|,|y|\}$, and $|x|=0$ if and only if $x=0$.  Write $R$ for the valuation ring 
of $K$ and $\M$ for its maximal ideal. The residue field is $k=R/\M$, endowed with the usual ring language, and 
the `residue map' from $K\times K$ to $k$ is defined by $\Res(x,y)$: this is the residue of $x/y$ in $k$ if 
$|x|\le |y|$, and is $0$ otherwise. For $x\in R$, we denote by $\res(x)$ its residue in $k$. 

The sorts $S_n$ and $T_n$ (for $n\geq 1$) are defined as follows. First, $S_n$\index{sn@$S_n$} is the collection of all 
codes for free $R$-submodules of $K^n$ on $n$ generators (we will also call these {\em lattices}, or {\em $R$-lattices}). 
\index{lattice} Thus, if we identify $\GL_n(K)$ with the set of all ordered bases of the $K$-vector space $K^n$, then there is a map $\rho_n:\GL_n(K)\rightarrow S_n$ taking each basis to a code for the $R$-lattice spanned by it; the maps $\rho_n$ are
part of our language $\L_{\GG}$ (viewing $\GL_n(K)$ as a subset of $K^{n^2}$). In particular, elements of $S_n$ act 
as codes for the elements of a $\emptyset$-definable quotient of a subset of $K^{n^2}$.  For $s\in S_n$, we shall often 
write $\Lambda(s)$\index{lambda@$\Lambda(s)$} for the lattice coded by $s$. Observe that $R^n$ is a rank $n$-lattice, 
so has a ($\emptyset$-definable) code in $S_n$. Also, for each $\gamma\in \Gamma\setminus \{0\}$, 
$\gamma R:=B_{\leq \gamma}(0)$  is a rank one lattice, so has a code in $S_1$; the latter is interdefinable with $\gamma$.

For $s\in S_n$, write $\red(s)$ for $\Lambda(s)/\M\Lambda(s)$. This has $\emptyset$-definably the structure of an $n$-dimensional vector space over $k$. Let $T_n$\index{tn@$T_n$} be the set of codes for elements of
$\bigcup \{\red(s): s\in S_n\}$; that is, each $t\in T_n$ is a coset for some member of $\red(s)$ for some $s\in S_n$, so is a code for a coset of $\M \Lambda(s)$ in $\Lambda(s)$.  For each $n\ge 1$, we have the functions 
$\tau_n:T_n\to S_n$ defined by $\tau_n(t)=s$ if and only if $t$ codes an element of $\red(s)$. 
We shall put $\S:=\bigcup_{n\geq 1}S_n$ and $\T:=\bigcup_{n\geq 1}T_n$.

As shown in Section 2.4 of \cite{hhm}, the sorts $S_n$ and $T_n$ can be described as codes for members of coset spaces of matrix groups. This both helps to give
an intuition about them, and also supports the proofs, for example in Chapter~11. We sketch the details.

First, observe that $\GL_n(K)$ has a $\emptyset$-definable transitive action on $S_n$: for each $A\in \GL_n(K)$, if $s\in S_n$
 and $B$ is a matrix whose columns are a basis for $\Lambda(s)$, then $A(s)$ is a code for the lattice with basis the columns 
of $AB$. The stabiliser of the code for $R^n$ is just
$\GL_n(R)$, the group of $n\times n$ matrices over $R$ which are invertible over $R$. Thus, $S_n$ can be regarded as a set of 
codes for the coset space $\GL_n(K)/\GL_n(R)$. However, we find it more useful to work with the following upper triangular 
representation.

It is noted in \cite[Lemma 2.4.8]{hhm} that every $R$-lattice $A$ in $K^n$ has a basis such that the corresponding
 matrix is upper
 triangular. Indeed, let $\pi_i:K^n\rightarrow K^{n-i}$ be the projection to the last $(n-i)$-coordinates, and 
let $A_i:=\ker(\pi_i)$. 
Then for each $i$, $A_i\cong R^i$ and $A_{i+1}/A_i \cong R$ (see the proof of \cite[Proposition 2.3.10]{hhm}).
It is possible to choose a basis $(u_1,\ldots,u_n)$ of $A$ such that $(u_1,\ldots,u_i)$ 
is a basis of the free $R$-module $A_i$
 for each $i$: choose $u_{i+1}$ so that its image generates $A_{i+1}/A_i$. 
The matrix whose $i^{\th}$ column is $u_i$ for each $i$ is upper triangular.

Let $B_n(K)\subset \GL_n(K)$ be the group of
invertible upper triangular matrices over $K$, and $B_n(R)$ be the
corresponding subgroup of $\GL_n(R)$ (where inverses are required to
be over $R$). Let $\TB(K)$ be the set of triangular bases of $K^n$, 
that is, bases $(v_1,\ldots,v_n)$ where $v_i \in K^i \times (0)$ 
(i.e., the last $n-i$ entries of $v_i$ are zero). An element
$a=(v_1,\ldots,v_n)\in \TB(K)$ can be identified with an element of $B_n(K)$, 
with $v_i$ as the $i^{\th}$ column. 
Now $B_n(R)$ acts on $B_n(K)=\TB(K)$ on the right. Two elements 
$M,M'$ of $\TB(K)$ generate the same rank $n$ lattice precisely if 
 there is some $N \in \GL_n(R)$ with $MN=M'$, and as $M,M' \in B_n(K)$, 
we must have $N \in \GL_n(R) \cap B_n(K)=B_n(R)$. Using the last paragraph, this gives an identification 
of $S_n$ with  the 
set of orbits of $B_n(R)$ on  $\TB(K)$. Equivalently, $S_n$ can be identified 
with the set of (codes for) left cosets of $B_n(R)$ in $B_n(K)$. This is a natural way of
 regarding $S_n$ as a quotient of a power of $K$ by a $\emptyset$-definable 
equivalence relation.

We 
can also treat $T_n$ as a set of codes for a finite union of coset spaces. For each 
$m=1,\ldots,n$, let $B_{n,m}(k)$ be the set of elements of $B_n(k)$ whose 
$m^{\th}$ column has a 1 in the $m^{\th}$ entry and other entries zero. Let
$B_{n,m}(R)$ be the set of matrices in $B_n(R)$ which reduce (coefficientwise) modulo 
$\M$ to an element of $B_{n,m}(k)$. Also define $B_{n,0}(R)=B_n(R)$. Then $B_{n,m}(k)$ and $B_{n,m}(R)$ are groups. Let $e\in S_n$, and put $V:=\red(e)$. We 
may put $\Lambda(e)=aB_n(R)$ for some $a=(a_1,\ldots,a_n)\in \TB(K)$. So $\Lambda(e)$ is the orbit 
of $a$ under $B_n(R)$, or the left coset $aB_n(R)$ where $a$ is regarded as a 
member of $B_n(K)$, and $(a_1,\ldots,a_n)$ is a triangular basis of the lattice $\Lambda(e)$. There is a  filtration 
$$  \{0\}=V_0 < V_1 <\ldots < V_{n-1}<V_n
$$
of $V$, where $V_i$ is the $k$-subspace of $\red(e)$ spanned by
 $\{\red(a_1),\ldots,\red(a_i)\}$ (here $\red(a_j)=a_j+\M e$).
The filtration is canonical, in that if also $\Lambda(e)=a'B_n(R)$ then $V_i$ is spanned by
$\{\red(a_1'),\ldots,\red(a_i')\}$ for each $i$.

Let $\TB(V)$ be the set of triangular bases of $V$,
that is, bases $(v_1,\ldots, v_n)$ where $v_i \in V_i \setminus V_{i-1}$. 
Now $B_n(k)$ acts sharply transitively on $\TB(V)$ on the right, with the action defined by
$$  (v_1,\ldots,v_n)(a_{ij})=(a_{11}v_1, a_{12}v_1+a_{22}v_2,\ldots,
		   \Sigma_{i=1}^n a_{in}v_i).
$$
For each $i=0,\ldots,n$, put $O_i(V)=V_i\setminus V_{i-1}$ (so $O_0(V)=\{0\}$). 
It is easily verified that two elements of $\TB(V)$ are in the same orbit 
under $B_{n,m}(k)$ precisely if they agree in the $m^{\th}$ entry. Thus,
$O_m(V)$ (the set of $m^{\th}$ entries of triangular bases) can be identified with the left coset space
$\TB(V)/B_{n,m}(k)$, and 
$V\setminus \{0\}$ with $\bigcup_{m=1}^n \TB(V)/B_{n,m}(k)$. Now, if $M$ is the triangular basis
$(a_1,\ldots,a_n)$ of the lattice with code $e\in S_n$, put 
$\RED(M):=(a_1+\M\Lambda(e),\ldots,a_n+\M\Lambda(e))$, a triangular basis for $V$. 

\medskip

{\em Claim.} If $M,M'\in \TB(K)$, then they 
are $B_{n,m}(R)$-conjugate (i.e. there is $N\in B_{n,m}(R)$ with $MN=M'$)
 precisely if they generate the same lattice, and their reductions
$\RED(M)$ and $\RED(M')$ are $B_{n,m}(k)$-conjugate. 

{\em Proof of Claim.} Suppose $MN=M'$, where $N\in B_{n,m}(R)$. Then $MB_n(R)=M'B_n(R)$, so $M,M'$ generate the same lattice. Also, reducing mod $\M$, we have $\RED(M)\red(N)=\RED(M')$ where $\red(N)\in B_{n,m}(k)$ is obtained by reducing each entry mod $\M$. Conversely, suppose $M,M'$ generate the same lattice. Then there is $a\in B_n(R)$ such that $MA=M'$. Thus,
$\RED(M)\red(A)=\RED(M')$. Suppose also there is $B\in B_{n,m}(k)$ with $\RED(M)B=\RED(M')$. Then $\red(A)=B$, so $A\in B_{n,m}(R)$.

By the claim and the paragraph before it, $M,M'$ are $B_{n,m}(R)$-conjugate
precisely if they generate the same 
lattice $A$, and $\RED(M)$, $\RED(M')$ have the same element of $\red(A)$ in the $m^{\th}$ entry. Via 
the identification of $\TB(K)$ with $B_n(K)$, 
we now obtain an $\emptyset$-definable map  $\phi:\bigcup_{m=0}^n B_{n}(K)/B_{n,m}(R)\rightarrow T_n$.
For $m\in \{1,\ldots,n\}$ and $M\in B_n(K)$ let $\phi(MB_{n,m}(R))= v+\M A$, where $A$ is the lattice 
spanned by the columns of $M$, and
$v$ is the $m^{\th}$ entry of $\RED(M)$. Also put $\phi(MB_{n,0}(R)):=\M \Lambda$, where $\Lambda$ is the lattice spanned by the columns of $M$.
Thus, we may identify $T_n$ with  $\bigcup_{m=0}^n B_{n}(K)/B_{n,m}(R)$.

From now on, we view  $S_n$ as a set of codes for members of  $B_n(K)/B_n(R)$, and $T_n$ as a set of codes for elements of 
$\bigcup_{m=0}^n B_n(K)/B_{n,m}(R)$.

Formally, in the language $\L_{\GG}$, in addition to the maps $\rho_n:\GL_n(K)\rightarrow S_n$ and $\tau_n:T_n\rightarrow S_n$, there are,
for each $n\geq 1$ and $m=0,\ldots,n$, maps $\sigma_{n,m}: B_n(K)\rightarrow T_n$: for $A\in B_n(K)$, let
$\sigma_{n,m}(A)$ be a code for the coset $AB_{n,m}(R)$. 
As mentioned above, $\L_{\GG}$ also has the ring language on $K$ and $k$, the ordered group language on $\Gamma \setminus \{0\}$,
and the value map $|-|:K\rightarrow \Gamma$ and map $\Res:K^2\rightarrow k$.
 Any completion of ACVF admits elimination of quantifiers in a specific definitional expansion of $\L_{\GG}$, as shown in \cite[Section 3.1]{hhm}.
 We will not
need the precise details of the expansion here. The completions of ACVF are determined by the characteristics of $K$ and $k$.
 The main theorem of \cite{hhm} is the following.

\begin{theorem}
ACVF admits elimination of imaginaries in the sorts  of $\GG$.
\end{theorem}

We work inside a large, homogeneous, sufficiently saturated model $\U$ of ACVF, in the sorts of $\GG$, namely
 $K,\Gamma, k,S_n,T_n$ (for $n>0$). By {\em substructure of $\U$}, \index{substructure} 
we generally mean `definably closed subset of $\U$'.
We occasionally refer to elements  of $K^{\eq}$, but by elimination of imaginaries, these will have codes in $\U$.
Sometimes realisations of types over
 $\U$ are considered, but we will
assume that all sets of parameters come from $\U$, are small relative to
the size of $\U$, and can contain elements of all the geometric sorts unless 
specifically excluded. If $C\subset \U$, we write $\Gamma(C)=\dcl(C)\cap\Gamma$ 
\index{value group!gammaC@$\Gamma(C)$}and
$k(C)=\dcl(C)\cap k$. If $C$ is a subfield of $\U$, not necessarily algebraically closed, 
we write $\Gamma_C$ \index{value group!value group of C@$\Gamma_C$} for the value group of $C$, and $k_C$ for the residue field. Thus, $\Gamma(C)= \Qq \otimes \Gamma_C$.
Suppose $C\subseteq A$ are sets,  with
$A=(a_\alpha:\alpha<\lambda)$. As in Part~I, when we refer to $\tp(A/C)$, 
we mean the type of the infinite tuple listing $A$, indexed by $\lambda$; the 
particular enumeration is not important, and is often omitted. 
Likewise, if $h$ is an automorphism of some model then $h(A)$ will denote the
tuple $(h(a_\alpha):\alpha<\lambda)$, and the statement $h(A)=g(A)$ means that 
the corresponding tuples are equal.

\begin{remark}\label{torsormodule} {\rm 
By Lemma 2.2.6(ii) of \cite{hhm}, any closed ball $u$ of non-zero 
radius has code interdefinable with an element of $S_1 \cup S_2$. If the ball contains 0 
then it is already a 1-dimensional lattice; and otherwise, $u$ is $\emptyset$-interdefinable with the code for the 
$R$-submodule of $K^2$ generated by $\{1\}\times u$. }
\end{remark}

\section{The sorts internal to the residue field.}

The following lemma is used repeatedly. Part (i) is Lemma 2.1.7 of \cite{hhm}, and (ii) is an easy adaptation.
We emphasise that for each $s\in S_n$, $\Lambda(s)$ is definably $R$-module isomorphic to $R^n$, but this isomorphism 
is in general not canonical, as $\Lambda(s)$ has no canonical basis. The point below is that over any algebraically closed 
base $C$ {\em in the field sort}, any $C$-definable lattice has a $C$-definable basis. If the valuation on $C$ is trivial, this 
holds by (iii), and if it is non-trivial, it holds as $C$ is the field sort of a model of ACVF.

\begin{lemma}\cite[Lemma 2.1.7]{hhm} \label{2.1.7}
Let $C$ be an algebraically closed valued field (or more generally suppose  
$\acl_K(C\cap K) \subseteq C\subseteq\acl(C\cap K)$). Suppose $s\in S_n\cap \dcl(C)$

(i) $\Lambda(s)$ is $C$-definably isomorphic to $R^n$, and there 
are $a_1,\ldots,a_n \in C^n$ which form an $R$-basis for $\Lambda(s)$.

(ii) $\red(s)$ is $C$-definably isomorphic to $k^n$.

(iii) If the valuation on $C$ is trivial,
$\Lambda(s)=R^n$.
\end{lemma}

For each $s\in S_n$, $\red(s)$ is a finite-dimensional vector space over $k$. As the residue field is
a stable, stably embedded subset of the structure, so is $\red(s)$. As we have seen  in Part~I, the 
stable, stably embedded sets can play an important role for the independence theory of a structure. In 
an algebraically closed valued field we give the following definition. 

\begin{definition}\label{intkc} \index{vs@$\VS_{k,C}$} \index{int@$\Int_{k,C}$} \rm
 For any parameter set $C$, let $\VS_{k,C}$ be the many-sorted structure whose
sorts are the $k$-vector spaces $\red(s)$ where $s\in\dcl(C)\cap\S$. Each sort $\red(s)$ is equipped
with its $k$-vector space structure. In addition, $\VS_{k,C}$ has, as its $\emptyset$-definable 
relations, any $C$-definable relations on products of the sorts.
\end{definition}

In \cite{hhm} we used the notation ${\rm Int}$ rather than $\VS$; we have changed the notation to $\VS_{k,C}$ 
to emphasise that the structure consists of vector spaces, and not {\em all} of the $k$-internal sets.

\begin{proposition}\cite[Proposition 2.6.5]{hhm} \label{2.6.5}
For any parameter set $C$, $\VS_{k,C}$ has elimination of imaginaries.
\end{proposition}

Clearly, $\VS_{k,C}$ is contained in $\St_C$.
By the following proposition from \cite{hhm}, they are essentially the same. Recall that a $C$-definable set $D$ is
{\em $k$-internal}\index{internal!to the residue field} if there is finite $F\subset \U$ such that $D \subset \dcl(k \cup F)$. 
It is clear that if $s\in S_n$ is $C$-definable, then $\red(s)$ is $k$-internal; for if $B$ is a basis of $\red(s)$, then 
$\red(s) \subset \dcl(k \cup B)$.

\begin{proposition}\cite[Proposition 3.4.11]{hhm}  \label{263} 
Let $D$ be a $C$-definable subset of $K^{\eq}$. Then

(i) $D$ is $k$-internal if and only if $D$ is stable and stably embedded.

(ii) If $D$ is $k$-internal then $D\subset \dcl(C \cup \VS_{k,C})$.
\end{proposition}

\begin{remark}\label{*C}\rm
It follows  by the last proposition and Remark~\ref{*C1} that condition $(*_C)$ of Lemma~\ref{*C0} holds in ACVF.
\end{remark}

In \cite[Lemma 2.6.2]{hhm} several other conditions equivalent to $k$-internality
are given. In particular, we have 

\begin{lemma}\cite[Lemma 2.6.2]{hhm}\label{2.6.2}
Let $D$ be a $C$-definable set. Then the following are equivalent.

(i) $D$ is $k$-internal.

(ii) $D$ is finite or (after permutation of coordinates) contained in a finite union of sets of the form
$\red(s_1) \times \ldots \times \red(s_m)\times F$, where $s_1,\ldots,s_m$ are 
$\acl(C)$-definable elements of $\S$ and $F$ is a $C$-definable finite set of tuples.
\end{lemma}

\section{Unary sets, 1-torsors, and generic 1-types.}

In strongly minimal and o-minimal contexts, one often argues by induction on dimension,
fibering an $n$-dimensional sets over an $n-1$-dimensional set with $1$-dimensional
fibers, thus reducing many questions to the one-dimensional case over parameters.
This can also be done for definable subsets of $K^n$, when $K \models ACVF$, but
is less clear for the lattice sorts $S_n$ and for $T_n$.  It turns out however that a good
substitute exists, and we proceed to describe it.  

In particular, we will use this process to define {\em sequential independence}.  As noted earlier, 
one variable definable sets in the field sorts are finite unions of balls. This yields
a natural notion of independence: if $C\subseteq B$ and $a\in K$, then $a$ is independent from $B$ over 
$C$ if any $B$-definable ball containing $a$ contains a $C$-definable ball (not necessarily properly) containing 
$a$. We shall develop the theory in the more general setting of `unary sets', to obtain a form of independence 
for the $S_n$ and $T_n$ sorts.     

First, recall that a {\em torsor} \index{torsor} $U$ is an $R$-module $A$ is a set equipped with a regular 
i.e. sharply 1-transitive)
action of $A$ on $U$. A {\em subtorsor} \index{subtorsor} is a subset of $U$ of the form $u+B$, where $B$ 
is an $R$-submodule of $A$.

For $\gamma\in \Gamma$, let $\gamma R:=\{x\in K: |x|\leq \gamma\}$ and $\gamma \M=\{x\in K:|x|<\gamma\}$. 
These are both $R$-submodules of $K$.
A {\em definable 1-module}\index{one-module@1-module} is an $R$-module in $\Keq$ which is definably isomorphic
to $A/\gamma B$, where $A$ is one of $K$, $R$ or $\M$, $B$ is one of $R$
or $\M$ and $\gamma\in\Gamma$ with $0\le\gamma\le 1$. The 1-module is {\em open}
if $A$ is $K$ or $\M$, and {\em closed} if $A$ is $R$. A {\em definable 1-torsor}\index{one-torsor@1-torsor} is a set 
equipped with a definable regular  action of a definable 1-module on it, and an 
$\infty$-definable 1-torsor is an intersection of a chain (ordered under inclusion) of definable subtorsors of a definable 
1-torsor. Typically,
($\infty$-) definable 1-torsors arise as cosets of $(\infty$-) definable 1-modules, in a larger 1-module.
 We will sometimes call a definable 
1-torsor a {\em ball} if $\gamma=0$. A {\em 1-torsor} is either a  definable
1-torsor or an $\infty$-definable 1-torsor, and we call it a {\em $C$-1-torsor} if 
the parameters used to define it come from $C$; we do not here require that the isomorphism to $R/\gamma \M$, etc, 
is $C$-definable. Finally, a {\em $C$-unary set}\index{unary!set} is a $C$-1-torsor
or an interval $[0,\alpha)$ in $\Gamma$. A {\em unary type}\index{type!unary} \index{unary!type} over $C$ is the type
of an element of a  $C$-unary set. The most natural examples of unary sets are just  balls. 

 \begin{lemma} \label{incbyone}  
Suppose $a$ lies in a  $C$-unary set $U$.
Then $\rk_{{\mathbb Q}}(\Gamma(Ca) / \Gamma(C)) +  \trdeg(k(Ca) / k(C) )\leq  1$.
\end{lemma}

{\em Proof.}  Let $M$ be any model containing
$C$.  Suppose the lemma is false; 
 Then there exist two elements $b,c \in \Gamma(Ca) \union k(Ca)$ with $b \notin \acl(C)$ and $c \notin \acl(Cb)$. 
 Find $b' \models \tp(b/C)$ with  $b' \notin M$.  Conjugating by
 an element of $\Aut(\U/C)$   we may assume 
$b \notin M$.  Find $c' \models \tp(c/Cb)$ with $c' \notin \acl(Cb)$.  Conjugating 
 by an element of $\Aut(\U/Cb)$ we may assume $c \notin \acl(Mb)$.  
   But over $M$, $U$ is definably isomorphic to a subset of a quotient of $K$.  
So there exists an element $a' \in K$ with 
$\rk_{{\mathbb Q}}(\Gamma(Ma') / \Gamma(M)) +  \trdeg(k(Ma') / k(M) ) >  1$, a
contradiction.  \qed  
 
 \medskip

In \cite[Section 2.3]{hhm} a notion of {\em relative radius} of a definable subtorsor \index{radius of subtorsor}
 is defined. We repeat
 it, as it is occasionally used here.  Let $U$ be a torsor of the module $A$; a
{\em  subtorsor} \index{torsor!subtorsor} is a subset of the form $V=u+B$, where $B$ is a  submodule
of $A$.  Note that $B = \{v-v': v,v' \in V \}$. 
 Let $\Sub(U)$ be the set of definable subtorsors of $U$.  

We shall say that a definable 1-torsor $U$ is {\em special} \index{1-torsor!special} if it is a torsor of a 1-module $A$ which is definably 
isomorphic to a quotient of $R$ or $\M$ or a proper quotient of $K$, or equals $K$ itself (rather than just being definably isomorphic to $K$). A {\em special unary set} \index{unary set!special}
 is a subset of $\Gamma$ of form $[0,\alpha)$ or a special 1-torsor.

\begin{lemma} \label{radius}  Let $U$ be a $C$-definable special 1-torsor. 
Then there exists a $C$-definable function   $\rad: \Sub(U) \to \G$ such that for any 
$V \in \Sub(U)$, 
$$V' \mapsto \rad(V')$$
is a bijective, order-preserving map between $\{V' \in \Sub(U): V \subseteq V' \}$
and an interval in $[0,\infty]$.  \end{lemma}

{\em Proof.}  
(i) Suppose $U$ is a 1-torsor of the 1-module $A$, with definable subtorsor $V$. This means that $V$, a subset of $U$, is a
 torsor of a definable submodule
$B$ of $A$.
 We put $\rad(V):=\rad(B)$, so have to define $\rad(B)$. Suppose first $A$ is closed.
Then for some unique $\gamma$, $B=\gamma R A$
or $\gamma \M A$. 
Then
$\rad(V):=\gamma$. 
If $A$ is open (but not definably isomorphic to a quotient of $K$),
 then a definable submodule has radius $\gamma$ if it has the  form
$\gamma R A$ or $\bigcap (\delta R A: \delta>\gamma)$.

The definition of  radius for a subtorsor of a torsor arising from a proper quotient of $K$ is clear. 
First, any such quotient, if non-trivial,  is definably isomprphic to $K/R$ of $K/\M$, so we only consider these cases.
If $A$ is definably isomorphic to $K/R$, then a definable submodule $D$ has radius 
$\rad(D):=\gamma$ if $\gamma$ is greatest such that $\gamma R D =\{0\}$; if $A$ is definably isomorphic
to $K/\M$, then $D$ has radius $\gamma$ where  is greatest such that
$\gamma R D$ is isomorphic to $\{0\}$ or $k$. 

Finally, if $A=K$, the definition of radius for subtorsors of $U$ is clear.
 \qed

\begin{remark}\label{radius2} \rm
We extend the definition of $\rad$ to the case when $U$ is a $C-\infty$-definable 1-torsor which is the 
intersection of a chain $(U_i:i\in I)$
of definable subtorsors of some 1-torsor $V$, where $V$ is definably isomorphic to a quotient of $R$ or $M$. 
In this case, we arbitrarily fix some
$i_0\in I$, and for any definable subtorsor $W$ of $U$, define $\rad(W)$ with respect to $U_{i_0}$, i.e.
 by regarding $W$ as a subtorsor of the {\em definable}
1-torsor $U_{i_0}$. This device ensures that the radius lies in $\Gamma$ rather than in its Dedekind completion.
\end{remark}

If $T$ is a closed 1-torsor, say   of the closed 1-module $A$, we   write
$\red(T)$\index{redT@$\red(T)$} for its reduction $T/ \M A$,  the quotient of $T$ by the
action of $\M A$.   This has definably, without extra parameters, the structure
of a 1-dimensional affine space over $k$, so is strongly minimal.

The following result of \cite{hhm} shows that any element of the geometric sorts
can be thought of as a sequence of realisations of unary types. If $a=(a_1,\ldots,a_m)$ is a sequence of elements from the sorts
$\GG$, we say {\em $a$ is unary}\index{unary!sequence}
 if, for each $i=1,\ldots,m$, $a_i$ is an element of a unary set defined over $\dcl(a_j:j<i)$.
Trivially, any finite sequence of elements of  $K$ is unary, since $K$ itself is a unary set, as $K=K/0.R$.

\begin{proposition} \cite[Proposition 2.3.10]{hhm}\label{getun}
Let $s\in \U$. There is a unary sequence (called a {\em unary code}\index{unary!code}) $(a_1,\ldots,a_m)$ 
such that $\dcl(s)=\dcl(a_1,\ldots,a_m)$.
\end{proposition}

The proof of Proposition~\ref{getun} exploits the triangular form of the matrices in the identification of $S_n$ with $B_n(K)/B_n(R)$ 
and $T_n$ with $\bigcup_{i=0}^n B_n(K)/B_{n,m}(R)$.  
Take the case of $S_n$.  The group
$B_n(K)$ has a normal unipotent
subgroup $U_n(K)$ (the strictly upper triangular matrices) with quotient $D_n(K)$ isomorphic to 
$(K^*)^n$.  Now $D_n(K)/D_n(R)$ is isomorphic to $\G^n$, and leads to 1-torsors
of the type $\G$.  On the other hand   $U_n(K)$ has a sequence of normal 
subgroups $N_j$ with successive quotients     isomorphic to the additive   group  of $K$.
This leads to fibers of the form $g N_{j+1} U_n(R)    / N_j U_n(R) $;
these are torsors for 
$$N_{j+1} U_n(R) / N_j U_n(R)  \cong N_{j+1} / N_j (N_{j+1} \cap U_n(R))$$
Since $N_{j+1}/N_j \cong (K,+)$ and since $N_{j+1} \cap U_n(R) \neq (0)$, 
this leads to 1-torsors for modules isomorphic to $K/R$.  Similarly $T_n$
can be analyzed by elements of $\G$ and of 1-torsors for modules
isomorphic to $K/\M$.

This shows that the statement can be improved somewhat:

\begin{proposition}  \label{getun+}
Let $s\in \U$. There is a   sequence $(a_1,\ldots,a_m)$ and $C_0,\ldots,C_{m-1}$
such that $\dcl(s)=\dcl(a_1,\ldots,a_m)$, $a_i \in C_{i-1}$, 
 and $C_i$ is either $\G$ or $K$  or an $(a_1,\ldots,a_{i})$-definable torsor of an
$R$-module isomorphic to $K/R$ or $K/\M$.  \end{proposition}
 
In particular, the unary sets involved can be chosen to be special. 
Alternatively, using the proof of Proposition~\ref{getun} given in \cite{hhm},
we can obtain a decomposition with unaries of the form $\G,K$ or closed  1-torsors
(belonging to  {\em non-free}
1-generated R-modules); again, these are all special.

In \cite{hhm} we developed a theory of independence for the unary types, including results on 
definability of  types and orthogonality to the value group. Part of the goal in this monograph is to develop
 these results more thoroughly for $n$-types. We reproduce the definitions and elementary results  from
\cite{hhm}.

First, recall from Section 2 of \cite{hhm} that a {\em Swiss cheese} \index{Swiss cheese} is a subset of $K$ of the f
orm $t\setminus (t_1 \cup\ldots \cup t_m)$, where $t$ (the {\em block}) \index{Swiss cheese!block} is a ball of $K$ or 
the whole of $K$, and the $t_i$ (the {\em holes})\index{Swiss cheese!hole} are distinct proper sub-balls of $t$
(which could be field elements). More generally, if $U$ is a 1-torsor, then a subset of $U$ of the form
$t\setminus (t_1 \cup \ldots \cup t_m)$ (where $t,t_1,\ldots,t_m$ are definable subtorsors) is regarded as a Swiss cheese 
of $U$, with corresponding notions of `hole' and `block'. 
If $U_1,U_2$ are Swiss cheeses of $U$, we say that they are {\em trivially nested} \index{Swiss cheese!trivially nested} if the hole of one 
is equal to the block of another;
in this case, $U_1 \cup U_2$ can be represented as a {\em single} Swiss cheese. 
The following lemma, partly due to Holly \cite{holly}, is a consequence 
of quantifier elimination in ACVF.

\begin{lemma} \cite[2.1.2 and 2.3.3]{hhm}\label{2.3.3}
Let  $U$ be a $C$-1-torsor.

(i) Let $X$ be a definable subset of $U$. Then $X$ is uniquely expressible as
the union of a  finite set $\{A_1,\ldots,A_m\}$ of Swiss
cheeses, no two trivially nested.  

(ii) Suppose in addition that $C=\acl(C)$. If $a,b\in U$ and neither of $a,b$ lie in a $C$-definable proper subtorsor of $U$,
 then $a\equiv_C b$.
\end{lemma}

The proof uses the fact that any definable subset of $K$ is a Boolean combination of balls, and the observation that any 
definable unary set is in definable bijection with $K$, an interval of $\Gamma$, or a ball.

\begin{definition}\label{generictype} \rm
Let  $U$ be an $\acl(C)$-unary set  and $a\in U$. Then $a$ is 
{\em generic in $U$ over $C$}\index{generic}  if $a$ lies in no $\acl(C)$-unary proper subset of $U$.
\end{definition}

In the motivating case when $a\in K$, $U$ will be a $C$-definable ball (possibly $K$ itself) or the intersection of a 
chain of $C$-definable balls. Then $a$ is generic in $U$ over $C$ if and only if there is no sub-ball of $U$ which is 
algebraic over $C$ and contains $a$.

\begin{remark}\label{uniquegen}\rm
(i) By Lemma~\ref{2.3.3}, if $a,b$ are generic over $C$ in a $C$-unary set $U$, 
then $a\equiv_{\acl(C)} b$. Thus, we may talk of {\em the generic type of $U$}\index{generic!type} 
(over $C$) as the type of an element of $U$ which is generic over
$C$. This uniqueness, like Lemma~\ref{realisegen},  follows from Lemma~\ref{2.3.3}.

(ii) If $T$ is a closed 1-torsor then the above notion of genericity for the strongly minimal 1-torsor 
$\red(T)$ agrees with that from  stability theory. 
That is, if $T$ is $C$-definable then $t\in \red(T)$ is generic over $C$
if it does not lie in any $C$-definable finite subset of $T$.
Also, suppose $T$ is a $C$-definable closed 1-torsor, 
and $a$ is generic
in $\red(T)$ over $C$. Then all elements of $a$ have the same type
over $C$; for otherwise, some $C$-definable subset of $T$ intersects
infinitely many elements of $\red(T)$ in a proper non-empty subset,
contradicting Lemma~\ref{2.3.3}(i). 

(iii) We adapt slightly the above language, by saying that if $\gamma_0\in\Gamma(C)$, then 
$\gamma$ is {\em generic over $C$ below $\gamma_0$}\index{generic!in value group}
if for any $\epsilon \in \Gamma(C)$, if $\epsilon<\gamma_0$ then $\epsilon<\gamma$. That is, 
$\gamma$ is generic in the unary set $[0,\gamma_0)$. 
\end{remark}

\begin{lemma}\label{realisegen} \cite[Lemma~2.3.6]{hhm}
Suppose $C=\acl(C)$, and $a$ is an element of 
a $C$-unary set $U$. Then $a$ realises the generic 
type over $C$ of a unique $C$-unary subset $V$ of $U$.
\end{lemma}

To see this in the case when $U$ is a 1-torsor, let $V$ be the intersection of the set of $C$-definable subtorsors of $U$ containing $a$.

It follows from the lemma that if $C=\acl(C)$ and $a$ is a field element, then $\tp(a/C)$ is the generic type over $C$ 
of a unary set. Likewise if $\gamma\in \Gamma(C)$ and $a$ is a ball (say closed) of radius $\gamma$, then
$a$ lies in the $C$-unary set $K/\gamma R$, so realises the generic type over $C$ of some unary set; namely,
the intersection of the $C$-definable subtorsors of $K/\gamma R$ which contain $a$.

\begin{lemma}\label{defintype} \cite[Lemma~2.3.8]{hhm}
Let $C$ be any set of parameters.

(i) If $p$ is the generic type over $C$ of a $C$-definable unary set, then
$p$ is definable over $C$.

(ii) Let $\{U_i:i\in I\}$ be a descending sequence of $C$-definable subtorsors 
of some $C$-1-torsor $U$, with no least element, and let $p$ be the generic type over ${\cal U}$ of field 
elements of ~$\bigcap(U_i:i\in I)$. Then $p$ is not definable.
\end{lemma}

For example, in (i), suppose $U$ is the closed ball $R$, and $\phi(x,y)$ is some formula. There is some $n_\phi$ such that for each $c$, $\phi(x,c)\in p$ if and only if, for all but at most $n_\phi$ elements $\alpha$ of $k$,
$\phi(x,c)$ holds for all elements of $R$ with residue $\alpha$. Thus,
then $(d_px) (\phi(x,y))$ is just 
$$\exists \xi_1 \ldots \exists \xi_{n_\phi+1}(\bigwedge_{i=1}^{n_\phi+1} \xi_i\neq \xi_j 
         \wedge (\forall x\in R)(\bigvee_{i=1}^{n_\phi+1}\res(x)=\xi_i \rightarrow \phi(x,y))).
$$

\begin{definition}\label{genindsingleton} \rm
Let $a$ be an element of a unary set, and $C, B$ be sets
of parameters with $C=\acl(C)\subset \dcl(B)$. We say that $a$ is {\em
generically independent \index{generically independent} \index{independence!generic, $\dnf^g$}
from $B$ over $C$}, and write $a\dnf^g_C B$,
if either $a \in \acl(C)$, or, if $a$ is generic over $C$ in a
$C$-unary set $U$, it remains generic in $U$ over $B$. \end{definition}

Without extra assumptions, this notion is not symmetric. If $a\dnf^g_C B$ and $b\in B$, we might not have
$b\dnf^g_C Ca$. See Example~\ref{exsym}  for a counterexample. However, the following easy result gives a kind 
of stationarity principle for generic extensions of unary types, and yields that they have invariant extensions. It will be 
extended in the next chapter to arbitrary types.

\begin{proposition} \label{exist1gen} \cite[Proposition 2.5.2]{hhm}\index{invariant type!and generic independence} 
Let $B,C$ be sets of parameters with $C=\acl(C)\subseteq \dcl(B)$, and let $p$ be the type  of an element 
of a $C$-unary set $U$. Then there is a
unique unary type $q$ over $B$ extending $p$ such that if $\tp(a/B)=q$
then $a\dnf^g_C B$. \end{proposition}

\section{One-types orthogonal to $\Gamma$.}

In Section 2.4 of \cite{hhm} there is a complete description of definable functions from $\Gamma$ to $\U$. 
We shall not need these in their general form, but we quote a result which underpins some of the orthogonality 
discussion below.  

\begin{proposition}\cite[Proposition 2.4.4]{hhm} \label{newclassifyfun}
Let $B$ be a set of parameters, $U$ a $B$-1-torsor, 
$\alpha,\gamma \in \Gamma$, and $t$ be a subtorsor of $U$ of
radius $\gamma$ (possibly 0) with $t\in \acl(B\alpha)\setminus \acl(B)$.
Then $\gamma\in \dcl(B\alpha)$ and there is an  $s\in \acl(B)$ (a
subtorsor of $U$) with $\rad(s)<\gamma$, such that $t\in\{B_{\leq
\gamma}(s),B_{<\gamma}(s)\}$. \end{proposition}

\begin{definition}\label{orth1} \rm
Let $C=\acl(C)$, and $a$ be an element of a $C$-unary set.
We write $\tp(a/C)\perp \Gamma$, and say $\tp(a/C)$ is {\em orthogonal} \index{orthogonal!to value group} to
$\Gamma$ if, for any algebraically closed valued field $M$ such that
$C\subseteq \dcl(M)$ and $a\dnf^g_C M$, we have $\Gamma(M)=\Gamma(Ma)$.
\end{definition}

This definition will be extended to arbitrary types in Chapter~\ref{orthogonal}. To see the reason for going up to a model 
$M$, suppose $a\in R$ is not algebraic over $\emptyset$, and $s:=\lceil B_{<1}(a)\rceil$. Then 
$\Gamma(s)=\Gamma(sa)=\{0,1\}$. However, if $M$ is any model with $s\in \dcl(M)$, then $M$ 
contains a field element $b$ in the ball coded by $s$. 
Now if $a\dnf^g_b M$ then $\gamma:=|b-a|\in \Gamma (Ma)\setminus \Gamma(M)$. Indeed if $\gamma \in\dcl(M)$ then
$B_{\leq \gamma}(b)$ is a proper sub-ball of $B_{<1}(a)$ containing $a$, contradicting genericity. This argument shows 
that the generic type of an open ball, or the intersection of a chain of balls with no least element, cannot be
 orthogonal to $\Gamma$. In fact, we have the following.

\begin{lemma}\label{eqorthog} \cite[Lemma~2.5.5]{hhm}
Let $C=\acl(C)$ and  $a\not\in C$ lie in  a $C$-unary set 
$U$. Then  the following are equivalent:

(i) $a$ is generic over $C$ in a closed subtorsor of $U$ defined over  $C$.

(ii) $\tp(a/C)\perp \Gamma$.

\medskip\noindent
Furthermore, if $A:=\acl(Ca)$ then condition 

(iii) $\trdeg(k(A)/k(C)) =1$ 

implies both (i) and (ii). If in addition $C=\acl(C\cap K)$, then (i), (ii) are equivalent to (iii).
\end{lemma}

The following lemmas follow in \cite{hhm} from the analysis of definable functions with domain $\Gamma$.

\begin{lemma}\cite[Lemma 3.4.12]{hhm}\label{3.4.12}
If $C= \acl(C)$, and $\alpha\in \Gamma$, then $\acl(C\alpha)  =\dcl(C\alpha)$.
\end{lemma}

\begin{lemma}\cite[Lemma 2.5.6]{hhm} \label{2.5.6}
If  $T$ is a $C$-1-torsor which is not a closed 1-torsor, then the following are equivalent:

(i) no proper subtorsor $T'$ of $T$ is algebraic over $C$;

(ii) for all $a$ generic in $T$, $\Gamma(C)=\Gamma(Ca)$.
\end{lemma}

To see the direction (i) $\Rightarrow$ (ii), suppose that $a$ is generic in $T$ and $\delta\in \Gamma(Ca)\setminus \Gamma(C)$. 
Then there is a $C$-definable function $T\rightarrow \Gamma$ with $f(a)=\delta$, and $f^{-1}(\delta)$ is
 a proper $C\delta$-definable subset of $T$. 
It follows easily from Lemma~\ref{2.3.3} that there is a $C\delta$-definable proper subtorsor $T_\delta$ of $T$. By 
Proposition~\ref{newclassifyfun}, $T_\delta$ is a neighbourhood of a subtorsor $T'$ of $T$ which is definable over $C$.

\begin{definition}\cite[Definition 2.5.9]{hhm} \label{2.5.9} \rm
If $C=\acl(C) $ and $a$ lies in some unary set, we say that
 $\tp(a/C)$ is {\em order-like} \index{order-like} if $a$ is generic over $C$ in a $C$-unary set which is either 
(i) contained in $\Gamma$, or (ii) an open 1-torsor, or (iii) the
 intersection of a chain of $C$-definable 1-torsors with no least element, such that this
 intersection contains some proper $C$-definable subtorsor.
\end{definition}

\begin{lemma}\cite[Remark 2.5.10]{hhm} \label{orderlike}
(i) Suppose that $\tp(a/C)$ is order-like (and in case (ii) of the last definition, assume  also $C=\acl(C\cap K)$). 
Then $\Gamma(C)\neq \Gamma(Ca)$. 

(ii) Suppose that $a$ lies in some unary 
set but $\tp(a/C)$ is not order-like. Then $\Gamma(C)=\Gamma(Ca)$.
\end{lemma}

{\em Proof.} (i) This follows from Lemma~\ref{2.5.6} (ii)$\Rightarrow$ (i).

(ii) Apply Lemma~\ref{eqorthog} if $\tp(a/C)$ is the generic type of a closed ball, and Lemma~\ref{2.5.6} otherwise. \qed

\medskip 

The next lemma gives a symmetry property of $\dnf^g$ to be extended in
 Propositions~\ref{fieldsym} and \ref{noembed} below. A special easy case (which plays a role in the proof)
is when $a$ and $b$ are each generic over $C$ in a $C$-definable closed ball.

\begin{lemma}\cite[Lemma 2.5.11]{hhm}\label{2.5.11}
Suppose $C=\acl(C)$ and $a$ and $b$ are respectively elements of the $C$-1-torsors $U$ and $V$. 
Assume that at least one of $\tp(a/C)$, $\tp(b/C)$ is not order-like. Then
$a\dnf^g_C \acl(Cb)$ if and only if $b\dnf^g_C \acl(Ca)$.
\end{lemma}

\section{Generic bases of lattices.}

We shall need repeatedly a notion of {\em generic basis} \index{generic!basis} for  a lattice, introduced in 
Section 3.1 of \cite{hhm}.
For $s \in S_n$, $B(s):=\{a\in (K^n)^n:a=(a_1,\ldots,a_n), \Lambda(s)=Ra_1+\ldots+Ra_n\}$, the
set of all bases of $\Lambda(s)$. We shall describe an $\Aut(\U/C)$-invariant extension $q_s$ of
the partial type $B(s)$ over $C$, where $s\in \dcl(C)$. As $\red(s)^n$ is a definable set of
Morley rank $n^2$ and degree 1 in the structure $\VS_{k,C}$, it has a unique generic
type (in the sense of stability theory) $q_{\red(s)^n}$ over $\U$. Now
 $a=(a_1,\ldots,a_n)\models q_s$ if and only if $(\red(a_1),\ldots,\red(a_n))\models q_{\red(s)^n}.$ To 
show that $q_s$ is complete, observe that there is a $\U$-definable
isomorphism $\Lambda(s)\rightarrow R^n$. Thus we may suppose $\Lambda(s)=R^n$. Now
$q_{\red(s)^n}$ is just the type over $\U$ of a generic element
$(\beta_1,\ldots,\beta_{n^2})$ of $k^{n^2}$. It follows easily from
Remark~\ref{uniquegen}(ii) that for such a sequence, any two tuples
$(b_1,\ldots,b_{n^2})$, where $\res(b_i)=\beta_i$ for each $i$, have the
same type. This gives completeness of $q_s$, and, along with the
invariance of $q_{\red(s)^n}$, yields invariance of $q_s$. 

We call a realisation of $q_s|C$ a {\em generic basis} of $\Lambda(s)$ over $C$ or 
{\em generic resolution}\index{generic!resolution} of $s$ over $C$, and also talk of a
generic resolution over $C$ of a sequence $s_1,\ldots,s_m$ of codes of $C$-definable lattices;
 the latter is a sequence $b_1,\ldots,b_m$, where each $b_i$ is a generic basis of $\Lambda(s_i)$
over $Cb_1\ldots b_{i-1}$. The order of the sequence is irrelevant to genericity.

\begin{lemma}\cite[Remark 3.1.1]{hhm} \label{3.1.1}
Let $s\in S_n \cap \dcl(C)$, let $B\supseteq C$, and suppose that $a\models q_s|_B$. Then
$\Gamma(B)=\Gamma(Ba)$.
\end{lemma}

We omit the proof, but for example, if $B$ is a model then $\Lambda(s)$ is interdefinable with $R^n$ and a
 realisation of $q_s|B$. The latter is just a generic sequence of realisations of the closed ball $R$,
 so by Lemma~\ref{eqorthog} does not extend the value group .



\chapter{Sequential independence}\label{gindependence}


In this chapter, we extend Definition~\ref{genindsingleton} of generic independence for an element of a unary type
 to a definition of sequential  independence for arbitrary tuples.  This yields 
 invariant extensions of arbitrary types over algebraically closed sets.  
 For stably dominated types, we observe as a consequence that sequential 
 independence coincides with independence defined in Part 1.

Sequential independence can equally well be deduced from generic
independence for one-types  in an o-minimal or weakly 
o-minimal theory: see Example~\ref{additive}.  However as discussed
in the introduction, the need to work over algebraically
closed sets makes it impossible here to reduce to   ambient dimension one:  
one cannot stay within the family of   algebraically closed sets while adding one point
at a time.  It is necessary to add the whole algebraic closure along with the new point;
the reduction is only to  pro-finite covers of unary sets, hence by compactness
to finite covers of unary sets.  It is for the same reason that pseudo-finite fields
do not admit quantifier-elimination, but require quantifiers over algebraically bounded sets.

As Examples~\ref{exsym} and \ref{exfun} show, sequential independence is not preserved 
under permutations of the variables, even in the  field sort; the same issues would arise in the o-minimal case. However, we do 
obtain uniqueness of sequentially independent extensions, and in particular, this serves to show the 
existence of invariant extensions of {\em any} type over an algebraically closed set (Corollary \ref{extendable}  below). 
Sequential independence will be used in Chapter~\ref{orthogonal} in the definition of orthogonality to $\Gamma$. The 
latter turns out to be the same as stable domination, and thus leads to a more 
symmetric form of independence. Using sequential independence we also show, at the end of this chapter, that
 over an algebraically closed set $C$ the collection of $n$-types which extend to a $C$-definable type over $\U$ is dense in $S_n(C)$.

For more examples involving sequential independence, see Chapter~\ref{examples} and the end of 
Chapter~\ref{Jindependence}.

\begin{definition}\index{acl-generated@$\acl$-generated} \index{dcl-generated@$\dcl$-generated} \rm 
If $A\subseteq\acl(Ca)$ for some finite tuple $a$, we say 
that $A$ is {\em finitely $\acl$-generated over $C$}, and call $a$ an
{\em $\acl$-generating sequence}. Likewise, if $A\subseteq\dcl(Ca)$ 
then $A$ is {\em finitely $\dcl$-generated over $C$}, with {\em $\dcl$-generating 
sequence} $a$. 
\end{definition}

\begin{definition}\index{independence!sequential, $\dnf^g$} \index{sequential independence}
 \label{g-indep}\rm 
Let $A, B,C$ be sets. For $a =(a_1,\ldots,a_n)$, and $U=(U_1,\ldots,U_n)$, 
define $a\dnf^g_C B$ {\em via} $U$ (`$a$ is sequentially independent from  $B$ over $C$ via $U$') 
to hold if for each $i\leq n$, $U_i$ is an $\acl(Ca_1\ldots a_{i-1})-\infty$-definable unary set, and 
$a_i$ is a generic element of $U_i$ over $\acl(BCa_1\ldots a_{i-1})$. We allow here the degenerate case 
when $a_i\in \acl(Ca_1\ldots a_{i-1})$, formally putting $U_i=\{a_i\}$.

We shall say $A\dnf ^g_C B$ via $a,U$ if 
$a$ is an $\acl$-generating sequence for $A$ over $C$
and $a\dnf^g_C B$ via $U$ (and we sometimes 
omit the reference to $U$). We say $A\dnf^g_C B$ if 
$A\dnf^g_C B$ via some $a,U$. Finally, we say $A\dnf^g_C B$
{\em via any generating sequence} if for any unary $\acl$-generating sequence $a$ 
for $A$ over $C$, with $a$ chosen from $A$,  we have $a\dnf^g_C B$  via $U$ for some $U$. 
\end{definition}

\begin{remark} \label{infinite-g}\rm
Definition~\ref{g-indep} also makes sense  for transfinite tuples $a=(a_i:i<\lambda)$. If $U=(U_i:i<\lambda)$,
we say that $a\dnf^g_C B$ if for each $i<\lambda$, $U_i$ is an
$\acl(Ca_j:j<i)-\infty$-definable unary set, and $a_i$ is generic in it over $\acl(BCa_j:j<i)$. We shall refer to such $a$ as a 
{\em unary transfinite sequence}. \index{unary transfinite sequence} For concreteness, we 
generally work 
with finite sequences,
and make occasional comments indicating how most of the  results lift to the transfinite version. The infinite version is 
used in the proof of 
Theorem~\ref{moregtoJ}.
\end{remark}

The following cautionary examples show that $\dnf^g$ is 
not symmetric in general, and that the relation $a \dnf^g_C B$ 
can depend on the order of $a=(a_1\ldots,a_n)$.

\begin{example} \label{exsym} \rm
(Failure of symmetry) Suppose that $b$ is generic over $\emptyset$ in the open ball $\M=B_{<1}(0)$ 
and  $a$ is generic over $b$, also in  $B_{<1}(0)$.  Then $a\dnf^g_\emptyset b$. But
$|b|<|a|$, so $b\df^g_\emptyset a$, as the  $a$-definable ball $B_{\leq |a|}(0)$, which is smaller than $B_{<1}(0)$, contains $b$.
\end{example}

\begin{example}\label{exfun} \rm
(Dependence on the enumeration) Let $b_1\in {\cal M}$ and $b_2\in 1+{\cal M}$. Choose $a_1\in {\cal M}$ 
generic over $b_1b_2$, and $a_2$ generic in  $1+{\cal M}$ over $a_1b_1b_2$.
Then $a_1a_2 \dnf^g_C b_1b_2$. However, by genericity of $a_2$, 
$|a_1-b_1|<|a_2-b_2|<1$, so $a_1  \df^g_{Ca_2} b_1b_2$ as $a_1\in B_{\leq |a_2-b_2|}(0)$, and 
hence $a_2a_1  \df^g_C b_1b_2$.
\end{example}
For more complicated examples of sequential independence, see Chapter~\ref{examples}.

Despite these examples, sequential independence behaves in some ways as one would
expect of an independence relation, and has useful properties:  there is transitivity and monotonicity
for the parameters on the right (by the next lemma),
and sequentially independent extensions always exist (by 
Proposition~\ref{existence} (i)).

\begin{lemma}\label{transitivity}
(i) If $C \subseteq B \subseteq D$, then
$A\dnf^g_C D$ via $a,U$ if and only if $A\dnf^g_C B$ 
via $a,U$ and $A\dnf_B D$ via $a,U$.

(ii) If $a_1,U_1$ are $k$-tuples and $a_2,U_2$ are $(n-k)$-tuples, 
and $a=a_1a_2$, $U=U_1U_2$, then $a\dnf^g_C B$ via $U$ if and only if
$a_1\dnf^g_C B$ via $U_1$ and $a_2\dnf_{Ca_1} B$ via $U_2$. 
\end{lemma}

{\em Proof.} This is immediate from the definition. \qed

\bigskip

\begin{lemma}\label{exist1}
For any $A,B$ with $B\subseteq A$ and $A$ finitely $\acl$-generated over $B$, 
we have $A\dnf^g_B B$.
\end{lemma}

{\em Proof.} By transitivity, it suffices to do this when $A=\acl(Ba)$ for 
some single $a\in \U$, and this is immediate except when 
$a$ is a code for an element of $\S \cup \T$. However, by Proposition~\ref{getun}, in this 
case $L$ has a unary code
$(e_1,\ldots,e_n)$ with each $e_i$ in a unary set $V_i$ definable over 
$\acl(e_j:j<i)$.
If $V_i$ is a 1-torsor, choose $U_i$ to be the intersection of the 
$\acl(Be_j:j<i)$-definable subtorsors of $V_i$ which contain $e_i$, 
and argue similarly if 
$V_i$ is an interval of $\Gamma$. Then $U_i$ is an 
$\acl(Be_j:j<i)$-$\infty$-definable unary set, and $e_i$ is generic in $U_i$ 
over these parameters. \qed

\begin{proposition} \label{existence}
Let $A,B,C$ be sets, and $a=(a_1,\ldots,a_n)$ be 
in $A$ such that $A\subseteq \acl(Ca_1\ldots a_n)$, and each $a_i$ 
lies in a unary
set defined over $\acl(Ca_j:j<i)$. 

(i) There are $A',a'$ such that $A'a'\equiv_C Aa$ and
$A' \dnf^g_C B$ via $a'$.

(ii) There is $B'\equiv_C B$ such that $A \dnf^g_C B'$ via $a$.
\end{proposition}

{\em Proof.} 
(i) This is by induction on $i$. The case $i=1$  follows from 
Proposition~\ref{exist1gen}. For  any $1\le j\le n$ write 
$A_j=\acl(Ca_1\ldots a_j)$, $A_j'=\acl(Ca_1'\ldots a_j')$. Suppose for 
the inductive hypothesis that 
for some $i$ we have $A_i' \dnf^g_C B$ via $(a_1',\ldots,a_i')$ and 
$\sigma_i \in \Aut(\U/C)$ with $\sigma_i(A_i)=A_i'$, $\sigma_i(a_j)=a_j'$ for 
$j=1,\ldots,i$. By Proposition~\ref{exist1gen}, there is 
$a_{i+1}'$ such that $A_ia_{i+1} \equiv_C A_i'a_{i+1}'$ and
$a_{i+1}' \dnf^g_{A_{i}'} B$. Now $\sigma_i$ extends to
$\sigma_i':A_i \cup\{a_{i+1}\} \rightarrow A_i' \cup\{a_{i+1}'\}$,
which extends to $\sigma_{i+1}:A_{i+1} \rightarrow A_{i+1}'$. Finally, put
$A':=A_n'$.

(ii) First find $A'a'$ as in (i). Let $\tau\in \Aut(\U/C)$ 
be such that $\tau(a)=a'$ (so $\tau(A)=A'$). Then as 
$A'\dnf^g_C B$ via
$a'$, we have $A \dnf^g_C \tau^{-1}(B)$ via $a$. 
Now put $B':=\tau^{-1}(B)$.
\qed   

\medskip

Note that the transfinite analogues of \ref{transitivity}, \ref{exist1} and \ref{existence} also hold.                               

\medskip
\index{invariant extensions!in ACVF|(}
We now use sequential independence to show that all types over $C$ 
have $\Aut(\U/C)$-invariant extensions. The essential point is that
the uniqueness or `stationary' statement of Proposition~\ref{exist1gen}
also extends to sequentially independent extensions of $n$-types, but requires more work because
of problems with algebraic closure. First, we need three  lemmas. 
The first is immediate, the second is a restatement of Proposition 3.4.13 of \cite{hhm}, and the third
is the key to the uniqueness result.

\begin{lemma} \label{conjgerms}
Let $p$ be a type over $\U$, and let $g_1,\ldots,g_r$ be  
definable functions which agree on $p$. Then there is a 
$\lceil\{g_1,\ldots,g_r\}\rceil$-definable function 
$g$ which agrees with the $g_i$ on $p$.
\end{lemma}

{\em Proof.} Define $g$ by putting $g(x)=y$ if $g_1(x)=\ldots =g_r(x)=y$, 
and $g(x)=0$
(or some $\emptyset$-definable element of the appropriate sort) otherwise. \qed
\medskip

The following lemma is taken from \cite{hhm}.  Recalling that generic types of closed unary sets are
stably dominated,  an independent proof of the closed case
can be found in Chapter~\ref{strongcodes}.   The other cases reduce to the closed case using the classification
in \cite{hhm} of definable maps from $\G$.

\begin{lemma}\cite[Proposition 3.4.13]{hhm} \label{endofhhm}
Let  $\acl(C)=C$,  and let $U$ be a $C$-unary set. Let $f$ be a definable function (not necessarily $C$-definable) with 
range in
$\U$ such that for all $x\in U$ we have $f(x)\in \acl(Cx)$. Then there is a $C$-definable function $h$ with the same germ on $U$ as $f$.
\end{lemma}

\begin{lemma} \label{moreinv}
Suppose $C \subseteq B$, with $\acl(C)=\dcl(C)$, and 
$a:=(a_1,\ldots,a_n)$, with $a\dnf^g_C B$. 
Suppose also $s\in \U$ with 
$s\in \dcl(CaB) \cap \acl(Ca)$. Then $s \in \dcl(Ca)$.
\end{lemma}

{\em Proof.} 
We shall prove the lemma by induction (over all $a,s, B, C$) 
on the least $m$ such that there is a sequence $d=(d_1,\ldots,d_m)$
such that
$d\dnf^g_C B$, $d \in \dcl(Ca)$, and 
$a \in \acl(Cd)$.
Clearly the result holds if $m=0$.

So suppose we have $d=(d_1,\ldots,d_m)$ satisfying the 
above with respect to $a$.
We first show that $a$ may be replaced by $d$. Let 
$F_0$ be the finite set of conjugates of $a$ over $CBd$. 
Since $s\in \dcl(CBa)$, there is a $CBd$-definable 
function $h$ on $F_0$ with
$h(a)=s$. The graph of $h$ is the finite set of conjugates 
of $(a,s)$ over $CBd$,
and these are also conjugate over the smaller set $Cd$. 
Hence, as $(a,s)$ is algebraic over
$Cd$, all its conjugates over $CBd$ are algebraic 
over $Cd$, so 
$\lceil h \rceil \in \acl(Cd)$. Let $(t_1,\ldots,t_r)$ be 
a code in $\U$ for $\lceil h\rceil$ (which 
exists by  elimination of imaginaries). Now
$(d,t_i)$ satisfies the  hypotheses in the statement of the lemma
for $(a,s)$. If we could prove that $t_i\in \dcl(Cd)$ for 
each $i$, then 
$\lceil h\rceil \in \dcl(Cd)$, and hence $s=h(a)\in 
\dcl(Cad)=\dcl(Ca)$, as required. Thus, in the 
assumptions we now replace $a$ by $d$, and assume
$s\in \dcl(CBd) \cap \acl(Cd)$, and must show
$s\in \dcl(Cd)$. 

As $s\in \dcl(CBd)$, there is a $CBd_1\ldots d_{m-1}$-definable 
partial function $g$ with $g(d_m)=s$. Put $C':=\acl(Cd_1\ldots d_{m-1})$. By minimality of $m$, $d_m\not\in C'$. 
Let $p:=\tp(d_m/C')$, so $p$ is the generic type  of a
unary set $U$ which is $\infty$-definable over $C'$. 
We shall assume $U$ is not a unary subset of $\Gamma$,
for otherwise we could adjust the argument below, replacing $g$ 
by a function 
$g^*:K\rightarrow \U$, where $g^*(x)=g(|x|)$.  Thus, we may suppose that $U$ is the intersection of a chain
$(U_i:i\in I)$ of $C'$-definable 1-torsors which are all subtorsors of $U_0$, say; possibly $U_i=U_0$ for each $i$. 
Let $A':=C'\cap \dcl(Cd_1\ldots d_m)$.
For each $i$, $U_i$ has just finitely many conjugates over $Cd_1\ldots d_{m-1}$, so if $U_i'$ is such a conjugate, we
 cannot have one of $U_i,U_i'$ containing the other. Given two subtorsors of $U_0$, either one contains the other or
 they are incomparable.
 Hence, as $d_m\in U_i$ for each $i$, it follows that each $U_i$ is $A'$-definable. That is, 
$p$ is the generic type 
of a 1-torsor (namely $U$) which is $\infty$-defined over $A'$.

By Lemma~\ref{endofhhm}
there is a $C'$-definable function $g'$ with the same $p$-germ as $g$, 
so $g'(d_m)=s$. 
Let $g_1=g',\ldots,g_r$ be 
the conjugates of $g'$ over $BA'$. 
As $g$ is defined over $BCd_1\ldots d_{m-1} \subseteq BA'$,
and $g$ and $g_1$ have the same $p$-germ, also $g$ and  $g_i$ have 
the same $p$-germ for all $i$. By Lemma~\ref{conjgerms} 
that there is a 
$\lceil\{g_1,\ldots,g_r\}\rceil$-definable function 
$h$ which has 
the same $p$-germ as all the $g_i$. As $\{g_1,\ldots,g_r\}$ is a set of 
conjugates over $BA'$,
$\lceil h\rceil \in \dcl(BA')$. Also, as $\lceil g_1\rceil \in \dcl(C')$
and $g_i \equiv_{Cd_1\ldots d_{m-1}} g_1$, $\lceil g_i \rceil 
\in \dcl(C')$ for each $i$, so $\{g_1,\ldots,g_r\}$ is $C'$-definable, 
and hence
$\lceil h\rceil \in C'=\acl(Cd_1\ldots d_{m-1})$. Choose finite unary  $e \in
A'$ with $(d_1,\ldots,d_{m-1})$ as an initial segment
so that 
$\lceil h\rceil \in \dcl(CBe)$; then $e\dnf^g_C B$, as $e\in \acl(Cd_1\ldots d_{m-1})$. 
Let $(t_1,\ldots,t_{\ell})$ be a code in $\U$ for
$\lceil h\rceil$. Then, by the induction hypothesis (with
$e$ replacing $a$, the $t_i$ replacing $s$, and
$(d_1,\ldots,d_{m-1})$ replacing $d$), we have 
$t_i\in \dcl(Ce)$ for each $i$, so   
$\lceil h\rceil
\in \dcl(Ce)$. Since $Ce \subseteq \dcl(Cd_1\ldots d_m)$, and
$h(d_m)=g(d_m)$,
$s=h(d_m) \in \dcl(Cd_1\ldots d_m)$, as required.\qed

\begin{theorem} \label{moreinv2}
Suppose $C \subseteq B $ with $C=\acl(C)$, and 
$a=(a_1,\ldots,a_{n})$,
$a'=(a_1',\ldots,a_{n}')$ with $a\equiv_C a'$, $a\dnf^g_C B$
and $a'\dnf^g_C B$. Then $a\equiv_B a'$.
\end{theorem}

{\em Proof.} This is by induction on $n$, and the case $n=1$ comes from
 Proposition~\ref{exist1gen}. We assume the result holds for $i-1$, so by 
applying an automorphism we may assume
$a_1=a_1',\ldots, a_{i-1}=a_{i-1}'$.
Then $a_i$  is   generic in a unary set $U = \bigcap(U_j:j\in J)$ defined over $\acl(Ca_1\ldots a_{i-1})$.
Here the $U_j$ are $\acl(Ca_1\ldots a_{i-1})$-definable unary sets; we allow the case $J$ infinite
as well as $|J|=1$, and in the latter case $U=U_1$ is allowed to be a point.  
Let $U'=\bigcap(U_j':j\in J)$  be the conjugate of $U$ over $\acl(Ca_1\ldots a_{i-1})$ containing $a_i'$.
For each $j\in J$, let $c_j$ be a code
in $\U$ for the set of conjugates of $U_j$ over $Ca_1\ldots a_{i-1}B$. Then
$c_j\in \dcl(Ca_1\ldots a_{i-1}B) \cap \acl(Ca_1\ldots a_{i-1})$, so 
by Lemma~\ref{moreinv},
$c_j\in \dcl(Ca_1\ldots a_{i-1})$. It follows that since $U_j'$ is a conjugate
of $U_j$ over $Ca_1\ldots a_{i-1}$, then
$U_j\equiv_{Ba_1\ldots a_{i-1}} U_j'$. Hence by compactness and saturation, applying an automorphism over
$Ba_1\ldots a_{i-1}$, we may suppose $U=U'$. 
Now by Proposition~\ref{exist1gen} we have 
$a_i \equiv_{Ba_1\ldots a_{i-1}} a_i'$, as required. 

\qed

\bigskip

We list a number of corollaries. We shall refer to them collectively 
as `uniqueness' results (for sequentially independent  extensions). In particular, we obtain
 the existence of invariant 
extensions of types over algebraically closed sets.

\begin{corollary} \label{leftright}
Suppose $C=\acl(C)$, and $A,A',B,B' $ all contain $C$ with
 $A\dnf^g_C B$ via $a$, $A' \dnf^g_C B'$ via $a'$,
$Aa\equiv_C A'a'$ and $B \equiv_C B'$. Then 
$AaB\equiv_C A'a'B'$.
\end{corollary}

{\em Proof}. Using  an automorphism over $C$, we may suppose 
that $B=B'$. 
Let $a_1, a_1'$ be any finite sequences from $A,A'$ respectively such that
$aa_1\equiv_C a'a_1'$. Then $aa_1\dnf^g_C B$ and $a'a_1'\dnf_C^g B$.
Thus, by Theorem~\ref{moreinv2}, $aa_1\equiv_B a'a_1'$.   \qed

\begin{corollary}\label{infgstat}
The analogue of Corollary~\ref{leftright} also holds when $a=(a_i:i<\lambda)$ is a transfinite sequence.
\end{corollary}

{\em Proof.} Again, we may suppose $B=B'$. For each $i<\lambda$, let
$A(i)=\acl(Ca_j:j<i)$. As usual, we view the $A(i)$ as sequences. We prove inductively that
$A(i)a_i\equiv_{B} A(i)'a_i'$. 

Assume that this holds for all $j<i$. If $i$ is a limit, it follows that $A(i)\equiv_B A(i)'$.
Thus, we may suppose $A(i)=A(i)'$. Then as $a_i\equiv_{A(i)} a_i'$ and $A(i)$ is algebraically closed,
$A(i)a_i \equiv_B A(i)a_i'$ by Proposition~\ref{exist1gen}. 

Suppose now $i=k+1$ is a successor ordinal. Again, we may suppose $A(k)=A(k)'$. Now $a_ka_{k+1}\equiv_{A(k)} a_k'a_{k+1}'$ and 
$a_ka_{k+1}\dnf^g_{A(k)} B$ and $a_k'a_{k+1}' \dnf^g_{A(k)} B$. Hence, by Corollary~\ref{leftright},
$A(k+1)a_{k+1}\equiv_B A(k+1)'a_{k+1}'$, as required. \qed

\medskip

The next two corollaries also have analogues where $a$ is indexed by an infinite ordinal.

\begin{corollary}\label{invariance}
Let $C=\acl(C)$, let $M$ be a model containing $C$, 
and suppose $A\dnf^g_C M$ via $(a_0,\ldots,a_{n-1})$. Then
any automorphism of $M$ over $C$ is elementary over $A$.
\end{corollary}

{\em Proof.} Let $g$ be an automorphism of $M$ over $C$, 
and let $(a_\alpha:\alpha<\lambda)$ be an enumeration of $A$, 
with $A\dnf^g_C M$ via $(a_0,\ldots,a_{n-1})$.
There is a sequence $(a_\alpha':\alpha<\lambda)$ such that $g$ 
extends to an elementary map $a_\alpha\mapsto a_\alpha'$. Now 
$(a_\alpha:\alpha<\lambda)\equiv_C (a_\alpha':\alpha<\lambda)$, 
and $(a_\alpha':\alpha<\lambda)\dnf^g_C M$ 
via $(a_0',\ldots,a_{n-1}')$. It follows by Theorem~\ref{moreinv2} 
that the map $h$ fixing $M$ pointwise and taking each
$a_\alpha'$ to $a_\alpha$ is elementary. Now the elementary map $hg$ fixes
$A$ pointwise and induces $g$ on $M$, as required.
\qed

\begin{corollary}\label{extendable}
Let $C=\acl(C)$, let $a$ be a finite sequence, and put 
$A=\acl(Ca)$. Let $M$ be a model containing $C$. Then 
there is $A'a'\equiv_C Aa$ such that $\tp(A'/M)$ is invariant under 
$\Aut(M/C)$. 
\end{corollary}

{\em Proof.} Let $d$ be a unary code for $a$.
By Proposition~\ref{existence} there is $d'\equiv_C d$ with
$d'\dnf^g_C M$. Now apply Corollary~\ref{invariance}.
\qed
\index{invariant extensions!in ACVF|)}
\medskip

Next, we obtain some initial results on $\St_C$ and stable domination in ACVF, using the
above uniqueness results.

\begin{proposition}\label{gdom}
In ACVF, assume $C\subseteq B $. Suppose $\tp(A/C)$ is stably dominated, and 
let $a$ be a unary $\acl$-generating sequence for $A$ over $C$ (possibly transfinite). 
Then $A\domind_C B$ if and only if $A\dnf^g_C B$ via $a$.
\end{proposition}

{\em Proof.} By Corollary~\ref{stab1.1}, $\tp(A/\acl(C))$ is also stably dominated.
Thus, $A\domind_C B \Leftrightarrow A\domind_{\acl(C)} B$. But
$A\domind_{\acl(C)} B$ if and only if $\tp(A/B)$ has an extension to an 
$\Aut(\U/\acl(C))$-invariant type over $\U$ 
(by Lemma~\ref{stab0.1}). Also, $A\dnf^g_C B$ via $a$ 
if and only if $A\dnf^g_{\acl(C)} B$ via $a$, and  by 
Corollary~\ref{invariance} the latter implies
that $\tp(A/B)$ extends to an
$\Aut(\U/\acl(C))$-invariant type over $\U$. The result now 
follows from the uniqueness part of Proposition~\ref{definable} (ii). \qed  
\index{stable domination!in ACVF|)}

\begin{lemma} \label{notclosed} 
Let $a$ be generic over $C$ in a $C$-1-torsor $U$ which is not closed. Then $\acl(Ca)^{\st}=\acl(C)$.
\end{lemma}

{\em Proof.}
Suppose there is $c\in \dcl(Ca) \cap \St_C$ with $c\not\in \acl(C)$.
Then there is a definable function $f:U\rightarrow \St_C$ with $f(a)=c$. 
There is also an $\acl(C)$-definable (i.e. not just $\infty$-definable) 1-torsor $U_0$ with $U$ as a subtorsor.
Choose
parameters $b$ so that $U_0$ is $b$-definably isomorphic 
with a 1-torsor of the form $A/\gamma B$ (for $A\in \{K,R,\M\}$, $B\in \{R,\M\}$,
 and $\gamma \in\Gamma$ with $0\leq \gamma\leq 1$); so, working over $b$,  there is a natural meaning for $|x-y|$, where $x,y\in U_0$.
 Also pick $a'\in U$. For each $\gamma<\rad(U)$, let 
$\hat{f}(\gamma):=\lceil\{f(x):|x-a'|=\gamma\}\rceil$. Then $\hat f$ is a definable function 
from a totally ordered set to a stable stably embedded set, so we may 
suppose that $\hat f$ is constant.
Thus there is a set $V$ such that for all $\gamma \in \dom(\hat{f})$, 
$\{f(x):x\in U \wedge |x-a'|=\gamma\}=V$. As $c\not\in \acl(C)$, $V$ is infinite. Now pick $y\in V$. Then 
$f^{-1}(y)$ contains a proper non-empty subset of $\{x:x\in U \wedge |x-a'|=\gamma\}$ for each $\gamma \in \Gamma$.
Thus $f^{-1}(y)$  is a definable 
subset of $U$ which is not a finite union of Swiss cheeses, which 
is a contradiction.      \qed

\begin{proposition} \label{gtostab}
In ACVF, suppose $C\subseteq B $. Suppose  $a\dnf^g_C B$ for some (possibly transfinite) tuple $a$, and let $A:=\acl(Ca)$. Then

(i) $A^{\st} \dnf_C B^{\st}$.

(ii) $k(A)$ and $k(B)$ are linearly disjoint over $k(C)$
and $\Gamma(A)$ and $\Gamma(B)$ are $\Qq$-linearly independent over $\Gamma(C)$.
\end{proposition}

{\em Proof.}   By Corollary~\ref{extendable},
$\tp(A/B)$ extends to an $\Aut(\U/\acl(C))$-invariant type.
Parts (i) and the first assertion of (ii) follow immediately from Remark~\ref{C''}(ii). For (ii), note that
$\tp(\Gamma(A)/\Gamma(B))$ has an $\Aut(\Gamma(\U)/\Gamma(C))$-invariant extension over $\Gamma(\U)$, and 
this implies $\Gamma(A)\cap \Gamma(B)=\Gamma(C)$. \qed

\medskip

The following lemma, related to Lemma~\ref{notclosed} will be used in the proof of Theorem~\ref{fulldom}.

\begin{lemma}\label{stabdomgamma}
Suppose $\tp(a/C)$ is stably dominated. Then $\tp(a/C) \vdash \tp(a/C\Gamma(\U))$.
\end{lemma}

{\em Proof.} This is immediate from the definition of stable domination, as $\dcl(\Gamma(\U) \cup C) \cap \St_C=\dcl(C)$. \qed

\medskip

For subsets of the field sort, the results on generic (i.e. sequential) independence
 from \cite{hhm} for
one variable and the above uniqueness results have two useful applications  
for extensions of transcendence degree one, or more generally 
for extensions increasing $\G$ by rank at most one.    Here,  $\rk_{{\mathbb Q}}$ 
 refers to dimension as a vector space over ${\mathbb Q}$.

\begin{proposition} \label{fieldsym}
Let $C \leq A,B$ be algebraically closed valued fields, and suppose
$\rk_{{\mathbb Q}} (\Gamma(AB)/\Gamma(C)) \leq 1$. Then the condition $A\dnf^g_C B$  is symmetric in $A$ and $B$, and does not depend 
on the choice of $\acl$-generating sequences.
\end{proposition}

{\em Proof.}  For convenience, in the proof we suppose that $A$ and $B$ have finite transcendence degree over 
$C$, but this is not essential.
Let $a\in K^n$ and $b\in K^m$ with
$A=\acl_K(Ca)$, $B=\acl_K(Cb)$, and $A\dnf^g_C B$ via $a$.
For any $1\leq i \leq n$ and $1\leq j \leq m$,  
$$a_i{\dnf}^g_{\acl(Ca_1\ldots a_{i-1}b_1\ldots b_{j-1})} b_j.
$$
Since $\rk_{{\mathbb Q}}(\Gamma(Ca_1\ldots a_ib_1\ldots b_j)/\Gamma (Ca_1\ldots a_{i-1}b_1\ldots b_{j-1})\leq 1$,
 it follows from Lemma~\ref{orderlike}(i) that
at least one of the types $\tp(a_i/\acl(Ca_1\ldots a_{i-1}b_1\ldots b_{j-1}))$ and
$\tp(b_j/\acl(Ca_1\ldots a_{i-1}b_1\ldots b_{j-1}))$ is not 
{\em order-like}, in the sense of Definition~\ref{2.5.9}. 
Hence, by Lemma~\ref{2.5.11},
$b_j\dnf^g_{\acl(Ca_1\ldots a_{i-1}b_1\ldots b_{j-1})} a_i$.
Since this holds for each $j$, $b \dnf^g_{\acl(Ca_1\ldots a_{i-1})} a_i$,
and as this holds for all $i$, we obtain $b\dnf^g_C a$. \qed

\begin{proposition} \label{noembed}
Let $C \leq A,B$ be algebraically closed valued fields, and suppose
that $\Gamma(C)=\Gamma(A)$,
$\trdeg(B/C)=1$, and there is no embedding of $B$ into $A$ over 
$C$. Then

(i) $\tp(A/C) \cup \tp(B/C) \vdash \tp(AB/C)$;

(ii) $\Gamma(AB)=\Gamma(B)$.
\end{proposition}

\begin{remark} \label{symrem}\rm
In the proposition, if $A=\acl_K(Ca)$ where 
$a \in K^n$, then by Proposition~\ref{existence} there is some $A'a'\equiv_C  Aa$ with
$A'\dnf^g_C B$ via $a'$. Hence  (i) implies that 
$A\dnf^g_C B$ via $a$. Likewise, if 
$B=\acl(Cb)$ then $B\dnf^g_C A$ via $b$.
\end{remark}

{\em Proof of Proposition~\ref{noembed}.}
(i) Let $b \in B\setminus C$. It suffices to show $b\dnf^g_C A$. For then $B\dnf^g_C A$ via $b$. Hence, if $A'\equiv_C A$ then
 also $B\dnf^g_C A'$ via $b$. So by Corollary~\ref{leftright}, $A'B\equiv_C AB$ 
as required. 

 We may suppose that $b$ is generic over $C$ in the
$C-\infty$-definable 1-torsor $U$. Suppose $b\df^g_C A$. Then there is an
$A$-definable proper subtorsor $V$ of $U$, with $b\in V$. If $\gamma$ denotes the radius of $V$, then
$\gamma\in \Gamma(A)=\Gamma(C)$. By Lemma~\ref{2.1.7}, there is $a\in A$ with $a\in V$. As $b$ does not embed into
$A$ over $C$, $b\not\equiv_C a$, so there is a $C$-definable proper subtorsor $V'$ of $U$ containing just one of $a,b$, and it must contain $a$ and not $b$. Since $a\in V$, $V\cap V'\neq\emptyset$; hence as $b\not\in V'$, $V'\subset V$. Now as $V$ has radius $\gamma \in\Gamma(C)$ and contains $V'$, $V$ is $C$-definable (it is defined as a ball of radius $\gamma$ containing $V'$). This contradicts that $b$ is generic in $U$ over $C$. So $b\dnf^g_C A$.

(ii) Again, let $b\in B\setminus C$, and suppose that $b$ is generic over $C$ in the $C-\infty$-definable 1-torsor $U$.
It suffices to show that if $d$ is in the field $A(b)$, then $|d| \in \Gamma(B)$; for $\Gamma(AB)$ is the definable closure
in $\Gamma$ 
of the set of such values. We may suppose that $d=\Sigma_{i=0}^n a_ib^i=(a_1'-b)\ldots (a_n'-b)$ (as $A$ is algebraically closed). 
Thus, it suffices to show that $|b-a|\in 
\Gamma(B)$ for any $a\in A$. If $a\not\in U$, then there is a $C$-definable ball $U'$ containing $U$ but not
$a$, and $|b-a|=|x-a|$ for any $x\in U'$, so $|b-a|\in \Gamma(A)=\Gamma(C)$. So suppose that $a\in U$. As in (i), there is a $C$-definable 1-torsor $V$ 
containing $a$ and excluding $b$. Thus, $|b-a|$ is the value of $|b-x|$ for 
any $x\in V$, so $|b-a|\in \Gamma(B)$, as required. \qed

\bigskip

As a further corollary of the uniqueness results, we show that in the geometric 
sorts, Lascar strong types and strong 
types coincide, so the Lascar group over any parameter set is profinite (see \cite{kp} 
and \cite{las} for background). The proof works in any theory in which every type
has an invariant extension. This has also been noted independently by A. Ivanov \cite{iv}, where 
there is further information on connections between invariant extensions, Lascar strong types, and 
Kim-Pillay strong types.

\begin{corollary}\label{lascar}\index{type!Lascar strong!in ACVF}
Let $a,b$ be finite tuples, and $C$ be a parameter set. Then $a$, $b$ have the same 
Lascar strong type over $C$ if and only if they have the same type over $\acl(C)$.
\end{corollary}

{\em Proof.} By elimination of imaginaries, we may suppose $a,b$ are from $\U$.
The direction $\Rightarrow$ is immediate, so suppose 
$\tp(a/\acl(C))=\tp(b/\acl(C))$. 
Choose $c\equiv_{\acl(C)}a$ 
with $c\dnf^g_C ab$. Then choose a model $M\supset C$ with 
$a\dnf^g_C M$ and $c\dnf^g_{Ca} M$. Then $c\dnf^g_C Ma$ (transitivity), so
$c\dnf^g_C M$, so by Theorem~\ref{moreinv2}, $c\equiv_M a$. 
Thus $a$ and $c$ have the same Lascar strong type, and similarly $b$ and 
$c$ have the same Lascar strong type, so $a$ and $b$ have the same Lascar strong
type (all over $C$). \qed

\bigskip

We conclude this chapter with a result on definability of types which extends Lemma~\ref{defintype}. We show that, if
$C =\acl(C)$, the set of $n$-types over $C$ which have an extension to 
a $C$-definable type over $\U$ is dense in the Stone space $S_n(C)$. 
 We begin with a general lemma.

\begin{lemma} \label{dense2}
Suppose that $\U$ is a sufficiently saturated model of some complete theory $T$, and  $C\subset \U$,
with $\dcl(C)=\acl(C)$. Suppose also $p=\tp(b/\U)$ is 
$C$-definable, 
and $b'\in\acl(Cb)$. Then $\tp(bb'/\U)$ is $C$-definable. 
\end{lemma}

{\em Proof.} Clearly $\tp(b'/\U b)$ is definable. It follows that
$q:=\tp(bb'/\U)$ is definable. Let $\phi(yy',x)$ be a formula, and let 
its $q$-definition be
$(d_qyy')\phi(yy',x)=\theta(x,c)$, with canonical parameter $c\in \U^{\eq}$.
If $\{c_i:i\in I\}$ is a complete set of conjugates of $c$ over $C$, there is a corresponding set
$\{q_i:i\in I\}$ of conjugates of $q$ under $\Aut(\U/C)$. These all extend
$\tp(b/\U)$, so for each $i\in I$ there is $b_i'$ such that
$bb_i'\models q_i|\U$. As the $q_i$ are distinct, the $b_i'$ are also distinct.
However, there are at most $\Mult(b'/Cb)$ such $b_i'$, so 
$|I|\leq \Mult(b'/Cb)$.
Thus, $c\in \acl(C)$. \qed

\begin{remark} \label{inftup} \rm
In the last lemma we may allow $b'$ to be a tuple of infinite length
(with the same proof).
\end{remark}

In ACVF, if $C=\acl(C)$, given an $n$-type $p(x)$ over $C$ and a formula $\phi(x,y)$ over $C$ with 
$n=l(x)$, $m:=l(y)$,  write
$$
  S_px\phi(x,y):=\{b:\mbox{~there is~} 
              a \models p \mbox{~such that~} a{\dnf}^g_{C} 
                      b\mbox{~and~} \phi(a,b)\}.
$$
If the set $S_px\phi(x,y)$ is definable over $C$, let
$d_px\phi(x,y)$ be a formula over $C$ defining it.
\begin{theorem} \label{densestone}
Let $C =\acl(C)$. 

\noindent 
(a) Let  $\theta(x)$ be an $\L_{\GG}$-formula over $C$ with $l(x)=n$, 
and suppose $\models\exists x\theta(x)$. 
Suppose also that $\theta(x)$ implies that each $x_i$ 
is algebraic over $Cx_1\ldots x_{i-1}$ or lies in a unary set 
defined over $Cx_1\ldots x_{i-1}$. Then 
there is an 
$n$-type $p$ over $C$ containing $\theta$
such that: 

(i) for all $\phi(x,y)$ over $C$,
$S_px\phi(x,y)$ is $C$-definable;

(ii) $p$ has an 
extension to a $C$-definable type $p'$ over $\U$ such that for any $a$ 
realising $p'$ we have
$a\dnf_C^g \U$.

\noindent
(b) The set of $n$-types over 
$C$ which have an extension to a $C$-definable type over $\U$
is dense in the Stone space $S_n(C)$.

\end{theorem}

{\em Proof.} 
(a) 
(i) 
We use induction on $n=l(x)$. Suppose first that 
$n=1$. If $\theta(x_1)$ has finitely many solutions, then choose
$p$ to be any type over $C$ containing $\theta$.
Since $p$ is isolated, in this case the result is immediate. 
So suppose $\theta(x_1)$ is non-algebraic. 
Then the solution set of $\theta(x_1)$ lies in a $C$-definable unary set $U_1$.
We assume $U_1$ is a 1-torsor; the case when $U_1 \subset \Gamma$ 
is easier, using o-minimality of $\Gamma$.
By Lemma~\ref{2.3.3}, the  solution set of $\theta(x_1)$ is a finite 
union
of disjoint Swiss cheeses $s_1\setminus (t_{11} \cup \ldots \cup 
t_{1\ell_1}),
\ldots, s_k \setminus (t_{k1}\cup \ldots \cup t_{k\ell_k})$. 
We may assume no two of these are trivially nested,  so this representation is unique. Hence, as 
$C=\acl(C)$,  $\lceil s_i \rceil \in C$ for each $i$.
Let $p$ be the generic type of $s_1$ over $C$. Then $\theta(x)\in p$. Also, $p$ is $C$-definable by 
Lemma~\ref{defintype}, and likewise the generic type of $s_1$ over $\U$ extends $p$ and is $C$-definable, by the same schema. It follows   that  $S_px\phi(x,y)$ is $C$-definable for each $\phi$.

Suppose now $n>1$, and let $\theta'(x_1,\ldots,x_{n-1})$ be the formula
$\exists x_n\theta(x_1,\ldots,x_n)$. By induction, there is a type
$p_{n-1}$ over $C$ containing $\theta'$ such that for each 
$\phi(x_1,\ldots,x_{n-1},y)$,
$S_{p_{n-1}}(x_1,\ldots,x_{n-1})\phi(x_1,\ldots,x_{n-1},y)$ is 
$C$-definable.
Let $(e_1,\ldots,e_{n-1})$ realise $p_{n-1}$. 
By the $n=1$ case, there is $s\in \acl(Ce_1\ldots e_{n-1})$ and  a definable 1-type $q$ over $Ce_1\ldots e_{n-1}s$
containing
$\theta(e_1,\ldots,e_{n-1},x_n)$ so that for any formula 
$\phi(x_1,\ldots,x_{n-1},z,x_n, y)$ over $C$,
$S_qx_n\phi(e_1,\ldots,e_{n-1},s, x_n,y)$ is $Ce_1\ldots e_{n-1}s$-definable, 
by the formula $d_qx_n\phi(e_1,\ldots,e_{n-1}, s,x_n,y)$. Let $e_n$ realise $q$, and
$p=\tp(e_1,\ldots,e_n/C)$. Also, let $p^*=\tp(e_1,\ldots,e_{n-1},s/C)$, a definable type by 
Lemma~\ref{dense2}. Put $x':=(x_1,\ldots,x_{n-1})$. Then for any $\phi(x',z,x_n,y)$ over 
$C$, we have $b \in S_px\phi(x',z,x_n,y)$
$$
\begin{array}{rl}
\Longleftrightarrow &
\exists x'z \models p^*\!(x'z
\dnf^g_C b
\wedge \exists x_n \models q\ (x_n \dnf^g_{Cx'} 
b \wedge \phi(x',z,x_n,b)))\\
\Longleftrightarrow &
\exists x'z \models p^*\!(x'z
\dnf^g_C b \wedge d_qx_n \phi(x',z,x_n,b))\\
\Longleftrightarrow & d_{p^*}x'z(d_q x_n)
\phi(x',z,x_n,b).
\end{array}
$$

(ii) 
Choose $a$ realising $p$
with
$a\dnf^g_C \U$. Put $p':=\tp(a/\U)$. Then for any formula 
$\phi(x_1,\ldots,x_n,y_1,\ldots,y_m)$,
if $b\in \U^m$ then $a\dnf^g_C b$, so 
$\phi(x,b)\in p$ if 
and only if 
$b\in S_px\phi(x,b)$. 

(b) Let $\theta(x)$ be a formula over $C$. 
By Proposition~\ref{getun}, we may assume that $\theta(x)$ asserts that $x$ is a unary sequence.
By (a)(ii), there is a $C$-definable type $p$ over $\U$ containing 
$\theta(x)$.
\qed

\chapter{Growth of the stable part.}\label{dominationinacvf}
\index{stable domination!in ACVF|(}

The stable part $\St_C$ of ACVF over a parameter set $C$ grows with $C$.
Later on we will need to study with exactitude how $\St_C$ grows when,
for instance, one adds to $C$ an element of $\G$.   Here we fix a tuple $a$
and consider the stable part $\St_C(a)$ of $a$ over $C$, as $C$ grows.

Recall  (Definition~\ref{intkc}) that in ACVF, if $C$ is any parameter set,
then 
$\VS_{k,C}$ is the
 many sorted structure consisting of 
a sort $\red(s)$ for each $C$-definable lattice with code $s$, with the induced $C$-definable 
relations. This is essentially the same as  $\St_C$, for by Proposition~\ref{263}, every element
 of any $C$-definable, stable, stably embedded set is coded in $\VS_{k,C}$. In particular, if the 
parameter set $C$ is a model  $M$  of ACVF, and $s\in S_n\cap \dcl(M)$, then $\red(s)$ is
 in $M$-definable bijection with $k^n$ (Lemma~\ref{2.1.7}). Hence, over a model $M$, the structure $\VS_{k,M}$, and 
hence also $\St_M$, is essentially 
just $k$, with the $M$-definable relations.  On the other hand when $C$ is 
not a model, the essential number of sorts can be large. We shall later 
obtain some results by `resolving' $C$, i.e. finding an appropriate a model $M$ with $C\subset \dcl(M)$,
 to study this situation. 

 We will 
show in this section that $\St_C(a)$ is nevertheless always countably generated.  
First observe the following.

\begin{example} \rm  The residue field of a finitely generated extension of a valued field $L$ need
not be finitely generated over the residue field of $L$.   

Indeed, let $F$ be a field, 
 $F((t))$ the Laurent series field with the  usual valuation, trivial on $F$.  
 Let   $a_n \in F^{\alg}$ be any sequence
of elements, $s = \sum a_n t^n \in F((t))$, $L=F(s,t)$, $L^h$ the Henselization.
Then $a_0 \in L^h$; indeed if $P$ is the minimal monic polynomial for $a_0$ over $F$,
then $a_0$ is the unique solution of $F$ in the neighborhood $s+\M$ of $s$.  
But then it follows that $\sum a_{n+1} t^n = (s-a_0)/t \in L^h$, and inductively
each $a_i \in L^h$.  So the residue field of $L^h$ contains $F(a_0,a_1,\ldots)$.   
But this is also the residue field of $L$.

Similarly, using the generalized power series $s'=\sum t^{n+\frac{1}{n}}$, we see that there exists a valuation on the rational function field in two variables whose value group
is $\Qq$.  
\end{example}

Indeed, for any valued field $F$, any countably generated algebraic extension of 
the residue field of $F$ is contained in the residue field of some finitely generated valued
field extension of $F$.   However, by Corollary~\ref{ctbl2}  below,  this is the worst that can happen.

Call a type $\tp(a/C)$ {\em stationary} \index{type!stationary} if $\dcl(Ca) \meet \acl(C) = \dcl(C)$.
Equivalently, $\tp(a/C) \vdash \tp(a/\acl(C))$.  

\begin{lemma} \label{stationary-t}  If $\tp(a/C)$ and $\tp(b/C(a))$ are stationary,
then so is $\tp(ab/C)$. \end{lemma}

{\em Proof.}  $\dcl(Cab) \meet \acl(C) \subseteq \dcl(Cab) \meet \acl(Ca) \subseteq
\dcl(Ca)$ by the stationarity of $b/Ca$; hence by the stationarity of $a/C$,
   $\dcl(Cab) \meet \acl(C) \subseteq \dcl(Ca) \meet \acl(C) = \dcl(C)$.  \qed
   
   \medskip
   
For the rest of the chapter we work in ACVF.

\begin{example}    If $a \in \Gamma$ then $\tp(a/C)$ is always stationary. \end{example}

{\em Proof.}  Any $\phi \in \tp(a/\acl(C))$  is a finite union 
of intervals in $\G$.  The endpoints of these intervals, being $\acl(C)$-definable,
are actually $C$-definable using the linear ordering.  Hence $\phi$ is $C$-definable,
so $\phi \in \tp(a/C)$.  Thus $\tp(a/C) \vdash \tp(a/\acl(C))$. \qed

\medskip

An extension $C'$ of $C$ is said to be {\em finite} if $C' \subseteq \acl(C)$
and $C'$ is a finitely $\dcl$-generated extension of $C$.  In this case, $C'=\dcl(Ce)$
for some finite $e$, and the number of realizations $\Mult(e/C)$ of $\tp(e/C)$
does not depend on the choice of $e$.  We write  $\Mult(C'/C) =\Mult(e/C)$ .

\begin{lemma} \label{1type-dichotomy} Let $X$ be a $C$-definable special unary set, and 
let $a\in X$.
  Then one of the following holds:
\begin{enumerate}  
\item For some finite extension $C' $ of $C$, $\tp(a/C')$ is stationary.
\item  For some countable $C_0 \subseteq C$, $\tp(a/\acl(C_0)) \vdash \tp(a/\acl(C))$.
\end{enumerate}
\end{lemma} 

{\em Proof.} 
\def\Oo{R}
Since every type in $\G$ is stationary, 
we may therefore suppose that $X$ is a torsor
for a unary $\Oo$-module $M$, or $X=K$.
 It follows (Lemma~\ref{radius}) that any proper sub-torsor $Y$ of
$X$ can be assigned a {\em radius} $\rad(Y)$; and $Y \in \dcl(C,Y',\rad(Y))$ 
for any subtorsor $Y'$ of $Y$.  
 
Let $W$ be the collection of $\acl(C)$-definable   subtorsors of $X$.  For $w \in W$,
let $C(w):=\dcl(Ce)$, where $e$ is a code in $\U$ for $w$;
 let $m(w) = \Mult (C(w) / C)$.   Let $W_a = \{w \in W: a \in W\}$.  

Note that if $Y' \in W$ and $Y' \subset Y$,
then $Y \in \dcl(C,Y', \rad(Y))$; but $\rad(Y)$ is an $\acl(C)$-definable element of $\G$,
hence is $C$-definable; so $Y \in \dcl(C,Y')$, and $C(Y) \subseteq C(Y')$.  Thus
$m(Y) \leq m(Y')$, and if equality holds then also $C(Y)=C(Y')$.  

Note also
that if $w,w' \in W_a$ then so is $w \cap w'$, and $m(w \cap w') :=\Max\{ m(w),m(w')\}$.
In particular, if $m(w)$ is unbounded on $W_a$ then for all 
   $w \in W_a$ there exists $w' \in W_a$ with $w \supset w'$ and $m(w) < m(w')$. 
   
{\em Case 1.}   $m(w)$ is bounded on $W_a$.

 Let $Y$ be chosen in $W_a$
and  $m(Y)$ maximal possible.   If $Y' \in W_a$ with  $Y' \subset Y$,
then   $m(Y') \geq m(Y)$, so by maximality $m(Y')=m(Y)$, and $C(Y')=C(Y)$.  
 This shows that any element of $W_a$ is $C(Y)$-definable.  Let $B$ be the intersection
 of all elements of $W_a$; then $a$ realizes over $\acl(C)$ the generic type of $B$;
 it follows that   $\tp(a/C') \vdash \tp(a/C)$.  Thus (1) holds, with $C':=C(Y)$.
 
 {\em Case 2.}  $m(w)$ is unbounded on $W_a$.  Then we can find $w_1 \supset w_2 \supset \ldots  \in W_a$
 with $m(w_1)< m(w_2) < \ldots$.   Let $B = \meet_{n=1}^\infty w_n$.  
 
 We claim that $B$ contains no proper $\acl(C)$-definable subtorsors.
 Indeed, suppose $B$ contains $w^*\in W$.  Then each $w_n$ contains $w^*$, and it follows
 as above that $w_n \in \dcl(C,w^*)$.  So $m(w_n) \leq m(w^*)$ for each $n$, a contradiction.

 Let $C_0$ be a countable substructure of $C$ such that each $w_i$ is $\acl(C_0)$-definable.
 Then $B$ is $\infty$-definable over $\acl(C_0)$, and $\tp(a/\acl(C_0)) \vdash \tp(a/\acl(C))$.  \qed

 \begin{corollary}\label{ctbl}  Let $a$ be any finite tuple, and $C$ any base structure. 
 Then there exists a countably $\dcl$-generated algebraic extension $C'$ of $C$
 such that $\tp(a/C')$ is stationary.  \end{corollary} 
 
 {\em Proof.}  
We prove this by induction on the length of a unary code for $a$, so may assume that
$a=(a_1,\ldots,a_n)$ is unary. If $n=1$, then either option in  Lemma~\ref{1type-dichotomy} clearly implies the conclusion.  

So put $b=(a_1,\ldots,a_{n-1})$, and by induction let $C'$ be a countably $\dcl$-generated (over $C$)
subset of $\acl(C)$ such that $\tp(b/C')$ is stationary. By Lemma~\ref{1type-dichotomy} there
is countably generated $D=C'(b,b_1,b_2,\ldots) \subseteq \acl(C'b)$ such that $\tp(a_n/D)$ is stationary. Put
$b':=(b,b_1,b_2,\ldots)$.
Now construct a tower $C_0=C' \subseteq C_1 \subseteq C_2\subseteq\ldots$ of finitely generated extensions of $C'$ within
$\acl(C')$ such that for each $m$, $\tp(b_m/bb_1\ldots b_{m-1}, C_m)\vdash \tp(b_m/b b_1\ldots b_{m-1},\acl(C'))$; this
 is possible for as
$\tp(b_m/C'bb_1\ldots b_{m-1})$ is algebraic, it is isolated. Put $C'':= \bigcup_{i\geq 0} C_i$.
Now $\tp(b'/C'')\vdash \tp(b'/\acl(C'))\vdash \tp(b'/\acl(C''))$, so $\tp(b'/C'')$ is  stationary. Also, $\tp(a_n/C''b')$
 is stationary, since $C'b'\subseteq C''b'\subseteq \acl(C'b')$ and $\tp(a_n/C'b')$ is stationary.
Thus, by Lemma~\ref{stationary-t}, $\tp(b'a_n/C'')$ is stationary, and hence so is
$\tp(ba_n/C'')$. 
 \qed

\begin{corollary}\label{ctbl2} Let $a$ be any finite tuple, and $C$ any base structure. 
  Then $\St_C(a)$ is countably generated over $C$.  \end{corollary}
 
{\em Proof.}  
{\em Claim.}  There exists a finite $b \in \St_C(a)$ with
$\St_C(a) \subseteq \acl(Cb)$.  

{\em Proof of Claim.}  Using Proposition~\ref{getun}, it suffices to show that if
  $d$ lies in a unary torsor over $Cb$ then $\St_C(bd) \subseteq 
\acl(\St_C(b),e)$ for some finite $e$.  If $d \in \G$ or if $\tp(d/C(b))$ is generic
in some open   or properly $\infty$-definable $C(b)$-torsor, then $\St_C(bd) \subseteq 
\acl(\St_C(b))$ by Lemma~\ref{notclosed}.  If    $d$ is generic in a closed $C(b)$-definable torsor $Y$, let
$\bar{d}$ be the the image of $d$ in $\res(Y)$.   Then $d$ is generic in an open 
$C(b,\bar{d})$-torsor (the 
inverse image of $\bar{d}$), so  $\St_C(bd) \subseteq 
\acl(\St_C(b),\bar{d})$ by the previous case; let $e = \bar{d}$. 

Now let $b$ be as in the claim, and replace $C$ by $\dcl(Cb)$.  In this way
we may assume $\St_C(a) \subseteq \acl(C)$.  On the other hand 
by Corollary~\ref{ctbl}, $\dcl(C(a)) \meet \acl(C) \subseteq \dcl(C,Z)$
for some countable $Z \subseteq \acl(C)$.  
So $\St_C(a) \subseteq \dcl(C,Z)$.  
For any tuple $z$ from $Z$,
$\tp(z/\St_C(a))$ is isolated, so there exists a finite $z' \in \St_C(a)$ such that
$\tp(z/z') \vdash \tp(z/C,\St_C(a))$.  Let $Z' = \{z': z \in Z^m, m=1,2,\ldots \}$.  Then
$Z'$ is countable, and $\tp(Z/C,Z') \vdash \tp(Z/C,\St_C(a))$.  So $\St_C(a) \subseteq \dcl(C,Z')$.
\qed

\begin{proposition}\label{ACVFBS}
  Condition (BS) holds in ACVF, so the conclusion of
Proposition~\ref{strongcode2} holds.\index{BS!for ACVF}
\end{proposition}

{\em Proof.} Immediate from Corollary~\ref{ctbl2} and  Lemma~\ref{boundedm} \qed
 
\medskip

\chapter{Types orthogonal to $\Gamma$}\label{orthogonal}

In this chapter, we develop a theory 
of orthogonality to $\Gamma$ for $n$-types in the sorts of $\GG$. This extends 
Definition~\ref{orth1}, which gives a notion of orthogonality to $\Gamma$   for unary types.
At first sight, Definition~\ref{ortho2} appears, as with sequential independence,  to be
dependent on the choice of a generating sequence.  We shall show that it is independent
of the choice, by proving that orthogonality  to $\Gamma$ is 
 the same as stable domination, 
in ACVF. We also extend the notion of the {\em resolution} of a lattice. Over a parameter set $C$, it is possible
that a lattice is defined, but does not contain any $C$-definable elements. We show how to add a basis of the lattice to $C$, without 
increasing the definable closure in either the value group or the residue field. We conclude the chapter by proving that, over the value
group, an indiscernible sequence is an indiscernible set  (Proposition~\ref{indiscseq}).

Some of  the main results of this chapter will be superseded later on.  These resolution
statements, in the finitely generated case, will be majorized in Chapter 11 by a theorem yielding
a canonical minimal resolution.  The equivalence of 
orthogonality to $\G$ with stable domination (Theorem~\ref{invarstabdom}) is a special case of 
 Theorem~\ref{fulldom}, which will be proved independently. The present hands-on proof has
 the merit of giving more information, 
 when an invariant type is {\em not} orthogonal to $\Gamma$, of
 the nature of the base change needed to see the non-orthogonality.  
 Theorem
 ~\ref{indiscseq} can also be viewed in this light, since the same result over a 
 large model will become, in Chapter 12, an immediate consequence of the theory of stable domination.

  \begin{definition}\label{ortho2}\index{orthogonal!to value group} \rm
If $a $ is a unary sequence (possibly transfinite), we say that
$\tp(a/C)$ is {\em orthogonal} to $\Gamma$ (written
$\tp(a/C) \perp \Gamma$) if, for 
any model $M$ with $a \dnf^g_C M$ and $C \subseteq \dcl(M)$,
we have $\Gamma(M)=\Gamma(Ma)$.

\end{definition}

Note that for any $C\subseteq C'\subseteq \acl(C)$, $\tp(a/C)\perp \Gamma$ if and only if $\tp(a/C')\perp \Gamma$.
The above definition is the natural notion of orthogonality which generalises that of stability, and is based on sequentially independent extensions (for which we have existence and uniqueness).
We will slightly extend this definition in \ref{moreorthog}, where it is given as a condition on $\acl(Ca)$ ($a$ any unary sequence) and is independent of the choice of $a$.  Recall from 
Lemma~\ref{eqorthog} that a {\em unary} type is orthogonal to $\Gamma$ precisely if it is the generic type of a closed 1-torsor.

\begin{lemma}\label{ortheasy} 
(i) Let $a=a_1a_2$ where $a_1,a_2$ are tuples, and suppose 
$\tp(a/C)\perp\Gamma$.
Then $\tp(a_1/C)\perp \Gamma$.

(ii) Suppose $\tp(a_1/C)\perp \Gamma$ and 
$\tp(a_2/Ca_1)\perp \Gamma$.
Then $\tp(a_1a_2/C) \perp \Gamma$.

(iii) Suppose $C \subseteq B$ and $a\dnf^g_C B$. Then
$\tp(a/C) \perp \Gamma$ if and only if $\tp(a/B) \perp \Gamma$.

(iv) Suppose that for each $i=1,\ldots n$, $a_i$ is algebraic or 
generic in a closed unary set over $Ca_1\ldots a_{i-1}$. Then 
$\tp(a/C) \perp \Gamma$.

(v) If $\tp(a/C) \perp \Gamma$ then $\Gamma(C)=\Gamma(Ca)$. 
\end{lemma}

{\em Proof.} 
(i)  Suppose
$a_1\dnf^g_C M$, and choose $a_2' 
\equiv_{Ca_1} a_2$
with $a_2'\dnf^g_{Ca_1} M$. Then as 
$a_1a_2'\dnf^g_C M$
we have $\Gamma(Ma_1a_2')=\Gamma(M)$, so 
$\Gamma(Ma_1)=\Gamma(M)$.

(ii) Choose $M$ so that $a_1a_2 \dnf^g_C M$. Then
$a_1 \dnf^g_C M$ and $a_2\dnf^g_{Ca_1} M$,
so $\Gamma(M)=\Gamma(Ma_1)=\Gamma(Ma_1a_2)$.

(iii) We may suppose that $C=\acl(C)$ and $B=\acl(B)$. If $\tp(a/C)\perp \Gamma$ and $a\dnf^g_B M$ with $B\subseteq M$, then
as $a\dnf ^g_C B$ we have 
$a\dnf^g_C M$, so $\Gamma(M)=\Gamma(Ma)$. Thus,
$\tp(a/B)\perp \Gamma$.
Conversely, suppose $\tp(a/B)\perp \Gamma$ and 
$a\dnf^g_C M$, with $C\subseteq M$. There is $M'\equiv_B M$ with
$a\dnf^g_B M'$, so $a\dnf^g_C M'$, and by 
Corollary~\ref{leftright} we have
$aM\equiv_C aM'$. We have $\Gamma(M'a)=\Gamma(M')$, 
so $\Gamma(Ma)=\Gamma(M)$.

(iv) This is immediate from Lemma~\ref{eqorthog} applied stepwise, 
together with (ii).

 (v) Suppose $\gamma\not\in 
 \Gamma(C)$, with
$\gamma=f(a)$, say, where $f$ is a $C$-definable function. Let 
$\gamma'\equiv_C \gamma$. We must show $\gamma'=\gamma$. 
Choose a model $M\supset C$ with $a\dnf^g_C M$. As 
$\tp(a/C)\perp \Gamma$, 
$\Gamma(M)=\Gamma(Ma)$, so $\gamma \in M$, and hence
$a\dnf^g_C \gamma$. Choose $a'\equiv_C a$ with 
$a'\dnf^g_C \gamma \gamma'$.  Then $a\gamma\equiv_C 
a'\gamma'\equiv_C a'\gamma$,
so $f(a')= \gamma=\gamma'$. Hence $\gamma \in \Gamma(C)$. \qed

\bigskip

Part (iv) above and Lemma~\ref{eqorthog} yield in particular that if $C\leq A$ are 
algebraically closed valued fields with $A=\acl_K(Ca)$, and 
$\trdeg(A/C)=\trdeg(k(A)/k(C))$ (informally, all of the extension is in the 
residue field), then $\tp(a/C)\perp \Gamma$.
The converse is false, as shown in Example~\ref{promex}. It is also shown,
 in Example~\ref{orthsing}, that it is possible that $\tp(a/C)\perp \Gamma$ for
 all singletons 
$a\in \acl(A)\cap K$, even though $\tp(A/C)\not\perp\Gamma$ (in the sense of Definition~\ref{moreorthog} below).

\bigskip

Recall from Chapter~\ref{acvfbackground}  the notion of a {\em generic basis} of the lattice with code 
$s\in S_n\cap \dcl(A)$. We say that a lattice $\Lambda(s)$ defined over $A$ is {\em resolved}\index{resolved} 
in $B \supset A$ (or just that $s$ is {\em resolved} in $B$) if $B$ contains a basis for $\Lambda(s)$. If 
$t$ codes an element $U(t)\in \red(s)$, then $t$ is {\em resolved} in $A$ if $s$ is resolved in $A$ and 
$A\cap K^n$ contains an element of $U(t)$. We say $A$ is resolved if all elements of 
$A$ are resolved.

\begin{definition}\index{resolution!generic closed} \rm Let $A$ be an $\L_{\GG}$-structure. Then a 
{\em generic closed resolution} of $A$ is a structure 
$B=\acl(A\cup\{b_i:i<\lambda\})$, where, for each $i<\lambda$, $b_i$ is a 
generic basis of some lattice defined over 
$\acl(Ab_j:j<i)$, and each lattice of $B$ has a basis in $B$.
\end{definition}

\begin{remark} \rm \label{resexist}

(i) If $A=\acl(A)$, then $A$ is resolved if and only if $A=\dcl(A\cap K)$. For example, if $A$ is resolved
and $t,s\in A$ with $t$ a code for $U(t)\in \red(s)$, then $U(t)$ contains some $(a_1,\ldots,a_n)\in K^n$, so is defined as
$U(t)=(a_1,\ldots,a_n)+\M \Lambda(s)$.

(ii)  If $A$ is an $\L_G$-structure then $A$ 
has a generic closed resolution
$B=\acl(A\cup\{b_i:i<\lambda\})$, and $B$ may be chosen so that $\lambda\leq |A|$. 
Furthermore, $\Gamma(B)=\Gamma(A)$, by Lemma~\ref{3.1.1}.

(iii) More generally, suppose $(s_1,\ldots,s_m)$ is a tuple of codes for lattices in $A$. 
We definably identify $(s_1,\ldots,s_m)$ with a single $s$ which is a code for the lattice
$\Lambda(s_1) \times \ldots \times \Lambda(s_m)$. Let 
$b$ be a generic basis for $\Lambda(s)$ over $A$ in the sense of Chapter~\ref{acvfbackground}. Then
$\tp(b/A) \perp \Gamma$ (treating $b$ as a tuple of field elements). For suppose 
$b\dnf^g_A M$ for some model $M$. By Proposition~\ref{gtostab}, $\St_A(Ab) \dnf_A \St_A(M)$. 
It follows (from the definition of generic basis) that $b$ is a 
generic basis of $\Lambda(s)$ over $M$, and hence, by (ii) above, 
$\Gamma(Mb)=\Gamma(M)$.

\end{remark}

Let $A=\acl(A)$, and suppose that all lattices of $A$ are resolved in $A$. Then for any
lattice $\Lambda(s)$ of $A$, the $k$-vector space $\red(s)$ has a basis in $A$; namely, 
the set $\{\red(b_1),\ldots,\red(b_n)\}$, where $\{b_1,\ldots,b_n\}$ is a basis 
of $\Lambda(s)$. Hence, there is an $A$-definable bijection $\red(s) \rightarrow k^n$. 
Let $k_{\rep}(A)$ 
be the sub-field of $k(A)$ consisting of $\res(0)$ and elements 
$\res(a)$ where $a\in A$ and $|a|=1$. Let ${\cal B}$ be a transcendence basis of
$k(A)$ over $k_{\rep}(A)$, and let $J\subset K$  with $|J|=|{\cal B}|$ and 
${\cal B}=\{\res(c):c\in J\}$. 
Then $\acl(A\cup J)$ is called a
{\em canonical open resolution of $A$}. \index{resolution!canonical open}

\begin{lemma} \label{canopen}
Let $A=\acl(A)$, and suppose all lattices of $A$ are resolved.
 Then $A$ has a canonical open resolution $B$, which is unique up to 
$A$-isomorphism, and
$\Gamma(A)=\Gamma(B)$. We have $B \subseteq \dcl(B \cap K)$. 
\end{lemma}

{\em Proof.} 
We adopt the notation above, with $B=\acl(A\cup J)$. The existence part is immediate.

Let ${\cal B}=\{\alpha_j:j\in J\}$.
Put $E:=A \cap K$. Then $A=\acl(E \cup {\cal B})$:
 for elements of $S_n \cap A$ lie in $\dcl(E)$ by assumption, and elements of
 $T_n \cap A$ lie 
in
$\dcl(E \cup {\cal B})$ by the argument before the lemma. Also, 
if $d_j$ is a field element chosen in $\alpha_j$ for each $j\in J$, then
for each $j$, $\alpha_j$ is a generic element of $k$ over
$E\cup ({\cal B}\setminus \{\alpha_j\}) \cup \{d_i:i<j\}$. It follows that 
$\Gamma(A \cup\{d_i:i<j\})=\Gamma(A\cup \{d_i:i\leq j\})$, so
$\Gamma(A)=\Gamma(B)$. 

The uniqueness assertion is also clear: the type of $(d_j:j\in J)$
over $A$ does not depend on the choice of ${\cal B}$, essentially because the 
type of a generic sequence from $R$ is uniquely determined.
\qed

\bigskip

Before making the general connection, for invariant types, between stable domination
and orthogonality to $\Gamma$, we work out an example. It will also be used in the proof.

\begin{example} \label{invarstabdomex}  Let $p=\tp(a/\U)$ be an $\Aut(\U/C)$-invariant type
of a single field element $a$. Suppose that
 $\Gamma(\U)=\Gamma(\U a)$. Then $p$ is the generic type of a closed ball, defined over $C$. \end{example}
 
{\em Proof.} Let $M$ be a small model, $C \subseteq M $.  Now $a$    is generic over $M$ in the intersection $V$ of a 
chain $(U_i:i\in I)$ of 
$M$-definable balls. We first show that $V$ is closed.

Suppose first that $V$ contains no proper $M$-definable ball.  If  
$a \in B$, where $B$ is a $\U$-definable ball, and $B$ is a proper sub-ball of $V$,
then by invariance $a \in B'$ for any $\Aut(\U/C)$-conjugate $B'$ of $B$.  But it is easy to find a $C$-conjugate
$B'$ of $B$ with $B \cap B' = \emptyset$.  This contradiction shows that there is no 
such ball $B$, so that $a$ is generic in $V$ over $\U$.  
This implies that $V$ is closed; for by saturation $\U$ contains some field element $d\in V$, and if $V$ is not closed then 
$|d-a|\in \G(\U a)\setminus \G(\U)$, which is impossible.

In the other case, $V$ contains a
proper $M$-definable ball, hence a point
$e \in M$. Let $\gamma = |e-a|$, and let $B_{\leq \gamma}(e) = \{x: |x-e| \leq \gamma \}$.  
Then $\gamma \in \G(Ma) \subset \G(\U)$, and $\gamma$ is fixed by $\Aut(\U/M)$, so $\gamma \in M$. 
 Thus $B_{\leq \gamma}(e)$ is $M$-definable;
since $a \in B_{\leq \gamma}(e)$ we must have $V=B_{\leq \gamma}(e)$, so again, $V$ is closed.  

We now show $V$ is $C$-definable.
 If $a \in V'$ for some maximal 
open sub-ball $V'$ of $V$, with $V'$ defined over $\U$
 then by definition of $V$, $V'$ cannot be defined over $M$; so $V'$ has an $\Aut(\U/M)$-conjugate 
$V'' \neq V$; but then $V' \cap V'' = \emptyset$, a contradiction as $a\in V''$ by invariance.  Thus
 $a$ is generic in the closed ball $V=B_{\leq \gamma}(e)$ over $\U$.   Finally, $V$ is unique
 with this property, so $V$ is $\Aut(\U/C)$-invariant and hence $C$-definable.
\qed

\bigskip

In the next theorem, the base change results for stable domination in Chapter 4 are crucial.

\begin{theorem} \label{invarstabdom}
Let $p$ be an $\Aut(\U/C)$-invariant $*$-type (so possibly in infinitely many variables). Suppose that for any model $M\supseteq C$ and $a\models p|M$ we have 
$\Gamma(M)=\Gamma(Ma)$.
 with  Then $p|C$ is stably dominated.
\end{theorem}

{\em Proof.} 
We shall show that if $C\subseteq M$ and $a\models p|M$  then $\tp(a/M)$ is
 stably dominated. For then by 
Theorem~\ref{stab7},  $\tp(a/C)$ is stably dominated.

We first show how to replace $a$ by a tuple of field elements.  
Choose a small model $M$ with $C\subseteq M$.
Let $a':=aa_1$ be a  sequence such that $\acl(M a')$ is
 a generic closed resolution of $a$ 
over $M$. We choose $a_1$ to contain, for each lattice with code $s$ in the sequence $a$, a basis 
of $\Lambda(s)$ which is generic over $M$. Note that $\tp(a'/M)$ has a canonical
$\Aut(\U/M)$-invariant extension $p'$ over $\U$ whose restriction to the $a$-variables is $p$.
Put $A:=\acl(M a')$. Now, as in the paragraph before 
Lemma~\ref{canopen},
 let $k_{\rep}(A)$ be the subfield of $k(A)$ consisting of $\res(0)$ and elements 
$\res(x)$ where $x\in A$ and $|x|=1$. Let $\B,J$ be as in the definition of canonical 
open resolution of $A$. Let $a_2$ list $J$, and put $a^*:=a'a_2$.
Then $\acl(M a^*)$ is a canonical open resolution of $A$. 
All elements in $\S$ in 
$\acl(M a')$ are resolved, and $\Gamma(M)=\Gamma(M a')$. Also,  $\acl(M a^*)$ is resolved.
Let $p^*$ be any $\Aut(\U/M)$-invariant extension over $\U$ of $\tp(a^*/M)$ whose
 restriction to the $aa'$-variables is $p'$. By omitting some elements of $a^*$ if necessary, we may suppose that
$a^*$
 is a sequence of field elements (with $a\in \dcl(Ma^*)$). 
 Now suppose $a^*\models p^*|M'$ for some model $M'\supseteq M$.
Then $k(Ma')$ and $k(M')$ are independent in $\VS_{k,M}$. 
It follows that $\acl(M' a^*)$ is a canonical open resolution of $\acl(M' a')$, and that 
$\Gamma(M' a^*)=\Gamma(M' a')=\Gamma(M' )$.

 We shall show that $\tp(a^*/M)$ is 
stably dominated. For then $\tp(a/M)$ 
is stably dominated by Proposition~\ref{bits}(ii), and this suffices as described above.
 Put $A^*:=\acl(Ma^*) \cap K$.

To show that $\tp(a^*/M)$ is stably dominated, we must show that if $N$ is a model of finite
 transcendence degree over $M$,
and $k(A^*)$ and $k(N)$ are independent over $k(M)$, then $a^*\dnf^g_M N$; for this
 ensures, by Theorem~\ref{moreinv2},
that $\tp(N/Mk(A^*))\vdash \tp(N/A^*)$.
We prove this by induction on the transcendence degree $n$ of $N$ over $M$.

Suppose first that $n=1$, that is, $N=\acl(Mb)$ where $b\in K$. Then by
 Proposition~\ref{fieldsym}, 
we have $A^*\dnf^g_M N \leftrightarrow b\dnf^g_M MA^*$, and there is no 
dependence on the choice 
of generating sequence. By Proposition~\ref{noembed}, we may suppose $b$ 
 embeds into $A^*$ over $M$.  By Example \ref{invarstabdomex}, $b$ is generic in a closed ball over $M$. This closed ball is in $M$-definable bijection with $R$, so we may suppose that $b\in R$.
Since $k(Mb)$ and $k(A^*)$ are 
independent over $k(M)$, we have
$b\dnf^g_M A^*$, as required.

For the inductive step, suppose that $M\subset M'\subset N$ with $\trdeg(N/M')=1$. By 
induction,
$a^*\dnf^g_M M'$, and by Proposition~\ref{gtostab}, $k(A^*)$ and $k(M')$ are independent 
(in $\St_M$ over $k(M))$. 
Hence, $\trdeg(k(M'a^*)/k(M'))=\trdeg(k(Ma^*)/k(M))$. The latter equals
$\trdeg(k(Na^*)/k(N))$, by the assumption that $k(A^*)$ and $k(N)$ are independent over $k(M)$. Thus, $k(M'a^*)$ and $k(N)$ are independent over $k(M')$.
It follows that, by the $n=1$ case, 
$a^*\dnf^g_{M'} N$. Hence, by transitivity of $\dnf^g$, we have
$a^*\dnf^g_M N$ as required. \qed

\begin{corollary} \label{equivorthstabdom}\index{orthogonal!and stable domination}
Let  $a$ be a unary sequence over $C$ (possibly transfinite). Then $p:=\tp(a/C)$ is stably dominated if and only 
if it is orthogonal to $\Gamma$.
\end{corollary}

{\em Proof.} We may suppose $C=\acl(C)$, as both conditions are unaffected by this.  
Suppose first that $p$ is stably dominated.
Then by Proposition~\ref{definable}, $p$ 
has a unique $\Aut(\U/C)$-invariant
 extension $p|\U$. Let $M\supseteq C$ be a model.
If
$a\models p|M$, then $a\dnf^g_C M$ (see Proposition~\ref{gdom}). 
Thus, by uniqueness of sequential extensions, if  $a\dnf^g_C M$, then $a\models p|M$, a stably dominated type by
 Proposition~\ref{dropstabdom}. Hence, 
by Lemma~\ref{stabdomgamma}, $\tp(a/M)\vdash \tp(a/M\Gamma(\U))$. It follows that 
$\Gamma(M)=\Gamma(Ma)$.
Thus $\tp(a/C)\perp \Gamma$.

Conversely, suppose $p\perp \Gamma$. Let $p|\U$ be the $\Aut(\U/C)$-invariant extension given by sequential independence.
Then for any model $M\supseteq C$, if $a\models p|M$
we have $a\dnf^g_C M$, so $\Gamma(M)=\Gamma(Ma)$. It follows by 
Theorem~\ref{invarstabdom} that $p$ is stably dominated. \qed

\bigskip

Observe that in the above setting, $\tp(a/C)$ is stably dominated if and only if
$\tp(\acl(Ca)/C)$ is stably dominated (by Proposition~\ref{bits}(iii) and Corollary~\ref{uptoacl2}). 
We use this freely.

\begin{lemma} \label{orthogwelldefined}
Let  $\acl(Ca)=\acl(Ca')$ with $a,a'$ both unary sequences over $C$. Suppose that $\tp(a/C)\perp \Gamma$. 
Then 

(i) $a\dnf^g_C B$ if and only if $a'\dnf^g_C B$,

(ii) $\tp(a'/C)\perp \Gamma$.
\end{lemma}

{\em Proof.} We may suppose $C=\acl(C)$. Let $A:=\acl(Ca)$. By Corollary~\ref{equivorthstabdom},
$\tp(A/C)$ is stably dominated.
Hence, by Proposition~\ref{definable}, for any saturated model $M$ there is a unique 
$\Aut(M/C)$-invariant extension of $\tp(A/C)$ to $M$. 

(i) Suppose  $a\dnf^g_C B$. By Proposition~\ref{gtostab},  $A^{\st} \dnf_C B^{\st}$ 
and hence by stable domination, the unique $\Aut(\U/C)$-invariant  extension of $\tp(A/C)$ extends 
$\tp(A/B)$. 
There is $a''\equiv_C a'$ with $a''\dnf^g_C B$ such that $\tp(\acl(Ca'')/B)$ extends to some $\Aut(\U/C)$-invariant 
extension of $\tp(A/C)$ (Corollary~\ref{invariance}). Hence, by the above uniqueness, $a'\dnf^g_C B$.

Now suppose $a'\dnf^g_C B$. Since $A=\acl(Ca)=\acl(Ca')$, the same argument shows that
$A \dnf^g_C B$ via any generating sequence.

(ii) Suppose $a'\dnf^g_C M$. Then by (i) $a\dnf^g_C M$, so $\Gamma(Ma)=\Gamma(M)$.
But $\Gamma(Ma)=\Gamma(Ma')$. Hence $\Gamma(Ma')=\Gamma(M)$.
\qed

\bigskip

\begin{definition}\label{moreorthog} \rm
Let
 $A=\acl(Ca)$, with $\tp(a/C)$ unary. We write  $\tp(A/C)\perp \Gamma$  if $\tp(a/C)\perp\Gamma$. 
By the last lemma, this definition is independent of the choice of $\acl$-generating sequence $a$.
\end{definition}

\begin{proposition}  \label{perpgamma}
Suppose  
$\tp(A/C)\perp \Gamma$. Then the following are equivalent.

(i) $\St_C(A) \dnf_C \St_C(B)$.

(ii) $A\domind_C B$.

(iii) $A\dnf^g_C B$ via some generating sequence.

(iv)  $A\dnf^g_C B$.

(v)  $B\domind_C A$.

(vi) $B\dnf^g_C A$ via some generating sequence.

(vii)  $B\dnf^g_C A$.

\end{proposition}

{\em Proof.} We may suppose $C=\acl(C)$. First, by Corollary~\ref{equivorthstabdom}, $\tp(A/C)$ is stably dominated.

The equivalence (i) $\Leftrightarrow$ (ii) is by the definition of stable domination, and 
(iii) $\Leftrightarrow$ (iv) comes from Lemma~\ref{orthogwelldefined}.  
We get (ii) $\Rightarrow$ (iii) from Proposition~\ref{gdom} (as there is a unary 
sequence $a$
with $A=\acl(Ca)$); (iii) $\Rightarrow$ (i) is by Proposition~\ref{gtostab}.

The equivalence (ii) $\Leftrightarrow$ (v) is Proposition~\ref{symmetric}(i). To prove 
(v) $\Rightarrow$ (vi), suppose $B\domind_C A$, and  $B=\acl(Cb)$, $b$ a unary sequence. 
There is an automorphism $\sigma$
over $C$ with $\sigma(b) \dnf^g_C A$. By Lemma~\ref{gtostab}, 
we have $\St_C (\acl(C\sigma(b))) \dnf_C \St_C(A)$, which gives
$\sigma(B) \domind_C A$. It follows that $B\equiv_A \sigma(B)$, so
$b\dnf^g_C A$, as required. Next, we show (vi) $\Rightarrow$ (i). For
if $b\dnf^g_C A$ then by Proposition~\ref{gtostab} we have 
$\St_C(B) \dnf_C \St_C(A)$, which gives (i) by symmetry of stable forking.
Finally, (vii) $\Rightarrow$ (vi) trivially, and (vi) $\Rightarrow$ (vii) is immediate, 
for the first 6 conditions are 
equivalent, and in the implication
(v) $\Rightarrow$ (vi), the choice of $b$ was arbitrary.
\qed

\bigskip

We now extend the resolution results further, in order to obtain 
 \ref{indiscseq} below.

\begin{lemma}\label{kint}
Suppose $A=\acl(A) \cap K$, and ${\cal B} \subset k$. Then
$\acl(A \cup {\cal B}) \cap K=A$.
\end{lemma}

{\em Proof.} If $a\in \acl(A\cup {\cal B}) \cap K$, then by elimination of imaginaries 
in ACF, the set of conjugates 
of $a$ over $A \cup {\cal B}$ is coded by some tuple $a'$ from $K$. Now $a'$ lies in 
some $A$-definable $k$-internal subset
of $K^m$ for some $m$, and it follows from Lemma~\ref{2.6.2} that $a'\in A$.
\qed

\begin{lemma}\label{9.12ii}
 Suppose $C =\acl(C)\subseteq A=\acl(A)$, $\tp(A/C)\perp \Gamma$, and $C$ is resolved, and that 
 all elements of $S_n$ in $A$ are resolved and that $A$ is finitely $\acl$-generated over $C$. Let $B$ be a
canonical open resolution of $A$. Then $\tp(B/C)\perp \Gamma$.
\end{lemma}

{\em Proof.} 
In the notation before Lemma~\ref{canopen}, let ${\cal B}$ be a transcendence 
basis of $k(A)$ over
$k_{\rep}(A)$, and suppose ${\cal B}=\{\res(c):c\in J\}$. Let $B:=\acl(A\cup J)$.
We may suppose that ${\cal B}$ is finite. There is a finite sequence $a_1$ from $A_K$ with
$A_K=\acl(Ca_1) \cap K$. Let $a_2$ be an enumeration of ${\cal B}$. Then as in Lemma~\ref{canopen},
$A=\acl(Ca_1a_2)$. Clearly $a_1a_2$ is unary. By Lemma~\ref{orthogwelldefined}, $\tp(a_1a_2/C)\perp \Gamma$.
Let $M$ be a model, and suppose $a_1a_2J\dnf^g_C M$.
It suffices to show 
$\Gamma(M)=\Gamma(MB)$. Now $a_1a_2\dnf^g_C M$, so
$\Gamma(M)=\Gamma(MA)$ since $\tp(a_1a_2/C)\perp \Gamma$.

We claim that $\acl(BM)$ is a canonical open resolution of $\acl(AM)$; 
for then, by Lemma~\ref{canopen}, $\Gamma(MB)=\Gamma(MA)$, so $\Gamma(MB)=\Gamma(M)$.
We require that ${\cal B}$ is a 
transcendence basis of 
$k(AM)$ over $k_{\rep}(AM)$; that is, ${\cal B}$ is algebraically independent over
 $\{\res(x):x\in K \cap \acl(AM)\}$. By Lemma~\ref{kint}, 
$\acl(AM \cap K)=\acl(A_K\cup M_K \cup{\cal B}) \cap K=
\acl(A_KM_K) \cap K$. Also, $a_2\dnf^g_{A_K} A_KM_K$, and it follows that
the elements of $a_2$ are algebraically independent over $k_{\rep}(AM)$. This yields the 
claim, and hence the lemma.
\qed

\begin{lemma} \label{orthres}
Suppose $C=\acl(C)$ is resolved, and suppose
$\tp(a/C)\perp \Gamma$.
Then there
is a resolved $\L_{\GG}$-structure $B$ containing $Ca$, with 
$\tp(B/C)\perp \Gamma$.
\end{lemma}

{\em Proof.} Choose $A$ to be a generic closed resolution of $Ca$, 
and $B$ to be a canonical open resolution of $A$, and apply Remark~\ref{resexist}(iii) and
Lemma~\ref{9.12ii}. \qed

\bigskip

\begin{proposition}\label{twores}
Let $C$ be an $\L_{\GG}$-structure. Then

(i) there is a resolved $\L_{\GG}$-structure $C'\supseteq C$ with $\Gamma(C)
=\Gamma(C')$,

(ii) there is a resolved $\L_{\GG}$-structure $C''\supseteq C$ with
$k(\acl(C))=k(C'')$.
\end{proposition}

{\em Proof.} We may suppose that $C=\acl(C)$.

(i) Apply Remark~\ref{resexist}(ii) and Lemma~\ref{canopen}.

(ii) It suffices to show that we may resolve any set coded by an element  
of $S_n \cup T_n$ without increasing  $k(C)$. 

Let $s$ be a code for a lattice  of $S_n$. 
Then $\red(s)$ is an $n$-dimensional vector space over $k$.
We show that we can elements to $C$ to get a basis of $\red(s)$ to $C$ without increasing $k(C)$.
Let $U_1,\ldots,U_r$ be a maximal $k$-linearly  independent subset of  
$\red(s)$ with codes  in  $C$, and suppose that $r<n$. 
Choose $U$ generic in $\red(s)$ over $C$. We may suppose that $k(C) \neq k(C\lceil U\rceil)$, 
since otherwise we may add $\lceil U\rceil$ to $C$ and so get a larger $k$-linearly independent 
set in $\red(s)$. Thus there is a 
$C$-definable partial function $f$ with $f(\lceil U\rceil)\in k\setminus k(C)$. Let $X$ equal
$$
   \{U'\in \red(s): \lceil U'\rceil \in \dom(f) \mbox{~and~} U_1,\ldots,U_r,U'\mbox{~are linearly 
                             independent over~} k\}.
$$ 
Then $\{f(\lceil U'\rceil):U'\in X\}$ is an infinite definable 
subset of $k$, so as $k$ is strongly minimal, there is $U'\in X$ with 
$f(\lceil U'\rceil)\in k(C)$. Thus, we may add $\lceil U'\rceil$ to $C$ without increasing $k(C)$.

It remains to show that if $U\in \red(s)$ with $\lceil U\rceil\in \dcl(C)$, then there is
$a\in U \cap K^n$ with $k(C)=k(Ca)$. We prove, by induction on $n$, that this 
holds for any coset  in $K^n$ of an element with code in $T_n$ (that is, any $U+x$ where $\lceil U\rceil\in T_n$ and $x\in K^n$).  
So let $U$ be such a coset, and let $U_1$ be its projection to the first coordinate. 
Then $U_1$ is an open ball, so if $a_1$ is chosen generically in $U_1$ over $C$, 
then $k(C)=k(Ca_1)$ by Lemma~\ref{eqorthog}. Now $U^1:=\{x\in K^{n-1}:(a_1,x)\in U\}$ is a coset of an 
element with code in $T_{n-1}$. Hence, by induction, there is $b\in U^1$ with $k(Ca_1)=k(Ca_1b)$. 
Now $(a_1,b)\in U$ and $k(C)=k(Ca_1b)$, as required.     \qed

\bigskip

These results on resolutions will be extended further in Chapter~\ref{primemodels},
 where prime resolutions are considered.

The next proposition indicates another sense in which all instability arises from 
interaction with the value group.

\begin{proposition}\label{indiscseq}\index{indiscernible!in ACVF}
Let  $(a_i:i<\omega)$ be an indiscernible sequence  over  
$C\cup\Gamma(\U)$. Then $\{ a_i:i<\omega\}$ is an indiscernible set over $C\cup\Gamma(\U)$.
\end{proposition}

{\em Proof.} We may assume $C=\acl(C)$.
In the proof below, we often just write $\Gamma$ for $\Gamma(\U)$.

Since $(a_i:i < \omega)$ is indiscernible over $C \cup\Gamma$, 
we have
$$
  \Gamma(Ca_{i_1}\ldots a_{i_n})  = \Gamma(Ca_{j_1}\ldots a_{j_n}) $$ 
for any $n$ 
                   and any $i_1 < \ldots < i_n, j_1<\ldots < j_n$ .  
By adding all such sets to $C$, we may arrange that for all such $n$ and
$i_1<\ldots< i_n$ we have 
$$
\Gamma(Ca_{i_1}\ldots a_{i_n})=\Gamma(C).  \eqno(1)
$$
Also, for each $n$ there is a fixed $C_n\supset C$ such that
for all $ i_1<\ldots< i_n <j_1<\ldots <j_n$ (chosen in $\omega$) we have
$\acl(Ca_{i_1}\ldots a_{i_n})\cap \acl(Ca_{j_1}\ldots a_{j_n})=C_n$. By 
expanding $C$ we may suppose that $C_n=C$ for each $n\in \omega$; that is, for all
$n<\omega$ and $i_1<\ldots<i_n <j_1<\ldots<j_n$,
$$\acl(Ca_{i_1}\ldots a_{i_n})\cap \acl(Ca_{j_1}\ldots a_{j_n})=C. \eqno(2)$$
We preserve also
$$ (a_i:i< \omega) \mbox{~is an indiscernible sequence over~} C \cup \Gamma.\eqno(3)$$

We will show that $\tp(a/C)\perp  \Gamma$ for each finite subsequence $a$ of
$(a_i:i<\omega)$. 
We first argue that this suffices, so suppose it holds. 
Let $a'_i:=\VS_{k,C}(a_i)$ for each $i<\omega$. Then
$(a'_i:i<\omega)$ is an indiscernible sequence in the stable structure
$\VS_{k,C}$, so is an indiscernible set.  
It follows from this and $(2)$ that, for any $n$ and $i_1<\ldots < i_n$ and $j_1<\ldots 
< j_n$ with $\{i_1,\ldots i_n\} \cap \{j_1,\ldots,j_n\}=\emptyset$,
we have $\acl(Ca'_{i_1}\ldots a'_{i_n}) \cap \acl(Ca'_{j_1}\ldots a'_{j_n})=C$. 
It follows (as in Proposition~\ref{stab0.7}) that $\{a'_i:i\in \omega\}$ is an independent 
set in $\VS_{k,C}$.
Using orthogonality to $\Gamma$ and  Proposition~\ref{simplifications}(ii),  we obtain
$$\dcl(a'_{i_1},\ldots a'_{i_n}) =\VS_{k,C}(a_{i_1}\ldots a_{i_n})=
\VS_{k,C}(a_{i_1}\ldots a_{i_n}\Gamma),$$
(the second equality holds as any definable function from $\Gamma$ to the stable structure 
$\VS_{k,C}$ has finite range).
Hence, by indiscernibility,  for any 
$i,j\in \omega\setminus \{i_1,\ldots,i_n\}$,
$$
  \tp(a'_i/\VS_{k,C}(a_{i_1}\ldots a_{i_n}\Gamma)) 
=\tp(a'_j/\VS_{k,C}(a_{i_1}\ldots a_{i_n}\Gamma)).
$$
 It follows as in Remark~\ref{symmetry}  that
$$
   \tp(a_i/\VS_{k,C}(a_{i_1}\ldots a_{i_n}\Gamma)) =\tp(a_j/\VS_{k,C}(a_{i_1}\ldots 
a_{i_n}\Gamma)),
$$
and hence by Proposition~\ref{perpgamma} that
$$\tp(a_i/Ca_{i_1}\ldots a_{i_n}\Gamma)=\tp(a_j/Ca_{i_1}\ldots a_{i_n}\Gamma).$$ This yields the 
indiscernibility claimed in the theorem.

It remains to prove that $\tp(a/C)\perp  \Gamma$ for each such $a$. For this,
 it suffices to show 

\medskip
{\em Claim.} 
Suppose that $C'=\acl(C')\supseteq C$, and $(1')$, $(2')$ and $(3')$ hold, where these are 
(1), (2), and (3) with the base $C$ replaced by $C'$.
Suppose also  
$b$ is chosen in $K$ so that for each $n$ and $i_1<\ldots <i_n$ we have
$a_{i_1}\ldots a_{i_n}\dnf^g_{C'} b$. Then $(1'')$, $(2'')$, and $(3'')$ hold, where 
these are (1), (2), and (3) respectively, with the base $C$  
replaced by $C''=\acl(C'b)$. 

\medskip

Given the claim, suppose $M$ is a model with $a\dnf^g_C M$ for 
each finite subsequence  $a$ of $(a_i:i\in \omega)$. We may suppose $M=\acl(Cb)$ for some 
finite sequence 
$b$ of field elements. By applying the claim repeatedly, we find 
$\Gamma(M)=\Gamma(Ma)$ for each finite 
subsequence $a$ of $(a_i:i<\omega)$. This proves $\tp(a/C) \perp \Gamma$.

\medskip

{\em Proof of Claim.}

We first prove $(1'')$.
This certainly holds if $b$ is a generic element of a closed 1-torsor over $C'$. 
For then $\tp(b/C')\perp \Gamma$, so $b\dnf^g_{C'}\{a_i:i< \omega\}$, by 
Proposition~\ref{perpgamma}. Hence,
 $$
  \Gamma(C')=\Gamma(C'a_i:i<\omega)=\Gamma(C'ba_i:i<\omega)=\Gamma(C'b). \eqno(4)
$$
It remains to consider the case when $b$ is a generic element of an open 1-torsor
$U$ over $C'$ (or of an $\infty$-definable 1-torsor --- this case is similar). 
First, if not all elements of $U$ have the same type over $C'$, then, as any $C'$-definable subset of $U$ is a finite union
of $C'$-definable Swiss cheeses, some proper subtorsor of $U$
is algebraic over $C'$. Hence
$\Gamma(C'b)\neq\Gamma(C')$ by Lemma~\ref{2.5.6}, so $\rk_{{\mathbb Q}}(\Gamma(C'b))=\rk_{{\mathbb Q}}(\Gamma(C'))+1$ by 
Lemma~\ref{incbyone}.
In this case $(1'')$ follows from $(1')$ by rank considerations (Lemma~\ref{incbyone}). Thus,
 we assume
 all elements of $U$ have the same type over $C'$. We may suppose they do not all have the 
same type over
$C' \cup \{a_i: i < \omega\}$, since otherwise (4) is again valid. Hence,
there are $i_1,\ldots,i_n$ and an
$a_{i_1}\ldots a_{i_n}$-definable proper subtorsor $U(a_{i_1},\ldots,a_{i_n})$
of $U$ which intersects the $\tp(b/C')$ non-trivially. By $(1')$, the
radius of $U(a_{i_1},\ldots,a_{i_n})$ is in $C'$, and so also is the
radius of $U(a_{i'_1},\ldots,a_{i'_n})$ for any other $a_{i'_1},\ldots,a_{i'_n}$
from the sequence. So we have a uniform family of definable subtorsors of an
open 1-torsor, all of the same radius. The distance between two such torsors is
a $C'$-definable function of $2n$ variables, hence by $(1')$ has 
constant value, say $\gamma$, and $\gamma\in C'$. But then all the 
torsors $U(a_{i'_1},\ldots,a_{i'_n})$, and in particular, $U(a_{i_1},\ldots,a_{i_n})$
are contained in a $C'$-definable closed 1-torsor properly contained in $U$.
 This contradicts the assumption that all elements of $U$ have the same 1-type over $C'$.

To see $(2'')$, suppose $e\in \acl(C'ba_{i_1}\ldots a_{i_n}) \cap
\acl(C'ba_{j_1}\ldots a_{j_n})$, where $i_1<\ldots <i_n<j_1<\ldots <j_n$.
 By increasing $n$ if necessary we may suppose
the number $m$ of conjugates of $e$ over $C'ba_{i_1}\ldots a_{i_n}$ is as small as possible.  
Then there is a definable function $f$ on
$\tp(b/C')$ taking $b$ to a code for the $m$-set $X$ consisting of realisations of 
$\tp(e/C'ba_{i_1}\ldots a_{i_n})$. By minimality of $m$, $X$ is also a complete
 type over $C'ba_{j_1}\ldots a_{j_n}$; otherwise $e$ has fewer than $m$ conjugates over $C'ba_{i_1}\ldots a_{i_n}a_{j_1}\ldots a_{j_n}$.
Hence by indiscernibility (over $C'b$) $X$ is a complete type over any such 
$C'ba_{\ell_1}\ldots a_{\ell_n}$.
Thus, $\lceil X \rceil \in\dcl (C'ba_{\ell_1}\ldots a_{\ell_n})$ for any $\ell_1<\ldots < \ell_n$.
Thus, there is a definable function $f$ on $\tp(b/C')$ with $f(b)=\lceil X\rceil$, 
and $\lceil f\rceil\in
\acl(C'a_{i_1}\ldots a_{i_n}) \cap \acl(C'a_{j_1}\ldots a_{j_n})=\acl(C')$,
so $e\in \acl(C'b)$, as required.

Finally, for $(3'')$, first observe that by uniqueness of  sequentially independent extensions
 (Corollary~\ref{leftright}), $(a_i:i<\omega)$ is an indiscernible sequence over $C''$. Condition $(3'')$ now follows easily from $(1'')$.
\qed

\chapter{Opacity and prime resolutions}\label{primemodels}

In this chapter we extend the results on resolutions from Chapter~\ref{orthogonal}. In particular, 
we show (Theorem~\ref{R7}) that over an algebraically closed base in the field sorts, any finite tuple 
has a unique minimal atomic `prime' resolution, which does not extend either the value group or the residue field. 
We talk of `prime resolutions' rather than `prime models' since the resolution may not be a model of ACVF, as the 
valuation may be trivial on it.

 The first few lemmas (\ref{R1}--\ref{R8}) are proved in complete generality.
We suppose $M$ is a sufficiently saturated structure, and emphasise that $M$ consists of  elements of the home sort,
 that is, we are not identifying $M$ with $M^{\eq}$. Unless otherwise specified, parameter sets and tuples below 
live in $M^{\eq}$. In the application to ACVF, $M$ will be the sort $K$ of field elements. We work
 in the model $M$, over a parameter base $C$.

\begin{definition} \rm \label{resolution} Let $A$  be a substructure of $M^{\eq}$.
A {\em pre-resolution} \index{pre-resolution} of $A$ is a substructure $D$ of $M^{\eq}$ with 
$A \subseteq \acl(D \meet M)$. 

A {\em resolution} \index{resolution} of   $A$ is an
algebraically closed 
structure   $D \subseteq M$ such that $A \subseteq \dcl(D)$. 
A resolution $D$  is {\em prime} \index{resolution!prime} over $A$
if $D$ embeds over $A$ into any other resolution.
\end{definition}

Usually, $A$ will have the form $Ca$ where $C \subseteq M$, and $a$ is a finite sequence
of imaginaries.  

In the intended application to ACVF, the sequence $a$ may be taken from $\S \cup \T$. We shall exploit the fact 
that $S_n$ and $T_n$ can be viewed a coset spaces of groups of upper triangular matrices. Such groups have 
sequences of normal subgroups whose quotients are strongly minimal (the additive or multiplicative group of the 
field, in the case of $S_n$). This leads us to the following sequence of definitions and lemmas.

\begin{definition} \rm
Let $E$ be a definable equivalence relation on a definable set $D$
(in $M$). We say $E$ is {\em opaque} \index{equivalence relation!opaque}\index{opaque!equivalence relation} 
if, for any definable $Z\subset D$, there are $Z'$ (a union of $E$-classes) and $Z''$ (contained in the union of
 finitely many $E$-classes) such that $Z=Z' \cup Z''$. 
\end{definition}

For example, in ACVF, the partition of $R$ into cosets of $\M$ is an opaque equivalence relation; the same holds for 
the partition of any closed ball into open sub-balls of the same radius.

\begin{lemma}\label{R1}
Let $E$ be a $C$-definable opaque equivalence relation on $D$, and 
let $a=b/E$ be a class of $E$. Then either $a\in \acl(C)$, or 
$\tp(b/Ca)$ is isolated, where the isolating formula $\phi(y)$
states that $y/E=a$.
\end{lemma}

The terms `isolated type' and likewise, below, `atomic' are interpreted with respect to $\Th(M)$.
Thus, $\tp(b/Ca)$ might be isolated but not realised in some resolution of 
$Ca$ (e.g. if the latter is trivially valued, in the ACVF context).

\medskip

{\em Proof.} If $\phi(y)$ does not isolate a complete type over $Ca$, 
there is a formula $\psi(y,x)$ over $C$ such that $\psi(y,a)$ defines a 
proper non-empty subset of $a$. Let $Z$ be defined by $\psi(y,y/E)$. 
Then $Z$ meets the class $a$ in a proper non-empty subset. Since $E$ is 
opaque, there are just finitely many such classes. Hence $a\in \acl(C)$.
\qed

\begin{lemma} \label{R2}
Let $a_0,\ldots,a_{N-1}$ be a sequence of (imaginary) elements. Assume that $a_n=b_n/E_n$ for each 
$n\in \{0,\ldots,N-1\}$, where $E_n$ is an opaque equivalence relation defined over
$A_n:=C \cup \{a_j:j<n\}$. Define $B_n$ by downward recursion:   $B_N=\emptyset$
and $B_n:=B_{n+1} \cup\{b_n\}$ if  $a_n\not\in \acl(B_{n+1} \cup A_n)$,
and $B_{n}:=B_{n+1}$ otherwise. Let $I:=\{n: b_n\in B_0\}$
(which equals $\{n: a_n\not\in \acl(B_{n+1} \cup A_n)\}$). Then

(i) $A_N \subset \acl(B_0\cup C)$,

(ii) if $n\in I$ then $\tp(b_n/A_N \cup B_{n+1})$ is isolated.

(iii) $B_0$ is atomic over $A_N$.
\end{lemma}

{\em Proof.} (i) We have for each $n$
$$  a_n\in \acl(B_n \cup A_n).\eqno(1)$$
If $n\not\in I$, this follows by definition as 
$a_n\in \acl(B_{n+1}\cup A_n)\subset \acl(B_n \cup A_n)$. If $n\in I$,
then $a_n =b_n/E_n\in \acl(B_n)$.

Also, for each $k$, we have
$$~~\mbox{~if~} n\geq k \mbox{~then~} a_n \in \acl(B_k \cup A_k).\eqno(2)$$
This is proved by induction on $n \geq k$. If $n\geq k$ then by (1) 
(as $B_k \supseteq B_n$)
we have
$a_n\in \acl(B_k\cup A_n)$. By induction, 
$A_n \subset \acl(B_k \cup A_k)$, so
$a_n \in \acl(B_k \cup A_k)$.

Applying (2) with $k=0$ we have (i).

(ii) Let $n\in I$. We apply Lemma~\ref{R1} with $C$ replaced 
by $A_n \cup B_{n+1}$
to obtain that $\tp(b_n/B_{n+1} \cup A_n \cup\{a_n\})$ is isolated. 
By (2), $A_N \subset \acl(B_{n+1} \cup A_n \cup\{a_n\})$. Hence (ii) holds.

(iii) This follows immediately from (ii). \qed

\begin{definition} \rm (i) We say that $\tp(a/C)$ is {\em opaquely layered} \index{opaquely layered}
(or {\em $a$ is opaquely layered over $C$}) if there exist
$a_0,\ldots,a_{N-1}, E_0,\ldots, E_{N-1}$ satisfying the hypotheses of Lemma~\ref{R2} such that
$\dcl(C,a)=\dcl(C,a_0,\ldots,a_{N-1})$.

(ii) Let $E$ be a definable equivalence relation on a definable set $D$. \index{equivalence relation!opaquely layered}
Then $(D,E)$ is {\em opaquely layered over $C$} if $D,E$ are defined over $C$ 
and for each $a\in D$,
the imaginary $a/E$ is opaquely layered over $C$. We say $(D,E)$ is {\em opaquely layered everywhere}
if for any $C'$ and any pair $(D',E')$ definable over $C'$ and definably 
isomorphic (over some $C''$) to $(D,E)$, the pair $(D',E')$ is opaquely layered over 
$C'$.

(iii) Let $G$ be a definable group, and $F$ a definable subgroup. \index{opaque!group} Then 
$G/F$ is {\em opaquely layered over $C$} (respectively, {\em opaquely layered everywhere}, 
respectively, {\em opaque}) if the pair $(G,E)$ is, where $Exy$ is the equivalence relation $xF=yF$.
\end{definition}

Note that opacity is preserved under definable bijections, so `opaque' implies
`opaquely layered everywhere'. Also, if $\tp(a/C)$ is opaquely layered and $C \subset B$, 
then $\tp(a/B)$ is opaquely layered.

\begin{lemma} \label{extra}
Suppose $\tp(a/C)$ and $\tp(b/Ca)$ are opaquely layered. Then $\tp(ab/C)$ is opaquely layered.
\end{lemma}

{\em Proof.} Opaque layering of $\tp(ab/C)$ is witnessed by concatenating a
 sequence for $\tp(a/C)$ (as in Lemma~\ref{R2}) and a sequence for $\tp(b/Ca)$. \qed

\begin{proposition}\label{R3}
Suppose that $C\subset M$ and $a$ is a finite sequence of imaginaries, and suppose that
$\tp(a/C)$ is opaquely layered. Let $a_0,E_0,\ldots, a_{N-1},E_{N-1}$ witness the opaque layering.
Then $Ca$  has a pre-resolution  $D$ 
 which is atomic over $Ca$ and embeds 
into $D'$ for any    $D'$   which contains $C$ and an element of each
 class $E_i$-class $a_i$. The embedding can be taken to be an elementary map.
\end{proposition}

{\em Proof.} Let $I$ and 
$b_i$  ($i\in I$) and $B_n$ ($n<N$) be as in Lemma~\ref{R2}. Put
$D: = B_0\cup C$.
By Lemma~\ref{R2}, $D$ is atomic over $Ca$, and $Ca \subseteq \acl(D)$.  

For the second part, let $D'$ be a resolution of $a$ over $C$ which contains an element $b_i'$ of each $E_i$-class $a_i$.
By part (ii) of Lemma~\ref{R2}, arguing
 inductively, the embedding property holds. More precisely, adopt the notation $A_n$, $B_n$ of ~\ref{R2}, and suppose an embedding
 $f_{n+1}:B_{n+1}\rightarrow D'$ has been defined over $Ca$. If $a_n\in \acl(B_{n+1}\cup A_n)$, then $B_n=B_{n+1}$ and 
we may put $f_n:=f_{n+1}$. Otherwise, all elements of the equivalence class $a_n$ have the same type over $CaB_{n+1}$, so 
we may extend $f_{n+1}$ to $f_n$ by putting $f_n(b_n)=b_n'$. Finally, $f_0$ is elementary, so extends to an embedding
$D\rightarrow D'$. \qed

\medskip

The following corollary applies for instance to $\Th(\Qq_p)$ (taken in the $\GG$-sorts.)

\begin{corollary}\label{R3.1}  Assume any definably closed  subset of $M$ is
an elementary submodel of $M$.  Let  $C\subset M$ and $a$  a finite sequence of imaginaries, and suppose that
$\tp(a/C)$ is opaquely layered.  Then $Ca$ has an atomic pre-resolution
that embeds  into any pre-resolution; the embedding can be taken to be an elementary map.
\qed
 \end{corollary}

Define a {\em $\dcl$-resolution} \index{resolution!dcl} of $Ca$ to be a subset $D$ of $M$ with $Ca 
\subseteq \dcl(D)$. Note that (ii) below is true for any field; (i) holds
for $\Th(\Qq_p)$  for instance, as well as for ACVF if one works over a
nontrivially valued field.

\begin{corollary}\label{R3.2}  Assume:

(i)   any algebraically   closed subset of $M$ is
an elementary submodel of $M$.

(ii) $M$ admits elimination of finite imaginaries; i.e. any finite set  of
tuples of $M$ is coded by elements of $M$.

Let  $C\subset M$ and let $a$ be  a finite sequence of imaginaries, and  
suppose that
$\tp(a/C)$ is opaquely layered.    Let $a_0,E_0,\ldots, a_{N-1},E_ 
{N-1}$ witness the opaque layering. Then $Ca$ has an atomic 
$\dcl$-resolution  $B$.
Furthermore, $B$ embeds into any $B' = \dcl(B')$ which contains $C$ and an element  of
each class $E_i$-class $a_i$. 

Assume in addition that  whenever an
$E_i$-class $Q$ is defined over  a set $C'$ and isolates a type $q$ over
$C'$, this type is stationary.  Then  $\tp (B/C) \vdash \tp(B/\acl(C))$.
  
\end{corollary}
{\em Proof.}  Let $D$ be as in Proposition \ref{R3};
so    $Ca \subseteq \acl(D \meet M)$, and $D/Ca$ is atomic.
By construction $D$ is obtained from $C$ by adding finitely many  
elements of $M$;
at each stage one adds an element of an $E_i$-class which isolates a  type
over $Ca$ and the preceding elements; by assumption, this type is
stationary.   It follows that $\tp(D/Ca)$ is stationary.

Now $\acl(D) \meet M$ is an elementary submodel of $M$; so $a \in \dcl
(e)$ for some $e \in \acl(D) \meet M$.  The orbit of $e$ over $Da$ is
finite, hence coded by some tuple $e' $ from $M $.  We have $e' \in
\dcl(Da)$, and $a \in \dcl(De')$.  Let $B = De'$.  Then $B$  is a
$\dcl$-resolution of $Ca$, and all the conditions hold:  $\tp(e'/Da)$ is
isolated and  stationary, since $e' \in \dcl(Da)$, so $\tp(B/Ca)$ is
isolated and stationary;  and if $B' = \dcl(B')$  contains $C$ and an
element of each class $E_i$-class $a_i$, then $D \subseteq B'$ and $a \in
E'$, so $e'  \in \dcl(Da) \subseteq B'$. \qed 


\begin{lemma} \label{R4}
Let $G$ be a $C$-definable group, and $N,H,F$  be $C$-definable subgroups. Assume 
that $N$ is normal, $N \cap H =\{1\}$,
$NH=G$. Also suppose $F=N_FH_F$, where 
$K_F:=K\cap F$ for any $K\leq G$.
Suppose that the coset space $H/H_F$ is opaquely layered over $C$, and for each $h\in H$, 
$N/(N\cap F^h)$ is opaquely layered over $C\cup\{\lceil hH_F\rceil\}$. Then  
$G/F$ is opaquely layered over $C$.
\end{lemma}

{\em Proof.} Let $g\in G$. Then there are unique $n\in N$ and $h\in H$ with $g=nh$. By
Lemma~\ref{extra} and our assumptions, it suffices to show that
$$\dcl(C,\lceil nhF\rceil)=\dcl(C,\lceil hH_F\rceil,\lceil n(N\cap F^h)\rceil),$$
where $F^h=hFh^{-1}$.

For the containment $\supseteq$, suppose that $nhF=n'h'F$, with $n,n'\in N$ and $h,h'\in H$. 
We must show $hH_F=h'H_F$ and $n(N\cap F^h)$
Then, first, $h'=hf''$ for some $f''\in H_F$: indeed, $n'h'=nhf$ for some $f\in F$, and
$f=n''f''$ where $n''\in N_F$ and $f''\in H_F$, hence $n'h'=nhn''f''=nhn''h^{-1}hf''$,
so $Nh'=Nhf''$. As $h',hf''\in H$ and reduction 
modulo $N$ is injective on $H$, we have $h'=hf''$. Thus, $hH_F=h'H_F$. 
Hence $\lceil hH_F\rceil \in \dcl(C,\lceil nhF\rceil)$. Also, 
$nhF=n'h'F=n'hf''F=n'hF$, so $nF^{h}=n'F^h$, and hence $n'(N\cap F^{h})=n(N\cap F^h)$. As $h'=hf''$, $F^h=F^{h'}$.
So
$\lceil n(N\cap F^h)\rceil \in \dcl(C, \lceil nhF\rceil)$.

For the containment $\subseteq$, suppose that $hH_F=h'H_F$ and $n(N\cap F^h)=n'(N\cap F^{h'})$, where
$n,n'\in N$ and $h,h'\in H$. The first equality yields that $hF=h'F$ so $n'hF=n'h'F$ and also $Fh^{-1}= Fh'^{-1}$; so
$F^h=F^{h'}$, and
$n(N\cap F^h)=n'(N\cap F^h)$ (by the second inequality). Hence, $n^{-1}n'\in F^h$, so $n'hF=nhF$, as required. \qed

\begin{corollary}\label{R5}
Let $\{1\}=G_0\leq G_1\leq \ldots \leq G_N=G$ and $H_i \leq G_{i+1}$ (for $i=0,\ldots N-1$)
be $\emptyset$-definable groups, with $G_i$ normal  in $G$, $G_i\cap H_i=\{1\}$ and
$G_iH_i=G_{i+1}$ for each $i$. Let $F$ be a definable subgroup of $G$ such that for each $i$, 
$G_{i+1}\cap F=(G_i \cap F)(H_i \cap F)$ and
$H_i/(H_i \cap F)$ is opaquely layered everywhere. Then $G/F$ is opaquely layered everywhere.
\end{corollary}

{\em Proof.} This follows from Lemma~\ref{R4} and induction on $i$. Observe that if 
$G_i/(G_i \cap F)$ is opaquely layered everywhere, then by normality of $G_i$, and since opaque layeredness everywhere is preserved by definable
 bijections,
$G_i/(G_i \cap F^h)$ is opaquely layered  everywhere for each $h\in G$. \qed

\bigskip

In the next lemma, by `group scheme' \index{group scheme} we mean a system of polynomial equations, viewed as 
a formula, which for any integral domain $A$ defines a subgroup of $B_n(A)$. Here, $B_n$ 
is the group scheme such that $B_n(A)$ is the group of upper triangular matrices over $A$ which are invertible 
in $\GL_n(A)$. As usual, $K$ denotes a large algebraically closed field. We denote by 
$(e_1,\ldots,e_n)$ the standard basis for $K^n$.

\begin{lemma} \label{R8}
Fix $n$, and let $N={n+1 \choose 2}$. 

(i) The group scheme $B_n$ has, for each $i=0,\ldots,N$, 
 group subschemes $G_i,H_i$ over ${\mathbb Z}$ such that for any
subring $A$ of $K$, the following hold:
$G_{i+1}(A)=G_i(A)H_i(A)$, $G_i(A)\cap H_i(A)=\{1\}$, $H_i$ is isomorphic to 
$(A,+)$ or the group of multiplicative units of $A$, 
$G_i(A)$ is normal in $B_n(A)$, $G_0=\{1\}$, and $G_N=B_n$. 

(ii) For each $i$, if $a\in G_i$, $b\in H_i$, and $ab$ fixes $e_n$ (acting by left
 multiplication)
then $a$ and $b$ both fix $e_n$.
\end{lemma}

{\em Proof.} The scheme $B_n$ is defined so that $B_n(K)$ is the stabiliser of a maximal flag
$$\{0\}=V_0 \subset V_1 \subset \ldots \subset V_{n-1} \subset V_n=K^n;$$ here $V_i$
is the space spanned by $\{e_1,\ldots,e_i\}$.
For $0\leq j\leq i\leq n$, let $G_{ij}(K)$ be the subgroup of $B_n(K)$ consisting of
 elements $x$ such that $x$ acts as the identity on $V_{i-1}$, and $x-I$
maps $V_i$ into $V_j$. Then
$G_{1,1}=B_n$, $G_{i,0}=G_{i+1,i+1}$, $G_{n,0}=\{1\}$, and
$$G_{n,0}<G_{n,1}<\ldots <G_{n,n}=G_{n-1,0}<\ldots <G_{i,0}<\ldots  
< G_{i,i}=G_{i-1,0}<\ldots<G_{1,1}=B_n.$$
The $G_{ij}$ are all normal in $B_n$ since they are defined in terms of the flag. Also, 
it is easy to see that
$G_{i,i}=G_{i,i-1}\unlhd H_{i,i-1}$, where $H_{i,i-1}$ is the group of maps fixing each 
$e_j$ ($j\neq i$) with eigenvector $e_i$
(so $H_{i,i-1}(A)$ is isomorphic to $G_m(A)$, the group of units of $A$). Likewise,
if $j<i$ then
$G_{ij}=G_{i,j-1} H_{i,j-1}$, where $H_{i,j-1}(A)$ is isomorphic to 
$(A,+)$, and consists of elements fixing all $e_k$
with $k\neq i$ and with $e_i \mapsto e_i+te_j$ ($t\in A$). 

In terms of matrices, $G_{ii}=G_{i-1,0}$ and $G_{ij}$ (for $j<i$) is the subset of $B_n(A)$ consisting 
of matrices which have top left $(i-1)\times (i-1)$-minor equal to $I_{i-1}$, and its $i^{\th}$ column
is the transpose of $(a_1,\ldots,a_n)\in A^n$ where $a_i=1$ and $a_k =0$ for $k>j$ with $k\neq i$.  
Also, $H_{i,i-1}$ is the group of invertible diagonal matrices in $B_n(A)$ with
ones in all diagonal entries except the $(i,i)$-entry. For $j<i-1$, $H_{ij}$ has ones on the diagonal, 
any element of $A$ in the $(j,i)$-position, and zeros elsewhere. It is easily verified that these are all groups. 

(ii) If $i<n$ then all elements of $H_{ij}$ fix $e_n$. Thus, we may suppose that
$c=ab \in G_{nj}$ for some $j$. Now if $c$ fixes $e_n$, then $c=1$, 
so $a=b=1$. \qed

\bigskip

From now on, we revert to ACVF, so $K$ denotes an algebraically closed valued field.
Recall, the notation $B_{n,m}(R)$ from Chapter~\ref{acvfbackground}. In particular, $B_{n,n}(R)$
is the inverse image under the natural map $B_n(R) \rightarrow B_n(k)$ of the 
stabiliser in $B_n(k)$ of the standard basis vector $e_n$
of $k^n$.

\begin{lemma} \label{R6}\index{opaque!group}
Let $K$ be an algebraically closed valued field.

(i) The groups $K/R$ and $K/\M$ (under addition) are opaque.

(ii) If $V$ is the group of units of $R$, then the multiplicative groups $K^*/V$ and
$K^*/(1+\M)$ are opaque.

(iii) Let $G=B_n(K)$, $F=B_n(R)$ and $F'=B_{n,n}(R)$. Then $G/F $ and $G/F'$ are both opaquely layered everywhere.
\end{lemma}

{\em Proof.} (i) Every definable subset of $K$ is a finite Boolean combination of balls. Thus,
it suffices to verify that if $U$ is a ball, then, for all but finitely many cosets $a+R$ of $R$, $a+R$ is  
either contained in $U$ or disjoint from $U$ (and likewise for $\M$).
This is obvious.

(ii) This is similar to (i).

(iii) This follows by (i), (ii), Corollary~\ref{R5} and Lemma~\ref{R8}. 
For the case of $G/F$, \ref{R8} provides a sequence 
$\{1\}=G_0<\ldots <G_N=G$ as in \ref{R5}, and the corresponding groups $H_i$. 
By the description above in terms of matrices of the $H_{ij}$, 
the groups $H_i/(H_i \cap F)$ and $H_i/(H_i \cap F')$ have the form described in (i) and (ii). 
Since 
\ref{R8} applies both to $G_i(K)$ and $G_i(R)$, the semidirect product decomposition
 required in \ref{R5} is clear for $F$. 

We must also verify that for each $i$,
$G_{i+1} \cap F'=(G_i \cap F')(H_i \cap F')$. Suppose $c\in G_{i+1} \cap F'$. Then $c=ab$ for some
$a\in G_i \cap F$ and $b\in H_i \cap F$. Applying \ref{R8} for the field $k$, we see that this
 remains true in the reduction modulo $\M$. Hence 
the reduction of $ab$ fixes the reduction of $e_n$. Thus the reductions of 
$a$ and $b$ both fix the reduction of $e_n$, by the last part of \ref{R8}. It follows that
$a\in G_i \cap F'$ and $b\in H_i \cap F'$, as required.
\qed

\begin{theorem}\label{R7} \index{resolution!prime}
Let $C$ be a subfield of the algebraically closed valued field $K$, and let $e$ be a
 finite set of imaginaries. Then $Ce$
admits a resolution $D$ which is minimal, prime and atomic over $Ce$. 
Up to isomorphism over $Ce$, $D$ is the unique prime resolution of $Ce$.
Also, $k(D)=k(\acl(Ce))$ and $\Gamma(D)=\Gamma(Ce)$.
\end{theorem}

{\em Proof.} 
We first prove existence of an atomic prime resolution. The easy case is 
when $e$ has a trivially valued resolution $B$ over $C$. 
In this case, by Lemma~\ref{2.1.7},  the only $B$-definable lattice
 is $R^n$. Since $R^n$ is resolved in $C$, 
we reduce to the case when $e=(e_1,\ldots,e_n)$ is a sequence of elements of $k$. 
Now $e$ is opaquely layered over $C$ via $e_i,E_i$, where
$E_i$ is just the equivalence relation $x+\M=y+\M$. Let $D$ be the atomic 
resolution of $Ce$ over $C$
given by Proposition~\ref{R3}. Let $B'$ be any resolution of $e$ over $C$. To prove 
primality, we must show that $B'$ satisfies the last condition in Proposition~\ref{R3}.
The quantifier 
elimination of \cite[Theorem 2.1.1(iii)]{hhm}
with sorts $K,k,\Gamma$  yields the following, which is what is required: if $e_i\in k$ and 
$e_i\in \acl(B')$, then 
$\acl(B')$ contains a field element $x$ with $\res(x)=e_i$. 

Now suppose that every resolution of $e$ over $C$ is non-trivially valued. 
The tuple of imaginaries $e$ has the same definable closure over $C$ as some pair
$(a,b)$ where $a\in S_n$ and $b\in T_m$ for some $n,m$. We identify $S_n$ with 
$B_n(K)/B_n(R)$. In Chapter~\ref{acvfbackground} we identified $T_m$ with $\bigcup_{i=1}^{\ell}B_m(K)/B_{m,\ell}(R)$.
However, if $t\in B_m(K)/B_{m,\ell}(R)$ for $\ell<m$, then $t$ may be identified with a pair
$(s',t')$ where $s'\in B_m(K)/B_m(R)=S_m$, and $t'\in B_{\ell}(K)/B_{\ell,\ell}(R)$;
here, $t\in \red(s')$, and the coset $t'$ consists of the $\ell \times \ell$ top left
 minors of the matrices in the coset $t$. Thus, after adjusting $a$ and $b$ if necessary, 
we may suppose that
$a\in B_n(K)/B_n(R)$, and $b\in B_m(K)/B_{m,m}(R)$.
By Lemma~\ref{R6}, $B_n(K)/B_n(R)$ and $B_m(K)/B_{m,m}(R)$ are opaquely layered, so 
$\tp(a/C)$ and $\tp(b/C)$ are opaquely layered (so $\tp(b/Ca)$ is). Hence,
 by Lemma~\ref{extra},
$\tp(ab/C)$ is opaquely layered.  

By Proposition~\ref{R3}, $Ce$ has an atomic resolution $D$. 
Again, to prove primality, we must show that any resolution $B'$ of 
$Ce$ satisfies the final hypothesis in that proposition. Now $B'$ is non-trivially 
valued, so  is a model of ACVF,
so contains an element of any $B'$-definable equivalence class, as required.

The field $D$ (obtained in either of the above cases) is a minimal 
resolution of $e$ over $C$: for if $D'$ is an 
algebraically closed subfield of $D$ with
$C \subset D'$ and $e\in \acl(D')$, then by primality $D$ embeds into $D'$ over $Ce$; hence 
as $D$ has finite transcendence degree over $C$, $D=D'$. 

To see uniqueness, suppose that
$D'$ is another prime resolution of $Ce$ over $C$. Then $D'$ embeds into $D$ over $Ce$. 
By minimality of $D$ the embedding is surjective, so $D'$ is $Ce$-isomorphic to $D$. This 
gives uniqueness. 

The last two assertions follow from Proposition~\ref{twores}.  \qed

\begin{corollary} \label{R7.1}
Let $E=\{e_i:i\in \omega\}$ be a countable set of imaginaries, and $C\subset K$.  Then 
$CE$ admits a resolution $D$ over
$C$ such that $k(D)=k(\acl(CE))$ and $\Gamma(D)=\Gamma(CE)$.
\end{corollary}

{\em Proof.} Put $C_n=Ce_0\ldots e_{n-1}$. Then each $C_n$ admits a resolution $D_n$ over 
$C$ as in the theorem. The inclusions $C_n \rightarrow C_{n+1}$ 
yield embeddings $D_n \rightarrow D_{n+1}$. Let $D$ be the direct limit of this system of maps.
\qed

\begin{corollary} \label{R7.2}  Let $C\subset K$, and $e$ a finite  
tuple of imaginaries.
Then $Ce$ admits a $\dcl$-resolution $D$ with  $k(D)=k( Ce)$ and $ 
\Gamma(D)=\Gamma(Ce)$.
\end{corollary}

{\em Proof.}   We may assume that $C$ is non-trivially valued, to
ensure (working over $C$) the first assumption of 
 Corollary~\ref{R3.2}. By \ref{R3.2}, there exists a $\dcl$-resolution  $D$
with $\tp(D/Ce)$ stationary, and embedding into any resolution.  By the 
latter fact, and using Corollary~\ref{R7.1}, we have
$\Gamma(D)=\Gamma(Ce)$ and $k(D) \subseteq k(\acl(Ce))$.  But $k(D) / Ce$
is stationary, so by  stable embeddedness $k(D) / k(Ce)$ is stationary,
and it follows that $k(D) = k(Ce)$.  \qed

%
%

\chapter{Maximally complete fields and domination}\label{maximallycomplete}

We now focus on independence relations $A\dnf_C B$ when all the sets $A,B,C$
 are in the field sort. Recall that a valued field is {\em maximally complete}\index{maximally complete} if it has 
no proper immediate extension. In this chapter we prove that over a maximally complete base, without any assumption 
of stable domination, still a field is dominated by its value group and residue field (Proposition~\ref{dom1}).
A slightly harder but more powerful result is Theorem~\ref{fulldom}, that, again over a maximally complete base, the 
type of a field is stably dominated over its definable closure in $\Gamma$. 

The following result  will be used to give a criterion for 
orthogonality to $\Gamma$ (Proposition~\ref{maxim}), and for these domination results. We do not
know its origin, but part (i) is Lemma 3 of \cite{baur}.
It will also be  important in 
Chapter~\ref{maximummodulus}. In the next proposition, if $A$ is a subfield of $K$, then $R_A$ denotes $R\cap A$ and $\Gamma_A$
 denotes the value group of $A$. 

\begin{proposition}\label{subs}
Let $C<A$ be an extension of non-trivially valued fields, 
and suppose that $C$ is maximally complete. 

(i) Let $V$ be a finite dimensional $C$-vector subspace of $A$. Then 
there is a basis
$\{v_1,\ldots,v_k\}$ of $V$ such that for any $c_1,\ldots,c_k\in C$,
$|\sum_{i=1}^k c_iv_i|=\max\{|c_iv_i|:1\leq i\leq k\}$.

(ii) Assume also that $\Gamma_C=\Gamma_A$, and let $S=V\cap R_A$.
Then there are generators $d_1,\ldots,d_k$
of $S$ as a free $R_C$-module such that
$|d_i|=1$ for each $i=1,\ldots,k$, and for any $f_1,\ldots,f_k \in R_C$,
$|\Sigma_{i=1}^kf_id_i|= \max\{|f_1|,\ldots,|f_k|\}$.
\end{proposition}

Note that (ii) implies that $\{\res(d_1),\ldots,\res(d_k)\}$ 
is linearly independent over $k(C)$: for if $r_1,\ldots,r_k \in R_C$, 
then
$$\begin{array}{rl}
 \sum_{j=1}^k \res(r_j)\res(d_j)=0 \Longleftrightarrow 
                                     |\sum_{j=1}^k r_j d_j|<1
        &  \Longleftrightarrow \max\{|r_1|,\ldots,|r_k|\}<1  \\
        &  \Longleftrightarrow \res(r_j)=0 \mbox{~for all~} j.
\end{array}
$$
We shall say that ${\cal B}= \{v_1,\ldots,v_k\}$ is a {\em separated} basis of $V$ if (i) holds. 
We say in addition that it is {\em good}\index{good separated basis} if,
whenever $b,b'\in {\cal B}$ and $|b|\Gamma_C=|b'|\Gamma_C$, we have $|b|=|b'|$.

\bigskip

{\em Proof.} 
(i) 
The separated basis is built inductively.
Suppose $\{v_1,\ldots,v_{\ell}\}$ is a separated basis of a $C$-subspace $U$ of $V$, and 
$v\in V\setminus U$. 

{\em Claim.} There is $w\in U$ 
such that $|v-w|=\inf\{|v-u|:u\in U\}$. 

{\em Proof of Claim.} We construct a transfinite sequence
$(w^\nu)$ ($\nu$ an ordinal),
 where $w^\nu= \sum_{i=1}^{\ell}a_i^\nu v_i$, such that if
$\gamma^\nu=|v-w^\nu|$ then the sequence $(\gamma^\nu)$ is decreasing.
At stage $\nu+1$, if $|v-w^\nu| \leq |v-x|$ for all $x \in U$, then stop.
 Otherwise, there is $w^{\nu+1}\in U$ with 
$\gamma^{\nu+1}:=|v-w^{\nu+1}| <\gamma^\nu$. 
Now consider a limit ordinal $\lambda$.
For any $\nu <\nu'<\lambda$,
$$|(v -w^\nu)-(v-w^{\nu'})|=\max\{\gamma^\nu,\gamma^{\nu'}\}=\gamma^\nu,$$
so $|\sum_{i=1}^{\ell}(a_i^\nu-a_i^{\nu'})v_i|=\gamma^\nu$. It follows 
by the
inductive hypothesis that $|a_i^\nu-a_i^{\nu'}| \leq \gamma^\nu/|v_i|$ for 
each $i=1,\ldots,\ell$. For each $i$, by choosing a pseudoconvergent subsequence 
of $(a_i^\nu)$ (not necessarily cofinal) and by using the maximal completeness of $C$, we 
may find some
 $a_i^\lambda$ such that 
$|a_i^\lambda-a_i^\nu| \leq \gamma^\nu/|v_i|$ for 
all $\nu\leq \lambda$. Hence, for all $\nu<\lambda$, 
$$|v-\sum_{i=1}^{\ell}a_i^\lambda v_i| 
\leq \max\{|v-\sum_{i=1}^{\ell}a_i^\nu v_i|, 
|\sum_{i=1}^{\ell}(a_i^\nu-a_i^{\lambda})v_i|\} \leq \gamma^\nu.$$
Thus, the induction proceeds at limit stages, where we put
$w^\lambda=\sum_{i=1}^{\ell}a_i^\lambda v_i$. As $C$ is maximally
complete, the pseudoconvergent sequence
$(w^\lambda)$ has a 
limit $w=\sum_{i=1}^{\ell}a_iv_i$, and this $w$ 
satisfies the claim.

Given the claim, put $v_{\ell+1}:=v-w$. Then 
$\{v_1,\ldots,v_{\ell+1}\}$ is a separated basis of the 
$C$-subspace of $V$ which it spans: indeed, if $|\Sigma_{i=1}^{\ell+1} c_iv_i|<
\max\{|\Sigma_{i=1}^{\ell} c_iv_i|,|c_{\ell +1}v_{\ell +1}|\}$, then 
$|v_{\ell +1}+\Sigma_{i=1}^\ell c_{l+1}^{-1}c_iv_i|<|v_{\ell+1}|$, so
$|v-(w-\Sigma_{i=1}^\ell c_{\ell+1}^{-1}c_i v_i)|<|v-w|$, contradicting the choice of $w$.

(ii) Let $\{v_1\ldots,v_k\}$ be as in (i). Since $\Gamma_C=\Gamma_A$, for each 
$i=1,\ldots,k$ there is $e_i\in C$ with $|e_iv_i|=1$. Now put $d_i:=e_iv_i$ for 
each $i$.  \qed

\begin{lemma} \label{MVFnew}
Let $C<A$ be valued fields, and suppose that $C$ is maximally complete.
Let $V$ be a finite or countable dimensional subspace of $A$ (as a vector space over $C$). Then 
$V$ has a  good separated basis ${\cal B}$ over $C$.
\end{lemma}

{\em Proof.} 
 Suppose that ${\cal B'}$ is a good separated basis of the
 finite-dimensional subspace $U$ of $V$,
 and let
$v\in V\setminus U$. We shall show that ${\cal B'}$ extends to a good separated basis of $U+Cv$. 

By the claim in the last proof, there is $w\in U$ such that
$|v-w|=\inf \{|v-u|:u\in U\}$. Put $b:=v-w$. Then by the last proof, ${\cal B'} \cup \{b\}$ is a separated basis for the
subspace (over $C$) which it spans.

Suppose first that for any $u\in U$, $|b|\neq |u|$. In this case, $|b|\neq \gamma|u|$
for any $\gamma \in \Gamma(C)$. Now ${\cal B'}\cup \{b\}$ is a separated basis 
of the $C$-space it spans (by the proof of Proposition~\ref{subs}) and still 
is good.

Now suppose that there is $u\in U$ with $|b|=|u|$. Put $u=\Sigma_{i=1}^m a_ib_i$, 
where ${\cal B'}=\{b_1,\ldots,b_m\}$, and $a_1,\ldots,a_m \in C$.  As ${\cal B'}$ is separated,
$|b|=|u|=|a_ib_i|$ for some $i$, say with $i=1$. We claim that $|b-a_1b_1|=|b|$.
For otherwise, $|b-a_1b_1|<|b|=|a_1b_1|$. But by the choice of $w$,
$|b-a_1b_1|=|v-(w+a_1b_1)| \geq |v-w|=|b|$, a contradiction.
Now put $b_{m+1}:=a_1^{-1}b$. We have $|b_{m+1}|=|b_1|=|b_{m+1}-b_1|$, the latter equality as
$|a_1b_1|=|b|=|b-a_1b_1|$, so $|b_1|=|ba_1^{-1}-b_1|=|b_{m+1}-b_1|$. 
Since ${\cal B}' \cup \{b\}$ is a separated basis, so is ${\cal B'} \cup \{b_{m+1}\}$. 
The latter also is good, since ${\cal B}'$ is good and $|b_{m+1}|=|b_1|$. \qed

\begin{remark}\label{linindepsep} \rm
Let $C<A$ be an extension of non-trivially valued fields, with $C$ maximally complete, and
suppose that
 ${\cal B}$ is a separated basis for a subspace $V$ of $A$. Then 
if $b_1,\ldots, b_{\ell}\in {\cal B}$ and $|b_1|=\ldots =|b_{\ell}|$, then
$1, \res(b_2/b_1),\ldots,\res(b_{\ell}/b_1)$ are linearly independent over $k(C)$.

Indeed, suppose $r_1,\ldots, r_{\ell}\in C$ with $|r_i|\leq 1$ for each $i$, and
$\Sigma_{i=1}^{\ell} \res(r_i)\res(b_i/b_1)=0$.  Then
$|\Sigma_{i=1}^{\ell} r_ib_i/b_1|<1$, so $|\Sigma_{i=1}^{\ell} r_ib_i|<|b_1|$. 
As ${\cal B}$ is separated, $|r_jb_j|<|b_1|$ for each $j$, and 
it follows that
$|r_j|<1$ for each $j$, so each $\res(r_j)=0$, as required. 
\end{remark}

\begin{lemma}\label{Cmax}
Let $C\leq A,B$ be algebraically closed valued fields, and suppose 
that $C$ is maximally complete. 
Assume that $\Gamma(C) = \Gamma(A)$ and that $k(A)$ and $k(B)$ are 
linearly disjoint over $k(C)$. Let $E$ be the subring of $K$ generated by $A \cup B$.
Then

(i)  Let $a_1,\ldots,a_n \in A, b_1,\ldots,b_n \in B$.  
Then
there exist $d_1,\ldots,d_k\in A$, $b_1',\ldots,b_k' \in B$   such that in $A \tensor_C B$ we have
 $ \sum_{i=1}^n a_i \tensor b_i =  \sum_{j=1}^k d_j \tensor b_j'$
while in $E$ we have $|e| = \max_{j=1}^k |d_j| |b_j'|$.  

(ii) $\Gamma(E)=\Gamma(B)$.

(iii) Let $a=(a_1,\ldots,a_n) \in A^n$, $b_1,\ldots,b_n \in B$, $e = \sum_{i=1}^n a_ib_i$.  
Let $a'=(a_1',\ldots,a_n') \models \tp(a/C)$.  Let $e' = \sum _{i=1}^n a_i'b_i$.
Then $|e'| \leq |e|$.

(iv) If $A'\equiv_C A$ and  $k(A')$ is also linearly disjoint from
 $k(B)$ over $k(C)$, then $A\equiv_B A'$.
\end{lemma}

{\em Proof.}

(i)  We may suppose
$|a_i|\leq 1$ for each $i$: for if $\gamma=\max\{|a_1|,\ldots,|a_n|\}>1$,
choose $c\in C$ with $|c|=\gamma$, and replace each $a_i$ by $a_ic^{-1}$ and
$b_i$ by $b_ic$. 

Now by Proposition~\ref{subs}(ii) and the remark following it, there are $d_1,\ldots,d_k\in A$ such that 
$|d_1|=\ldots = |d_k|=1$, $\res(d_1),\ldots,\res(d_k)$ are linearly independent 
over $k(C)$,
and $a_1,\ldots,a_n$ are in the $R_C$-module generated by $d_1,\ldots,d_k$.
Thus, there are $c_{ij}\in R_C$ ($i=1,\ldots, n$, $j=1,\ldots, k$), so that
$a_i=\sum^k_{j=1}c_{ij}d_j$ for each $i=1,\ldots,n$. In particular,
$$\sum_{i=1}^n a_i \tensor b_i=\sum_{i=1}^n\sum_{j=1}^k c_{ij}d_j \tensor b_i=
\sum_{j=1}^k d_j \tensor (\sum_{i=1}^n c_{ij}b_i) =\sum_{j=1}^k d_j \tensor b_j',$$
where $b_j'=\sum_{i=1}^n c_{ij}b_i \in B$ for each $j$.    

Let  $e = \sum_{i=1}^n a_ib_i$; then  in particular $e = \sum_{j=1}^k d_j b_j'$.

We may partition $\{1,\ldots,k\}$ as $I_1 \cup \ldots \cup I_m$, where
for each $j$,
$I_j:=\{i:|b_i'|=\gamma_i\}$. For each $j=1,\ldots,m$, put
$e_j:=\sum_{i\in I_j} d_ib_i'$. Then $e=e_1+\ldots +e_m$. We claim
that $|e_i|=\gamma_i$ for each $i=1,\ldots,m$. So fix $i\in \{1,\ldots,m\}$,
and $\ell\in I_i$. Put $f_j= b_{\ell}'^{-1}b_j'$ for each $j\in I_i$. Then
$|f_j|=1$ for each $j\in I_i$. As $\res(d_1),\ldots,\res(d_k)$ are 
linearly independent over $k(C)$, and $k(A)$ and $k(B)$ 
are linearly disjoint over
$k(C)$, $\res(d_1),\ldots, \res(d_k)$ are linearly independent over
$k(B)$. It follows that $|\sum_{j\in I_i} d_jf_j|=1$, yielding the 
claim.  (i) clearly follows.

 (ii) Let $e\in E$. Since $\Gamma(E) \subseteq \dcl(\{|x|:x\in E\}$, it suffices to 
show $|e|\in \Gamma(B)$.   By (i) we have $e=\sum_{j=1}^k d_j b_j'$ 
with $d_j \in A, b_j' \in B$, and  $|e| = \max_{j=1}^k |d_j| |b_j'|$; so $|e| \in \G(B)$.

(iii)  Let $d_j,b_j'$ be as in (i).  Let $d'=(d_1',\ldots,d_k')$ be such that $\tp(a',d'/C)=\tp(a,d/C)$.
 Then $\sum_{i=1}^n a_i'b_i = \sum_{j=1}^k d_j' b_j'$
 (since the equality already holds in the tensor product.)  
  We have $|d_j'| = |d_j|$ using $\G(A)=\G(C)$.  Thus
  $|e'|  = | \sum_{j=1}^k  d_j' b_j' | \leq \max_{j=1}^k |d_j| |b_j'| = |e|$.

(iv)   In the notation of (ii), by quantifier elimination (e.g. in $\L_{\div}$) we must show that
 if $(a_1',\ldots,a_n')\equiv_C
 (a_1,\ldots,a_n)$ and $e':=\sum_{i=1}^n a_i'b_i$, then $|e'|=|e|$. The proof 
of (i) yields this.
\qed

\bigskip

We deduce a criterion for orthogonality to $\Gamma$ which will be used in the 
next chapter. It implies in particular that if $C$ is an algebraically closed valued field
with no proper immediate extensions and $a\in K^n$, then $\tp(a/C) \perp \Gamma$
if and only if $\Gamma(C)=\Gamma(Ca)$. Recall the notation $\St_C$ from Part~I and
Chapter~\ref{acvfbackground}. If $C$ is a field, then $\St_C\subset \dcl(C\cup k(C))$; in general,
$\St_C\subset \dcl(C\cup \VS_{k,C})$.

\begin{proposition}\label{maxim}
Let $C=\acl_K(C)$, and let $F$ be a maximally complete immediate extension of $C$, 
and $a\in K^n$. Then the
following are equivalent.

(i) $\tp(a/C)\perp \Gamma$.

(ii) $\tp(a/C)\vdash \tp(a/F)$ and $\Gamma(Ca)=\Gamma(C)$.
\end{proposition}

{\em Proof.}
Assume (i). Then $\Gamma(Ca)=\Gamma(C)$ by Lemma~\ref{ortheasy} (v). 
Furthermore, since the extension $C\leq F$ is immediate, $k(A)\dnf^g_C F$, where
$A=\acl_K(Ca)$. 
As $C=\acl_K(C)$, $\St_C \subset \dcl(C \cup k(C))$. Thus, $\St_C(A) \dnf_C \St_C(F)$, 
so
by
Proposition~\ref{perpgamma}, $A\dnf^g_C F$ via $a$. Hence
$\tp(a/C)\vdash \tp(a/F)$.

Conversely, suppose (ii). We first show that $\Gamma(F)=\Gamma(Fa)$. 
For this, it suffices to show that if $b=(b_1,\ldots,b_m)\in F^m$ then
$\Gamma(Cb)=\Gamma(Cba)$. 
Suppose this is false, and
let $i$ be least such that
$\Gamma(Cb_1\ldots b_i)\neq \Gamma(Cb_1\ldots b_i a)$.
Then as $\tp(a/C)\vdash \tp(a/F)$, we have $\tp(a/Cb_1\ldots b_{i-1})\vdash \tp(a/F)$. 
It follows that
$\tp(b_i/Cb_1\ldots b_{i-1})$ implies $\tp(b_i/Cb_1\ldots b_{i-1}a)$.
As $\Gamma(Cb_1\ldots b_{i-1}a)\neq \Gamma(Cb_1\ldots b_i a)$, 
$b_i$ is not in $\acl(Cb_1\ldots b_{i-1}a)$, so
$\tp(b_i/\acl(Cb_1\ldots b_{i-1}))$ is not realised in 
$\acl(Cb_1\ldots b_{i-1}a)$. Hence, by
Proposition~\ref{noembed} (ii), $\Gamma(Cb_1\ldots b_i)=\Gamma(Cb_1\ldots b_i a)$,
 contradicting the choice of $i$.

By (ii), $a\dnf^g_C F$, so to show $\tp(a/C)\perp \Gamma$
it suffices by Lemma~\ref{ortheasy}(iii) to prove $\tp(a/F)\perp \Gamma$. Put 
$A:=\acl_K(Fa)$, and let $B$ be an algebraically closed field extending $F$, with
$A\dnf^g_F B$. By Lemma~\ref{gtostab}, $k(A)$ and $k(B)$ are 
linearly disjoint over
$k(F)$. Since $\Gamma(F)=\Gamma(A)$ (by the last paragraph), it follows from Lemma~\ref{Cmax} 
that $\Gamma(A \cup B)=\Gamma(B)$, 
as required.
\qed

\begin{remark} \rm
Proposition~\ref{maxim} actually characterises maximally complete fields. For  suppose
 $C=F$ is not maximally complete.
Then there is a chain $(U_i:i\in I)$ of $C$-definable balls with no least element such that 
no  element of $C$ lies in $U:=\bigcap(U_i:i\in I)$.
Let $a\in K$ lie in $U$. Then $\Gamma(C)=\Gamma(Ca)$ 
but $\tp(a/C)$ is not orthogonal to $\Gamma$ by Lemma~\ref{eqorthog}.
\end{remark}

In Chapter~\ref{dominationinacvf} we studied domination of a type by its stable part. Here, we
examine domination of a field by its  value group and residue field. For these
results, we do not need to assume orthogonality to $\Gamma$, but do need the 
assumption that the base is a maximally complete  valued field. 
As a consequence, we obtain the domination results (\ref{fulldom2}, \ref{fulldom}) which are the goal of this chapter.
First we need several technical lemmas. Analogous to the notion of orthogonality to the
value group, we have a notion of orthogonality to the residue field.

\begin{definition}\index{orthogonal!to residue field} \rm
Let $a=(a_1,\ldots,a_n)$.
We shall say that $\tp(a/C)$ {\em is orthogonal to $k$} (written $\tp(a/C)\perp k$) 
if for any model $M$ with $a\dnf^g_C M$, we have $k(M)=k(Ma)$. 
\end{definition}

For orthogonality to $k$, the obvious analogues of Lemma~\ref{ortheasy} (i), (ii), (iii), (v) hold. Observe also that if $\tp(a/C)\perp k$ then $\St_C(a)=\dcl(C)$ (see the proof of Lemma~\ref{ortheasy}(v)).

The following is well known, but for want of a reference we give a proof.

\begin{remark}\label{somewhere}
Let $F<L=\acl(F(a))$ be an extension of  valued fields (where $a$ is a finite 
sequence). Then 
$$\trdeg(L/F)\geq  \trdeg (k(L)/k(F))+ \rk_{{\mathbb Q}}(\Gamma(L)/\Gamma(F)).$$
\end{remark}

{\em Proof.}
We may suppose that $a$ is a singleton. Then
$\trdeg (k(L)/k(F))\leq 1$ and  $\rk_{{\mathbb Q}}(\Gamma(L)/\Gamma(F)) \leq 1$.
Suppose that $k(L)\neq k(F)$. Then certainly
$a\not\in \acl(F)$.
We may
 assume that $|a|=1$ with $\res(a)\not\in k(F)$. Thus, $a$ is generic in the 
closed ball $R$, and it follows from Lemma~\ref{eqorthog} that 
$\Gamma(L)=\Gamma(F)$. \qed

\bigskip

\begin{lemma}\label{intdesc}
Let $C,L$ be algebraically closed valued fields with $C \subset L$. Let 
 $a_1,\ldots,a_m\in L$, $\gamma_i:=|a_i|$ and let
$d_i$ be a code for the open ball $B_{<\gamma_i}(a_i)$
(so $d_i \in \red(\gamma_iR$)). Put
$C':=\acl(C\gamma_1\ldots \gamma_m)$, $C'':=\acl(Cd_1\ldots d_m)$ and
$C''':=\acl(Ca_1\ldots a_m)$. Then
$\St_{C'}(L)=\acl(C''  k(L))\cap \St_{C'}$.
\end{lemma}

{\em Proof.} First, consider the case $m=0$. In this case, the conclusion of the lemma
is that
$$\St_{C}(L)=\acl(C k(L))\cap \St_C.$$
This is clear by Lemma~\ref{2.6.2}(ii), since if $s\in \acl(C) \cap S_n$
 then $\red(s)$ is $C$-definably isomorphic to $k^n$.

In general, we argue by induction on $m$. First observe that $d_i\in \St_{C'}(L)$ 
for each $i$. Hence $C''k(L) \subset \St_{C'}(L)$. As $L=\acl(L)$, $\St_{C'}(L)$ is algebraically closed in
$\St_{C'}$, so we have the containment $\supseteq$ of the statement. 

In the other direction, suppose first
$C'=\acl(C\gamma_1\ldots \gamma_{\ell})$ for some $\ell <m$.
Then by induction, we obtain 
$$\St_{C'}(L)=\acl(\acl(Cd_1\ldots d_{\ell})  k(L)) \cap \St_{C'}
          \subset \acl(C''  k(L)).
$$ 
Thus, we may assume that the $\gamma_i$ are ${\mathbb Q}$-linearly independent in 
$\Gamma(L)$ over $\Gamma(C)$. It follows easily 
that each $a_i$ is generic in the open ball $d_i$ over $C''a_1\ldots a_{i-1}$. 
In particular, $\tp(a_i/C''a_1\ldots a_{i-1})\perp k$: indeed, if
$a_i\dnf^g_{C''a_1\ldots a_{i-1}} N$, for any model $N$, then as $N$ contains 
a field element of $d_i$,  $\Gamma(N)\neq \Gamma(Na_i)$, and hence $k(N)=k(Na_i)$ by 
Remark~\ref{somewhere}. 
 From this, we obtain that $a_i\dnf^g_{C''a_1\ldots a_{i-1}}
k(L)$; for example if $\beta$ is a finite sequence from $k(L)$ then the above orthogonality gives
$\beta \dnf^g_{C''a_1\ldots a_{i-1}} a_i$, and then Proposition~\ref{perpgamma} applies. Thus,
 $\tp(a_i/C''a_1\ldots a_{i-1} k(L))\perp k$. 
It follows that $\tp(a_1\ldots a_m/C''k(L))\perp k$.
Hence, $\St_{C''k(L)}(C''k(L)a_1\ldots a_m)\subseteq \dcl(C''k(L))$, so
$$
  \St_{C'' k(L)}(\acl(C''' k(L))=\acl(C''  k(L)) \cap \St_{C''k(L)}.
$$
However, applying the case $m=0$ to the field $C'''$, we have
$$
  \St_{C'''}(L)=\acl(C'''  k(L))\cap \St_{C'''}.
$$
A fortiori, $\St_{C'}(L) \subset \acl(C'''  k(L)).$ So
$$\St_{C'}(L)\subset\St_{C'}(\acl(C''' k(L))\subset \St_{C'' k(L)}(\acl(C''' k(L))
                 \subset\acl(C''  k(L)).
$$
This proves the lemma. \qed

\begin{lemma}\label{intdesc2}
Let $C$ be an algebraically closed valued field. Let $L,M$ be extension fields of $C$,
with $\Gamma(L) \subset \dcl(M)$. Consider sequences $a_1,\ldots,a_r$ and $b_1,\ldots,b_s$ 
from  $L$ such that  $|a_1|,\ldots,|a_r|$ is a ${\mathbb Q}$-basis of $\Gamma(L)$ over 
$\Gamma(C)$ and $\res(b_1),\ldots,\res(b_s)$ form a transcendence 
basis of $k(L)$ over $k(C)$ (we do not here assume that $r,s$ are finite). Let $e_i\in M$ with $|a_i|=|e_i|$ for each $i$. Let
$C^+:=\dcl(C \cup \Gamma(L))$. Then the following conditions on $C,L,M$ are equivalent.

(1) For some $a_1,\ldots,a_r,b_1,\ldots,b_s,e_1,\ldots,e_r$ as above, 
$$\res(a_1/e_1),\ldots,\res(a_r/e_r),\res(b_1),\ldots,\res(b_s)$$ are algebraically
 independent over $k(M)$.

(2) 
For all $a_1,\ldots,a_r,b_1,\ldots,b_s,e_1,\ldots,e_r$ as above, 
$$\res(a_1/e_1),\ldots,\res(a_r/e_r),\res(b_1),\ldots,\res(b_s)$$ are 
algebraically independent over $k(M)$.

(3) $\St_{C^+}(L)$ and $\St_{C^+}(M)$ are independent in the stable structure 
$\St_{C^+}$.
\end{lemma}

{\em Proof.} Let $a_1,\ldots,a_r,b_1,\ldots,b_s,e_1,\ldots,e_r$ satisfy the 
hypotheses of the lemma, and let $d_i$  be the open ball of radius $|a_i|$ around $a_i$.
By Lemma~\ref{intdesc}, 
$$\St_{C^+}(L)=\acl(Cd_1\ldots d_r \res(b_1)\ldots \res(b_s)) \cap \St_{C^+}.
$$
For notational convenience we now assume $r,s$ are finite, but the argument below, 
applied to subtuples, yields the lemma without this assumption.
Since $\Gamma(Cb_1\ldots b_s)=\Gamma(C)$, it follows that $d_1,\ldots,d_r,\res(b_1),\ldots,\res(b_s)$ are independent in $\St_{C^+}$,
so the Morley rank of any tuple enumerating $\St_{C^+}(L)$ over $\St_{C^+}$ is $r+s$.
We have $\dcl(Md_i)=\dcl(M\res(a_i/e_i))$ for each $i$, as over the model $M$, $k$ is in definable bijection with $\red(B_{\leq |a_i|}(a_i))$. Thus, (3) holds if and only if
 $$\RM(\tp(\res(a_1/e_1),\ldots,\res(a_r/e_r),\res(b_1),\ldots,\res(b_s)/M))=r+s,
$$
which holds if and only if 
$\res(a_1/e_1),\ldots,\res(a_r/e_r),\res(b_1),\ldots,\res(b_s)$ 
are algebraically independent over $k(M)$. Thus (1) and (3) are equivalent. Since this 
equivalence is valid for any choice of
$a_1,\ldots,a_r,b_1,\ldots,b_s,e_1,\ldots,e_r$ and (3) does not mention the choice, 
(1) and (2) are also equivalent.
\qed

\bigskip

We begin with a field-theoretic version of our  domination results (Proposition~\ref{dom1});
and deduce a   statement allowing for imaginaries (\ref{fulldom2}).  We will 
then  deduce what turns out to be a stronger model theoretic fact, 
Theorem~\ref{fulldom}.  This is be the basis of the notion of metastability.
Lemma~\ref{intdesc2} provides an algebraic rendering of the statement of ~\ref{fulldom}.

 In the following, we say that a function $f$ on a field induces the function $h$ on the value group
if $h(|x|)=|f(x)|$ for all $x$ in the domain of $f$. Similarly, $f$ induces $h'$ on the residue field if
$h'(\res(x))=\res(f(x))$ for all $x\in \dom(f) \cap R$.
If $A,B,C$ are structures, we say that a map $f:A\rightarrow B$ (or, formally, the pair $(f,B)$)
is unique up to conjugacy over $C$ (subject to certain conditions) if, whenever 
$f_i:A\rightarrow B_i$ ($i=1,2$) are two maps satisfying these conditions, 
there is an isomorphism $\ell:B_1\rightarrow B_2$ over $C$ such that $f_2=\ell \circ f_1$.
Below, we write $LM$ (or sometimes $\langle L,M\rangle$)
for the field generated by $L \cup M$.

\begin{proposition}\label{dom1}
Let $C$ be a maximally complete 
 valued field, and let $L,M$ be 
 valued fields containing $C$. 
 Let $h:\Gamma(L) \rightarrow \Gamma$ and $h':k(L)\rightarrow k$ be 
embeddings, with $h(\Gamma(L)) \cap \Gamma(M)=\Gamma(C)$, and with $h'(k(L))$ and $k(M)$ linearly disjoint over $k(C)$. 

(i) Up to conjugacy over 
$M \cup h(\Gamma(L)) \cup h'(k(L))$, there is a unique pair $(f,N)$ such that
$f$ is a valued field embedding  $L\rightarrow N$ over  $C$
which induces $h$ and $h'$,
and 
$\langle f(L),M\rangle =N$ (as fields). 

(ii) With $f, N$ as in (i),
$\Gamma(N) = \langle h(\Gamma(L)), \Gamma(M)\rangle$ (as subgroups of $\Gamma$),
and $k(N) = \langle k(f(L)),k(M)\rangle$ (as fields).
\end{proposition}

{\em Proof.}  %
Without loss of generality, we may assume that
$h$ and $h'$ are the identity maps, so $\Gamma(L)$ and $\Gamma(M)$ 
are independent over $\Gamma(C)$ in the sense of ordered groups, 
and $k(L)$ is linearly disjoint from $k(M)$ over $k(C)$. 
To prove (i),
 we must show that if $L'$ is isomorphic to $L$ over $C \cup \Gamma(L)\cup k(L)$, then
the isomorphism extends to a valued field isomorphism $LM \rightarrow L'M$ over $M$.

We first show that $L$ and $M$ are linearly disjoint 
over $C$.
So suppose $u_1,\ldots,u_n \in L$ are linearly independent over $C$, and span 
a $C$-subspace $U$ of $L$.
Now $U$ has a good separated basis ${\cal B}=\{b_1,\ldots, b_n\}$ by
Lemma~\ref{MVFnew}. By Remark~\ref{linindepsep},  if $b_1,\ldots,b_{\ell}\in {\cal B}$ with 
$|b_1|=\ldots =|b_{\ell}|$ then the elements $1,\res(b_2/b_1),\ldots,\res(b_{\ell}/b_1)$
are linearly independent over $k(C)$. 
We must show that $u_1,\ldots, u_n$ are linearly independent over $M$, so it 
suffices to show that $b_1,\ldots,b_n$ are linearly independent over $M$.
This will follow from the following claim.
 
\medskip

{\em Claim.} Let $x=\Sigma_{i=1}^n b_im_i$ where 
$m_1,\ldots,m_n\in M$. Then
$|x|=\max\{|b_i||m_i|:1\leq i\leq n\}$. 

{\em Proof of Claim.} Suppose that the claim
 is false. Put $\gamma:=\max\{|b_i||m_i|:1\leq i\leq n\}$ and $
J:=\{i:|b_im_i|=\gamma\}$. Then $|\Sigma_{i\in J} b_im_i|<\gamma$. For any distinct 
$i,j\in J$ we have $|b_i/b_j|=|m_j/m_i|\in \Gamma(L)\cap \Gamma(M)=\Gamma(C)$; so 
$|b_i|=|b_j|$ (as ${\cal B}$ is good), and hence $|m_i|=|m_j|$. Fix an element (say 1) in $J$ and let 
$J':=J\setminus \{1\}$. So $|b_i/b_1|=|m_i/m_1|=1$. Now as 
$|\Sigma_{i\in J} b_im_i|<\gamma$, we have
$|1+\Sigma_{i\in J'}(b_im_i)/(b_1m_1)|<1$. Thus, in the residue field,
$1+\Sigma_{i\in J'}\res(b_i/b_1)\res(m_i/m_1)=0$.  Hence, the elements
$1,\res(b_i/b_1)$ (for $i\in J')$ are linearly dependent over $k(M)$. 
Since $k(L)$ and $k(M)$ are linearly disjoint over $k(C)$, 
these elements are also linearly dependent over $k(C)$, 
contradicting the choice of the basis ${\cal B}$.

\medskip

Now suppose that $f:L\rightarrow L'$ is a valued field isomorphism 
inducing the identity on $\Gamma(L)\cup k(L)$. Then, by the above argument,
 $L$ and also $L'$ are 
linearly disjoint from $M$ over $C$ (so independent in the sense for pure algebraically closed fields). Hence,
we may extend $f$ by the identity on $M$ to a field isomorphism, also denoted $f$, from
$LM$ to $       L'M$. If $x\in LM$ then $x\in UM$ for
 some finite dimensional $C$-subspace $U$ of $L$ with a separated basis
$\{b_1,\ldots,b_n\}$, say. We may write $x=\Sigma_{i=1}^n m_ib_i$. Then, as in the claim
$|x|=\max\{|b_i||m_i|:1\leq i\leq n\}$. Likewise,
$|f(x)| =\max\{|f(b_i)||m_i|:1\leq i\leq n\}$. Since $|b_i|=|f(b_i)|$ for each $i$, 
it follows that $|x|=|f(x)|$, so $f$ is an isomorphism of valued fields.
This proof also gives that $|x|=|b||m|$ for some $b\in L$ and $m\in M$. This yields both assertions of (ii).
\qed

\bigskip

  In the following corollary, the conditions on residue fields and on value groups can
be viewed as independence in the theories of algebraically closed fields and divisible
(ordered) Abelian groups, respectively. The fields are inside $\U$, and types are in the sense of ACVF. 

\begin{corollary} \label{fulldom2}  Let $F$ be a maximally complete 
valued field,
$F \subset A=\dcl(A)$.  
Let $M$ be an extension of $F$ with $k(A),k(M)$ linearly disjoint over $k(F)$,
and $\G(A) \meet \G(M) = \G(F)$.  
   Then $\tp(M/F,k(A),\Gamma(A))\vdash \tp(M/A)$.
   \end{corollary}

{\em Proof.} Without loss of generality 
we may assume $A = \acl(F \union F')$, $F'$ finite.  By Theorem~\ref{R7.2} a prime $\dcl$-resolution $L$ 
of $A$ exists and satisfies:    $k(L) = k(A), \G(L)= \G(A)$.  So
the corollary for $L$ implies, a fortiori, the same for $A$.  Thus
we may assume $A=L$ is resolved.   In this case the corollary 
is immediate from   Proposition~\ref{dom1}.
\qed

\medskip

\begin{remark} \rm  
The assumption that $C$ is maximally complete is needed in these results,
and thus the stronger ones that follow.  
 For otherwise $C$ has a proper
 extension $L$ with $k(C)=k(L)$ and $\Gamma(C)=\Gamma(L)$. Then the hypotheses 
of Corollary~\ref{fulldom2} are vacuously true, but taking $M=A=L$, the conclusion is not. \end{remark}

\medskip

Recall the general notion of domination of invariant types from Definition~\ref{invariant},
along with   the semigroup $\overline{\Inv}(\U)$.  
We have natural embeddings of $\overline{\Inv}(k)$ and of $\overline{\Inv}(\G)$ into $\overline{\Inv}(\U)$.  With respect
to these, we have:

\begin{corollary} \label{inv-dom}  $\overline{\Inv}(\U) \cong \overline{\Inv}(k) \times \overline{\Inv}(\G)$. 
\end{corollary}

{\em Proof.}   Let $p$ be an invariant type of $\U$.  
 Let $C$ be a maximally complete algebraically closed valued field, such that 
 $\dcl(C)$ is a base for $p$.  Let $a \models p | C$, $b = k(Ca), d = \G(Ca)$.
 Then $b,d$ are definable functions of $a$, and $p$ yields invariant types
 $p_k,p_\G$ obtained by applying these definable functions to $p$.  
 It follows from Theorem  ~\ref{fulldom2} (ii) that $p$ is domination-equivalent
 to $p_k \tensor p_\G$.  This shows that the natural homomorphism 
  $\overline{\Inv}(k) \times \overline{\Inv}(\G) \to \overline{\Inv}(\U)$ is surjective.
 Injectivity is clear from the orthogonality of $k$ and $\G$.  \qed 

\medskip

Note that $\overline{\Inv}(k) \cong {{\mathbb N}}$; the unique generator is the generic type $p_{RES}$ of the residue field.     In  the next chapter we describe $\overline{\Inv}(\G)$.

Proposition~\ref{dom1} is used in the proof of the following strengthening, 
leading up to 
Theorem~\ref{fulldom}.
 
\begin{proposition} \label{domovergamma}
Let $C$ be a maximally complete algebraically closed valued field, and let $L,M$ be
algebraically closed valued fields containing $C$, with $L$ of finite transcendence degree over $C$.
Let $g:\Gamma(L)\rightarrow \Gamma(M)$ be an embedding over $C$. Put
 $C^+:=\dcl(C\cup g(\Gamma(L)))$. 
Then, up to conjugacy over $M$, there is a unique pair $(f,N)$ such that $f$ is a valued field embedding
$L\rightarrow N$ over $C$  inducing $g$, 
$\St_{C^+}(f(L))$ and $\St_{C^+}(M)$ are independent in the stable structure $\St_{C^+}$, and
$\langle f(L),M\rangle =N$ (as fields).
\end{proposition}

In the next lemma, and the proof of Proposition~\ref{domovergamma}, we shall write $\Gamma$
 additively and use a valuation $v$ rather than a norm,  
as we exploit the ${\mathbb Q}$-linear structure.
If $\gamma,\delta$ are in the value group $\Gamma$ we write
$\gamma<<\delta$ if $n\gamma<\delta$ for all $n \in \omega$.
Below, when we consider places, formally $\infty$ is also in the range. Recall from Chapter~\ref{acvfbackground} 
the connection between places and valuations.

\begin{lemma}\label{place}
 Let $v:L\rightarrow \Gamma$ be a valuation on a field $L$ with 
corresponding place $p:L\rightarrow \res(L)$, let $F$ be a subfield of $\res(L)$, 
and let $p':\res(L)\rightarrow F$ be a place which is the identity on $F$. Let 
$p^*=p'\circ p:L\rightarrow F$
be the composed place, with induced valuation $v^*:L\rightarrow \Gamma^*$.  
Suppose $a\in L$ with $p(a)\in \res(L)\setminus \{0\}$, and
$p^*(a)=0$.  

(i) If $b\in L$ with $v(b)>0$, then $0<v^*(a)<<v^*(b)$,

(ii) If $$\Delta:= 
 \{v^*(x), -v^*(x): x\in L, p(x)\not\in\{\infty, 0_{\res(L)}\}, p^*(x)=0\}\cup\{0_{\Gamma^*}\},
$$
then $\Delta$ is a convex subgroup of $\Gamma^*$ and there is an isomorphism 
$g:\Gamma^*/\Delta\rightarrow \Gamma$ such that $g\circ v^*=v$.
\end{lemma}

{\em Proof of Lemma~\ref{place}.} See (8.4)--(8.7) of \cite{endler}, and the 
discussion at the end of that chapter. \qed

\bigskip

{\em Proof of Proposition~\ref{domovergamma}.} 
The existence of $f$ is clear, and we focus on uniqueness.

So suppose $ f_i:L\rightarrow N_i$ (for $i=1,2$) are two maps inducing $g$ and satisfying 
the hypotheses of the proposition, with $v_i$ denoting the valuation on $N_i$. We will construct an isomorphism
 $\ell:N_1\rightarrow N_2$
over $M\cup g(\Gamma(L))$ with $f_2=\ell\circ f_1$. The idea is to perturb the valuations $v_i$ to obtain valuations 
$v_i'$ satisfying the additional hypothesis of Proposition~\ref{dom1} (value group independence over $\Gamma(C)$). 
This will yield that the $v_i'$ are conjugate, and it will follow that the $v_i$ are conjugate too.

Choose $a_1,\ldots,a_r\in L$ so that if $e_1=f(a_1),\ldots,e_r=f(a_r)$, then 
$v(e_1),\ldots,v(e_r)$ form a ${\mathbb Q}$-basis of  $g(\Gamma(L))$ over $\Gamma(C)$.  Also choose 
$b_1,\ldots,b_s\in L$ so that
$\res(b_1),\ldots,\res(b_s)$ form a transcendence basis of $k(L)$ over $k(C)$.
By Lemma~\ref{intdesc2} applied to $f_i(L)$ and $M$, for $i=1,2$ the elements
$$\res(f_i(a_1)/e_1),\ldots,\res(f_i(a_r)/e_r),\res(f_i(b_1)),\ldots,\res(f_i(b_s))$$
are algebraically independent over $k(M)$. Let $m: k(f_1(L))\rightarrow k(f_2(L))$ be any isomorphism over $k(C)$.

For $j=0,\ldots,r$ and $i=1,2$ let
$$  R^{(j)}_i:=\acl(k(M),\res(f_i(a_1)/e_1),\ldots, \res(f_i(a_j)/e_j),
 \res(f_i(b_1)),\ldots,\res(f_i(b_s)))
$$
(so $R^{(0)}_i:=\acl(k(M),\res(f_i(b_1)),\ldots,\res(f_i(b_s)))$).
For each $i,j$ choose a place $p_i^{j}:R_i^{(j+1)}\rightarrow R_i^{(j)}$ over $R_i^{(j)}$
 (i.e.
which are the identity on $R_i^{(j)}$), and a place $p_i^*:k(N_i)\rightarrow R_i^{(r)}$ over
 $R_i^{(r)}$. 
(It will eventually turn out that $p_i^*$ is the identity and $R_i^{(r)}=k(N_i)$.)
Let $p_{v_i}:N_i\rightarrow k(N_i)$ be the place corresponding to the given valuation $v_i$.
Also let $p_{v_i'}:N_i\rightarrow R_i^{(0)}$ be the composed place 
$p_i^0\circ p_i^1\circ \ldots \circ p_i^{r-1}\circ p_i^*\circ p_{v_i}$
with corresponding valuation $v_i'$ on $N_i$ (determined up to value group isomorphism).
 Let $N_i'$ be the valued field $(N_i,v_i')$,
(so $N_i,N_i'$ have the same field structure, but possibly different valuations).

Let $i\in \{1,2\}$. Then all the $p_i^j$ and $p_i^*$ are the identity on 
$k(M)$, so $p_{v_i}$ and $p_{v_i'}$ 
agree on $M$; so $(M,v_i)$ and $(M,v_i')$ are isomorphic, as a valuation is determined up to isomorphism of value groups 
by the corresponding place. We therefore
identify $(M,v_i)$ and $(M,v_i')$, identifying also the value groups. 

For $i=1,2$ define the valued field embedding $f_i':L\rightarrow N_i'$ by $f_i'(x)=f_i(x)$.
Observe here that since $p_{v_i'}$ is the identity on $\res(f_i(L))$, it determines (up to isomorphism of value groups)
a valuation on $f_i(L)$ which agrees with $v_i|_{f_i(L)}$. That is, for $x,y\in f_i(L)$,
$v_i(x)\leq v_i(y)\Leftrightarrow v_i'(x)\leq v_i'(y)$. Since we have identified the value groups of $(M,v_i)$ and $(M,v_i')$,
we should not identify those of $(f_i(L),v_i)$ and $(f_i(L),v_i')$, but they are isomorphic.  
The map $f_i'$ induces some $h_i':\Gamma(L)\rightarrow \Gamma(N_i')$ by
$h_i'(v(x))=v_i'(x)$. 

Successively applying Lemma~\ref{place} (with $p^*$ replaced by
$p_i^j\circ \ldots \circ p_i^{r-1}\circ p_i^*\circ p_{v_i}$), we find that
$$0<v_i'(f_i(a_1)/e_1)<<\ldots <<v_i'(f_i(a_r)/e_r)<<v_i'(b)$$ for any $b\in M$ with $v(b)>0$. 
Let $\Delta_i$ be the subspace of the ${\mathbb Q}$ vector space  $\Gamma(N_i')$ generated  
 by $v_i'(f_i(a_1)/e_1),\ldots,v_i'(f_i(a_r)/e_r)$.
By Lemma~\ref{place}, $\Delta_i$ is a convex subgroup of $\Gamma(N_i')$.
Since $\Gamma(L)$ and $h_i'(\Gamma(L))$ both have ${\mathbb Q}$-rank $r$ over $\Gamma(C)$, 
$\Gamma(N_i')=\Delta_i \oplus \Gamma(M)$. (This yields that $p_i^*$ is the identity, since 
otherwise $\Gamma(N_i')$ would have additional elements infinitesimal with respect to
 $\Gamma(M)$.)
In particular, $h_i'(\Gamma(L))\cap \Gamma(M)=\Gamma(C)$. Thus, there is an isomorphism
$\ell^*:\Gamma(N_1')\rightarrow \Gamma(N_2')$ which is the identity on $\Gamma(M)$ and
 satisfies
$\ell^*(v_1'(a_j))=v_2'(a_j)$ for $j=1,\ldots,r$. Likewise, there is a map
$m^*: k(N_1')\rightarrow k(N_2')$ which is the identity on $k(M)$ and extends $m$.

By Proposition~\ref{dom1} there is a valued field isomorphism $\ell:N_1'\rightarrow N_2'$ 
which 
is 
the identity on $M$, extends $\ell^*$ and $m^*$,
and satisfies
$\ell(f_1(x))=f_2(x)$ for all $x\in L$. We must verify that $\ell$ is also an 
isomorphism
$(N_1,v_1) \rightarrow (N_2,v_2)$. However, we have 
$\Gamma(N_i'):=\Delta_i \oplus \Gamma(M)$.
Let $\pi$ denote the projection map to the second coordinate.
Then, by the last assertion of Lemma~\ref{place}, for $x\in N_i$ we have 
$v_i(x)=\pi(v_i'(x))$. It follows that
$$\ell^*(v_1(x))=\ell^*(\pi(v_1'(x)))=\pi(\ell^*(v_1'(x)))=\pi(v_2'(x))=v_2(x)$$
for all $x\in N_1$, as required. \qed

\begin{remark}\rm
Our proof gave a slightly stronger statement than Proposition~\ref{domovergamma}. We showed that the isomorphism 
$\ell:N_1\rightarrow N_2$ can be chosen to extend any given isomorphism
$k(f_1(L))\rightarrow k(f_2(L))$ which is the identity on $k(C)$.
\end{remark}

\begin{theorem}\label{fulldom}
(i) Suppose that $C \leq L$ are   valued fields with $C$
maximally complete,   $k(L)$ a regular extension of $k(C)$,
and $\G_L/\G_C$ torsion free. Let $a$ be a  sequence from $\U$, $a \in \dcl(L)$.
   Then
$\tp(a/C \cup \Gamma(Ca))$ is stably dominated. 

(ii)  Let $C$ be a maximally complete algebraically closed valued field, $a$ a  sequence from $\U$,
and $A=\acl(Ca)$.   Then $\tp(A/C,\Gamma(Ca))$ is stably dominated.
\end{theorem}

{\em Proof.} (i) When $C$ is algebraically closed, this  is immediate from 
Proposition~\ref{domovergamma}. In general, 
let $\bar{C}$ be a maximally complete immediate
extension of $C^{\alg}$.  The assumptions imply that 
$\tp(L/C) \vdash \tp(L/C,k(C)^{\alg},\Qq \tensor \G(C)) = \tp(L / C, k(\bar{C}), 
\G(\bar{C}))$.  By Corollary~\ref{fulldom2}, 
$\tp(L/C) \vdash \tp(L/ \bar{C})$.   It follows that $\tp(L/C,\G(L)) \vdash \tp(L/\bar{C},\G(L))$.
But $\G(L \bar{C}) = \G(L)$.   
Since $\tp(L/ \bar{C},\G(L \bar{C}))$ is stably dominated, so is $\tp(L/C, \G(L))$.
By Corollary~\ref{st-dom-eq-0} (ii), $\tp(a/C \cup \Gamma(Ca))$ is stably dominated. 

(ii) This is a special case of (i); let $C \leq L$ be any valued field with $a \in \dcl(L)$.  
 
\qed

\medskip

 \begin{remark} \label{fulldom-r} \rm Suppose $C,L$ are as in Theorem~\ref{fulldom}.
  Let $a_1,\ldots,a_m \in L$ be such that $|a_1|,\ldots,|a_m|$ 
 is a $\Qq$-basis for $\Qq \tensor \G(L)$ over $\Qq \tensor \G(C)$.  Then
 $\tp(L/C \cup \Gamma(L))$ is  stably dominated by $k(L) \union \{r(a_1),\ldots,r(a_m) \}$
 over $C \cup \Gamma(L))$, where $r(a_i)$ denotes the image of $a_i$ under the reduction map
$|a_i|R\rightarrow \red(|a_i|R)$.  \end{remark}
 
{\em Proof.}  Let $f$ enumerate $k(L) \union \{r(a_1),\ldots,r(a_m) \}$, and let 
$C' = C \union \G(L)$.  By Lemma~\ref{intdesc}, $\St_{C'}(L) \subseteq \acl(C,f)$.
This, together with Theorem~\ref{fulldom}, gives (ii) of Proposition ~\ref{stab1},
and it remains to show (i), i.e. that $\tp(L/C',f) \vdash \tp(L/ \acl(C'),f)$.  
 
The valued field $C$ is Henselian.  Replacing $L$ by the Henselian hull $L^h$
will not change $\G(L)$ or $\k(L)$, so we may assume $L$ is Henselian too.
Thus  $\Aut_v(L^{\alg}/L)=\Aut(L^{\alg}/L)$ (any
field automorphism of $L^{\alg}$ over $L$ is a  valued field automorphism.)
Since $L$ is a regular extension of $C$, the homomorphism
$\Aut(L^{\alg} /L) \to \Aut(C^{\alg}/C)$ is surjective.  Equivalently,
$\tp(C^{\alg} / C) \vdash \tp(C^{\alg}/L)$.   But   $\dcl(C^{\alg}) = \acl(C)$;
so $\tp(\acl(C) / C) \vdash \tp(\acl(C)/L)$.  Since $C \subseteq \dcl(C,\Gamma(L),f) \subseteq \dcl(L)$, it follows a fortiori that
$$\tp(\acl(C) / C,\G(L),f) \vdash \tp(\acl(C)/L).$$
Thus,
$$\tp (L / C, \G(L),f) \vdash \tp(L / \acl(C),\G(L),f)$$
But by Lemma~\ref{3.4.12}, $\acl(C') = \acl(C \cup \G(L)) = \dcl(\acl(C), \G(L))$.  So 
$\tp(L / C',f) \vdash \tp(L/ \acl(C'), f)$ as required.  \qed

\chapter{Invariant types}\label{sequential2}
\section{Examples of  sequential independence} \label{examples}
We give here four examples concerning issues with sequential independence. We work in the field sort throughout 
the chapter. In the first two examples, $\acl$ denotes  field-theoretic algebraic independence (which is model-theoretic 
algebraic independence in ACVF in the field sort). At the end of the section, we state  a result 
giving a setting where forking and sequential independence coincide.

\begin{example} \label{promex}\rm
As promised after Lemma~\ref{ortheasy}, we give an example where, for $A=\acl(Ca)$, 
 $\tp(a/C)\perp \Gamma$ even though
$\trdeg(A/C)\neq \trdeg(k(A)/k(C))$. Let $M$ be a maximally complete 
algebraically closed non-trivially valued field with archimedean value group. 
Let $X$ be transcendental over $M$. There is a norm $|-|$ on $M(X)$  extending that on $M$,
defined by
$|\Sigma_{i=0}^m a_iX^i|=\Max\{|a_i|:0\leq i\leq m\}$. Then $X$ is generic in $R$ over $M$. 
Let $M\{X\}$ be the ring of convergent power series $\Sigma_{i=0}^{\infty} a_i X^i$ such that
$\lim_{i\rightarrow \infty}|a_i|=0$. Let $L$ be its field of fractions. 
The norm $|-|$ extends from
$M(X)$ to $M\{X\}$ by the same definition (maximum norm of coefficients), and thence to $L$. 
Choose a sequence $(n_j:j=1,2,\ldots)$ with $0<n_1<n_2<\ldots$ and $\lim n_{j+1}/n_j=\infty$.
Let $(a_i:i=1,2,\ldots)$ be a sequence from $M$ with $\lim |a_i|=0$. Put
$Y:= \Sigma_{i=1}^\infty a_i X^{n_i} \in L$. Then by the argument on p. 93 of \cite{rib},
$Y\not\in \acl(M(X))$. Now $Y$ is 
a limit of a pseudoconvergent sequence $(b_i:i\in I)$ (the polynomial truncations of $Y$) of elements of $M(X)$, so $Y$ 
is generic over $M(X)$ in a chain of $M(X)$-definable balls with no least element, namely the balls
$B_{\leq |Y-b_i|}(Y)$. Thus, by Lemma~\ref{eqorthog},
$k(\acl(M(X,Y)))=k(\acl(M(X)))$. By the definition of the norm on $L$, $\Gamma(L)=\Gamma(M)$, so 
$\Gamma(\acl(M(X,Y)))=\Gamma(M)$. It follows by Proposition~\ref{maxim} (with $C$ and $F$ of
 the proposition both equal to $M$) that
$\tp(XY/M)\perp \Gamma$. However $XY$ adds only 1 to the transcendence degree of 
the residue field, but $\trdeg (M(X,Y)/M)=2$. 

As a more general way to obtain such examples, let $\Delta$ be a countable divisible ordered 
abelian group, and $k,k'$ be algebraically closed fields
such that $k'$ is a transcendence degree 1 field extension of $k$. The field 
$C:=k((t))^\Delta$
of generalised power series of well-ordered support with exponents in $\Delta$ and 
coefficients in $k$ is maximally complete; and its extension $L:=k'((t))^\Delta$
has transcendence degree $2^{\aleph_0}$ over $C$. However, $\Gamma(L)=\Gamma(C)=\Delta$,
so by Proposition~\ref{maxim}, $\tp(C(a)/C)\perp \Gamma$ for any finite sequence $a$ from $L$.

\end{example}

\begin{example}\label{orthsing}\rm
We give an example where $\tp(A/C)\not\perp \Gamma$, but for every singleton 
$a\in \acl(A)$, $\tp(a/C)\perp \Gamma$. 

Let $F$ be a trivially valued algebraically closed field of characteristic 0, let $T$ be an indeterminate, 
and 
let the power series ring $F[[T]]$ have the usual valuation with $v(T)=1$. As usual, we 
shall write $F((T))$ for the corresponding field of Laurent series. Put
$C=\acl(F(T))$. Let $b,c \in F[[T]]$ 
be power series algebraically independent over $C$. Let $F'$ be a trivially valued 
algebraically closed extension of $F$, and $a\in F'\setminus F$. Put
$A=\acl(C(a,ab+c))$, and let $d\in A\setminus C$.

As $A\leq F'((T))^{{\mathbb Q}}$, $d$ has a Puiseux series expansion over $F'$, so $d\in F'((T^{1/n}))$ for some $n$. 
If $d\in F((T^{1/n}))$, then 
$k(F(T,b,c,d))=k(F)$, so $a\not\in \acl(F(T,b,c,d))$.
Hence, as $d\in \acl(F(T,a,b,c))$, we have $d\in F(T,b,c)$. 
This however is impossible, since
by an easy argument (see e.g. p.28 of \cite{mark}) $A\cap \acl(F(T,b,c))=C$. Thus, 
some term $a_mT^{m/n}$ of the Puiseux series for $d$ has $a_m \in F'\setminus F$. 
Choose the least such $m$, and put $f:=\sum_{m'<m} a_{m'}T^{m'/n}\in C$. Then 
$g:=(d-f)T^{-m/n}=a_m + e\in \acl(C(d))$ with $v(e)>0$. 
Thus, as $a_m\not\in F$,  $\res(g)\not\in k(C)$, so $\tp(g/C)\perp \Gamma$.
As $g$ and $d$ are interdefinable over $C$, $\tp(d/C)\perp \Gamma$.

Suppose for a contradiction that $\tp(A/C)\perp \Gamma$. Then by 
Corollary~\ref{equivorthstabdom}, $\tp(A/C)$ is stably dominated, that is, $\tp(A/C)$
is dominated by $A':=k(A)$; here we use that $\VS_{k,C}\subseteq \dcl(C \cup k(C))$, as $C$ is a field. However, as 
$k(C)$ is algebraic over $k(F((T)))$, we have
$A'\dnf^g_C F((T))$. Hence $A\dnf^d_C F((T))$, so 
$\tp(A/A'C)\vdash \tp(A/A' Cbc)$. Thus, the formula
$y=xb+c$ follows from $p(x,y)=\tp(a,ab+c/A' C)$, so $ab+c\in \acl(A'Ca)$. By quantifier 
elimination, as $A'$ consists of residue field elements, $ab+c \in \acl(Ca)$, which 
is impossible. 
\end{example}

We give two further examples concerning sequential extensions.  The
first takes
 place entirely in the value group (so in the theory of divisible ordered abelian groups),
but could be encoded into ACVF. It shows that a type of the form $\tp(A/C)$ can 
have many different sequentially independent extensions, via different sequences.
The second, in ACVF, shows how a type can have an invariant extension not arising 
through sequential independence
(at least via a sequence of field elements).

\begin{example}\label{additive} \rm
We suppose that $(\Gamma,<,+,0)$ is a large saturated divisible ordered abelian group, written additively, and 
with theory denoted DLO. If $C\subseteq B \subset \Gamma$ and $a\in \Gamma$, write
$a\dnf^g_C B$ if for any $b\in \dcl(B)$ with $b>a$ there is $c\in \dcl(C)$ with $a<c\leq b$. In other words, 
$a$ is chosen generically large over $B$ within its type over $C$. We extend the $\dnf^g$ notation sequentially 
as for ACVF; it has similar properties. In particular, every type over $C$ has an $\Aut(\Gamma/C)$-invariant 
extension over $\Gamma$.

Note that when $F\subseteq L $ are valued fields, $a,F,L$ contained in some $M \models \ACVF$ with 
$C = \{v(x): x \in F \}$, $B = \{v(x): x \in L \}$, and $v(a)\not\in F$, we have
$a\dnf^g_ F L$  if and only if $-v(a) \dnf^g _{C} B$.

Below, we write $B_n^+({\mathbb Q})$ for the  group of lower triangular
 $n\times n$ matrices over ${\mathbb Q}$ with strictly positive entries on the diagonal.

Let $n>1$. We show that there are $C\subset B\subset \Gamma$ and $a=(a_1,\ldots,a_n)\in \Gamma^n$
such that if $A=\dcl(Ca)$, then $A$ has infinitely many distinct
sequentially independent extensions over $B$. 

To see this, choose $C=\dcl(C)$ in $\Gamma$. Let $P_1,\ldots,P_n$ be cuts in $C$ (so
solution
 sets of complete 1-types over $C$), and pick $b_1\in P_1, \ldots, b_n\in
P_n$. We choose the $P_i$ so that
$0<b_1<\ldots<b_n$, and so that for each $i$, $\{x-y:y<P_i<x\}=\{x\in C:x>0\}$. We also
assume that the $P_i$ are independent, in the sense that if $x_1\in
P_1,\ldots,x_n\in P_n$ then $\{x_1,\ldots,x_n\}$ are ${\mathbb Q}$-linearly
independent over $C$; in particular, this choice of the $x_i$ determines a
complete $n$-type over $C$,
 denoted $q$.
Put $B:=\dcl(Cb)$ where $b=(b_1,\ldots,b_n)$. Choose $a=(a_1,\ldots,a_n)$
with $a_i\in P_i$ for each $i$, so that $a\dnf^g_C B$, so $a_i$ is
chosen generically large in $P_i$ over $Ba_1\ldots a_{i-1}$. Let
$A:=\dcl(Ca)$.

Now for each $D\in \GL_n({\mathbb Q})$, let $d:=Da$ (viewing $a,d$ as column
vectors). Then $A:=\dcl(Cd)$, and $D$ determines an extension over $B$ of
$\tp(A/C)$ which is sequentially independent via $d'$, where $d'\equiv_C
d$. Let $q_D$ be the corresponding extension over $B$ of $\tp(a/C)$.
Clearly $q_D$ has an $\Aut(\Gamma/C)$-invariant extension. 

We claim that if $D$ and $D'$ lie in distinct cosets
of $B_n^+({\mathbb Q})$, then $q_D\neq q_{D'}$.
 Clearly $\GL_n({\mathbb Q})$ acts transitively on the set of 
extensions over $B$ of $q$ of the above form, so it suffices to show that
the stabiliser of $\tp(a/B)$ is contained in $B_n^+({\mathbb Q})$.

Let $e_1:=a_1-b_1,\ldots,e_n:=a_n-b_n$. Then $e=(e_1,\ldots,e_n)$ is an
(inverse)
lexicographic basis for the ordered  ${\mathbb Q}$-vector space 
$$\Delta:=\{0\}\cup
\{x\in \dcl(Ba): 0<|x|<c \mbox{~for all~} c\in C^{>0}\}.$$ Now the flag
$$(0,{\mathbb Q} e_1,{\mathbb Q} e_1+{\mathbb Q} e_2,\ldots,{\mathbb Q} e_1+\ldots +{\mathbb Q} e_n)$$ 
is uniquely determined by the ordered vector space structure on
$\Delta$; the $k^{\th}$
 element is
the unique convex subspace of $\Delta$ of dimension $k$.  

We must show that if $D\in \GL_n({\mathbb Q})$, $d=Da$, and $a\dnf^g_C
B$ and $d\dnf_C^g B$ both hold, then $D\in B_n^+({\mathbb Q})$. So
suppose $d\dnf^g_C B$. If $b':=Db$ with $b':=(b_1',\ldots,b_n')$,
then $d_i$ is chosen in the cut over $C$ containing $b_i'$. Thus, arguing
as above, there will be $y_1,\ldots y_n \in B$ such that $$d_1-y_1\in {\mathbb Q}^+e_1, \ldots,
 d_n-y_n \in {\mathbb Q} e_1+\ldots +{\mathbb Q} e_{n-1}+{\mathbb Q}^+e_n.$$ It is
easily checked 
that if this holds then $D\in B_n^+({\mathbb Q})$.
\end{example}

In the rest of the chapter, the field valuation will be written additively.
In the next lemma, and the example which follows, we will use the
completion $\bF$ of a valued field $F$ with value group $\Delta$.  
For convenience, we assume that $\Delta$ is countable and Archimedean, so in particular has countable
 cofinality. 
The completion is a 
standard construction; it can be obtained as the field of limits of Cauchy
sequences from $F$, where $(a_n)$ is Cauchy if $v(a_{n+1}-a_n)$ is increasing and
cofinal in $\Delta$. The field  $\bF$ is characterized by the facts that $v(\bF) =
v(F)$, every Cauchy sequence in $\bF$ converges, and every element of
$\bF$ is the limit of a sequence from $F$.
It is easy to see that $\bF$ is  algebraically closed if $F$ is.   

We look at sequentially generic extensions (in the sense of ACVF) of types over $F$ of tuples
from $\bF$; this is 
 interesting for its own sake, and will also be used below.    We assume
$F$ is algebraically closed, $v(F)=\Delta \neq (0)$,  and $\trdeg_F (\bF)>1$.
 For any field $L$, $L^a$ denotes its algebraic closure.
Examples with $F \neq \bF$ (both algebraically closed) include
 Krasner's $F=({\mathbb Q}_p)^a$, with completion $C_p$; or   $F = F_0(t)^a$,
$F_0$ a countable field;
 note that $\bF$ is   uncountable.

 Suppose $F \subset F'  $ where $F',\bF \subset M$
 (a sufficiently saturated algebraically closed valued field), and $F'$ admits a valued
field embedding $\phi$ over $F$ into $\bF$. In this case, $\phi$ is
unique.  Working inside $M$, if $d\in F'$ then $\phi(d) $ is the unique element of $\bF$
with  $v(d-\phi(d)) > \Delta$.

Let $\rho(d) =  v(d-\phi(d))$, and  $\bar{\rho}(d) = \rho(d) +\Delta \in
\Gamma(M)/\Delta$.  Also, 
 let $P\br (x) = p^{{\mathbb Z}} \br(x)  $, where $p=1$ if $\char(F)=0$, and
$p=\char(F)$ otherwise.  
(So $P\br (x)$ gives less information than $\br (x)$ but a bit more than
the ${\mathbb Q}$-space ${\mathbb Q} \br(x)$.)

\begin{lemma}\label{lem1}  Assume $n$ is finite, $F' = F(d_1,\ldots,d_n)^a$ embeds into
$\bF$ over $F$ via $\phi$ as above, and $\Delta <
\rho(d_1) << \rho(d_2) << \ldots << \rho(d_n)$.
   Then $P\br$ takes just $n$  non-zero values on $F' $, namely  the
$P\br(d_i)$.

\end{lemma}

Before proving this, we mention a fact which is presumably well-known.

\begin{lemma}\label{genericzero}
 Let $h = h(x_0,....,x_n) $ be an irreducible polynomial over a
field $L$.  Let $h(a)=0$, where  $\trdeg(L(a)/L) = n$. Let $G(y_0,...,y_n) :=
h(a_0+y_0,...,a_n+y_n)$ be the Taylor expansion about $a$.    Fix $i$ between
$0,...,n$. Let $G_i(y_i) = G_i(0,...,0,y_i,0,..,0)$, and write $G_i(y_i) =
g_i(y_i^{q_i})$ where $q_i$ is 1 in characteristic 0, and in characteristic $p$, $q_i$ is a power of $p$
chosen so that
$g_i$ is a polynomial with some monomial of degree not divisible by
$p$. Then $g_i$  has nonzero  term of degree 1, and zero constant term.
\end{lemma}

{\em Proof of Lemma~\ref{genericzero}.}  First, $g$ has zero constant term as $h(a)=0)$.
We may suppose $i=0$.  Then by the genericity of the zero $a$,
$h(x,a_1,...,a_n)$ is irreducible over 
$L(a_1,...,a_n)$, and $G_0$ is the expansion of $h$
about $a_0$.  So we may assume $n=0$.
Write $h(x) = H(x^q)$, where  $q$ is  a power of $p$ and $H$ has a monomial of degree not divisible by
$p$ (in characteristic 0, $q=1$).  Then the derivative $H'$ is not identically $0$. Now $\deg_x H'(x^q) <
\deg_x H(x^q) = \deg h$, so $H'(x^q)$ cannot vanish at a root of $h$, by the irreducibility of $h$. So
$H'(a_0^q)\neq 0$.  However,  this is the linear coefficient of the expansion of $h$ about
$a_0$. \qed

\medskip

{\em Proof of Lemma~\ref{lem1}.}  Let $d_0 \in F' \setminus F$; we wish to compute
$\rho(d_0)$. The elements $d_0,\ldots,d_n$ are algebraically dependent over
$F$.  If some proper subset $d_0,d_{i_1},\ldots,d_{i_k}$ is algebraically
dependent,
 then we can use induction.    So assume they are not.  Let  $d=(d_0,\ldots,d_n)$. We have $h(d)=0$
for some  irreducible    $h \in F[X_0,\ldots,X_n]$; and $g(d) \neq 0$ for
all nonzero $g$ of smaller total degree. 

 Let $d'_i = \phi(d_i), d'=(d'_0,\ldots,d'_n),e_i=d'_i-d_i$.   Then  $h(d')=0$. 

Using multi-index notation, expand $h(d_0+y_0,\ldots,d_n+y_n) = \sum
h_{\nu}(d) y^\nu$, with $h_{\nu} \in F[X_0,\ldots,X_n]$.  So, since $h(d')=0$, 

$$   \sum h_{\nu}(d) e^\nu = 0.$$
In particular, the  lowest two  values $v (h_{\nu}(d) e^\nu)$   must be equal.

Now as $F'$ embeds into $\bar{F}$ over $F$, $v(F') = v (F) = \Delta$. 
 So $v(h_{\nu}(d))  \in
\Delta$.    Note
that $h_\nu(d) \neq 0$ whenever $h_\nu \neq 0$, since $h_\nu$ has smaller
degree than $h$.   

Consider the various monomials    $h_{\nu}(d)y^\nu$.  Since $d_0,d_2,\ldots,d_n$ are
independent, $y_1$ occurs in $h$; otherwise, $h(d)$ does not involve $d_1$.  
By Lemma~\ref{genericzero}, we may write $h(d_0+y_0,\ldots,d_n+y_n)=g(y_1^{q_1}) +f(d,y)$,
where: 
$g$ is a  polynomial in 1 variable (with coefficients polynomial
 in the $d_i$), $q_1=1$
in characteristic 0, and in characteristic $p$, $q_1$ is a power of $p$,
$g$ has zero constant term and non-zero linear term,  
and each monomial in $f$ involves some $y_i$ other than $y_1$. Thus, 
if
$\nu_1$ is the multi-index corresponding to $y^{\nu_1} = {y_1}^{q_1}$, 
then $h_{\nu_1}\neq 0$, so
$h_{\nu_1}(d) \neq 0$. 

Any  $\nu$ involving some $y_j$, $j>1$ has:  $\Delta< v(e^{\nu}) -
v(e^{\nu_1})$.  Thus in looking for the
summands of  smallest valuation, such $\nu$ can be ignored.
The same holds  for any $\nu$ involving $y_1$ to a higher power than $q_1$.
  Note
that by definition of $\phi$ we have $v(e_0) = v(d_0-d_0') = v(d_0 - \phi(d_0)) > \Delta$.
So $v(e_0e_1) - v(e_1) > \Delta$ and so  terms involving $y_0y_1$ must be of 
higher valuation than $e_1$ also.  
 
This leaves   indices $y^\nu$ of the form  $y_0^{m}$ only.  
Let 
$\nu_0$ be the least such;  $y^{\nu_0}=y_0^{q_0}$, $q_0$ least possible.
  (Again, $q_0=1$ in characteristic $0$, and is a power of $p$ in characteristic $p$.) 
As noted above, $v(e_0)> \Delta$.   
 Then terms $h_\nu(d) e^\nu$ where
$y^\nu=y_0^k$, $k>m_0$, again cannot have least value among the
summands, since they
have value greater than $h_{\nu_0}(d) e^{\nu_0}$.    We are
left with only two candidates for summands of least value; so they must
have equal value.

$$v ( h_{\nu_1}(d) e^{\nu_1} ) = v( h_{\nu_0}(d) e^{\nu_0}) $$

So $v(e_0) \in  (q_1/q_0) v(e_1) + \Delta$.   Thus $\rho(d_0) =
(q_1/q_0) v(e_1) + \Delta$,
and so $P\br(d_0) = P\br(d_1)$.  \qed

\bigskip

Working in the theory DLO, let $\Gamma$ be a large saturated model of DLO,
$\Delta$ a small subset of $\Gamma$ with no greatest element, and $B \subset \Gamma$.
Let $q^\Delta(u)$ denote the $\Aut(\Gamma/\Delta)$-invariant 1-type over $\Gamma$ containing the 
formulas $u>\gamma$ for all 
$\gamma \in \Delta$, and $u<\gamma$ for all $\gamma \in \Gamma$ with $\Delta<\gamma$. 
For $B\subset \Gamma$, let $q^\Delta_B$ denote the restriction $q^\Delta|B$. 

Let $F,M \models ACVF$ with $\bF \subset M$.
Let $\Delta:=\Gamma(F)$ and $c \in \bF
\setminus F$.   For any $L\supset F$, consider the generic extension  to $L$ (in the
$\dnf^g$-sense)
of $p_c(x)=\tp(c/F)$. 
Denote this by $p_c|L$.
Then, if $c\in L$,  $p_c|L$ contains all formulas of the DLO-type 
$q_L^\Delta(v(x-c))$.  Indeed, this determines $p_c|L$.   Note that $p_c|\U$ is
 $\Aut(\U/F)$-invariant.

We extend this notation for 2-types.  If $c=(c_1,c_2)$ with $c_i
\in \bF\setminus F$, and $L\supseteq F$,  let $p_c(x_1,x_2)|L$ be the sequentially 
generic extension
over $L$ of $\tp(c_2,c_1/F)$ (so we take the $x_2$ variable first).
 Then if $F\models \ACVF$ and $L \supseteq \bF$,
$p_c | L $ is implied by: $q^\Delta_L (v(x_2-c_2)),
q^\Delta_{Lx_2} (v(x_1-c_1))$.

\begin{example} \label{lotsinv}      \rm
There is $F\models \ACVF$, and $A=\acl(A)\supset F$, such that $\tp(A/F)$
 has an $\Aut({\cal
\U}/F)$-invariant
 extension 
$p$ over
${\U}$, such that if $A'$ realises $p$ then there is no $\acl$-generating
sequence $d$ of field elements from $A'$ with $d\dnf^g_F {\U}$.
\medskip

To construct this, we choose 
 $F$ as above,  let $\Delta:=\Gamma(F)$, and pick $c_1,c_2\in \bF$, algebraically
 independent over $F$.

For any $B\supseteq F$, let
$q_B :=q^\Delta_B$.
Then let
$s_B(u,v)$ be the 2-type over $B$ determined by specifying 
$q_B(u)$ and
$q_{Bu}(v)$. Let $D$ be any invertible matrix over ${\mathbb Q}$ which
is not lower-triangular, chosen so that the rows are positive elements of
${\mathbb Q}^2$ in the lexicographic order. Then let $r_B(u,v)$ be the 2-type
of DLO given as $$(\exists x,y)(s_B(x,y) \wedge (u,v)=D(x,y)).$$ Then $r_B$ is a
complete 2-type over $B$. Finally, $t^{c_1,c_2}_B(x_1,x_2)$ holds if we
have $r_B(v(x_1-c_1),v(x_2-c_2))$.

Now $t^{c_1,c_2}_B$ is a complete 2-type in the field sort over 
any $B\supseteq Fc_1c_2$.
Furthermore, by its definition,
$t^{c_1,c_2}_{{\U}}$ is $\Aut({\U}/Fc_1c_2)$-invariant.

\medskip

{\em Claim.} $t^{c_1,c_2}_{{\U}}$ is $\Aut({\U}/F)$-invariant.

\medskip

{\em Proof of Claim.} Let $\sigma \in \Aut({\U}/F)$, and put 
$(d_1,d_2):=\sigma((c_1,c_2))$. 
We must show $t^{c_1,c_2}_{{\U}}=t^{d_1,d_2}_{{\U}}$.
Thus, we must show $t^{c_1,c_2}_{{\U}}$ includes 
$r_{{\U}}(v(x_1-d_1),v(x_2-d_2))$. Suppose $N$ is a model 
containing $Fc_1c_2d_1d_2$,
and $(a_1,a_2)$ realises $t^{c_1,c_2}_N$. By the definition 
of $s$, $r$, and $t$,
we have $v(a_1-c_1)<v(b)$ and $ v(a_2-c_2) <v(b)$ for any $b\in N$ with $v(b)>\Delta$.
If $\gamma\in \Delta$ then by the condition on the $c_i$, there is for
each $i$ some $c_i'\in F$ with $v(c_i'-c_i)>\gamma$, and hence also
$v(c_i'-d_i)=v(c_i'-\sigma(c_i))>\gamma$. Thus, $v(c_i-d_i)>\gamma$ for all $\gamma \in
\Delta$. It follows that $v(a_i-c_i)<v(d_i-c_i)$ for each $i$. Thus,
$v(a_i-d_i)=v(a_i-c_i)$, so $r_{\U}(v(a_1-d_1),v(a_2-d_2))$ holds.

\medskip

Now suppose $a=(a_1,a_2)\models t^{c_1,c_2}_{\U}$, and
$A=(Fa)^{\alg}$,  $F':=(Fc_1c_2)^{\alg}$, $B:=(F'a)^{\alg}$. Let $p:=\tp(A/\U)$. 
Suppose that $A\dnf^g_F \bF$
via a sequence $d'=(d'_1,\ldots,d'_m)$ of field elements. By considering transcendence 
degree, $m=2$.
Now $A$ embeds into $\bF$ over $F$. Hence,
by the remarks before the example,  $d'\models p_d$ for some $d=(d_1,d_2)\in \bF^2$. 

 We  adopt the notation
$\phi,\rho,\bar{\rho}, P\br$ from Lemma~\ref{lem1}. Let $N:=\bF(d')^{\alg}$ and $F'':=F(d')^{\alg}$.
Clearly, $\Gamma(N)=\Delta+{\mathbb Q}\rho(d_1')+{\mathbb Q}\rho(d_2')$, with
$\Delta<<\rho(d_1')<<\rho(d_2')$. 
In Lemma~\ref{lem1}, we studied all possible sequentially
generic types of this form over $F$; and noticed that for such a type
$p_d$, 
 $P\br$ takes just
two values on $F'' \setminus F$. One of these values lies in the unique proper
nonzero convex subspace $S$ of $\Gamma(N) / \Delta$.  However, if $(a_1,a_2)
\models t_{\U}^{c_1,c_2} | \bF$,  then, by the choice of $D$,
$P \br (a_1) \neq P\br(a_2)$ and 
neither of
these two values corresponds to $S$. 
\qed

\end{example}

\section{Invariant types, dividing and sequential independence,}

We proceed towards a theorem asserting that invariant extensions, non-forking, non-dividing
and sequential independence are all essentially the same, over sufficiently good base structures.
Our previous results on stable domination reduce this  rather quickly to the case of $\G$.    
Most of the effort in this section is thus concentrated on $\G$.     

Let $A,B$ be divisible linearly ordered Abelian groups.  
By a {\em cut} in $A$ we mean a partition $A = L^- \union L^+$,
such that $a<b$ whenever $a \in L^-, b \in L^+$.  
We will sometimes describe a cut using a downward-closed subset $L' \subseteq A$ alone, 
having in mind $(L',L \setminus L')$.  

We let $A^{\geq 0} = \{a \in A: a \geq 0\}$.   Let $|h| = \pm h \geq 0$.
 If $A \leq B$ and $b \in B$, $b>0$,  the {\it   cut of } $b $ over $A$ is
$\ct_A(b) = \{a \in A: 0 \leq a \leq b\}$. If $(X^-,X^+)$ is any cut in $A$, the 
{\it corresponding convex subgroup}
is defined to be $H(X) := \{h \in A: h+X^- =X^-, h+X^+=X^+ \}$.      
We write $\Hct_A(b) = H(\ct_A(b))$.  

An embedding $A \subset B$ of divisible  ordered  Abelian groups
is called an   {\it i-extension} \index{i-extension} of $A$ if  there exists no $b \in B$ with $\ct_A(b)$
closed under addition.  Equivalently,  the map $H \mapsto (H \meet A)$
is a bijection between the convex subgroups of $A$ and of $B$; the inverse map then
takes a convex subgroup $H'$ of $A$ to the convex hull of $H'$ in $B$.  This shows
transitive of i-extensions.

A divisible ordered Abelian group $A$ is {\em i-complete} \index{i-complete} if it has no proper i-extensions.  
We will show  that any divisible ordered 
abelian group has an i-complete i-extension, unique up to isomorphism.  

The algebraically closed valued field $C$ will be called {\em vi-complete}\index{vi-complete} if it is maximally complete, and 
the value group is i-complete.

Let $C \leq A,B$ be subgroups of a divisible ordered Abelian group.  We  say that $A$  is i-free from
$B$ over $C$ if  there   are no $b_1 \leq  a \leq b_2$, $b_1,b_2 \in B $,
$a \in A \setminus C$, with $\ct_C(b_1)=\ct_C(b_2)$.   This implies in particular that $A \cap B=C$.
If $A = C(a)$, we will also say  
that $\tp(a/B)$ is an i-free extension of $\tp(a/C)$.

Below, we use the notion of sequential independence in $\Gamma$ from Definition~\ref{g-indep} and 
Remark~\ref{uniquegen}(iii),
 rather than that of Example~\ref{additive}.  In particular, $a\dnf^g_C B$ if for any $b\in \dcl(B)$ with $b\leq a$, there
 is $c\in \dcl(C)$ with $b\leq c\leq a$.

\begin{theorem}\label{vik}
Let  $C\leq B$ be algebraically closed valued fields, and $A$ be a field which is finitely generated over $C$. Suppose 
that $C$ is vi-complete. Then the following are equivalent.

(i) $\tp(A/B)$ does not fork over $C$.

(ii) $\tp(A/B)$ does not divide over $C$.

(iii) $k(A)$ and $k(B)$ are free over $k(C)$, and
$\G(A)$ is i-free from $\G(B)$ over $\G(C)$.  

(iv) $A$ is sequentially independent from $B$ over $C$ via some sequence of generators.

(v) $\tp(A/B)$ has an $\Aut(\U/C)$-invariant extension over $\U$.

The number of extensions  of $\tp(A/C)$ to $B$ satisfying these equivalent properties is
 at most $2^{\dim_{\Qq}(\G(A)/\G(C))}$.
 
\end{theorem}

Note that by Example~\ref{lotsinv}, the implication (v) $\Rightarrow$ (iv) fails without some assumption on $C$.

\begin{remark}\label{NIP}\rm
ACVF is an NIP theory, that is it does not have the {\em independence property} of \cite{shelah}; this is an easy
 consequence of quantifier elimination. 
By a recent observation of Shelah \cite{shelah4} and also Adler, for any complete NIP theory $T$ 
with sufficiently saturated model $\U$
and some elementary submodel $M$, if $p\in S(\U)$ then $p$ is $\Aut(\U/M)$-invariant if and only if it does not
 fork over $M$. Thus, the equivalence 
(i) $\Leftrightarrow$ (v) above does not require the vi-assumption on $M$.
\end{remark}

From here to the proof of the theorem, all lemmas 
(\ref{C1} to \ref{C9}) refer purely to the   category of divisible  ordered  Abelian groups.    
The material is in part  similar to \cite{shelah3}.

\begin{lemma} \label{C1}  Every divisible ordered Abelian group $A$ has an $i$-complete $i$-extension.  \end{lemma}

{\em Proof.}  Take a transfinite chain $A_\alpha$  of proper i-extensions of $A$; it is clear
that the union of a chain of i-extensions is again one.  Thus it suffices to show that the
chain stops at some ordinal, i.e. to bound the size of any i-extension $B$ of $A$.  
We show the bound $|B| \leq 2^{2^{|A|}}$, which can be proved rapidly and suffices
for our purposes, though a more careful analysis should give $2^{|A|}$.    For $b_1 \neq b_2 \in B$, let 
$c(\{b_1,b_2\}) = \ct_A(|b_1-b_2|)$.
Suppose
$|B| > 2^{2^{|A|}}$.  By the Erd\"os-Rado theorem
(\cite{jech} Theorem 69 and ex. 69.1), 
 there exist $y < y' < y'' \in  B$ and a
  $C \subseteq A$,
such that   $\ct_A(y'-y)=\ct_A(y''-y) = \ct_A(y''-y') =  C$.  
It follows that $C$ is closed under addition; this contradicts the definition of an i-extension.
\qed

\begin{lemma}  \label{C3} Let $X=(X^-,X^+)$ be a   cut in $A^{\geq 0}$, $H=H(X)$.
Let $X'=(X'^-,X'^+)$ be the image of $X$ in $A/H$.  
Then the following are equivalent:
\begin{enumerate}
\item There   exists an i-extension $B$ of $A$ and  $b \in B$ with $\ct_A(b)=X$.   

\item $X'^-$ has no maximal element, and $X'^+$ has no minimal element.

\end{enumerate}
\end{lemma}

{\em Proof.}    If $X'^-$ has a maximal element $a+H$, then the coset $a+H$ is cofinal in $X^-$;
so if $\ct_A(b) = X^-$ then $\ct_A(b-a) = H \union A^{< 0}$.   If $A/H$ has a minimal element 
$d+H$ above $X'^-$ then $\ct_A(d-a) = H^{\geq 0} \union A^{< 0}$.  Both cases   contradict the definition of an $i$-
extension. 

 In the converse direction, form
$B = A(b) = \{a + \alpha b: a \in A, \alpha \in \Qq\}$, with $a + \alpha b < a' + \alpha'  b$
iff $(a-a')/(\alpha' -\alpha) \in X$.

 To show that $B$ is an i-extension,
let $J \subset B$ be a convex subgroup, $J' = J \meet A$.
 Let $  \alpha b + a \in J$.
We must find $a' \in J \meet A$, $ \alpha b + a < a'$.
Multiplying all by $|1/\alpha|$, we may assume $\alpha = \pm 1$.  

Replacing $X$ by the dual cut $(-X^+, -X^-)$   
corresponds to replacing
$b$ by $-b$.  Hence it suffices to prove the claim for positive $\alpha$, so  we may assume $\alpha = 1$.

Note that $a+X^-$ is downward closed, and is not the cut of a convex subgroup (else
$ X^- $ is the cut of $H-a$, but then the image of $X^-$ in $A/H$ has a maximal element.)
Thus there exists $c \in a+X^-$ with $2c > a+X^-$.  Let $a' = 2c$.  Then $X^-< a'-a$
while $(1/2)(a'-2a) \in X^-$.  So $(1/2)(a'-2a) < b<a'-a$ and $b+a < a' < 2(b+a)$.
As $a+b \in J$, $2(a+b) \in J$, so $a' \in J$ as required.
\qed

\begin{corollary} \label{C7} Let $A$ be i-complete, $X$ a downward-closed subset of $A$. 
Let $H=H(X)$.  
 Then there exists $a \in A$
such that $$X= \{x: (\exists h \in H) ( x\leq h+a ) \}$$ or $$X=\{x: (\forall h \in H) ( x<h+a ) \}$$ 
\end{corollary}

{\em Proof.}       By Lemma~\ref{C3} , $X^-$ has a maximal element $a+H$, or else 
there exists a minimal element $a+H$
above $X^-$.       \qed

\medskip

Cuts of the form $ \{x: (\exists h \in H) ( x\leq h ) \}$ or $\{x: (\forall h \in H) ( x<h ) \}$
will be called  cuts of convex subgroups.  
The converse of Corollary~\ref{C7} is also easily seen to be true: if every cut of $A$ is a 
translate of a cut of a convex subgroup, then   $A$ is i-complete.

The following uniqueness statement will not be used but is included as background.

\begin{proposition}   Up to $A$-isomorphism, $A$ has a unique i-complete i-extension.
\end{proposition}
{\em Proof.}     Let $B_1,B_2$ be two i-complete i-extensions of $A$.  Let $f: B _1 \to B_2$
be an isomorphism with maximal domain $A' \supseteq A$.  If $A' = B_1, fA'=B_2$
we are done.  Otherwise say $A' \neq B_1$, $c_1 \in B_1 \setminus A'$.  We may assume
notationally that $f=\Id_{A'}$.
Let $X$ be the lower cut of $c_1$ over $A'$, and let $H$ be the corresponding
convex group.    By Lemma~\ref{C3},
as the cut is realized in the i-extension $B_1$,  the image of $X$ in $A'/H$ has
no least upper bound.  (In particular, $X$ has no maximal element.)  Now
$H_2 := H_{B_2}$, the convex hull of $H$ in $B_2$, is the group belonging to
 the convex closure $X_2$ of $X$ in $B_2$.  For if $h \in H$, then $h+X \subseteq X$,
and so $h+X_2 \subseteq X_2$.   Thus $H_2 \subseteq H(X_2)$.
 Conversely if $h \in H(X_2) \meet A'$, $h + X_2 \subseteq X_2$,
so $h+X \subseteq X_2 \meet A = X$, hence $h \in H$.  Thus $H(X_2) \meet A' = H$.
By the contrapositive of   Lemma~\ref{C3} and
by the i-completeness of $B_2$, the image of $X_2$ in $B_2/H_2$ does have a least
upper bound $c+H_2$.   Observe that this coset of $H_2$ cannot be represented
in $A'$; for if $(c+H_2) \meet A' \neq \emptyset$, then it is a least upper bound
for the image of $X$ in $A'/H$; but such a bound does not exist.

Now for $x \in X_2$, $c-x \geq h$ for some $h \in H_2$,
 while for $x > X_2$,
$c-x \leq h$ for some $h \in H_2$.  As $H$ is cofinal in $H_2$ (in both directions),
we may take $h \in H$.
Now if $x \in X$, then $c-x \geq h$ for some $h \in H$,
so $c \geq x+h$.  If $c \leq x$, then $x+h \leq c \leq x $ so $x + H_2 = c+H_2$, contradicting
the fact that $(c+H_2) \meet A' = \emptyset$.  Thus $X < c$.
If $c\geq  x$ , $x \in A'$, $x > X$, then $c-x \leq h $ for some $h \in H$.
So $x \leq c \leq x+h$  and again $c+H = x+H$, contradiction.  Thus the cut of $c$
over $A'$ is precisely $X$.  Hence $f$ may be extended to an isomorphism
$A'(c_1) \to A'(c)$, contradicting the maximality of $f$.

\

The following definition is phrased in terms of a universal domain $\U$ for 
the theory of divisible ordered Abelian groups; $C,A',A,B$ are substructures. 
Compare \ref{fulldom2}.  

\begin{definition} \index{i-domination} 
Let $C \leq A' \leq A$ be divisible ordered Abelian groups.  
$A$ is  {\em i-dominated} by $A' $
over $C$
if for any  $B \supset C$ with  $A'$   i-free from $B$ over $C$ we have:
   $\tp(B/A') \vdash \tp(B/A)$.
\end{definition}

Let $C \leq A$.  
Given a cut $P$ of $C$, let $P_A = \{a \in A: \ct_C(a) = P \}$.
 
\begin{lemma} \label{C5}    Let $C \subset A' \subset A$ be divisible ordered Abelian groups.
Assume:  for every   cut $P$ of $C$, $P_{A'}$ is (left) cofinal in $P_A$.
Then $A$ is
i-dominated by $A'$ over $C$.   \end{lemma}

Note that using the reflection $x \mapsto -x$, one sees that the right-sided condition is the same.

\medskip

{\em Proof.}   Given  $B \supset C$ with  $A'$   i-free from $B$ over $C$,  the isomorphism type 
over $A,B$
of the   group  
  $B+A$ is determined since $B \meet A'= C$, and the ordering is also determined as follows:
   $$b+a_1 \leq  b'+a_2 \iff b-b'  < a_2-a_1  \iff  
(\exists a' \in A')( b-b' \leq a' \leq a_2-a_1)$$     \qed

\medskip

If $H$ is a convex subgroup of $A$, let $q+H= q_H |A$ denote the type over $A$ of an element
$c$ with $\ct_c(A) = H\union A^{< 0}$.     If $b_1 \models q_H$, then $H$ remains a convex
subgroup in $A(b_1) = A + \Qq b_1$. 
 Let $b_2 \models q_H | (A(b_1))$, and let $q_H^2 | A = \tp((b_2,b_1) /A)$.
Proceed in this way to define $q_H^n$.    Thus $(c_1,\ldots,c_n) \models q_H^n$
iff $$H< c_n << c_{n-1} << \ldots << c_1 < (A^{\geq 0} \setminus H)$$
If $A \leq B$, $q_H(B)$ denotes the set of realizations
of $q_H$ in $B$.

\begin{lemma} \label{C6}  Let $H_1,\ldots,H_k$ be distinct convex subgroups of $A$. 
 Let $B$ be an extension
of $A$ generated by  tuples $b_i$ realizing $q_{H_i}^{m_i}$.  Then the  
isomorphism type of $B$ over $A$ is determined.  In particular, the coordinates of the elements
 $b_i$
 are $\Qq$-linearly
independent over $A$.  Thus if $B$ is finitely generated over $A$, then the cuts of only
finitely many convex subgroups of $A$ are realized in $B$.    
\end{lemma} 

{\em Proof.}     We have to show that $\union_{i=1}^k q_{H_i}^{m_i}(x_i)$ is a complete type.
Say $H_1 \subset \ldots \subset H_k$, and  use induction on   $k$. 
Let $A' = A(b_2,\ldots,b_k)$.  By induction, the $A$-isomorphism type of $A'$ is determined.
  It is easy to check that 
$H_1$ is a convex subgroup of $A'$, and   $q_{H_1}^{m_1}$ generates a complete type over $A'$, namely $q_{H_1} | A'$.  \qed

\begin{lemma} \label{C8} Let $A$ be i-complete, $A \subset B$, $B/A$ finitely generated.  Then there exist distinct
convex subgroups $H_1,\ldots,H_k$ of $A$, and $m_1,\ldots,m_k \in {{\mathbb N}}$,
such that $q _{H_i}^{m_i}$ is realized in $B$ by
some tuple $b_i$, $b=(b_1,\ldots,b_k)$, and $B$ is i-dominated by $A(b)$ over $A$.          
\end{lemma} 

{\em Proof.}   If $q _{H_1}^{m_1} \times \ldots \times q_{H_k}^{m_k}$ is realized in $B$ 
then by Lemma~\ref{C6} we have $\sum_{i=1}^{k} m_i \leq \rk_{\Qq} (B/A)$.  Thus
a maximal $b$ as in the statement of the lemma exists.  We claim that  $B$ is i-dominated by $A(b)$ over $A$.  
    
Let $P$ be a cut of $A$.   By Lemma~\ref{C5} it suffices to show that $P_{A(b)}$
is two-sided cofinal in $P_B$.

 Since $A$ is i-complete, possibly after inversion,  $P$ can be translated by an element
of $A$, to a cut of a convex group $H$; we may assume 
 $P =  \{x : (\exists h \in H) ( x<h  ) \}$
(Lemma~\ref{C7}).     Suppose for contradiction that
$a \in P_B$,  but all elements of $P_{A(b)}$  are below $a$ (or all above $a$).
Say $H=H_i$ (add a term with $m_i=0$ if necessary).  Then $(a,b_i)$
(respectively $(b_i,a)$) realizes $q _{H_i}^{m_i+1}$.  This contradicts the maximality
in the choice of $b$.  \qed

\begin{lemma} \label{C10}  Let $C \leq B$.
Let $H$ be a convex subgroup of $C$.

1)  $q_H $ has at most two i-free extensions to $B$.

2)  $q_H^m$ has at most $m+1$ i-free extensions to $B$.

3)  Let $H_1 < H_2 < \ldots < H_n$ be convex subgroups of $C$, $q=q_{H_1}^{m_1} \tensor \ldots \tensor q_{H_n}^{m_n}$.  
The  number of i-free extensions of $q$ to $B$ is at most $\Pi_{i=1}^n (1+m_i) \leq 2^{\sum_{i=1}^n m_i}$.
\end{lemma}  

{\em Proof.}  

(1) Let $a \models q_H$, with $C(a)$ i-free from $B$ over $C$.
By i-freeness, we have {\em either}  $r_1(a): \ H < a < q_H(B)$ {\em or} $r_2(a): \ q_H(B) < a < C^{\geq 0} \setminus H$. 
 Either of these two possibilities determines a complete type over $B$.
(Note that $q_H(B)$ is closed under multiplication by $\Qq^{>0}$.)

(2)  Let $(a_1,\ldots,a_m) \models q_H^m$.  By (1) for some $l$ we have $a_i \models 
r_1$ for $i \leq l$ and $a_i \models r_2$ for $i >l$.  It is easy to see that
$r_1(x_1) \union r_1(x_l)  \union x_l <<  b~ (\mbox{for~}b \in q_H(B) ) 
\union q_H^l(x_1,\ldots,x_l) \union r_2(x_{l+1}) \union \ldots 
\union r_2(x_m) \union q_H^{m-l}(x_{l+1},\ldots,x_m)$  generates a complete type.
This gives $m+1$ possibilities. 

(3)  For each $i$, there are $\leq m_i+1$ i-free extensions  $q_{H_i}^{m_i}$;
and it is easy to see that for any system $r_i$ of choices of such extensions,
$r_1(x_1) \union \ldots \union r_n(x_n)$ is a complete type.  \qed

\begin{lemma} \label{C9}   Let $A$ be generated over
$C$ by  realizations of $q_{H_i}^{m_i}$ for some finite set $H_1 \subset \ldots \subset H_n$ of 
distinct convex subgroups $C$.  Let $C \leq B$.
 If $A$ is i-free from $B$ over $C$, then   $A$ is
g-free from $B$ over $C$, via some sequence of generators, by definition.
 \end{lemma}

{\em Proof.}  
 This reduces to the case $n=1$ and $m_1=1$.  Let $a \in A$, $a \models q_H | C$.
 If $a < b$ for some $b \in B$ with $b \models q_{H}$, then by the i-freeness
 assumption, $a<b$ for all such $b$.  In this case, $\tp(a/B)$ is g-free over $C$
 (in the sense of Remark~\ref{uniquegen}(iii), where the valuation
is replaced by a norm).
 
 Otherwise, $a>b$ for all $b \in B$  with $b \models q_{H}$.  Then $-a$ is 
  g-free over $C$, and we use $-a$ as the generator.  
   \qed

\medskip


 {\em Proof of Theorem~\ref{vik}} 
 
(i) implies (ii)  by definition. 
 
(ii) implies (iii):    If $\tp(A/B)$ does not divide over $C$, then the same is true
of $\tp(k(A)/k(C))$ and $\tp(\G(A)/\G(C))$.  Hence we are reduced to showing
the implication inside $k$ and $\G$.  For $k$ it  is well-known
that non-dividing agrees with linear disjointness (over an algebraically closed field.)
As for $\G$, suppose $\G(A)$ is not i-free from $\G(B)$ over $\G(C)$.  Then there
exist $a \in \G(A) \setminus \G(C), b_0,b_0' \in \G(B)$ such that $b_0 \leq a \leq b_0'$
and $b,b'$ lie in the same cut over $\G(C)$.  Extend $(b_0,b_0')$
to an   indiscernible sequence $((b_i,b_i'): i < \omega)$ over $C$ with $b_0 < b_0' < b_1 < b_1' \ldots$.
 Clearly the formulas $b_i<x<b_i'$ are pairwise inconsistent, showing that $\tp(\Gamma(A)/\Gamma(B))$
 divides over $\Gamma(C)$.  

(iv) implies (v)  is Corollary~\ref{invariance}.
  
(v) implies (i) by definition of forking (see Remark~\ref{forking/inv}). 

 It remains to prove that (iii) implies (iv).   Observe that 
if $C \leq A' \leq A$ and $\tp(B/A') \vdash \tp(B/A)$, and if $A'$ is sequentially
independent from $B$ via some sequence of generators $a'$, then  $A$
is sequentially
independent from $B$ via any (unary) sequence of generators whatever extending $a'$.  
Moreover any extension of $\tp(A'/C)$ to $B$ extends uniquely to an extension of 
$\tp(A/C)$ to $B$, so the number of extensions with property (iii) is the same for $A'$ and for $A$.
  Hence in this case the theorem is reduced from $A$ to $A'$.

 Since $\tp(B/C, \k(A),\G(A)) \vdash \tp(B/A)$,  
we may assume $A = \k(A) \union \G(A)$, $B = \k(B) \union \G(B)$.
As $\k,\G$ are orthogonal, we can consider them separately.
Clearly $\k(A)$ is free from $k(B)$ over $k(C)$, and we are reduced to
$\G$.   Similarly, by Lemma~\ref{C8}, we may assume $A$ is generated over $C$
by realizations of $q_{H_i}^{m_i}$ for some distinct convex subgroups $C_i$ of $C$.
 This is given by Lemma~\ref{C9}.

The final assertion follows from Lemma~\ref{C10}.  \qed

\begin{remark}  \rm Forking is not symmetric even over a vi-complete $C$;
this is because i-freeness is not in general symmetric.  However one can find, given $\tp(A/C)$ and $\tp(B/C)$,
a 2-type that is non-forking in both directions.  \end{remark}

 For an $\Aut(\U/C)$-invariant convex subgroup $H$ we define an invariant type $q_H$ as follows:  if $B \leq G$, 
 $c \models q_H | B$ if $H \meet B < c < B^{\geq 0} \setminus H$.   Given a small subset
 $A$ of $\G$, let 
 $$H^-=H^-(A) = \{x \in \G: (\forall a \in A)(\forall n \in \Nn)(n |x| \leq a) \}$$
 $$H^+(A) =  \{x \in \G: (\exists a \in A)(\exists n \in \Nn)( |x| \leq n a) \}$$
  Also define $^-q_A=q_{H^-}$,   $^+q_A = q_{H^+}$.
   
If $A=\emptyset$ we set $^-q_A = q_\Gamma$, while $^+q_A$ is not defined.
 
 In this notation, we have $^-q_A = ^-q_{A'}$ if $A,A'$ are left-cofinal in each other,
 and similarly for $+$.  Let ${\mathcal Q}  = \{^-q_A,^+q_A: A \hbox{ small.} \}$.  
 
 Note that if $B$ is a small subgroup, a convex subgroup $H$ of $B$ corresponds
 to two invariant types, 
 $^+q_H$ and $^-q_J$ where $J = \{b \in B: b > H \}$.  Both   restrict to $q_H | B$.  
 When $B \leq B'$, and $c \models q_H |B$, we have:  $B(c)$ is i-free from $B'$ over
 $B$ iff $c$ realizes one of these types over $B'$.  Let $r$ be one of the
 types $^+q_H$ and $^-q_J$.  
 The notion `$c \models r | B'$ '  
 is thus finer than i-freeness, and correspondingly domination is 
easier than i-domination.  Indeed it is easy to see that $r^ 2$ is domination-equivalent to $r$.

 \begin{corollary}  $\overline{\Inv}(\G)$ is the free Abelian semigroup generated by 
 ${\mathcal Q}$, subject to the identity:  $x^2=x$. 
 \end{corollary} 

{\em Proof.} Let $q$ be an invariant type, with base $A$; we may take $A$ to be i-complete.
Note that if $A \leq B$ then $q|B$ is an i-free extension of $q|A$.
By Lemma~\ref{C8}, $q|B$ follows from $q|A$ together with $q' |B$, where
$q'$ is a certain type of the form $\tensor_i q _{H_i}^{m_i}$.  Thus we are reduced
to considering the case:  $q=q_H$, $H$ a convex subgroup of $A$.  In this 
case by Lemma~\ref{C10}, $q$ has at most two invariant extensions to $\U$, namely
  $^-q_H$ and $^+q_{A^{\geq 0} \setminus H}$.   Thus     $\overline{\Inv}(\G)$ 
  is generated by the said types.   The freeness is easy and left to the reader.
  \qed

\medskip

This should be combined with Corollary~\ref{inv-dom}.

%
%

\chapter{A maximum modulus principle}\label{maximummodulus}

In this chapter, we look at the relationship between sequential independence
and a form of independence, called {\em  modulus independence}, which is 
clearly symmetric.
If $A$ and $B$ are fields, and $C\cap K\leq A$ with $\Gamma(C)=\Gamma(A)$, we show that
sequential independence over $C$ implies modulus independence. This can be understood
as a maximum modulus principle: if a polynomial function 
takes a certain norm  on a realisation of a type, then it takes at least that norm 
on a generic point of the type (over parameters defining the 
function). Using the results of Chapter~\ref{orthogonal}, we conclude by 
showing that for fields, if there is some orthogonality to the value group, 
then sequential independence and modulus independence are equivalent. 

\begin{definition} \index{positive type, $\tp^+$} \rm 
Suppose $A\subset K$. Recall  that we  write $x_a$ for a tuple of variables 
corresponding to the elements  $a$ of $A$. Define
$$\tp^{<}(A/C):=\{|f(x_a)| < \gamma:
        f(X) \in {\mathbb Z}[X] \mbox{ and } \tp(A/C) \vdash 
                       |f(x_a)| < \gamma\}~,
$$
$$\tp^{\le}(A/C):=\{|f(x_a)| \le \gamma:
        f(X) \in {\mathbb Z}[X] \mbox{ and } \tp(A/C) \vdash 
                       |f(x_a)| \le \gamma\}~,
$$
and 
$$
   \tp^+(A/C):=\tp^{<}(A/C) \cup \tp^{\le}(A/C)~.
$$
\end{definition}

Because of the assumption that $C\cap K\leq A$, nothing is changed if we allow the polynomials $f$ in $\tp^+(A/C)$ 
to have coefficients in $C \meet K$ rather than ${\mathbb Z}$.

In the above, $\gamma$ ranges over $\Gamma$.  We have however:
 
\begin{lemma} \label{tp+}  Let $P(X)$ be a consistent partial type over $C$ in field variables, and
let $f \in (C \meet K) [X]$ be a polynomial.   Let $\g \in \G$.  Then: 

$$P \vdash |f(X)| \leq \gamma \hbox{ iff for some } \gamma \geq \gamma_0 \in \G(C), \  P \vdash |f(X)| \leq \gamma_0$$
$$P \vdash |f(X)| < \gamma  \hbox{ iff either } \gamma \in \G(C) \hbox{ or for some } 
\gamma > \gamma_0 \in \G(C), \  P \vdash |f(X)| \leq \gamma_0.$$    \end{lemma}

{\em Proof.}   If $P \vdash |f(X)| \leq \gamma$ then for some $\psi \in P$
we have $\U \models  (\forall x)(\psi(x) \to |f(X)| \leq \gamma )$.  Let
$I = \{\gamma' \in \Gamma:  \U \models (\forall x)(\psi(x) \to |f(X)| \leq \gamma' )$.
Then $I$ is $C$-definable, nonempty and closed upwards.  We cannot
have $I = \Gamma$, since if $c \models P$ and $\gamma<  |f(c)| $ then $\gamma \notin I$.
so  $I = [\gamma_0,\infty)$ or $I = (\gamma_0, \infty)$, with $\gamma_0 \in \G(C)$.
It suffices to show that $P \models |f(X)| \leq \gamma_0$.  Otherwise, 
let $a \models P, |f(a)| > \gamma_0$.   Choose $\gamma'$ with $\gamma_0 < \gamma' < |f(a)|$.  
Then $\gamma' \notin I$ according to the definition of $I$, but $\gamma' \in (\gamma_0,
\infty)$, a contradiction.  \qed

\medskip

Observe that $\tp^<$ determines $\tp^{\le}$. For suppose $\tp^<(a/C)=\tp^<(a'/C)$,
and $\tp(a/C)\proves |f(x)| \le \gamma$. Then for all $\delta>\gamma$, 
$\tp^<(a/C)\proves |f(x)| < \delta$, so $\tp^<(a'/C)\proves |f(x)| < \delta$,
and hence $\tp(a'/C)\proves |f(x)| \le \gamma$. 

Nonetheless, we use both strict and weak inequalities; this permits using
the above lemma to see that the type of   $\tp^+(A/C)$ consisting
of inequalities with parameters $\gamma \in \G(C)$ implies the full type.  

\begin{definition}\label{m-ind} \index{independence!modulus, $\dnf^m$} \rm
If $A$ and $B$ are subsets of $K$ such that
$$
     \tp(A/C) \cup \tp(B/C) \vdash \tp^+(AB/C)
$$
we will say that $A$ and $B$ are {\em modulus-independent} \index{modulus independent} over $C$ and
write $A\dnf^m_C B$.
\end{definition} 

The condition $A\dnf^m_C B$ says that if $h$ is any elementary map over $C$, then
$Ah(B)$ satisfies $\tp^+(AB/C)$.
Observe that $\dnf^m$ is symmetric, and is independent of any choice of
generating sequence for $A$ over $C$, so does not 
coincide with $\dnf^g$ (by Examples~\ref{exsym} and \ref{exfun}). 
The main theorem of the chapter shows that, under certain hypotheses, sequential independence implies modulus independence. This gives a consequence of $\dnf^g$-independence which is symmetric 
between left and right.

\begin{theorem}\label{maxmod}
Let $C=\acl((C\cap K) \cup H)$ where $H\subset \S=\bigcup_{m>0} S_m$, 
and let $A,B$ be valued fields, with $C \subseteq \dcl(A)\cap\dcl(B)$ and 
 $\Gamma(C)=\Gamma(A)$. Suppose $A\dnf^g_C B$, via a sequence $a=(a_1,\ldots,a_n)$
of field elements. Then $A\dnf^m_C B$.
\end{theorem}

The proof of Theorem~\ref{maxmod} requires a sequence of lemmas; some will also
be useful otherwise.   The case where $\tp(A/C) \perp \Gamma$ requires only a few
of these; see the proof of Theorem~\ref{maxmodapp}.

\begin{lemma} \label{assumingmax}
Let $A,B,C$ be valued fields, with $C \leq A\cap B$, and $C$ 
algebraically closed and maximally complete. Suppose $\Gamma(C)=\Gamma(A)$ and 
$k(A)$ and $k(B)$ are linearly disjoint over $k(C)$. 
Then, 

(i)  let $g \in B[X]$ be a polynomial over $B$,
  $a \in A^n$, and $\g \in \G(B)$, and suppose $|g(a)| \leq \gamma$ (respectively $|g(a)| < \gamma$).
Then $\tp(a/C) \vdash |g(x)| \leq \gamma$ (respectively $|g(x)| < \gamma$).

(ii) $A\dnf^m_C B$.   

\end{lemma}

{\em Proof.}  Note that by Lemma~\ref{Cmax}, $\Gamma(B) = \Gamma(AB)$.
 We argue just with $<$ as the argument is identical .

(i)  This is Lemma~\ref{Cmax}(iii).  

(ii) 
Let $f(X,Y) \in C[X,Y]$. We need 
to show that if $a,b$ are finite tuples from $A, B$ respectively,
and $\tp(AB/C) \vdash |f(x_{a},y_{b})| < \gamma$, then
$$\tp(A/C)\cup \tp(B/C) \vdash 
|f(x_{a},y_{b})| < \gamma,
$$
and similarly with $\le$ replacing $<$.  By Lemma~\ref{tp+}, we may assume 
$\gamma\in \dcl(C)$.   Letting $g(x)=f(x,e)$, we need to prove:
$$\tp(a/C) \vdash |g(x)| \leq \gamma$$
This follows from (i) (where the same is shown more generally for $\gamma \in \G(B)$). \qed

\begin{corollary} \label{Cmax+div}  Let $C=\acl(C) \leq \U$, and let $B$ be a valued
field with  $C \leq \dcl(B)$.  Let $g \in B[X]$.  Assume $p$ is a stably dominated 
$\Aut(\U/C)$-invariant type in the field sort,  $a \models p|B$, 
and $\gamma = |g(a)|$.  Then $p|C \vdash |g(x)| \leq \gamma$.   
\end{corollary}

{\em Proof.}   Note $\gamma \in \G(B)$, since $p$ is orthogonal to $\G$.  

Let $M$ be any algebraically closed, maximally complete valued field 
with $C \subseteq \dcl(M)$, and such that $\St_C(M) \dnf _C \St_C(B)$.
  Let $a \models p| M(b)$, where $b$ enumerates $B$.  By the stable domination of $p$, 
the hypotheses of  Lemma~\ref{assumingmax}
apply to $M(a),M(b), M$.  Hence $p| M \vdash  |g(x)| \leq \gamma$.    
Hence  $p|Ce \vdash |g(x)| \leq \gamma$ for some finite tuple $e$ from $M$. 
Let $I$ be a large index set, and let $(e_i: i \in I)$ be an indiscernible sequence
over $B$ with $e_1=e$ and $(\St_C(e_i): i\in I)$ a Morley sequence over $B$.
Since $\St_C(e) \dnf_C \St_C(B)$, this is also a Morley sequence over $C$.
Let $a' \models p |C$.  By Lemma~\ref{morley-div}, $\St_C(a') \dnf_C \St_C(e_i)$ 
for some $i$; so by stable domination, $a' \models p | Ce_i$.  Since
$\tp(e/B) = \tp(e_i/B)$, we have $p|Ce_i  \vdash |g(x)| \leq \gamma$.  
So $|g(a')| \leq \gamma$.  Since $a' \models p |C$ was arbitrary, 
$p|C  \vdash |g(x)| \leq \gamma$.  \qed

\begin{lemma}\label{toresolve}
Let $A,B,C,C'$ be fields with $C \leq A\cap B \cap C'$. Put 
$A':=\acl_K(AC')$ and $B':=\acl_K(BC')$. Assume 

(i) $\tp(A/C) \vdash \tp(A/C')$, and

(ii) $A'\dnf^m_{C'} B'$.

\noindent
Then $A\dnf^m_C B$.
\end{lemma}

{\em Proof.} By (ii), $\tp(A'/C') \cup \tp(B'/C') \vdash \tp^+(AB/C')$, so
$\tp(A/C') \cup \tp(B/C') \vdash \tp^+(AB/C')$. By (i), it follows that
$\tp(A/C) \cup \tp(B/C') \vdash \tp^+(AB/C')$, so
$\tp(A/C) \cup \tp(B/C') \vdash \tp^+(AB/C)$. 
It follows that $\tp(A/C) \cup \tp(B/C) \vdash \tp^+(AB/C)$.
For suppose $a \in A$, $b \in B$, $\gamma \in \Gamma(C)$,
and $|f(a,b)| < \gamma$ (where $f$ is over ${\mathbb Z}$). Let $\sigma$ be an automorphism 
(or elementary map) on $K$ over $C$. We must show that 
$|f(a,\sigma(b))| < \gamma$ (and similarly with $<$ replaced by $\le$). But clearly
$\tp(\sigma(A)/C) \cup \tp(\sigma(B)/\sigma(C')) \vdash \tp^+(\sigma(AB)/C)$,
and $|f(\sigma(a),\sigma(b))|\leq \gamma$. Hence, as
$a \equiv_C \sigma(a)$, we have $|f(a,\sigma(b))| < \gamma$. 
\qed

\begin{lemma}\label{maxmodorth}
Let $A,B,C$ be valued fields, with $C$ algebraically closed, 
$C\leq A,B$. Suppose $\tp(A/C) \perp \Gamma$ and $A\dnf^g_C B$. Then
$A\dnf^m_C B$.
\end{lemma}

{\em Proof.} Let $C'$ be a maximally complete immediate extension of $C$, chosen so 
that $A\dnf^g_C BC'$. Then by Proposition~\ref{maxim}, 
$\tp(A/C)\vdash \tp(A/C')$, and $\Gamma(AC')=\Gamma(C')$.
Put $A':=\acl(AC')$ and $B':=\acl(BC')$. Then $A'\dnf^g_{C'} B'$, 
so by Lemmas~\ref{gtostab} and  \ref{assumingmax}, 
$\tp(A'/C') \cup \tp(B'/C')\vdash \tp^+(A'B'/C')$. The lemma now follows from
Lemma~\ref{toresolve}. \qed

\begin{lemma}\label{interalgebraic}
Suppose $C=\acl(C)\cap K$, $a_1\in K$, and $C\subseteq B \subset K$. Assume 
also $\Gamma(C)=\Gamma(Ca_1)$ and $a_1\dnf^g_C B$. Let $a_2\in K$ be 
interalgebraic with $a_1$ over $C$.
Then $a_2\dnf^g_C B$.
\end{lemma}

{\em Proof.} Since $\Gamma(C)=\Gamma(Ca_1)$, $\tp(a_1/C)$ is not the generic 
type of an open ball.
By Lemma~\ref{orthogwelldefined}, the result is immediate if $a_1$ is generic in a closed
 ball over $C$, so we may suppose that
$a_1$ is generic in the intersection of a sequence of $C$-definable balls $(U_i:i \in I)$ with
 no least element, and that the $U_i$ are closed. Then $\tp(a_1/B)$ is
determined by a collection of formulas
$$\{x\in U_i:i \in I\} \cup \{x\not\in V_j:j\in J\},$$
where the $V_j$ are $B$-definable balls (possibly of radius zero)
all lying in $\bigcap(U_i:i\in I)$. Suppose that $a_2\df^g_C B$. 
Then there is a $B$-definable ball $V$ containing no $C$-definable balls, such that
$a_2\in V$. Let $\phi(x,y)$ be a formula over $C$ satisfied by $a_1a_2$ and 
implying that $x$ and $y$ are interalgebraic over $C$. Then
$$\{x\in U_i:i\in I\}\cup \{x\not\in V_j:j\in J\}\vdash 
\exists y(\phi(x,y)\wedge y\in V).$$
Hence, by compactness, there is $i_0\in I$ such that
$$\{x\in U_{i_0}\}\cup \{x\not\in V_j:j\in J\}\vdash 
\exists y(\phi(x,y)\wedge y\in V).\eqno{(*)}$$
Suppose that $a_1'$ is generic over $B$ in $U_{i_0}$. Then by $(*)$, there is 
$a_2'$ interalgebraic over $C$ with
$a_1'$ with $a_2'\in V$. Then $a_2'\df^g_C B$. Since 
$\tp(a_1'/C)\perp \Gamma$, this contradicts
Lemma~\ref{orthogwelldefined}. \qed

\medskip

{\em Proof of Theorem~\ref{maxmod}.}

Put $H=\{s_i:i<\lambda\}$. We may suppose $A=\acl_K(A)$.
For each $i$, let $A_i:=\acl_K(Ca_1\ldots a_{i-1})$, and let $U_i$ 
be the intersection (in $K$) of the $A_i$-definable balls containing $a_i$.
Let $\delta$ be the number of $i$ such that $U_i$ is not a 
single ball, that is, such that there is no smallest $A_i$-definable ball 
containing $a_i$. We argue by induction on $\delta$, over all 
possible $A,B,C$ satisfying the hypotheses of the theorem. The strategy is to reduce $\delta$ by replacing
a chain of balls with no least element by a closed ball chosen generically inside it.

First, we reduce to the case when $H=\emptyset$, that is, $C=\acl(C\cap K)$. So for fixed
$\delta$, we assume the result holds when $H=\emptyset$. Let $A^*\equiv_C A$ and
$B^*\equiv_C B$. We must show $\tp^+(AB/C)=\tp^+(A^*B^*/C)$. Let $(e_i:i<\lambda)$ be a 
generic closed resolution of
$(s_i:i\in \lambda)$ over $CABA^*B^*$ (so for each $i$, $e_i$ is a generic resolution of
 $s_i$ over
$CABA^*B^*\cup \{e_j:j<i\}$). Put $E:=\bigcup(e_i:i<\lambda)$ 
(treating each $e_i$ as a subset
of $K$). Then $e\dnf^g_C ABA^*B^*$ 
for each finite $e\subset E$, and $\tp(e/C)\perp\Gamma$ (by Remark~\ref{resexist}(iii)), so 
$ABA^*B^*\dnf^g_C e$, and
hence
$ABA^*B^* \dnf^g_C E$.
Put $C':=\acl(CE)$, $A':=\acl_K(AE)$, and 
$B':=\acl_K(BE)$. Now as $H\subset \dcl(B)$, $e\dnf^g_B A$ ($e$ as above), and  
$\tp(e/B) \perp \Gamma$. Hence, by 
Proposition~\ref{perpgamma}, $A\dnf^g_B Be$, so
$A\dnf^g_B B'$, both via $a$. Now as $A\dnf^g_C B$, we have $A\dnf^g_C B'$, so 
$A\dnf^g_{C'} B'$, and so $A'\dnf^g_{C'} B'$, all via $a$. Also, as the 
resolution is generic,
$\Gamma(A')=\Gamma(A)=\Gamma(C)$. 
We claim that the invariant $\delta$ for $(A',B',C')$ 
is no greater than for $(A,B,C)$. For suppose it is greater. Then for some $i$,
there is a chain of $\acl(C'a_1\ldots a_{i-1})$-definable closed balls which contain
$a_i$ but are not $\acl(Ca_1\ldots a_{i-1})$-definable. Let $s$ be (a code for) one such
closed ball. Then, as $\Gamma(A')=\Gamma(C)$,
$\gamma:=\rad(s)\in \Gamma(C)$, and for any generic closed resolution $E'$ of $H$ over $CA$, 
$s=B_{\leq \gamma}(a_i) \in \acl(Ca_1\ldots a_{i-1}E')$. Hence, as $H \subset C$,  
$s\in \acl(Ca_1\ldots a_{i-1})$, proving the claim. Since the elements of $H$ are
 definable over $E$,
$C'=\acl(C'\cap K)$, and it follows (by the case when $H=\emptyset$) that 
$\tp(A'/C') \cup \tp(B'/C') \vdash \tp^+(A'B'/C')$. 
Hence, $\tp(A/C') \cup \tp(B/C') \vdash \tp^+(AB/C')$. 
Also, as $AA^*\dnf^g_C E$ and $BB^*\dnf^g_C E$, it 
follows by Corollary~\ref{leftright} that $\tp(A/C')=\tp(A^*/C')$ and 
$\tp(B/C')=\tp(B^*/C')$. Hence, $A^*B^*$ satisfies $\tp^+(AB/C')$ and in particular 
$\tp^+(AB/C)$, as required.

Thus, we assume $C=\acl(C \cap K)$.  In particular,
each $A_i$ is resolved.
We may suppose there is  an immediate extension $C^*$ of $C$ such that
$\tp(A/C)$ does not imply a complete type over $C^*$. For otherwise, by
Proposition~\ref{maxim}(ii)$\Rightarrow$(i), $\tp(a/C)\perp \Gamma$, 
so Lemma~\ref{maxmodorth} is applicable. Choose $C^*$ so that $a\dnf^g_C C^*$.
We may now find
$C\leq C' \leq C_1 \leq C^*$ so that $\tp(A/C)$ implies a complete type 
over $C'$, but not over $C_1$, and so that
$\trdeg(C_1/C')=1$. 

Next, we reduce to the case when  $C=C'$. 
To justify this, put $A':=\acl_K(AC')$, 
and $B':=\acl_K(BC')$. Then $\Gamma(C')=\Gamma(A')$. Indeed, if not, then for some $i$,
$a_i$ is chosen generically over $\acl(C'a_1\ldots a_{i-1})$ in a $(C'a_1\ldots a_{i-1})-\infty$-definable 1-torsor $U$ 
which is not closed and contains a proper subtorsor $V$ also definable over $\acl(C'a_1\ldots a_{i-1})$.
As $a_i\dnf_{Ca_1\ldots a_{i-1}} C'$, $U$ is $\infty$-definable over $\acl(Ca_1\ldots a_{i-1})$.
Since $\Gamma(C)=\Gamma(Ca)$ and $a_i$ is generic in $U$ over $Ca_1\ldots a_{i-1}$, there is no $\acl(Ca_1\ldots a_{i-1})$-definable proper subtorsor of $U$. Hence, if $a_i'\in V$, then $a_i'\equiv_{Ca_1\ldots a_{i-1}} a_i$ but
$a_i'\not\equiv_{C'a_1\ldots a_{i-1}} a_i$. This contradicts the assumption that
 $\tp(a/C)\vdash \tp(a/C')$. Thus, we can apply
 the argument given below, with $(C',A',B')$ in place of $(C,A,B)$ to obtain:
$A'\dnf^m_{C'} B'$. Then
Lemma~\ref{toresolve} yields $A\dnf^m_C B$.

It follows that we may assume $C=C'$, so $\trdeg(C_1/C)=1$. We have $a\dnf^g_C C_1$. Pick
$c_1 \in C_1 \setminus C$, and let $V$ be the intersection of all $C$-definable balls 
containing  $c_1$. As $C^*$ is an immediate extension of $C$, $k(C_1)=k(C)$ and so $V$ is not closed. Also$\Gamma(C_1)=\Gamma(C)$,
so $V$ contains no $C$-definable proper subtorsors.

We now let $h$ be a $C$-elementary map, and must show $Ah(B)$ satisfies
$\tp^+(AB/C)$. Let $i$ be largest such that $\tp(A_i/C) \vdash \tp(A_i/C_1)$. Then $\tp(A_{i+1}/C)\not\vdash \tp(A_{i+1}/C_1)$ so
 by Proposition~\ref{noembed}(i) there is an embedding of $C_1$ into $A_{i+1}$ over $C$. Hence, there is
$a_i'\in A_{i+1}$ with $a_i'\in V$. Furthermore, as $\tp(A_i/C) \vdash \tp(A_i/C_1)$,
$\tp(c_1/C)\vdash \tp(c_1/A_i)$. Hence, 
$x\in V$ determines a complete type over $A_i$; for otherwise, different choices of 
$c_1$ in $V$ would have different types over $A_i$ whilst having the same type  over $C$. In particular, $a_i'\not\in A_i$, so $a_i$ and 
$a_i'$ are inter-algebraic over $A_i$.
Since $C$ is resolved and $C_1$ is an immediate extension of $C$, there is no smallest 
$C$-definable ball containing $V$; for if $V$ was an open ball, then $c_1$ would be 
generic over $C$ in an open ball, so $\Gamma(C_1)$ would properly contain $\Gamma(C)$. Hence, using Lemma~\ref{eqorthog} twice, 
 $k(A_i)=k(A_ia_i')=k(A_ia_i)$,
so $a_i$ is not generic in a closed ball over $A_i$. As $\Gamma(C)=\Gamma(A)$, $a_i$ 
is also not generic in an open ball over $A_i$, so $\tp(a_i/A_i)$ is generic in a chain 
of balls with no least element. In particular, $\delta>0$.
Since we assume the result holds for smaller values of $\delta$, our goal is to 
decrease $i$ by converting the choice of $a_i$ over $A_i$ to the 
choice of a generic in a closed ball.

Choose $\alpha \in \Gamma$ generic over $CABh(B)$ below the cut in $\Gamma(\U)$ made 
by  $\rad(V)$. Let $u:=B_{\leq \alpha}(a_i')$. Also, let $s$ be a code for the lattice in 
$S_2$ interdefinable over $\emptyset$ with $u$ (as in Remark~\ref{torsormodule}). 
Put $\hat{C}:=\acl(C\alpha s)=\acl(Cs)$ (as $\alpha \in \dcl(s)$). By Proposition~\ref{existence}(ii), there is an elementary map $f$ over $C_2$ such that
$A\dnf^g_{\hat{C}} f(B)$ via $a$. Put
$a':=(a_1,\ldots,a_{i-1},a_i',a_{i+1},\ldots,a_n)$.
\medskip

{\em Claim.} (i) 
$A\dnf^g_{C} f(B)$ via $a'$.

(ii) $A\dnf^g_{\hat{C}} f(B)$ via $a'$.

(iii) $A\dnf^g_C B$ via $a'$.

\medskip

{\em Proof of Claim.} (i) First, for $j<i$, note that $\tp(a_j/A_j)$ 
implies a complete type over $Cc_1$ and hence over $Ca_i'$ (as $a_i'\equiv_{A_j} c_1$). 
Hence, as $\Gamma(C) = \Gamma(A)$, $\tp(a_j/A_j)$ implies a complete type over 
$A_ja_i' \alpha$, and hence over $\acl(A_j s\alpha)$. Thus, 
$(a_1,\ldots,a_{i-1})\dnf^g_C \hat{C}$, so as
$A\dnf^g_{\hat{C}} f(B)$, we have $(a_1,\ldots,a_{i-1}) \dnf^g_C f(B)$. 

Next, we must check that $a_i'$ is generic in $V$ over $A_if(B)$. First, as $\alpha$ is 
generic below $\rad(V)$ over $AB$ and is fixed by $f$, it is also generic below $\rad(V)$ over $Af(B)$, so  
there is no $A_if(B)$-definable sub-ball of $V$ 
containing $u$. Likewise, $u\not\in \acl(A_if(B))$. Thus, it suffices to check 
  $$a_i'\mbox{~ is generic in~} u \mbox{~ over~} A_if(B)s.\eqno{(*)}$$ 
So suppose not. Then there is an open sub-ball $u'$ of $u$ of radius $\alpha=\rad(u)$, 
algebraic over $A_if(B)s$ and containing $a_i'$. 
Now $a_i'\in \acl(A_ia_i)$, so $u' =B_{<\alpha}(a_i')\in \acl(A_ia_is)$. Since 
$a_i \dnf^g_{\acl(A_is)} f(B)$, we have $a_i\dnf^g_{\acl(A_is)} u'$. The easy Lemma 2.5.3 of \cite{hhm}
now yields  that $u'\in \acl(A_is)$. 
Thus, if $b_1,b_2\in K$ are in the ball $s$ with $b_1\in u'$ and $b_2\not\in u'$, then $b_1\not\equiv_{A_is} b_2$. 
By 
\cite[Corollary 2.4.5(ii)]{hhm}, this contradicts 
the fact that $x\in V$ determines a complete type over $A_i$.

Finally, let $j>i$. We first show that $a_j$ is generic in $U_j$ over $A_j\alpha$. Suppose this is false. Observe 
 that $U_j$ is not closed, as otherwise there would be an element of the strongly minimal set $\red(U_j)$ which is algebraic over $A_j\alpha$ but not over $A_j$, which is impossible by the instability of $\Gamma$.
 There is an $A_j\alpha$-definable proper 
sub-ball of $U_j$, containing $a_j$; hence, by Proposition~\ref{newclassifyfun}, there is  a proper 
sub-ball $s_j$ of $U_j$ definable over $A_j$.
It follows by Lemma~\ref{orderlike} that $\Gamma(A_ja_j)\neq \Gamma(A_j)$, contradicting the fact that $\Gamma(C)=\Gamma(A)$.
 Hence $a_j\dnf^g_{A_j} \alpha$. Thus, as $s\in \dcl(A_j\alpha)$, $a_j\dnf^g_{A_j}\hat{C}$  so by the assumption on $f$, 
$a_j\dnf^g_{A_j} f(B)$, proving the claim.

(ii) Since $a\dnf^g_{\hat{C}} f(B)$, the only point to check is that 
$a_i'\dnf^g_{\hat{C}A_i} f(B)$. This was proved in (*) above.

(iii) Since $a\dnf^g_C B$, the only point to check is that $a_i'\dnf^g_{A_i} B$. 
Since $a_i\dnf^g_{A_i} B$, this follows from Lemma~\ref{interalgebraic}. 

\medskip

By parts (i) and (iii) of the claim,  $AB\equiv_C Af(B)$. Thus, if
 $A\dnf^m_C f(B)$ then also $A\dnf^m_C B$. Thus, we may replace 
$f(B)$ by $B$, that is, we may assume that
$A\dnf^g_C B$ via $a'$ and $A\dnf^g_{\hat{C}} B$ via $a'$.

Choose $a_{n+1},b_1,b_2$ generically in $u$ to ensure that  if
$\hat{A}:=\acl_K(Aa_{n+1})$, $\hat{B}:=\acl_K(Bb_1b_2)$, then $\hat{A}\dnf^g_{\hat{C}} \hat{B}$
via the sequence $(a_1,\ldots,a_{i-1},a_i',a_{i+1},\ldots,a_n,a_{n+1})$. Then $\hat{A}$ and $\hat{B}$
 are 
resolved. As $u$ is definable over
$a_i'a_{n+1}$ and over $b_1b_2$, we have $C_2 \subset \dcl(A_2) \cap \dcl(B_2)$. Also,
$\Gamma(C_2) = \Gamma(A_2)$: for $\Gamma(C)=\Gamma(A)$, $\Gamma(C)\neq \Gamma(C_2)$, and 
$\trdeg(A_2/A) \leq 1$ so $\rk_{{\mathbb Q}}(\Gamma(A_2)/\Gamma(A))=1$. 
Thus, the hypotheses of the theorem are satisfied with $C_2$ in place of $C$,
$(a_1,\ldots,a_{i-1},a_i',a_{i+1},\ldots,a_n, a_{n+1})$ in place of $a$, and $A_2,B_2$ 
in place of $A,B$. Furthermore, $\delta$ has been reduced: for $a_i',a_{n+1}$ were chosen 
generically in a closed ball (namely $u$), and for each $j\not\in \{i,n+1\}$, 
$a_j\dnf^g_{C_2A_i} B$. 
It follows that $\tp(A_2/C_2) \cup \tp(B_2/C_2) \vdash \tp^+(A_2B_2/C_2)$, so
$\tp(A/C_2) \cup \tp(B/C_2) \vdash \tp^+(AB/C_2)$.

It remains to replace $C_2$ by $C$ in this conclusion, so we must show
$Ah(B)$ satisfies $\tp^+(AB/C)$. Thus, we must show $B\equiv_{C_2} h(B)$.
Recall that $\alpha$ was chosen  generically below $\rad(V)$ over $CAB h(B)$). 
In particular, $B\equiv_{C\alpha} h(B)$.

Suppose first that $B$ (and hence also $h(B)$) has no element in $V$. 
Then $x\in V$ determines a complete type over $B$, and so
$a_i'\dnf^g_C B$, $a_i'\dnf^g_C h(B)$. By Corollary~\ref{leftright}, 
$B\equiv_{\acl_K(Ca_i')}h(B)$. As 
$\alpha$ is generic below $CABh(B)$,
 it is generic below $\rad(V)$ over
$\acl_K(Ca_i') B$ and over $\acl_K(Ca_i')h(B)$, so
$B\equiv_{\alpha \acl_K(Ca_i')} h(B)$. It follows by Lemma~\ref{3.4.12} that
$B\equiv_{C_2} h(B)$ in this case.

Finally, suppose  that there is $b\in B \cap V$. Then, as $a_i',b,h(b) \in V$,  
$\alpha>|a_i'-b|$ (so $b\in u$), and also $\alpha>|a_i'-h(b)|$. Choose a field element $d$
generically in the ball $u$ over $CBh(B)$. Then 
 $|d-b|=\alpha$ and $|d-h(b)|=\alpha$. Also, by the choice of $\alpha$,
$d\dnf^g_C B$ and $d\dnf^g_C h(B)$. Put $D:=\acl_K(Cd)$. Then by
 Corollary~\ref{leftright}, $B\equiv_D h(B)$. Since $\alpha=|d-b|=|d-h(b)|$,
$B\equiv_{D\alpha} h(B)$. Hence, by Lemma~\ref{3.4.12},
$B\equiv_{\acl(D\alpha)} h(B)$. Since $u=B_{\leq \alpha}(d)$, $B\equiv_{C_2} h(B)$, as 
required.
\qed

\bigskip

Using results from this and the last chapter, we  summarise how $\dnf^g$
and $\dnf^m$ behave, given some orthogonality to $\Gamma$.
First we need a lemma.

\begin{lemma} \label{ringbasis}
Let $a=(a_1,\ldots,a_n)\in K^n$, and let $C\leq A$ be 
algebraically closed valued fields, with $A=\acl_K(Ca)$ and $\Gamma(C) = \Gamma(A)$. 
Within the ring $C[a]$ there are $e_1,\ldots,e_m\in R$ 
algebraically independent over $C$ such that
$\{\res(e_1),\ldots,\res(e_m)\}$ is a transcendence basis 
of $k(A)$ over $k(C)$.
\end{lemma}

{\em Proof.} Choose $d_1,\ldots,d_m \in A$ such that $\{\res(d_1),
\ldots,\res(d_m)\}$
is a transcendence basis of $k(A)$ over $k(C)$.  For 
each $i=1,\ldots,m$,
as $d_i$ is algebraic over $C(a)$, there is a polynomial 
$f_i$ over $C[a]$ 
with $f_i(d_i)=0$. As $\Gamma(C)=\Gamma(Ca)$, we may 
arrange (by multiplying 
$f_i$ by a scalar from $C$) that the maximum of the absolute 
values of the coefficients of $f_i$ is $1$. Thus, the reduction of 
$f_i$ modulo
${\cal M}$ is a non-zero polynomial over $k(C(a))$, so 
$\res(d_i)$
is algebraic over the residues of the coefficients of $f_i$. Let
$E$ be the set of all the coefficients of all the $f_i$. Then 
each $\res(d_i)$
is algebraic over $\{\res(x):x\in E\}$. Hence, there are 
$e_1,\ldots,e_m\in E$
such that $\{\res(e_1),\ldots,\res(e_m)\}$ is  a transcendence basis of
$k(A)$ over $k(C)$. Then $e_1,\ldots,e_m$ are 
algebraically independent over $C$. \qed 

\begin{theorem}\label{allsame}
Let $C\leq A,B$ be algebraically closed valued fields with $A$, $B$ extensions of
$C$ of finite transcendence degree, and at least one of $\tp(A/C)$, $\tp(B/C)$
orthogonal to $\Gamma$. Then the following are equivalent.

(i) $A\dnf^m_C B$.

(ii) $A\dnf^g_C B$ via some generating sequence.

(iii) $A\dnf^g_C B$.

(iv) $B\dnf^g_C A$ via some generating sequence.

(v)  $B\dnf^g_C A$.

(vi) $A\domind_C B$.

\end{theorem}

{\em Proof.} We assume $\tp(A/C)\perp\Gamma$. Then (ii) $\Rightarrow$ (i) is
Theorem~\ref{maxmod}, and we prove (i) $\Rightarrow$ (ii). If instead
$\tp(B/C)\perp\Gamma$, then the same argument proves (i) $\Leftrightarrow$ (iv).
By Proposition~\ref{perpgamma}, we have the equivalence of (ii), (iii), (iv), (v), (vi).

So assume $A\dnf^m_C B$. By Lemma~\ref{ringbasis}, there are
$e_1,\ldots,e_m \in A$ algebraically independent over $C$ such
that $\res(e_1),\ldots,\res(e_m)$ form a transcendence basis of 
$k(A)$ over $k(C)$. As
$A\dnf^m_C B$, also $C[e_1,\ldots,e_m]\dnf^m_C B$, and 
it follows that $\res(e_1),\ldots,\res(e_m)$ are linearly independent (so
algebraically independent) over 
$k(B)$. For suppose that $\Sigma_{i=1}^m \res(b_i)\res(e_i)=0$. Then
$|\Sigma_{i=1}^m b_ie_i|<1$, so $\tp(b_1,\ldots,b_m/C) \cup \tp(e_1,\ldots,e_m/C) \vdash
|\Sigma_{i=1}^m x_iy_i|<1$. Hence $\res(b_1)=\ldots =\res(b_m)=0$: indeed,
 suppose say $\res(b_1)\neq 0$. Let $e_1'$ be chosen generically in $R$ 
over all other parameters.
Then $e_1'e_2\ldots e_m\equiv_C e_1\ldots e_m$, so $|e_1'b_1+\Sigma_{i=2}^m e_ib_i|<1$,
but  $$|e_1'b_1+\Sigma_{i=2}^m e_ib_i|=
|(e_1'-e_1)b_1+\Sigma_{i=1}^m e_ib_i|=|(e_1'-e_1)b_1|=1,$$
a contradiction.

 Hence, $k(A)\dnf^g_C B$. As $C$ is resolved, it follows
that $\St_C(A) \dnf_C \St_C(B)$. Hence, by Proposition~\ref{perpgamma},
$A\dnf^g_C B$ via some (in fact, every) generating sequence.
\qed

\medskip

\begin{theorem} \label{maxmodapp}  Let $C \subseteq \U$, and let $V$
be an  affine variety  defined over $C \meet K$.  Let 
$p|C$ be  a 
stably dominated  type over $C$ of elements of $V$, with $\Aut(\U/C)$-invariant extension $p$.  
    Let $F$ be
a regular function on $V$, defined over a field $L$ with $C  \subseteq \dcl(L)$.
   Then $|F(x)|$ has  
a maximum $\gamma_{\max}^F \in \G(C)$ on $\{x \in \U: x \models p |C\}$.    

 Moreover for 
$a  \models p|C$, we have: 
$a \models p|L$ if and only if $| F(a)|  = \gamma_{\max}^F$
for all such $F$.
\end{theorem}

{\em Proof.}  Since $p|C$ has an $\Aut(\U/C)$-invariant extension, it implies
$p | \acl(C)$; so with no loss of generality we may assume $C=\acl(C)$.  

Choose an embedding of $V$ in $\Aa^n$.  Then  $F$ lifts to a  
polynomial on $\Aa^{n}$.
Thus we may assume $V = \Aa^n$.  Write $F= F(x,b)$.
If $a \models p|L$ then $\G(La)=\G(L)$.  So for some $\alpha \in \G(L)$,
if $a \models p|L$ then $|F(a,b)| = \alpha$.  
The hypotheses of Corollary \ref{Cmax+div} hold.  Hence
 $p(x) | C  \vdash |F(x,b)| \leq \alpha $.    This shows that $\alpha $ is the maximum
 of $|F|$ on the realizations of $p|C$. 
 
    The `moreover'  follows from  quantifier
elimination, since  the set of formulas $\{|F(x)| = \gamma_{\max}^F: F[X]\in L[X] \}$
determines a  complete type over $L$. \qed

\medskip

Here is a two-sided version of the same result.  
If $p(x), q(y)$ are types, define
$p\times q(x,y):=p(x)\cup q(y)$.  If $p$ is stably dominated,    let 
$p\otimes q$ be the complete type with realisations 
$$\{(a,b): (a,b)\models p\otimes q \mbox{~and~} a\domind_C b\}.$$  
 
\begin{theorem}\label{maxmodap}
Let $F(x,y)$ be a polynomial over the algebraically closed valued field $C$, and $p,q$ be 
 stably dominated  types 
in the field sort  over $C$.  Assume $p$ is stably dominated.
Then $|F(x,y)|$ has a maximum $\gamma_{\max}\in \Gamma(C)$ on $p\times q$. Also,
$$(a,b)\models p\otimes q \Rightarrow |F(a,b)|=\gamma_{\max}.$$
\end{theorem}

{\em Proof.} First, if $(a,b)\models p\otimes q$ then $\Gamma(Cab)=\Gamma(b)=\Gamma(C)$, by stable domination (twice). 
Thus, $|F(x,y)|$ takes constant value in $\Gamma(C)$. Now apply Theorem~\ref{maxmod}. \qed

\medskip

Here is a   group-theoretic consequence of the maximum modulus principle.
If $G$ is an affine algebraic group over $K$, we have the ring of regular functions
$K[G]$.  We can also view $G$ as a definable group in ACVF.  
 Recall the discussion of translation invariant definable types above \ref{groupdom}.

\begin{theorem} \label{genchar2}  \index{group scheme}
Let $G$ be an affine algebraic group 
over a field $C$, and $H$ a  $C$-definable subgroup of $G(K)$.  Let $p$ be  a $C$-definable global 
type of
elements of $H$.    Then  the 
 following are equivalent: 

(i) $p$ is translation invariant, and stably dominated.

(ii) For any $f\in C[G]$, $b\models p$, and $a \in H$ then  
$|f(a)|\leq |f(b)|$.
\end{theorem}

{\em Proof.} 
(i) $\Rightarrow$ (ii) Since $p|C \perp \Gamma$, $|f(x)|$ takes a fixed value $\gamma_f$ on $p|C$. If
$(a,b)\models p\otimes p$ then $ab\models p$ (cf. Section 3 of \cite{hrush2}). In particular, $|f(ab)|=\gamma_f$. Also, it follows 
from Theorem~\ref{maxmodap} that for $(a,b)\models p\times p$, $|f(ab)|\leq \gamma_f$. 
Now $H=\{ab: (a,b)\models p\times p\}$. 

(ii) $\Rightarrow$ (i)  By quantifier elimination in ACVF in the language $L_{\div}$
 (see Theorem~\ref{robin}),
 condition (2) characterises $p$ uniquely. Also, condition (2) is translation-invariant,
 so the type characterised is translation-invariant. 
 Clearly $\gamma_f\in \Gamma(C)$ for each $f$. It follows by quantifier elimination and Proposition~\ref{maxim} that if $C'\supseteq C$ is a 
maximally complete 
 algebraically closed valued field then $p|C'$ is stably dominated, so $p|C$ is stably dominated by Theorem~\ref{stab7}.  \qed

\medskip

The equivalence implies in particular the uniqueness of stably dominated, 
translation-invariant types;
this can also be seen abstractly. The conditions (i), (ii) hold for example if $G=\SL_n(K)$, $H=\SL_n(R)$, and $p$
is the unique type of $H$ whose elements have reduction satisfying the generic type of $\SL_n(k)$.
See  \cite{hrush2} for more detail.

%
%
%

\chapter{Canonical bases and independence given by modules}\label{Jindependence}

Consider the
following situation in a pure 
algebraically closed field $K$.
Suppose $C\le A\cap B$ with $C$ algebraically closed. Let $p$ be a type, and put
$$I(p) = \{ f\in K[X]: \mbox{ for all } a\models p\ (f(a)=0) \}.$$
Then $A\dnf_C B$ in the sense of stability theory if and only if 
$I(a/C) = I(a/B)$ for all tuples $a$ from $A$. We refine
this to take the valuation into account, and use it to define  another notion
of independence. We will show that some orthogonality to $\Gamma$ is enough
to prove its equivalence with the other forms of independence.

In this chapter, for ease of notation, if $s\in S_n \cup T_n$ we identify $s$ with the subset of $K^n$
which it codes.
\begin{definition}  \label{Jinv} \index{independence!J-independence, $\dnf^J$} \rm
Let $p$ be a partial field type over some set of parameters. Define
$$
   J(p):=\{f(X)\in K[X]: p(x)\vdash |f(x)|< 1\}.
$$
More generally, if $p=p(x,v)$ is a partial type with 
$x=(x_1,\ldots,x_{\ell})$ field variables and 
$v=(v_1,\ldots,v_m)$ a tuple of  variables with $v_i$ ranging through 
$S_{n_i} \cup T_{n_i}$, then
$$
    J(p)=\{f(X,Y) \in K[X,Y]:
        p(x,v) \vdash (\forall y_1\in v_1)\ldots (\forall y_m \in v_m)
             |f(x,y)| < 1\}.
$$
Here, $y_i$ is an $n_i$-tuple of field variables for each $i=1,\ldots m$. The notation $y_i \in v_i$ is a slight abuse -- it means
 that $y_i$ lies in a lattice
(or in an element of  some $\red(s)$) {\em coded} by an element of $S_{n_i} \cup T_{n_i}$.
We shall often write $J(a/B)$ for $J(\tp(a/B))$.

Next, for a tuple $a$ in $\U$, and $B,C$ lying in arbitrary sorts, define
$a\dnf^J_C B$ to hold if and only if $J(a/C)=J(a/BC)$.
Finally, if $A \subseteq \U$, then $A \dnf^J_C B$ holds if $a\dnf^J_C B$ 
for all finite tuples $a$ from $A$. \index{J-independence}
\end{definition}

If $J_n(p):=\{f\in J(p): \deg(f)<n\}$ (where $\deg$ refers to total degree), then $J_n(p)$ can be regarded as an $R$-submodule 
of $K^N$ for some $N$: identify each polynomial with the correspondinq sequence of coefficients. 
Note that $J_n(p)$ is the union of a collection of {\em definable} $R$-modules (one module for each formula in $p$). 
For if $p=\{\phi_i:i\in I)$, let $J_n(\phi);=\{f(X)\in K[X]: \forall x(\phi(x)\rightarrow |f(x)|<1\}$, a definable $R$-module. 
Then an easy compactness argument shows that
$J_n(p)=\bigcup(J_n(\phi_i):i\in I)$.

In the above, we think of $K$ as the field sort of $\U$, so $K[X]$ denotes
$(K \cap \U)[X]$.  Observe that we do not need to consider
in addition the polynomials satisfying a weak inequality. 
For suppose $J(a/C)=J(a'/C)$, and $\tp(a/C)\vdash |f(x)|\leq 1$.  
Then for any $\delta\in \Gamma(K)$ with $\delta>1$ there is $d\in K$ 
with $|d|=\delta$. Now $d^{-1}f(x) \in J(a/C)=J(a'/C)$, so for all 
$a''\equiv_C a'$, $|f(a'')|<\delta$. It follows by saturation of $K$ that for
all such $a''$, $|f(a'')|\leq 1$, that is, $\tp(a'/C)\proves |f(x)|\le 1$.

For an invariant type $p$, $J(p)$ is defined as for types over a fixed set; here however
we can write more simply:  
 $J(p) =  \{f(X)\in K[X]:  (|f(x)|< 1) \in p\}$.  Note that if $p$ is a definable type,
 then $J(p)$ is {\em definable}, in the sense that the intersection with the space
 $K[X]_d$ of polynomials of degree $\leq d$ is definable, for any integer $d$.

If $p$ is a stably dominated type of field elements over $C=\acl(C) \subseteq \U$, and 
$p'$ is the corresponding invariant type,
the maximum modulus principle (in the form of Theorem~\ref{maxmodapp})
states precisely that $J(p)=J(p')$.  
 
  When $p$ is stably dominated, we will see that knowldege of $J(p)$ determines $p$.  Thus
 the sequence of codes for the submodules $(K[x]_d \meet J(p))$ can 
 be viewed as the {\em canonical base}  of $p$.   This can be
seen as giving geometric meaning to the canonical base of a stably dominated type,
analogous to the field of definition of a prime ideal in ACF.

For $J$-independence we have transitivity on the right, and unlike with $\dnf^g$,
there is no dependence on a generating sequence. On the other hand
 transitivity
 on the left fails (as follows from Example \ref{14.12}, along with Proposition \ref{equivsingles});
 hence symmetry can fail.  
Existence  of $J$-independent extensions of $\tp(A/C)$ 
can also fail (Example \ref{14.13}), though, by Theorem~\ref{moregtoJ} below, it holds if $\tp(A/C)\perp \Gamma$. 
 We show below that,
 given some orthogonality to $\Gamma$, $\dnf^J$ coincides with $\dnf^g$, and so has these properties.
There are further results and examples at the end of the chapter.

\begin{lemma} \label{Jimpliesw}
Suppose $C\leq B$ are valued fields with $C$ algebraically closed, and
$a \in K^n$ with $a\dnf^J_C B$. Then
$a\dnf^m_C B$.
\end{lemma}

{\em Proof.} Suppose that $\tp(ab/C)\vdash |f(x_a,y_b)|<\gamma$ where $b$ is a  tuple in $B$, and $f$ is over ${\mathbb Z}$.
We may suppose that $\gamma \in \Gamma(C)$.
 Choose $d\in K$ with $|d|=\gamma$.
We need to show that $\tp(a/C) \vdash |f(x,b)| < \gamma$.  Now $\tp(a/B) \vdash |f(x,b)|
<\gamma$, so
$d^{-1}f(x,b) \in J(a/B)= J(a/C)$. Hence, if $a'\equiv_C a$, then
$|f(a',b)|<\gamma$.

Next, suppose $\tp(ab/C)\vdash |f(x,y)|\leq\gamma$. Let $\delta\in \Gamma$ with $\gamma<\delta$.
Then $\tp(ab/C) \vdash |f(x,y)|<\delta$. So, arguing as above, if $a'\equiv_C a$ then
$|f(a',b)|<\delta$. Thus, if $a'\equiv_C a$ then $|f(a',b)|\leq \gamma$. 
\qed

\bigskip
Notice that the converse of the above lemma is not obviously true. This 
is because, in the definition of $J$-independence, the coefficients of 
the polynomials are allowed to be anywhere in $K$, whereas modulus 
independence only refers to polynomials with parameters from the given set. 
In the proof of the following lemma, we see a way of overcoming  this 
issue. Then we show that   sequential independence implies 
$J$-independence, given some orthogonality to $\Gamma$.

\begin{lemma}\label{Jres}
Suppose $C=\acl(C)\subseteq A\cap B$ with 
$\Gamma(C)=\Gamma(A)$, $A$ is resolved, 
and $A\cap K \dnf^J_C B$. Then $A\dnf^J_C B$.
\end{lemma}

{\em Proof.} For ease of notation, we just do the following special case: 
$t=(t_1\ldots,t_m)$ is a tuple of elements of $\red(s)$, where $s\in S_n$,
and $f(U)\in J(t/B)$. We must show that $f(U)\in J(t/C)$.
As $A$ is resolved, there is $a_i\in K^n \cap t_i\cap A$ for each $i=1,\ldots m$. 
Put $a=(a_1,\ldots,a_m)$. Furthermore, 
there is in $A$ a basis $b=(b_1,\ldots,b_n)$ for $s$ in $K$ (so $b_i \in K^n\cap A^n$ for 
each $i$). Note that 
$t_i:=\{a_i+\Sigma_{j=1}^n y_{ij}b_j: y_{ij}\in R\}$ for each $i$.
Let $q:=\tp(ab/B)$. Then, as $f(U)\in J(t/B)$,  $q(x_a,z_b)$ implies
\begin{eqnarray}\nonumber
  (\forall y_{ij}\in \M: 1\leq i\leq m, 1\leq j \leq n) & \\ \nonumber
(|f(x_1 +\Sigma_{j=1}^n y_{1j}z_j, \ldots,x_m+\Sigma_{j=1}^n y_{mj}z_j)|<1).\\ \nonumber
\end{eqnarray}
In particular, for any $y_{ij} \in \M_K$, 
$$
  q(x_a,z_b)\vdash 
    |f(x_1+\Sigma_{j=1}^n y_{1j}z_j,\ldots,x_m+\Sigma_{j=1}^n y_{mj}z_j)|<1.
$$
Thus, 
$f(x_1+\Sigma_{j=1}^n y_{1j}z_j,\ldots, x_m\Sigma_{j=1}^n y_{mj}z_j) \in J(q)$.
By assumption, $J(q)=J(q|_C)$, so 
$f(x_1+\Sigma_{j=1}^n y_{1j}z_j,\ldots, x_m+\Sigma_{j=1}^n y_{mj}z_j) \in J(q|_C)$.
By compactness, this forces that $q|_C(x_a,z_b)$ implies the first displayed inequality 
above, hence $f\in J(t/C)$. \qed

\medskip

It is in the proof of the next theorem that Chapter~\ref{maximummodulus} is used: $\dnf^m$ provides the link 
from $\dnf^g$ to $\dnf^J$.

\begin{theorem}\label{11.2}
Assume $C=\acl((C\cap K)\cup H)$, where $H$ is a subset of $\S=\bigcup_{n>0} S_n$. 
Suppose also that
$C\subseteq A\cap B$ and $\Gamma(C) = \Gamma(A)$.
Assume $A\dnf^g_C B$ via a sequence of field elements. Then $A\dnf^J_C B$.
\end{theorem}

{\em Proof.} 
Since, by Corollary~\ref{twores}, we may resolve $A$ generically over 
$B$ without losing the assumption
$\Gamma(C)=\Gamma(A)$, and may also resolve $B$ generically,
we may assume $A=\acl(A\cap K)$ and $B=\acl(B\cap K)$.
Let $f(x,e)\in J(A/B)$. Since $A$ is resolved, we may by Lemma~\ref{Jres} 
assume $x$ ranges through the field sorts. Notice that the tuple of 
coefficients $e$ may be anywhere in $K$.
We must show $f(x,e)\in J(A/C)$. Since this is a property of $\tp(A/C)$, 
and our assumption on $f(x,e)$ depends on $\tp(A/B)$, by 
translating $A$ over $B$ we may assume
$A\dnf^g_B e$. Hence $A\dnf^g_C B e$. It follows by
Theorem~\ref{maxmod} that
$\tp(A/C)\cup \tp(Be/C)\vdash \tp^+(A_KB_Ke/C)$, and 
so $\tp(A/C)\vdash |f(x,e)| < 1$.      \qed

Theorem~\ref{moregtoJ} and Corollary~\ref{orthright2} to follow can
also be proved very rapidly using the observation $J(p)=J(p|C)$
for stably dominated $p$ based on $C$, made at the beginning of the chapter.
However Lemma~\ref{orthright} covers cases that cannot be seen in this way.

\begin{theorem}\label{moregtoJ}
Assume $C=\acl(C)$, $A$, $B$ are ${\cal L}_{\GG}$-structures, finitely $\acl$-generated over $C$, and that 
$C\subseteq A\cap B$, with $\tp(A/C)\perp \Gamma$
and $A\dnf^g_C B$. Then $A\dnf^J_C B$.
\end{theorem}

{\em Proof.} We shall generically resolve $C$ and then use the last theorem.
Let $C'$ be a closed resolution of $C$, and $M$ a canonical open resolution of $C'$.
We may choose $C'$ and $M$  so that
$C'\dnf^g_C AB$ 
and $M\dnf^g_{C'} C'AB$. 
Note that here we make essential use of {\em transfinite} sequential independence.
In particular, as 
$\tp(A/C)\perp \Gamma$, $A\dnf^g_C C'$ by Proposition~\ref{perpgamma}, 
so $\tp(A/C')\perp \Gamma$ by Lemma~\ref{ortheasy}(iii).

Now $\tp(C'/C)\perp \Gamma$, so by Proposition~\ref{perpgamma}
$AB\dnf^g_C C'$ via any generating
sequence, and hence $A\dnf^g_{BC}BC'$. As $A\dnf^g_C B$, we obtain
$A\dnf^g_C BC'$, so $A\dnf^g_{C'} BC'$.
It follows that $\tp(A/BC')\perp \Gamma$, so as
$M\dnf^g_{C'B}A$ (since $M\dnf^g_{C'} AB$) we obtain
$A\dnf^g_{C'B} M$.  Since also $A\dnf^g_{C'} BC'$,
we have $A\dnf^g_{C'}BM$. Hence, $A\dnf^g_C BM$ (as $A\dnf^g_C C'$), so, finally,
$A\dnf^g_{M} B$.
In addition, as $M\dnf^g_{C'} A$, and $\tp(A/C')\perp \Gamma$, 
we have $A\dnf^g_{C'} M$. Thus, $A\dnf^g_C M$, so 
$\tp(A/M)\perp \Gamma$.

Put $A'':=\acl(AM)$ and $B'':=\acl(BM)$. Then
$A''\dnf^g_{M} B''$ and $\tp(A''/M)\perp \Gamma$. Furthermore, there 
is a resolution $A'''$ of $A''$
such that $A'''\dnf^g_M B''$ and $\tp(A'''/M)\perp \Gamma$ (by Lemma~\ref{orthres}). By 
Lemma~\ref{orthogwelldefined},
$A''' \dnf^g_M B''$ via a sequence of field elements.
It follows from Theorem~\ref{11.2} that $A'''\dnf^J_M B''$, so $A''\dnf^J_{M} B''$.

To prove the theorem, suppose $f(x,e)\in J(A/B)$; that is, 
for some subsequence $a$ of $A$, $\tp(A/B) \vdash |f(a,e)| < 1$ 
(for ease of notation -- not by
Lemma~\ref{Jres} --  we are assuming that $a$ 
is a tuple of field elements). We must show 
$\tp(A/C)\vdash |f(x,e)| < 1$, or equivalently, 
$\tp(e/C)\vdash |f(a,y)| < 1$. So suppose
$e'\equiv_C e$. Since we may translate $ee'$ 
over $AB$, we may suppose $M\dnf^g_{AB} ABee'$. 
Also, $M\dnf^g_C AB$, so
$M\dnf^g_C ABee'$. Thus $e\equiv_{M} e'$ by Corollary~\ref{infgstat}. 

We have $f(x,e)\in J(A''/B'')$, and $A''\dnf^J_{M} B''$, so 
$\tp(A''/M)\vdash |f(x,e)| < 1$. Since $x$ refers just to elements 
from $A$, $\tp(A/M)\vdash |f(x,e)|< 1$, so $\tp(e/M)\vdash |f(a,y)|\leq 1$. 
Hence, $|f(a,e')| < 1$. \qed

\bigskip

By the use of the following lemma, we can prove the analogue of Theorem~\ref{moregtoJ}, 
but with the orthogonality condition on the right.

\begin{lemma} \label{orthright}
Suppose $C=\acl(C\cap K)$ and   $C\subseteq A \cap B$. Suppose
 too that $A$ is resolved,  
$A\dnf^g_{C\cup \Gamma(A)} B$ via a sequence of field elements,
and that $\tp(\Gamma(A)/C)\vdash \tp(\Gamma(A)/C\Gamma(B))$. Then 
$A\dnf^J_C B$.
\end{lemma}

{\em Proof.} 
Let $A_1:=C\cup \Gamma(A)$. Since $\Gamma(A)$ consists of 
elements of $S_1$,
$A\dnf^J_{A_1} B\cup \Gamma(A)$ by Theorem~\ref{11.2}. Suppose now 
$f(X,Y)\in {\mathbb Z}[X,Y]$, and $e\in K^n$ with $f(x,e)\in J(A/B)$, so 
$\tp(A/B)\vdash |f(x,e)| < 1$ (so 
$x$ corresponds to some subtuple of $A$). We must show 
$\tp(A/C)\vdash |f(x,e)| < 1$. As usual,  we may take 
$x$ in the field sorts by Lemma~\ref{Jres}.

\medskip

{\em Claim.} $\tp(A/C)\cup \tp(A_1/B) \cup 
\tp(e/B)\vdash |f(x,y)| < 1$.

\medskip

{\em Proof.} Fix $A_1$ and $B$, and suppose $e'\equiv_B e$. We must show $\tp(A/A_1) \vdash |f(x,e')|<1$.
Suppose now that $A'\equiv_{A_1} A$, with $A'\dnf^g_{A_1} Be'$ via the corresponding sequence of field elements. 
Then $\tp(A'B)=\tp(AB)$ (since also $A\dnf_{A_1}B$), so $\tp(A'/B)\vdash |f(x,e')|< 1$
(as $\tp(A/B)\vdash |f(x,e)|<1$). By 
Theorem~\ref{11.2} (since $\Gamma(A)\subset S_1$), 
$A'\dnf_{A_1}^J Be'$, so $\tp(A'/A_1) \vdash |f(x,e')|< 1$. 
As $A'\equiv_{A_1} A$, this proves the claim.
 
\medskip

Now by assumption 
$\tp(A_1/C)\vdash \tp(A_1/C\Gamma(B))$, so $\tp(A_1/C)\vdash \tp(A_1/B)$ (as $\Gamma$ is stably embedded and $A_1 \subset C \cup \Gamma$).
 Hence, 
$\tp(A/C)\vdash |f(x,e)| < 1$, as required. \qed

\begin{corollary}\label{orthright2}
Suppose $C=\acl(C)\subseteq A\cap B$ with $A$, $B$ finitely $\acl$-generated 
over $C$, and that $\tp(B/C)\perp \Gamma$ and
$A\dnf^g_C B$. Then $A\dnf^J_C B$.
\end{corollary}

{\em Proof.} 
Since we may replace $A$ by a generic resolution over $B$, we may suppose that $A$ is resolved. 
By Proposition~\ref{perpgamma}, $B\dnf^g_C A$.
As in the proof of Theorem~\ref{moregtoJ}, let $C'$ be a generic closed 
resolution of $C$, and $M$ be a canonical open resolution of $C'$, with 
$C'\dnf^g_C AB$ and $M\dnf^g_{C'} C'AB$. Put $A'':=\acl(AM)$ 
and $B'':=\acl(BM)$, so as in \ref{moregtoJ} we have $B''\dnf^g_M A''$ 
and $\tp(B''/M)\perp \Gamma$. Then
 $B''\dnf^g_{M\cup \Gamma(A'')} A''$.
Hence, as $\tp(B''/M\cup \Gamma(A''))\perp \Gamma$, we have
$A''\dnf^g_{M\cup \Gamma(A'')} B''$ via a sequence of field elements.
 Now apply Lemma~\ref{orthright},
noting that $\Gamma(B'')=\Gamma(M)$ as $\tp(B/M)\perp \Gamma$, to obtain 
$A''\dnf^J_M B''$. The last part of the proof of Theorem~\ref{moregtoJ} 
now shows $A\dnf^J_C B$.
\qed

\bigskip

We have the following strong converse in any sorts, only assuming some
orthogonality to $\Gamma$.

\begin{theorem}\label{9.4+}
Let $A=\dcl(A)$, $B=\acl(B)$, $C=\acl(C)$. Suppose $A\dnf^J_C B$, and 
that $\tp(A/C)\perp \Gamma$, or $\tp(B/C)\perp \Gamma$. Then $A\dnf^g_C B$.
\end{theorem}

{\em Remark.} 1. Since one of $\tp(A/C),\tp(B/C)$ is orthogonal to $\Gamma$, 
the notion of independence in the conclusion is symmetric and independent of 
a choice of generating set. Furthermore, if $\tp(A/C)\perp \Gamma$, then 
the conclusion yields $\Gamma(AB)=\Gamma(B)$, and if $\tp(B/C)\perp \Gamma$, 
we obtain $\Gamma(AB)=\Gamma(A)$.

2. In the proof below, we use Srour's notion of an {\em equation} over a set $B$, 
that is, a formula $\phi(x)$ such that the intersection of any set of conjugates 
of $\phi(\U)$ over $B$ is equal to a finite sub-intersection (see e.g. \cite{sr}). We 
shall use results from \cite{ps} about equations in stable theories, applied within the
stable structure $\St_C$.

\medskip

{\em Proof.} It suffices to show $\St_C(A)\dnf \St_C(B)$. For 
by Proposition~\ref{perpgamma}, under either orthogonality hypothesis this yields 
$A\dnf^g_C B$.

Let $t=(t_1,\ldots,t_m)\in \St_C(A)$. We may suppose  $m=1$ and $t\in \red(s)$
where $s\in \dcl(C) \cap S_\ell$;
for if $t_i\in \res(s_i)$, where $s_i \in \dcl(C)$, then $(t_1,\ldots,t_m)$
 can be regarded as an element of $\red(\Lambda(s_1)\times \ldots \times \Lambda(s_m))$, and 
$\Lambda(s_1)\times \ldots \times \Lambda(s_m)$ 
is a $C$-definable lattice. 
 We must show $t\dnf_C \St_C(B)$, 
that is, $\RM(t/C)=\RM(t/\St_C(B))$ in the structure $\St_C$. Let $M$ be a model containing $B$, 
with $t\dnf^g_B M$. Then by Proposition~\ref{gtostab}, $t\dnf_B \St_B(M)$ in $\St_B$, so
$t\dnf_{\St_C(B)} \St_C(M)$ in $\St_C$ (compare Proposition~\ref{stab0.3},
 (i) $\Rightarrow$ (iii)); hence
$\RM(t/\St_C(B))=\RM(t/\St_C(M))$ in $\St_C$. Thus, our task is to show $\RM(t/C)=\RM(t/\St_C(M))$ in $\St_C$.
As $M$ is a model, 
there is an $M$-definable $R$-module isomorphism $\psi:\Lambda(s) \rightarrow R^\ell$ 
(an affine map) inducing an isomorphism (also denoted $\psi$) $\red(s) \rightarrow k^\ell$.  
Put $t'=(t_1',\ldots, t_{\ell}'):=\psi(t)\in k^{\ell}$. 

Let $P_C$ (respectively $P_M$) denote the solution sets of $\tp(t/C)$ (respectively 
$\tp(t/M)$), and let $Z_C$ (respectively $Z_M$) be the Zariski closure of 
$\psi(P_C)$ (respectively $\psi(P_M)$) in $k^{\ell }$. 
We have $P_C \supseteq P_M$ and $Z_C \supseteq Z_M$.
Now $\RM(P_C)=\RM(\psi(P_C))=\RM(Z_C)$, and similarly
$\RM(P_M)=\RM(Z_M)$. Thus, we must show $Z_C=Z_M$, as this gives $\RM(P_C))=\RM(P_M)$. 

So suppose $f(X)\in k(M)[X]$, with $f(t')=0$. We wish to show 
that $f$ vanishes on $Z_C$. Lift $f$ to a polynomial $F(X)$ over $R\cap M$. Then 
if $y_1\in t'_1, \ldots, y_\ell\in t'_\ell$ and $y=(y_1,\ldots y_\ell) \in K^{\ell}$ we 
have $|F(y)|<1$. Put $H(X):=F(\psi(X_1,\ldots,X_\ell))$. 
Then $|H(x_1,\ldots,x_\ell)|<1$ for any $x=(x_1,\ldots,x_{\ell})\in t$. 

{\em Claim.} If $u$ is a single variable ranging over elements of $T_\ell$, then
$\tp(t/B)\vdash f(\psi(u))=0$.

{\em Proof of Claim.} The formula $f(\psi(u))=0$ is over $M$: write it as $\delta(u,m)$.
Working over $B$, $\delta(u,m)$ is an equation
in the sense of Srour \cite{sr}. For if we fix
 an arbitrary affine bijection
$\red(s) \rightarrow k^\ell$ (not over $B$), then conjugates over
$B$ of the formula $f(\psi(u))=0$ become formulas $g(x)=0$ over $k$, and 
the ideal in $k[X]$ generated by any set of these is finitely generated.
Furthermore, the formula $f(\psi(u))=0$ 
can be interpreted in the stable structure 
$\St_B$, since this structure is stably embedded. That is, we may suppose $m$ is from $\Mst$.
Also, $t\dnf_B M$ and $f(\psi(t))=0$. Let
$X$ be the intersection of the solution sets of all the formulas $\delta
(u,m')$
such that $m'\in \Mst $ and $\models \delta(t,m')$ holds. Then by Proposition 4.2 of \cite{ps}, using that
$B=\acl(B)$, we find that $X$ is $B$-definable. The claim follows. 

It follows that
$$\tp(t/B) \vdash \forall x\in u(|H(x_1,\ldots,x_\ell)|<1).$$ 
Thus, $H\in J(t/B)$. But $J(t/B)=J(t/C)$, so
$$\tp(t/C)\vdash  \forall x\in u(|H(x_1,\ldots,x_\ell)|<1).$$ 
Hence, $\tp(t/C)\vdash f(\psi(u))=0$. 
It follows that $f$ vanishes on $Z_C$, as required.
 \qed

\bigskip

We summarise the  equivalences of different notions 
of independence, under an assumption of orthogonality to $\Gamma$.

\begin{theorem}\label{summary}

Let $A$, $B$, $C$ be algebraically closed ${\cal L}_{\GG}$-structures, with $C\le A\cap B$ and $A,B$ finitely $\acl$-generated over $C$,
and suppose that $\tp(A/C)$ is orthogonal to $\Gamma$.
Then the following are equivalent.

(i) $A\dnf^J_C B$

(ii)  $A\dnf^g_C B$

(iii) $B\dnf^g_C A$

(iv) $B\dnf^J_C A$

(v)  $A\domind_C B$.

\noindent
If all structures are valued fields, then all of the above are equivalent to

(vi) $A\dnf^m_C B$.
\end{theorem}

{\em Proof.} 
We have (ii) $\Leftrightarrow$ (iii) and (ii) $\Leftrightarrow$ (v) by Proposition~\ref{perpgamma}, and parts
(i) $\Rightarrow$ (ii) and  (iv) $\Rightarrow$ (iii) are from Theorem~\ref{9.4+}. For 
(ii) $\Rightarrow$ (i) see Theorem~\ref{moregtoJ}, and (iii) $\Rightarrow$ (iv) comes 
from Corollary~\ref{orthright2}. Finally, for the equivalence of (vi) in the case 
of fields, see Theorem~\ref{allsame}. \qed

\bigskip

We conclude with various examples. These illustrate that sequential independence
and $J$-independence are not equivalent in general.
In all such examples, stable domination of course fails. It is worth
noting that for singletons in the field sort, the notions are in fact
equivalent, as we first prove below. Thus the examples must all involve at least
two field elements, and there cannot be any orthogonality to the value
group.

\begin{proposition}\label{equivsingles}
Let $C=\acl(C) \subseteq B$ and $a\in K$ be a single element. 
Then $a\dnf^J_C B$ if and only if $a\dnf^g_C B$.
\end{proposition}

{\em Proof.} Suppose first that $a\df^g_C B$.
We suppose that $\tp(a/C)$ is the generic type of a unary set
$U$ (a ball, or the intersection of a sequence of balls, or $K$ itself). Then there is a 
$B$-definable proper sub-ball $V$ of $U$ containing $a$. If $U$ is a closed ball, 
and $V$ is an open ball of the same radius $\gamma$, 
choose  $d\in K$ with $|d|=\gamma$. Then $d^{-1}(x-a)\in J(a/B)\setminus 
J(a/C)$, so
$a\df^J_C B$; indeed, if $a'\equiv_B a$ then $a'\in V$, so $|x-a|<|d|$, but there is $a'\equiv_C a$ with $a'\not\in V$, and then $|x-a'|=|d|$. In all the other cases for $U$, $\rad(V)<\rad(U)$. Now choose $d\in
 K$ with $\rad(V)<|d|<\rad(U)$. Then again, $d^{-1}(x-a)\in J(a/B)\setminus J(a/C)$.

For the other direction, assume $a\dnf^g_C B$. Let $U$ be a $C$-unary set
such that $\tp(a/C)$ is the generic type of $U$. 
Let $f(x)\in J(a/B)$. We must show that $f\in J(a/C)$. So suppose that 
$a'\equiv_C a$ with $a'\not\equiv_B a$. We must show that $|f(a')|<1$. Choose 
$a''$ generically in $U$ over $B \lceil f\rceil a a'$. Then 
$a''\equiv_{B} a$,
so
$|f(a'')|<1$. Hence it suffices to 
show $|f(a')|\leq |f(a'')|$.

As $a'\not\equiv_B a$, there is a $B$-definable proper sub-ball $V$ of $U$ containing $a'$. 
Factorise $f$ as  $f(x)=d\prod(x-d_i)$.
 If $d_i\notin U$, then $|x-d_i|$ is constant for all
$x\in U$, so $|a''-d_i|=|a'-d_i|$. If $d_i\in V$, then $|a'-d_i| < |x-d_i|$ for all $x\in
U\setminus V$, so $|a'-d_i|<|a''-d_i|$. Suppose $d_i\in U\setminus V$. 
In this case, the generic choice of $a''$ forces $|a''-d_i|\geq |a'-d_i|$. Considered 
together, these three cases ensure $|f(a')|\leq |f(a'')|$, as required.
\qed

\begin{example} \rm
Here $J$-independence holds but sequential independence
fails.  Let $C$ be an algebraically closed valued field which is not
maximally complete. Let $\{ U_i: i\in I\}$ and $\{ V_j: j\in J\}$ be chains of
$C$-unary sets, chosen so that $U=\bigcap_{i\in I}
U_i$ and $V=\bigcap_{j\in J} V_j$ are complete types. We may suppose that
the cuts $\rad(U)$ and $\rad(V)$ in $\Gamma$ are equal;
also, that $U$ and $V$ are sufficiently independent that
if $a\in U$ then no element of $\acl(Ua)$ lies in $V$, and likewise with $U$
 and $V$ reversed.
Let $b_1 \in U$,
$b_2\in V$. Choose $a_1\in U$ generic over $Cb_1b_2$ and $a_2\in V$
generic over $Cb_1b_2a_1$. Notice that $\Gamma(Ca_1a_2)=\Gamma(C)\ne
\Gamma(Ca_1a_2b_1b_2)$. Then $a_1a_2 \dnf^g_C b_1b_2$, but $|a_2-b_2| >
|a_1-b_1|$ so $a_2a_1 \df^g_C b_1b_2$. By Theorem~\ref{11.2}, $a_1a_2
\dnf^J_C b_1b_2$ and hence also $a_2a_1 \dnf^J_C b_1b_2$.
 \end{example}

\begin{example} \rm \label{14.12}
This is an example where sequential independence holds but $J$-independence
fails. Let $C$ be any model. Let $p(x,y)$ be a type over $C$ such that
$p(x,y) \proves |x|< 1 \wedge |y|<1 \wedge |xy|=\gamma$, where $\gamma = |c|$ for
some $c\in C$. Choose $b$ generic in $\M$ and $a_1$, $a_2$ realising $p$
such that $a_1a_2 \dnf^g_C b$. The sequential independence implies that
$|a_1|> |b|$, and since $|a_1a_2|=\gamma$, we have $|a_2|<\gamma|b|^{-1}$.
Thus if $f(x) = c^{-1}bx_2$ then $f(x)\in J(a_1a_2/Cb)\setminus
J((a_1a_2/C)$. \end{example}

\begin{example} \rm \label{14.13}  \rm
We give an example of a type with no $J$-independent extension over a certain set.
Let $C$ be a valued field, let $\gamma$ lie in the Dedekind completion of $\Gamma(C)$ but not in $\Gamma(C)$,
and let $a\in K$ with $|a|=\gamma$. Put $A=C(a)$.
Let $B=C(b_1,b_2)$, where $|b_1|>|b_2|$ and both $|b_1|$, $|b_2|$ lie in the same cut of $\Gamma(C)$ as $\gamma$.
 If $A\dnf^J_C B$, then one cannot have $|a|<|b_1|$, for otherwise $b_1^{-1}x\in J(a/B)\setminus J(a/C)$.
 Thus $|a|\geq |b_1|$.  Similarly $|a ^{-1}| \geq |b_1|^{-1}$.  
Thus $|a|=|b_1|$.  But similarly $|a|=|b_2|$, a contradiction.  Thus, $\tp(A/C)$ has no $J$-independent extension over $B$.
\end{example}

\begin{remark} \rm 
 The proof in \cite{hhm} of elimination of imaginaries of \cite{hhm} can now be summarized as follows.
First, for any stably dominated
 type $q$, the $J$-invariant $J(q)$ shows that the canonical base of $q$ is coded.
 Indeed $J(q)$ is a submodule of a polynomial ring over $K$; hence it is canonically
 the union of $R$-submodules of finite-dimensional vector spaces, with canonical
 $K$-bases; one quickly reduces to sublattices, i.e. elements of $S_n$,
  and associated elements of $T_n$.

 Finite sets and unary sets are dealt with separately.  Given that, it remains to   code  
 functions $f$, on some neighborhood of each complete type over a set $C=\acl(C)$; in fact unary types suffice.   If the type $p$ is  stably dominated, we first code the germ of 
 $f$ on $p$.  Let $q$ be the stably dominated type obtained by applying the function
 $x \mapsto (x,f(x))$ to $p$.  Then the canonical base of $q$, coded above, serves
 also to code the germ of $f$.    Since the germ is strong, there exists a $C$-definable function agreeing, on a neighborhood of the given type, with the original one.  
 
 If $p$ is a limit of a $\G$-family of stably dominated
types, defined over $C$, we code the function restricted to each stably dominated type in the family using the previous step.  Over a rich enough base, any type is a $\G$-limit of stably dominated types; but such a family
may not be defined over $C$.  For 1-types however,   there exists such a family whose {\em germ} is defined over $C$.   
Using our understanding of definable functions from $\G$,
we are able to deal with germs of maps from $\G$ and conclude the proof.  

 The appropriate higher-dimensional generalization
of the     fact quoted above about 1-types, over an arbitrary base, is not yet clear.

\end{remark}


\chapter{Other Henselian fields} \label{henselian}

We give two examples of metastable theories other than ACVF.  This means that
much of the technology developed in the present manuscript applies to these theories.
The main intention here is to illustrate the way that theorems proved in the ACVF context
apply to other valued fields; all the real content of the lemmas below derives
from ACVF theorems in the main text.  At the end we make remarks concerning 
generalizations.

\begin{theorem} \label{hensel} $\Th(\Cc((t)))$ is metastable over the value group.  \end{theorem}

Let $T=\Th(\Cc((t)))$.    This theory admits elimination of valued field quantifiers, relative
to $\G$-quantifiers.  
It is simpler than other Hensel fields of residue characteristic $0$ in that it also eliminates
residue field quantifiers; this allows the statement below to have a simpler form,
involving $\G$ alone, but a similar quantifier elimination involving $k$ (or rather RV) is 
true in general.  These observations can be found in \cite{fvk}.

Let $L \models T$.  We view $L$ as a subfield of a model $M$ of ACVF.  As $L$
is Henselian, $\Aut(L^{\alg}/L) = \Aut_v(L^{\alg}/L)$, i.e. every field-theoretic automorphism
of $L^{\alg}/L$ preserves the valuation.  By 
Galois theory, $L=\dcl(L)$ within the field sort of $M$. 
We interpret the sorts $S_n$ by $(S_n)_{L} = \GL_n(L)/\GL_n(R_L)$.  We similarly
interpret the sorts $T_n$ though in the case of $T$ they are redundant. 
 In particular, $\G_L = \{v(a): a \in L \}$.  This is distinct from $\G(L)$,
 but $\G(L) = \G_L \tensor \Qq$ is the divisible hull of $\G_L$.  
 
Let  $\U^T$ be a universal domain for $T$.  We view it as a subset of a universal
domain $\U$ for ACVF, interpreting the $\GG$-sorts as above. 

It can be shown with methods similar to those used here that 
$T$ admits elimination of imaginaries in this language; but we will not require this fact.

{\noindent \bf Convention.}   By `formula', `definable function',
`type' we mean the quantifier-free ACVF notions; thus   $\tp(a/C)$   denote the
 quantifier-free (ACVF) type of $a$ over $C$.   The corresponding $T$-notions are
 marked with $T$, e.g.  
$\tp_T(a/C)$ denotes the $T$-type;   $\Inv^T_x(C)$ is the set of $\Aut(\U^T/C)$-invariant
 types of $\U^T$ in the  variable $x$.    We let $\G^T = \G_{\U^T}$, $S_n^T = (S_n)_{\U^T}$, 
 and denote the field
sort of $\U^T$ by  $K(\U^T)$.    For any $C \leq T$, let $\G^T(C) = \G(C) \meet \G^T$.  
When $d \in \G^T$,
$\tp_T(d)$ denotes the set of $\G$-formulas true of $d$ in $\U^T$.
\medskip

Let $\F$ be the collection of definable functions $f$ on a sort of $\GG$ into $\G^k$, with the property  that
 if $a \in \U^T$ then $f(a) \in \U^T$. Also let
  $\F(C)$ be the set of functions $f(x,c)$
with $f \in \F$ and $c \in C^m$.   We can think of $\F$ as generating all definable functions into $\G$, by virtue of:  

\begin{lemma} \label{ba0}  Let $f$ be a definable function into $\G$.  Then for some $n$,
 $nf \in \F$.  \end{lemma}
 
 {\em Proof.}  Since the algebraic closure of $\U^T$ is a  model of ACVF,
 and $\acl(\U^T)
 = \Qq \tensor \G^T$, it follows that any element of $\G$
 definable over $\U^T$ lies in $\Qq \tensor \G^T$.  Thus for any $a \in \U^T$,
 for some $n$, we have $n f(a) \in \G^T$. The graph of $n f$ is defined by a quantifier-free
 formula $\psi_n$, and we have shown that $\U^T \models (\forall x) \bigvee_{n \in \Nn} (\exists y) \psi_n(x,y)$.  By compactness, for some finite set $n_1,\ldots,n_k \in \Nn$,
 $\U^T  \models (\forall x) \bigvee_{i=1}^k   (\exists y) \psi_{n_i} (x,y)$, i.e. some
 $n_i f(x)  \in \U^T$ for any $x \in \U^T$.     Let $n= \Pi_{i=1}^m n_i$.
 Then $n f(x) \in \G^T$ for any $x \in \dom(f)$.  So $n f \in\F$.  \qed \medskip

Let $C \leq \U^T$,  $s \in S_n^T$ defined over $C$, and let $V = \red(s)$.
Write $V^T$ for the image of the $\U^T$-points under $\red$.

Let $x$ be a variable of one of the sorts $\GG$.  By a {\em basic formula} we
mean one of the form $\psi(g(x))$, where $g \in \F$, and $\psi$ is a  formula of $\G^T$.   
If $f,f' \in \F$ then $(f,f') \in \F$, so a Boolean combination of basic formulas is a basic formula.

\begin{lemma} \label{qe-h} 

1)  Let $c,c' \in \U^T$, and assume $\tp_T(g(c)) = \tp_T(g(c'))$
 for any $g \in \F$.  Then   $\tp_T(c)=\tp_T(c')$.  
 
2)  Every formula in the sorts $\GG$ is $T$-equivalent to a 
Boolean combination of  basic formulas. 

3)  For any $C \subseteq \U^T$ and $c,c' \in \U^T$,  if $\tp_T (g(c)/C) = \tp_T(g(c')/C)$ 
 for any $g \in \F(C)$,  then $\tp_T(c/C)=\tp_T(c'/C)$.  
 
 Equivalently,  let $C' = C \union \G^T(C(c))$, and let $T_{C'}$ be the elementary diagram
 of $C'$ in $\U^T$.   
Then  $\tp(c/C') \union T_{C'}  \vdash \tp_T(c/C')$.
\end{lemma}

{\em Proof.}  1)  By relative quantifier elimination, this is true for 
formulas in the field sort. We will reduce to this case using resolution.

Conjugating by an element of $\Aut(\U^T)$, we may assume $g(c)=g(c')$ for
any $g \in \F$.     By Lemma~\ref{ba0}, $\tp(c/\G)=\tp(c'/\G)$, and  in
particular $\tp(c)=\tp(c')$. We may assume $\G(c) \neq (0)$; otherwise the
hypothesis will apply  to the tuples $(c,t)$ and $(c',t)$, while the
conclusion for these tuples will be  stronger. By Corollary~\ref{R3.2}
(applied to ACVF) there exists  a sequence of  field elements $d $ such
that $\tp(d/c)$ is isolated, $\G(d) = \G(c)$, and $c \in \dcl(d) $.   Say
$c = F(d)$, $F$ an ACVF-definable function. Moreover $\tp(d/c)$ is realized
in any $F=\dcl(F)$ containing   representatives of the classes coded by
$c$ .  It follows that $d$  can be found in $\U^T$, and that there exists
$d' \in \U^T$ with $\tp(d'c')=\tp(dc)$.

As in Corollary  \ref{R7.2}, we have $\G(d)=\G(c)$.
Hence for any $g \in \F$ there exists a definable $h$ into $\G$ with
with $g(d)=h(c)$ and $g(d') = h(c')$.   Since $n h \in \F$ for some $n $,
it follows that $g(d) = g(d')$.  By the field case, we have $\tp_T(d) = 
\tp_T(d')$. Since $c=F(d), c'=F(d')$, we have $\tp_T(c) = \tp_T(c')$.   This
proves  (1).

 By (1),   any $T$-type is determined by
the basic formulas in it.  (2) is a standard consequence, using compactness.

(3) It follows from (2), replacing variables by constants, that every formula over $C$ is equivalent to a Boolean combination of  basic formulas over $C$.
   So if  if $\tp_T(g(c)/C) = \tp_T(g(c')/C)$ 
 for any $g \in \F(C)$,  then $\tp_T(c/C)=\tp_T(c'/C)$.  
  
  \qed
 \medskip

If $q \in \Inv_{xy}(C)$ and $p = q|x \in \Inv_{x}(C)$, say 
 $q$ is orthogonal to $\G$ relative to $p$ if for any $C' \supseteq C$, 
 if $(c,d) \models q | C'$ then $\G(C'c) = \G(C'cd)$.  
  \def\tP{{\tilde{P}}}
  
\begin{lemma} \label{ba2}  Let $P \in \Inv^T_x (C)$, $r \in S_{xy} (C)$,
 $p   \in \Inv_x (C)$. Let $\tilde{p} = p | \U^T$, and suppose $\tilde{p}\subseteq P$. 
Assume: $p \union r \vdash q$ for some $q \in \Inv_{xy}(C)$, and that
 $q$ is orthogonal to $\G$ relative to $p$. 
 Then $P \union r \vdash Q$ for a unique $Q \in \Inv_{xy}^T(C)$.
 \end{lemma} 

{\em Proof.}  Let $C \subseteq C' \subset \U^T$, 
$c \models P| C'$, $(c,d) \models r$.  We have to show that $\tp_T(d/C'c)$ is
determined.    But  by the relative orthogonality assumption, $\G^T(C'cd) = \G^T(C'c)$,
so by Lemma~\ref{qe-h}, $\tp(d/C'c)$ implies $\tp_T(d/C'c).$  \QED

\begin{corollary}  \label{ba4}  Let $C = \acl(C) \meet \U^T$.   Then every $T$ -type over $C$
extends to an $\Aut(\U^T/C)$-invariant $T$-type over $\U^T$.  
\end{corollary} 

{\em Proof.} First, it is an exercise (see also \cite{hmartin}) to check that 
any 1-type over $C$ of $T$ in the $\Gamma$-sort
has an invariant extension.

 We next argue that for a single field element $a$, $\tp_T(a/C)$ has an $\Aut(\U^T/C)$-invariant extension 
over $\U^T$. As for ACVF, if $(V_i:i\in I)$ is the set of $C$-definable closed balls
which contain $a$, we say $a$ is generic in $V:=\bigcap(V_i:i\in I)$. 

Suppose first 
that the chain has a least element, namely
 $V$, and that $V$ is a closed ball; note that as the value group is discrete, there 
is no open/closed distinction 
for balls. Let $W = \red(V)$ be the corresponding $k$-space.
Let $\red: V \to W$ be the natural map.   Consider $T$-types $Q(x,y)$ with  
$\red(y)=x$.  Then the
hypotheses of Lemma~\ref{ba2} apply, as $\{y: \red(y)=x\}$ is also a closed ball so
 has generic orthogonal to $\Gamma$.
Hence a complete $C$-invariant $T$-type is  
determined by
(i) genericity of $x$ in the affine $k$-space (ii) the formula $\red(y) 
=x$ (iii) a complete type over $C$ in $x$.

Suppose next that $V$ is not a ball (so $I$ has no least element) but has a $C$-definable
 point (or subtorsor) $x_0 $ inside. Let $\gamma:=|a-x_0|$.  Extend first $\tp(\gamma/C)$ to an invariant type, and then find the
 generic type of
$B_{\leq \gamma}(x_0)$
 using the closed ball case.

Finally, suppose that $I$ has no least element and $V$ has no $S$-definable ball. By
 the last case, for any
$d\in V$ we described a $C(d)$-invariant type consisting of elements 
of $V$; it did not depend on the choice of $d$, so is $C$-invariant.

To complete the proof of the lemma, it suffices by
 Lemma~\ref{st-dom-eq-0} (i) to prove it for types
of sequences of field elements.  We allow infinite sequences, and reduce immediately
to transcendence degree $1$, i.e. to $\tp(a,b / C)$ where $a$ is a singleton,
and $b$ enumerates a part of $\acl(C(a))$.  Let $r=\tp(a,b / C)$.  Let $P$ be an
$\Aut(\U^T/C)$-invariant extension of $\tp(a/C)$, and let $p$ be the restriction to a 
quantifier-free ACVF-type over $C$.  The conditions of Lemma~\ref{ba2} are met,
so $P \union r$ generates a complete $T$-type $Q$. 
 Since $r$ (being over $C$)
and $P | \U^T$ are $\Aut(\U^T/C)$-invariant, so is $Q$.  
\QED
 
\begin{lemma} \label{ba6}    Let $C$ be a maximally complete model of $T$, $C \prec M \models T$, 
and let $a$ be a tuple from $M$.  Let $C' = C \union \G^T({C(a)})$.   Then 
 $\tp_T(A/C')$ is stably dominated.  
\end{lemma}
\def\bL{{\mathbf L}}
{\em Proof.}   Let $C^+ = C \union \G(L))$, and let $\bL$ enumerate $L$.  
By Theorem~\ref{fulldom} (i), $\tp(\bL / C^+)$ is stably dominated;   by Remark~\ref{fulldom-r},
it is in fact stably dominated by a sequence $b \in \dcl(L) \meet \U^T$. 
Let $p$ (respectively $q$) be the $\Aut(\U/C^+)$-invariant type   extending
$\tp(b/ \acl(C^+))$ (respectively $\tp(\bL,b / \acl(C^+))$); let  $r=\tp(\bL,b/ C^+)$.  Then
stable domination via $b$ means that  $r \union p \vdash q$.  Noting
that $b$ lies in the stable part of $\U^T$, let  $P$ be the
$\Aut(\U^T /\acl^T(C^+))$-invariant type extending $\tp_T(b/\acl(C^+))$.  
By Lemma~\ref{ba2}, $P \union r$ generates an invariant type, extending
$\tp(\bL,b / \acl(C^+)$, which is therefore stably dominated.  Hence
$\tp_T(\bL / C^+)$ is stably dominated.  By Corollary~\ref{st-dom-eq-0} (ii),
   $\tp_T(a/C')$ is stably dominated.   \qed

\medskip

{\em Proof of Theorem~\ref{hensel}.}  Immediate from   Corollaries~\ref{ba4} and \ref{ba6},
taking into account   
Corollary~\ref{metastable-eq}.  \qed

\medskip

{\noindent \bf Generalizations.} 

\noindent
1)  Other value groups. \\
The only facts used about $\Th(\Zz)$ are:

a)   Every type over a set $C$ extends to an $\Aut(\U/C)$-invariant type.

b)  If $A \prec B \models \Th(\Zz)$ then $A/B$ is torsion free.  \\
The class of value groups satisfying include $\Th(\Zz^n)$, and intermediate groups 
between $\Zz$ and $\Qq$. 
 For applications the class seems comfortable.   

\noindent
 2)  Valued fields of mixed or positive characteristic $p$, with $p$-divisible value groups
 satisfying (a,b).    If the residue field remains algebraically closed, Theorem~\ref{hensel} goes through.

\noindent
3)  Non-algebraically closed residue fields.  This raises two  issues.
 
 a)  The quantifier elimination is relative to higher congruence groups.  In
 residue characterstic $0$, we have quantifier elimination  
  relative to $RV= K^*/ (1+ \M)$, rather than to $\G$.  
 Nevertheless the method of proof generalizes.  See \cite{hk} for related results,
 describing definable sets in terms of $RV$-definable families of ACVF-definable sets.

 b)  If the residue field is not stable, one cannot expect stably-dominated types;
 the notion needs to be refined.  
   
 Let $S$ be a partial type,    ${\mathcal D}$ be a collection of sorts.
Consider a *-definable map $f: S \to {\mathcal D}$; i.e. $f$ is represented by a sequence
 $f_i: S \to D_i$ with $D_i \in {\mathcal D}$.    
 Let   $I$ be an ideal on the (relatively) definable subsets of $f(S)$.    
 Say $S$ is dominated by ${\mathcal D}$ via $(f,I)$ if for any definable
 set $R \subseteq S$, for all $d \in f(S)$ outside an $I$-small set, 
  the fiber $f \inv(d)$ is disjoint from $R$ or contained in $R$.  
 
 In the case of Henselian fields of residue characteristic $0$, ${\mathcal D}$ will
 be the collection of definable sets internal to the residue field.  When $S$ is a complete stably
 dominated type, $f$ will be the restriction of a definable function of ACVF.
 Therefore the Zariski closure of $f(S)$ has (in ACVF) finite Morley rank. The ideal
 $I$ will consist of  definable sets whose Zariski closure has lower dimension than $f(S)$.

In mixed characteristic $(0,p)$ one has  quantifier elimination 
relative to the family of quotients $R/p^nR$.
The $p$-adics present an interesting case.      For a single such quotient, in the $p$-adic case,
the ideal $I$ is improper; but taken together the maps can be seen as having image in   $\Zz_p$, where the 
proper Zariski closed sets of $n$-dimensional space form 
a proper ideal.  For a partial type $S$, one has a dominating map not into $\Zz_p^n$ itself,
but into a pro-definable set internal to  $\Zz_p^n$ in an appropriate sense (uniformly
internal to the quotients $R/p^nR$).    This of course connects 
to Pillay's idea of compact domination, providing instances of the phenomenon outside
a group-theoretic setting.
 
 \medskip
 
\paragraph{Valued differential fields.}  
 Let ${\rm \widetilde{VDF}}$ be the model completion of the theory VDF of valued differential fields
(\cite{scanlon}).  ${\rm \widetilde{VDF}}$  extends ACVF and admits 
quantifier elimination.  The residue field is   stable as a sort in ${\rm \widetilde{VDF}}$; it has the induced structure
of a  model of the theory DCF
of differentially closed fields.
 
 \begin{theorem}\label{vdf}   ${\rm \widetilde{VDF}}$ is metastable.

 \end{theorem}

{\em Proof.} Let $\U$ be a universal domain for ${\rm \widetilde{VDF}}$.  Relying on 
Corollary~\ref{metastable-eq}, we will work without imaginaries.  We first show that
types extend to invariant types.  
Let $L_v$ be
  the language of valued fields and $\U_v$ the restriction of $\U$ to $L_v$. 
  Let $C \leq \U$, and let $p(x_0)$ be a type over $C$.  Add variables $x_i$ denoting $D^ix$, obtaining
a type $P(x_0,x_1,\ldots)$.  The new type $P$ is generated by the formulas $Dx_i = x_{i+1}$,
along with the restriction $P_v$ of $P$ to $L_v$.  
Let $Q_v$ be an $\Aut(\U_v / C)$-invariant $L_v$-type,  extending $P_v$.  
Given  a valued differential field   $C'$ extending
$C$, let $C''=C'(c_0,c_1,\ldots)$ be a valued field extension of $C'$, generated by 
a realization of $Q_v | C'$.  Then $C(c_0,c_1,\ldots)$  is linearly disjoint from $C'$ over $C$.
Therefore any derivation on $C(c_0,c_1,\ldots)$ extending the given derivation on $C$ extends
to a derivation on $C''$.  It follows that $Q=Q_v \union \{Dx_i=x_{i+1} \}$ is consistent.  By quantifier elimination 
$Q$  is a complete extendible type of ${\rm \widetilde{VDF}}$.  
This shows
that any type over $C$ extends to an $\Aut(\U/C)$-invariant type.

Similarly, let $M$ be a model of ${\rm \widetilde{VDF}}$, maximally complete as a valued field.  For any 
$a_0 \in \U$, let $a_0,a_1,\ldots$ and $P$ be as above.  Let $\g$ enumerate 
the value group of $M(a_0,a_1,\ldots)$, and let $c$ enumerate the residue field of $ M(a_0,a_1,\ldots)$.  Then $P$ is dominated
by $\tp(c/M,\g)$ over $M,\g$, in the sense of ACVF.  It follows as in Lemma~\ref{ba6} that $p$ is stably dominated.
\qed

%
%

\printindex

\noindent
Deirdre Haskell,\\Department of Mathematics and Statistics,\\McMaster University,\\Hamilton, Ontario L8S 4K1, Canada,\\haskell@math.mcmaster.ca

\bigskip

\noindent
Ehud Hrushovski,\\Department of Mathematics,\\The Hebrew University,\\
Jerusalem, Israel,\\
ehud@math.huji.ac.il

\bigskip

\noindent
Dugald Macpherson,\\
Department of Pure Mathematics,\\
University of Leeds.\\
Leeds LS2 9JT, England,\\
h.d.macpherson@leeds.ac.uk

\end{document}